\documentclass[review,3p]{elsarticle}

\usepackage{lineno,hyperref}
\modulolinenumbers[5]
\usepackage{subcaption}             
\usepackage{amsmath,amsfonts,amsthm}            

\usepackage{mathtools,amssymb}          

\usepackage{upgreek}
\usepackage[dvipsnames]{xcolor}
\usepackage{soul}
\usepackage{multirow}

\usepackage{array}
\newcolumntype{C}[1]{>{\centering\let\newline\\\arraybackslash\hspace{0pt}}m{#1}} 
\newcolumntype{L}[1]{>{\raggedright\let\newline\\\arraybackslash\hspace{0pt}}m{#1}} 
\newcolumntype{R}[1]{>{\raggedleft\let\newline\\\arraybackslash\hspace{0pt}}m{#1}} 

\usepackage{booktabs}       

\newcommand{\tabitem}{~~\llap{\textbullet}~~}           
\usepackage{makecell}       
\usepackage{pbox}           
\usepackage{empheq}

\usepackage{colortbl}
\usepackage{esvect}
\usepackage{spreadtab}
\usepackage{numprint}
\usepackage{xstring}

\usepackage[symbol]{footmisc}

\usepackage{siunitx}

\makeatletter       
\newcommand*{\rom}[1]{\expandafter\@slowromancap\romannumeral #1@}
\makeatother

\usepackage{enumitem}

\usepackage{cleveref}
\crefformat{section}{\S#2#1#3} 
\crefformat{subsection}{\S#2#1#3}
\crefformat{subsubsection}{\S#2#1#3}

\usepackage[labelformat=simple]{subcaption}	        	

\usepackage{graphicx}
\usepackage{wrapfig}
\usepackage{lipsum}

\usepackage{pgfplots}       

\usepackage[ruled,linesnumbered]{algorithm2e}		

\begin{document}

\begin{frontmatter}

\title{Balancing truncation and round-off errors in practical FEM: one-dimensional analysis}

 \author[1]{Jie Liu\corref{cor1}}
 \ead{j.liu-5@tudelft.nl}
 \author[1]{Matthias M\"oller}
 \ead{m.moller@tudelft.nl}
 \author[1]{Henk M. Schuttelaars}
 \ead{h.m.schuttelaars@tudelft.nl}
 
 \address[1]{Delft Institute of Applied Mathematics\\ Delft University of Technology\\ Van Mourik Broekmanweg 6, 2628 XE Delft, The Netherlands}
\cortext[cor1]{Corresponding author}

\begin{abstract}
In finite element methods (FEMs), the accuracy of the solution cannot increase indefinitely because the round-off error increases when the number of degrees of freedom (DoFs) is large enough. This means that the accuracy that can be reached is limited. A priori information of the highest attainable accuracy is therefore of great interest.   
In this paper, we devise an innovative method to obtain the highest attainable accuracy.
In this method, the truncation error is extrapolated when it converges at the analytical rate, for which only a few primary $h$-refinements are required, and the bound of the round-off error is provided through extensive numerical experiments. The highest attainable accuracy is obtained by minimizing the sum of these two types of errors.
We validate this method using a one-dimensional Helmholtz equation in space. 
It shows that the highest attainable accuracy can be accurately predicted, and the CPU time required is much less compared with that using the successive $h$-refinement. 
\end{abstract}


\end{frontmatter}

\section{Introduction}

Many problems in engineering sciences and industry are modelled mathematically by initial-boundary value problems comprising systems of coupled, nonlinear partial and/or ordinary differential equations. These problems often consider complex geometries, with initial and/or boundary conditions that depend on measured data \cite{Kumar2016}. 
In some applications, not only the solution, but also its derivatives are of interest \cite{Kumar2016,carey1982derivative}.
For many problems of practical interest, analytical or semi-analytical solutions are not available, and hence one has to resort to numerical solution methods, such as the finite difference, finite volume, and finite element methods. The latter will be adopted throughout this paper and applied to one-dimensional boundary value problems.

The accuracy of the numerically obtained solution is influenced by many sources of errors \cite{ferziger2012computational}: firstly, errors in the set-up of the models, such as the simplification of the domain and governing equations and the approximation of the initial and boundary conditions; next, truncation errors due to the discretization of the computational domain and the use of basis functions for the function spaces defined on it; then, the iteration error resulting from the artificially controlled tolerance of iterative solvers; finally, the round-off error due to the adoption of finite-precision computer arithmetics, rather than exact arithmetics.
One tacitly assumes that most errors are well-balanced and/or negligibly small.
In particular, the round-off error is often ignored based on the argument that it will be `sufficiently small' if just IEEE-754 double-precision floating-point arithmetics \cite{zuras2008ieee} are adopted.
In this paper, the focus is on the overall discretization error due to truncation and round-off. In particular, we will show that the latter might very well have a significant influence on the overall accuracy and propose a practical strategy to balance both error contributors.

The discretization error strongly depends on the number of degrees of freedom (``DoFs"), denoted by $N_h^{(p)}$, which is a function of the mesh width $h$ and the approximation order $p$. The truncation error, denoted by ${E}_{\rm {T}}$, dominates the discretization error only when $N_h^{(p)}$ is not too large, and it decreases with increasing mesh resolution and element degree as it can be expected from finite element theory \cite{gockenbach2006understanding}. Based on this, the commonly used approaches to reduce the truncation error are to reduce the mesh width ($h$-refinement), increase the approximation order ($p$-refinement), or apply both strategies simultaneously ($hp$-refinement) \cite{guo1986hp}. 
The round-off error, denoted by ${E}_{\rm R}$, is, however, only negligible for moderately small values of $N_h^{(p)}$ and dominates the overall discretization error if more and more DoFs are employed \cite{Babuska2018Roundoff}. 
Consequently, for a particular approximation order $p$, by performing $h$-refinement, the best accuracy is obtained at the break-even point where the discretization error is the smallest. We denote the highest accuracy by ${E}_{\rm {min}}^{(p)}$ and the optimal number of DoFs by $N_{\rm opt}^{(p)}$.

While $N_{\rm opt}^{(p)}$ is typically impractically large if low(est)-order approximations are used, it can be very small if high-order approximations are adopted, which are nowadays becoming more and more popular, and make the results more prone to be polluted by round-off errors.
Despite this alarming observation, to the authors’ best knowledge, only very few publications address the impact of accumulated round-off errors on the overall accuracy of the final solution \cite{ling1984numerical,mou2017example} or take them into account explicitly in the error-estimation procedure {\cite{ainsworth1992procedure,kelly1983posteriori}}.
The general rule of thumb is still to perform as many $h$-refinements as possible considering the available computer hardware.

The aim of this paper is to systematically analyze the influence of the round-off error on the discretization error, for the solution, and its first and second derivative, and propose a practical approach for obtaining ${E}_{\rm {min}}^{(p)}$.
The scope is restricted to one-dimensional model problems, i.e. Poisson, diffusion and Helmholtz equations, for which both the standard finite element method (FEM) and the mixed FEM\cite{boffi2013mixed} are considered.
To assess the general applicability of the aforementioned approach, the following factors are investigated: the element degree over a wide range, first and second derivative of the solution, type of boundary conditions and method of implementing them, choice and configuration of the linear system solver, order of magnitude of the solution and its derivatives, and equation type.

The paper is organized as follows. The model problem, finite element formulation and numerical implementation are described in Section \ref{section_model_problem_FEM_formulation_numerical_implementation}. The general behavior of the discretization error and the approach to predict ${E}_{\rm {min}}^{(p)}$ are discussed in Section \ref{section_behaviour_discretization_error_and_prediction}. Numerical results for determining the offset of the round-off error are shown in Section \ref{section_numerical_results_sensitivity}. The algorithm for realizing the approach is put forward in Section \ref{section_algorithm}, followed by its validation by a Helmholtz problem in Section \ref{section_validation}. The conclusions are drawn in Section \ref{paragraph on conclusion}.

\section{Model problem, finite element formulation and numerical implementation}	\label{section_model_problem_FEM_formulation_numerical_implementation}

\subsection{Model problem}

Consider the following one-dimensional second-order differential equation:
\begin{equation}
  \left(D(x) u_x \right)_x + r(x)u(x) = f(x),\qquad x \in I = (0,1),	\label{1D_general_Helmholtz_equation}
\end{equation}
with $u$ denoting the unknown variable, which can either be real or complex, $f(x) \in L^2 (I)$ a prescribed right-hand side, and $D(x)$ and $r(x)$ continuous coefficient functions.
By choosing $D(x)=1$ and $r(x)=0$, Eq. (\ref{1D_general_Helmholtz_equation}) reduces to the Poisson equation; for $D(x)>0$ and not constant, when $r(x)=0$, the diffusion equation is found, and when $r(x) \neq 0$, we obtain the Helmholtz equation. 
The boundary conditions are $u(x)=g(x)$ on $\Gamma_D$ and $u_x=h(x)$ on $\Gamma_N$. Here, $\Gamma_D$ and $\Gamma_N$ are the boundaries where, respectively, Dirichlet and Neumann boundary conditions are imposed.
In this paper, for all the equations investigated, the existence of the second derivative is guaranteed in the weak sense, i.e. $u \in H^2 (I)$, see \cite{boffi2013mixed}.

\subsection{Finite element formulation} 	\label{FE formulation}

For convenience, we introduce the two inner products:
\begin{subequations}
 \begin{align}
  (f_1(x), \,f_2(x) ) &= \int _I f_1(x) f_2(x) \, dx,	\\
  (g_1(x), \,g_2(x) )_{\Gamma} &= g_1(x_0) g_2(x_0),
 \end{align}
where $f_1(x)$, $f_2(x)$, $g_1(x)$ and $g_2(x)$ are continuous functions defined on the unit interval $I$, $\Gamma$ denotes the boundary of $I$, and $x_0$ denotes the value of $x$ on $\Gamma$.
\end{subequations}

\subsubsection{The standard FEM}

The weak form of Eq. (\ref{1D_general_Helmholtz_equation}) is derived in \ref{derivation_weak_form_SM}. Imposing the Dirichlet boundary conditions strongly, the weak form reads:
\begin{equation}
\centering
\boxed{ 
\begin{aligned}
&\text{Weak~form}~ 1 ~~~~~~~~~\\
&\text{Find $u \in H _D^1 (I)$ such that:} \\
&-({ \eta} _{ x }, \,D { u} _{ x }  ) + (\eta, \, ru) = (\eta, \, f ) - (\eta, \, D h {n} )_{\Gamma _N} \qquad \forall \eta \in H _{D0}^1 (I),\\
&\text{with} \\
&~~~~~~~~~~~~~~~~~~~~~~~~~H_{D} ^1 (I) = \{t \; | \; t \in H^1 (I), \; t = g \text{ on } \Gamma _D \},  \\
&~~~~~~~~~~~~~~~~~~~~~~~\,H_{D0} ^1 (I) = \{t \; | \; t \in H^1 (I), \; t = 0 \text{ on } \Gamma _D\},\\
&\text{where } {n} \text{ is 1 at $x=1$, and $-1$ at $x=0$.}
\end{aligned}		\label{1D_general_SM_weak_form_Diri_strong} 
}
\end{equation}
By imposing the Dirichlet boundary conditions in the weak sense, the weak form reads:
\begin{equation}
\centering
\boxed{
\begin{aligned}
&\text{Weak~form}~ 2 ~~~~~~~~~\\
&\text{Find } u \in H ^1 (I) \text{ such that:}\\
& - ( { \eta} _{ x }, \, D { u} _{ x } ) + (\eta, \, ru) + (\eta, \, D u_x n )_{\Gamma _D} - (\eta _x, \, u n )_{\Gamma _D} + (\eta, \, \rho u n )_{\Gamma _D} \\ 
&= (\eta, \, f ) - (\eta, \, D h n )_{\Gamma _N} - (\eta _x, \, g n )_{\Gamma _D} + (\eta, \, \rho g n )_{\Gamma _D} \qquad \forall \eta \in H ^1 (I), \\
&\text{where } \rho \text{ is a } \text{{positive}} \text{ value that serves as the penalty parameter}.
\end{aligned}	\label{1D_general_SM_weak_form_Diri_weak}
}
\end{equation}

\noindent Note that, the terms in the right-hand sides of Eqs. (\ref{1D_general_SM_weak_form_Diri_strong})--(\ref{1D_general_SM_weak_form_Diri_weak}) consist of information of the Neumann boundary conditions, and hence, if no Neumann boundary conditions are prescribed, these terms vanish. We use Weak form 1 if not stated otherwise.
Next, we approximate the exact solution $u_{\rm exc}$ by a linear combination of a finite number of basis functions:

\begin{equation}
 u_{\rm exc} \approx u_h^{(p)} = \sum _ {i=1} ^{m} u _{i} \varphi _{i}^{(p)}. \label{General_SM_u_approx}%
\end{equation}
Here, $\varphi _{i}^{(p)}$ are $C^0$-continuous Lagrange basis functions of degree $p$, denoted as $P_p$, with Gauss-Lobatto support points $x_j$, which feature the Kronecker-delta property, i.e. $\varphi _{i}^{(p)} (x_j)=\delta_{ij}$. The coefficients $u_i$ are the values of $u_h^{(p)}$ at the $\text{DoFs}$, as a direct consequence of the Kronecker-delta property of $\varphi _{i}^{(p)}$. The number of DoFs of $u_h^{(p)}$, denoted by $m$, equals $p \times t + 1$, where $t$ is the total number of the grid cells. 
Finally, taking the test function $\eta$ equal to $\varphi ^{(p)}_{k},~ k=1, \,2, \, \ldots , \, m$, we obtain
\begin{equation}
 A {U} = F,				\label{matrix equation std FEM}
\end{equation}
where $A$ is the stiffness matrix, $F$ the right-hand side and $U$ the discrete solution, i.e. the vector of the coefficients $u_i$.

\subsubsection{The mixed FEM}

As a first step, we introduce the auxiliary variable
\begin{subequations}
\begin{align}
   v(x) = - u _x, \label{Gene_MM_strong1} 
\end{align}  
allowing Eq. (\ref{1D_general_Helmholtz_equation}) to be rewritten as
\begin{align}
  -D_x v(x) - D(x) v_x + r(x)u(x) = f(x). \label{Gene_MM_strong2}
\end{align}	\label{1D_general_MM_2in1}%
\end{subequations}
Unlike the standard FEM, for the mixed FEM, the essential boundary conditions are imposed on $\Gamma _N$, and the natural boundary conditions on $\Gamma _D$.
The weak form of Eq. (\ref{1D_general_Helmholtz_equation}) using the mixed FEM, derived in \ref{derivation_weak_form_MM}, is given by:
\begin{subequations}
\begin{empheq}[box=\fbox]{align}
&\text{Weak~form}~ 3 ~~~~~~~~~\notag\\
&\text{Find $v \in H_{N}^1 (I)$ and $u \in L ^2 (I)$ such that:}	\notag\\
& ~~~~~~~~~~~~~~~~~~~~\,\,(w, \, v) - (w_x, \,  u ) = -(w, \, g n)_{\Gamma_D} \qquad \forall w \in H_{N0}^1 (I), \label{1D_General_MM_weak_1}\\ 
&  -(q, \, D_x v ) - (q, \, D v_x) + (q, \, ru) = (q, \, f) \qquad \forall q \in L ^2 (I), \label{1D_General_MM_weak_2}	\\
&    \text{with}\notag\\
& ~~~~~~~~~~~~~~~~~~~~~~~~~~~~~~ H_{N} ^1 (I) = \{t \; | \; t \in H^1 (I), \; t = -h \text{ on } \Gamma _N \},  \notag\\
& ~~~~~~~~~~~~~~~~~~~~~~~~~~~~\, H_{N0} ^1 (I) = \{t \; | \; t \in H^1 (I), \; t = 0 \text{ on } \Gamma _N\}.	\notag 
\end{empheq}
\label{1D_General_MM_weak_2in1}%
\end{subequations}
Next, we approximate the exact gradient $v_{\rm exc}$ and the exact solution $u_{\rm exc}$ by a linear combination of a finite number of basis functions:
\begin{subequations}
 \begin{align}
 v_{\rm exc} \approx v _h^{(p)} &= \sum _ {i=1} ^{n} v _{i} \varphi _{i}^{(p)},     \label{General_MM_var_approx1}  \\[3ex]
 u_{\rm exc} \approx u _h^{(p-1)} &= \sum\limits _ {j=1} ^{p} u _{cj} \psi _{j} ^{(p-1)} \text{ in cell }c, \text{ for } c=1,\,2, \, \ldots, \,t.  \label{General_MM_var_approx2}
\end{align}	\label{General_MM_var_approx}%
\end{subequations}
where $\varphi _{i}^{(p)}$ are of the same type of basis functions used in Eq. (\ref{General_SM_u_approx}), with coefficients $v_i$ the associated values of $v_h^{(p)}$ at the $\text{DoFs}$;
$\psi _{j} ^{(p-1)}$ are discontinuous Lagrange basis functions of degree $p-1$, denoted as $P_{p-1}^{\text{disc}}$, with coefficients $u_{c,j}$ the associated values of $u_h^{(p-1)}$ at the $\text{DoFs}$. 
This pair of elements will be referred to as $P_p/P_{p-1}^{\text{disc}}$.
Since the use of discontinuous basis functions, there are two independent $u_{c,j}$ at cell interfaces.
The number of DoFs for $v_h^{(p)}$, denoted by $n$, equals $p \times t + 1$, and the number of DoFs for $u_h^{(p-1)}$ equals $p \times t$. 
Finally, replacing the test functions $w$ and $q$ by $\varphi _{k}^{(p)} , ~{k} = 1, \,2, \, \ldots , \, p \times t + 1, $ and $ \psi _{e}^{(p-1)} ,~ {e} = 1, \,2, \, \ldots , \, p \times t$, respectively, the resulting coupled linear system of equations that has to be solved reads:
\begin{equation}
 \left[ \begin{array}{cc} M & B  \\ B^\top & 0 \end{array}\right] \left[ \begin{array}{cc} {V} \\ {U} \end{array}\right] =\left[ \begin{array}{cc} G \\ H \end{array}\right], \label{matrix equation mix FEM}
\end{equation}
where the mass matrix $M$, the discrete gradient operator $B$, and its transpose, the discrete divergence operator $B^\top$, are the components of the discrete left-hand side of Eqs. (\ref{1D_General_MM_weak_1})--(\ref{1D_General_MM_weak_2}), $G$ and $H$ are the components of the right-hand side, and $V$ and $U$ are the discrete first derivative and solution, i.e. the vectors of the coefficients $v_i$ and $u_{cj}$, respectively.

For the sake of readability, we will drop the superscript $(p)$, whenever the approximation order is clear from the context.

\subsection{Numerical implementation}

In what follows, we demonstrate how to obtain the numerical solution for Eq. (\ref{1D_general_Helmholtz_equation}) with specific coefficients and assess its quality. For the latter, both the error, obtained using the analytical solution or the finer numerical solution, and the order of convergence are investigated.

\subsubsection{Solution technique}

Unless stated otherwise, all results are computed in IEEE-754 double precision \cite{zuras2008ieee} using the deal.\rom{2} finite element code \cite{alzetta2018deal} that provides subroutines for creating the computational grid, building and solving the system of equations, and computing the error norms.

The computational mesh is obtained by globally refining a single element that covers the interval $I$, and the Dirichlet boundary conditions are imposed strongly unless stated otherwise.
The former means that, when the solution is real valued, using the standard FEM, the number of DoFs equals $2^{REF} \times p+1$ at the $REF$th refinement;
using the mixed FEM, the number of DoFs equals $2 \times 2^{REF} \times p+1$ at the $REF$th refinement.
For complex-valued problems, the above numbers double since deal.\rom{2} does not provide native support for complex-valued problems and, hence, all components need to be split into their real and imaginary parts.

To compute the occurring integrals, sufficiently accurate Gaussian quadrature formulas are used. 
Furthermore, unless stated otherwise, to solve the matrix equation, the UMFPACK solver \cite{davis2004algorithm}, which implements the multi-frontal LU factorization approach, is used as it results in relatively fast computations of the problems considered in this paper, and prevents the iteration errors for the iterative solvers. 

\subsubsection{Error estimation}

For the numerical results $var_h$, where $var$ can be $u$, $u_x$ and $u_{xx}$, the discretization error measured in the $L_2$ norm is used. This measure contains all types of errors, for example the truncation error, round-off error, etc. It is defined as
\begin{subequations}	\label{formula_abs_error}
\begin{align}		\label{formula_abs_error_analytical}
 E_{h} &= {\|var_{h}- {var}_{\rm exc}\|_{2}}
\end{align}
when the exact solution ${var}_{\rm exc}$ is available, or \cite{Runborg2012VerifyingNC}
\begin{align}		\label{formula_abs_error_numerical}
 {\widetilde {E_{h}}} &= {\|var_{h}- {var}_{h/2}\|_{2}}
\end{align}
otherwise,
\end{subequations}
where $var_{h/2}$ is the numerical solution computed on a mesh with grid size $h/2$. 
The derivatives, which are $u_{h,x}$ and $u_{h,xx}$ in the standard FEM and only $u_{h,xx}$ in the mixed FEM, are computed in the classical finite element manner, e.g. $u_{h,x} ^{(p-1)}=\sum\limits _{i=1}^m u_i\varphi_{i,x}^{(p)}$ yields an approximation to $u_x$ using standard FEM. Note that, each differentiation decreases the element degree by one.  

\subsubsection{Convergence of the solution}

When the number of DoFs is relatively large, but the round-off error does not exceed the truncation error, the discretization error converges at a fixed rate $\beta_{\rm T}$ theoretically\cite[Theorem~5.\rom{1}]{gockenbach2006understanding}. The value of $\beta_{\rm T}$ can be found in Table~\ref{convergence_order_sample_equations} for the element degree $p$ ranging from 1 to 5. In practice, it can be calculated from either 
\begin{subequations}	\label{formula_order_of_convergence}
\begin{align}
 \beta_{\rm T}=\log _2 \left( \frac{E_{h}}{E_{h/2}} \right)
\end{align}
using Eq. (\ref{formula_abs_error_analytical}), or
\begin{align}
 \beta_{\rm T}=\log _2 \left( \frac{\widetilde {E_{h}}}{\widetilde {E_{h/2}}} \right)		\label{formula_order_of_convergence_finer_solution}
\end{align}
using Eq. (\ref{formula_abs_error_numerical}).
\end{subequations}

\begin{table}[!ht]
\caption[sss]{Order of convergence for $u$, $u_{x}$ and $u_{xx}$.}
\label{convergence_order_sample_equations}
\hspace{0.7cm}
\begin{subtable}{0.5\textwidth}
\caption[sss]{The standard FEM}
\label{convergence_order_sample_equations_std}
\centering
 \begin{tabular}{c c c c} \hline
Elements & $u$ & $u_{x}$ & $u_{xx}$  \\	\hline
$P_1$ & 2 & 1 & n/a \\
$P_2$ & 3 & 2 & 1 \\ 
$P_3$ & 4 & 3 & 2 \\ 
$P_4$ & 5 & 4 & 3 \\ 
$P_5$ & 6 & 5 & 4 \\ \hline
\end{tabular}
\end{subtable}\hspace{-2.0cm}
\begin{subtable}{0.5\textwidth}
\centering
\caption[sss]{The mixed FEM}
\label{convergence_order_sample_equations_mix}
 \begin{tabular}{c c c c} \hline 
Elements & $u$ & $u_{x}$ & $u_{xx}$  \\	\hline
$P_1/P_0^{\rm disc}$ & 1 & 2 & 1 \\
$P_2/P_1^{\rm disc}$ & 2 & 3 & 2 \\
$P_3/P_2^{\rm disc}$ & 3 & 4 & 3 \\
$P_4/P_3^{\rm disc}$ & 4 & 5 & 4 \\
$P_5/P_4^{\rm disc}$ & 5 & 6 & 5 \\	\hline
\end{tabular}
\end{subtable}
\end{table}

\newpage

\section{General behaviour of the discretization error and approach to predict the highest attainable accuracy}      \label{section_behaviour_discretization_error_and_prediction}

In this section, based on \cite{Babuska2018Roundoff,WalterFrei}, we illustrate the general behaviour of the discretization error $E_h$  for Eq. (\ref{1D_general_Helmholtz_equation}) as a function of the number of DoFs $N_h$, and provide an approach to predict the highest attainable accuracy.

\begin{figure}[!ht]
 \centering
     \includegraphics[width=0.5\linewidth]{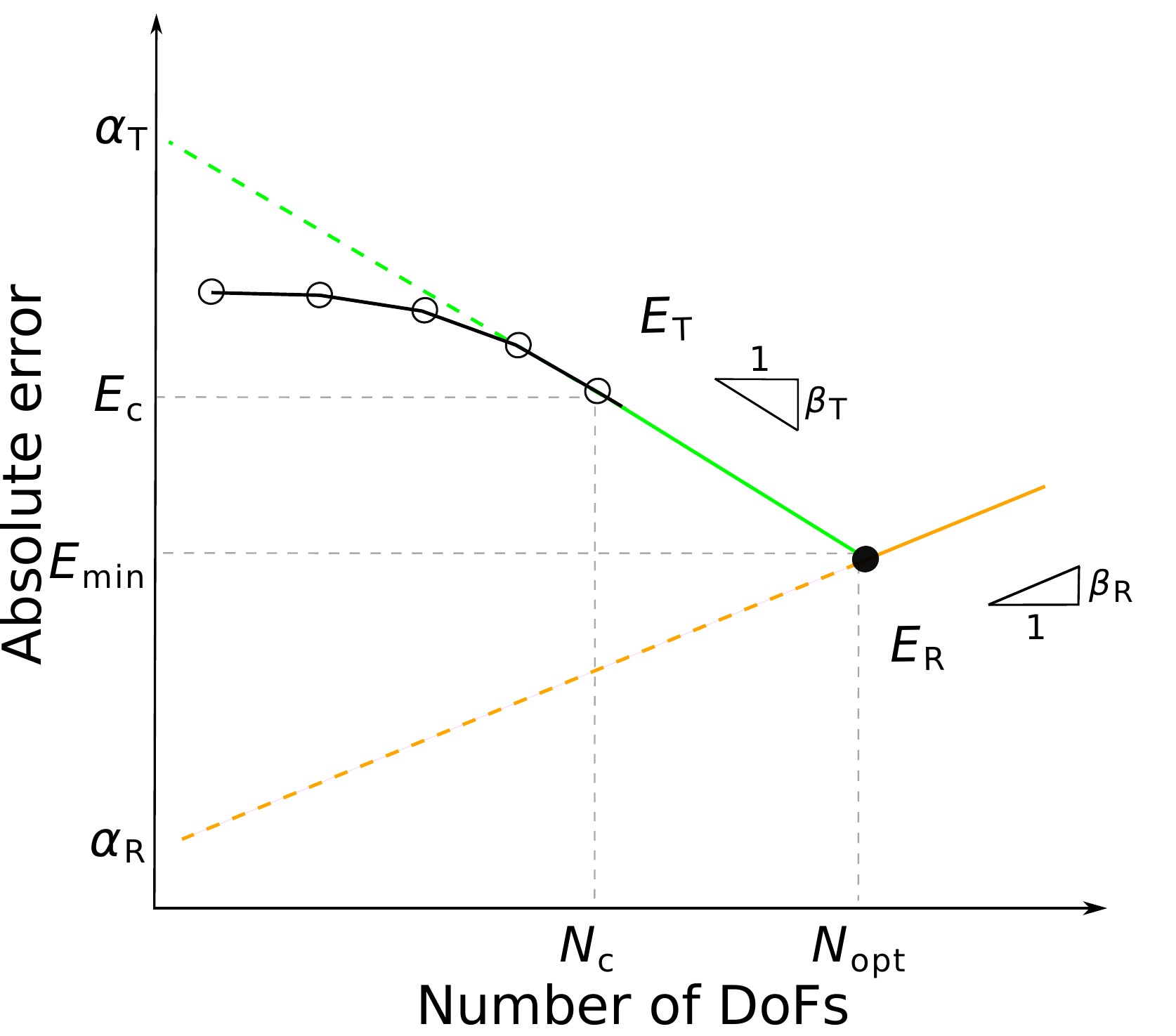}
     \caption{Conceptual sketch of the dependency of the discretization error on the number of $\text{DoFs}$.}
     \label{sketch_discretization_error_one_p}
 \end{figure}

The discretization error of one variable for one $p$ is illustrated in Fig.~\ref{sketch_discretization_error_one_p}, where log-log axes are used. 

As can be seen, the change of $E_h$ with $N_h$ can be divided into three phases according to $N_{\rm c}$ and $N_{\rm opt}$, for which the former is $N_h$ where $E_h$ begins showing the expected asymptotic convergence behavior, and the latter is $N_h$ where $E_h$ begins increasing, i.e. where the highest attainable accuracy is obtained. The features of $E_h$ in each phase are shown in Table~\ref{phase_discretization_error}.

\begin{table}[!ht]
\caption [sss] {Features of $E_h$ in different phases.}
\label{phase_discretization_error}
\centering
 \begin{tabular}{|l|C{3.2cm}|C{3.2cm}|C{4.0cm}|} \hline   
 & 1. {$N_h < N_{\rm c}$} & 2. {$N_{\rm c} \leqslant N_h < N_{\rm opt}$} & 3. {$N_{\rm opt} \leqslant N_h$} \\ \hline
Description & Decreasing but not converging at slope $\beta_{\rm T}$ & Decreasing and converging at slope $\beta_{\rm T}$, with the offset $\alpha_{\rm T}$ & Increasing and converging at slope $\beta_{\rm R}$, with the offset $\alpha_{\rm R}$ \\	\hline
Formula & - & $E_h=\alpha_{\rm T}{N_h}^{-\beta_{\rm T}}$ & $E_h=\alpha_{\rm R} {N_h}^{\beta _{\rm R}}$ \\	\hline
Dominant error & \multicolumn{2}{c|}{Truncation error} & Round-off error \\	\hline
\end{tabular}
\end{table}

As we will prove in section \ref{section_numerical_results_sensitivity}, the values of $\alpha_{\rm R}$ and $\beta_{\rm R}$ can be relatively fixed.
Therefore, the round-off error can be assessed before solving the problem.
Moreover, since $\alpha_{\rm T}$ can be inversed by using
\begin{equation}
 \alpha_{\rm T} = {E_{\rm c}}/{N_{\rm c}}^{- \beta_{\rm T}},		\label{formula_offset_truncation_error}
\end{equation}
at the beginning of phase 2, where ${E_{\rm c}}$ is the value of ${E_h}$ corresponding to ${N_{\rm c}}$, we can forecast $E_{\rm T}$ afterwards.


Obviously, $N_{\rm opt}$ happens when $E_{\rm T}+E_{\rm R}$ is the smallest. By solving
\begin{equation}
    \frac{d(E_{\rm T}+E_{\rm R})}{dN}=0,    \label{derivative_condition_N_opt}
\end{equation}
we can predict
\begin{subequations}
\begin{align}
 N_{\rm opt} = \left( \frac{\alpha_{\rm T} \beta_{\rm T}}{\alpha _{\rm R} \beta_{\rm R}} \right)^{\frac{1}{\beta_{\rm T} + \beta_{\rm R}}},
\end{align}
and hence, the highest attainable accuracy
\begin{align}
 E_{\rm min} = \alpha_{\rm T} {N_{\rm opt}}^{- {\beta _{\rm T}}}+\alpha_{\rm R} {N_{\rm opt}}^{{\beta _{\rm R}}}.
\end{align}
\end{subequations}

\section{Numerical quantification of the round-off error}  	\label{section_numerical_results_sensitivity}

In this section, we assess the general values for $\alpha_{\rm R}$ and $\beta_{\rm R}$ for variables $u$, $u_x$ and $u_{xx}$, using both the standard FEM and the mixed FEM.
We start with the preliminary results obtained from three benchmark equations, and then investigate the following factors: solution strategy, boundary condition and order of magnitude. 

\begin{table}[!ht]
\caption [sss] {Settings of the benchmark Poisson, diffusion and Helmholtz equations.}		
\label{benchmark one-dimensional equations} 
\centering
 \begin{tabular}{|C{2.5cm}|C{4cm}|C{4cm}|C{4cm}|} \hline   
{} & {``Poisson''} & {``diffusion''} & {``Helmholtz''} \\ \hline
{$D(x)$} & {$1$} & $1+x$ & $(1+i) e^{-x}$  \\	\hline
{$r(x)$} & {0} & 0 & $2 e^{-x}$ \\	\hline
{$f(x)$} & {$e^{- (x-1/2)^2} \left({4x^2 - 4x -1} \right)$}  & $2 \pi \cos (2 \pi x) - 4 {\pi}^2 \sin (2 \pi x)(x+1)$ & 0 \\ \hline
{$\|f(x)\|_2$} & {1.60} & {42.99} & {0.00} \\	\hline
\multirow{2}{2cm}{\centering Boundary conditions} & {$u(0) = e^{-1/4}$} & $u(0)=0$& $u (0) = 1$ \\	\cline{2-4}
&$u(1) = e^{-1/4}$ & $u_x(1)=2 \pi$  &$ u_x(1) = 0$ \\	\hline
Analytical solution $u_{\text{exc}}$ & {$e^{- (x-1/2)^2}$} & $\sin (2 \pi x)$ & $a e^{(1+i) x} + (1-a) e^{-i x}$, $a=1/{((1-i) e^{1+2i}+1)}$ \\	\hline
{$\|u_{\text{exc}}\|_2$} & {0.92} & 0.71 & 1.26 \\	\hline
\end{tabular}
\end{table}


\subsection{Preliminary results}		\label{section_preliminary_results}

We consider the benchmark equations given in Table~\ref{benchmark one-dimensional equations}, for which the $L_2$ norm of the analytical solution $u_{\rm exc}$ is of order 1. Element degrees $p$ range from 1 to 5.

\subsubsection{Benchmark Poisson equation}

\begin{figure}[!ht]
    \begin{subfigure}{5.5cm}
        \includegraphics[width=1.0\linewidth]{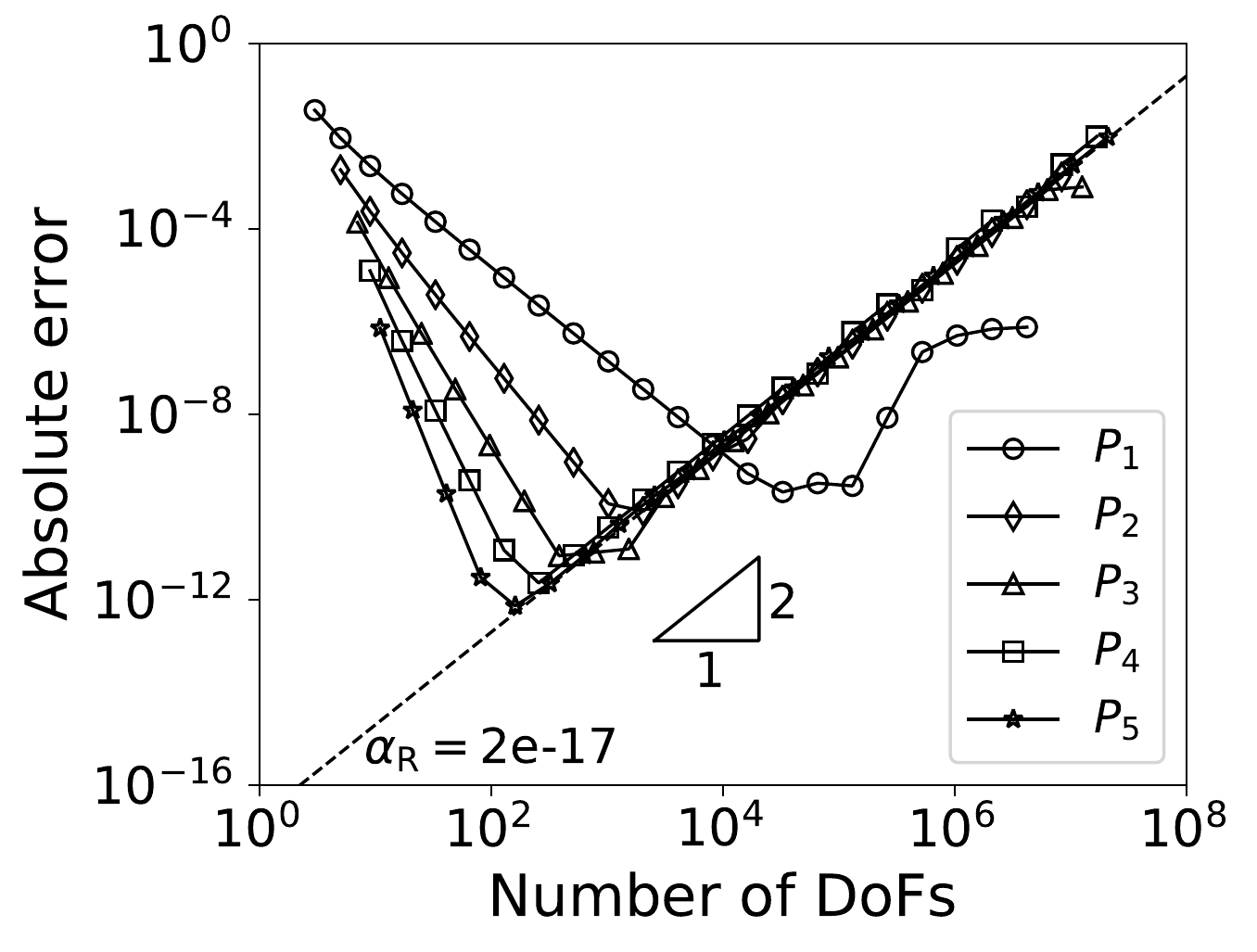}			
        \caption{Solution}
        \label{py_bench_Pois_SM_solu}
    \end{subfigure}
    \hspace{-0.2cm}
    \begin{subfigure}{5.5cm}
        \includegraphics[width=1.0\linewidth]{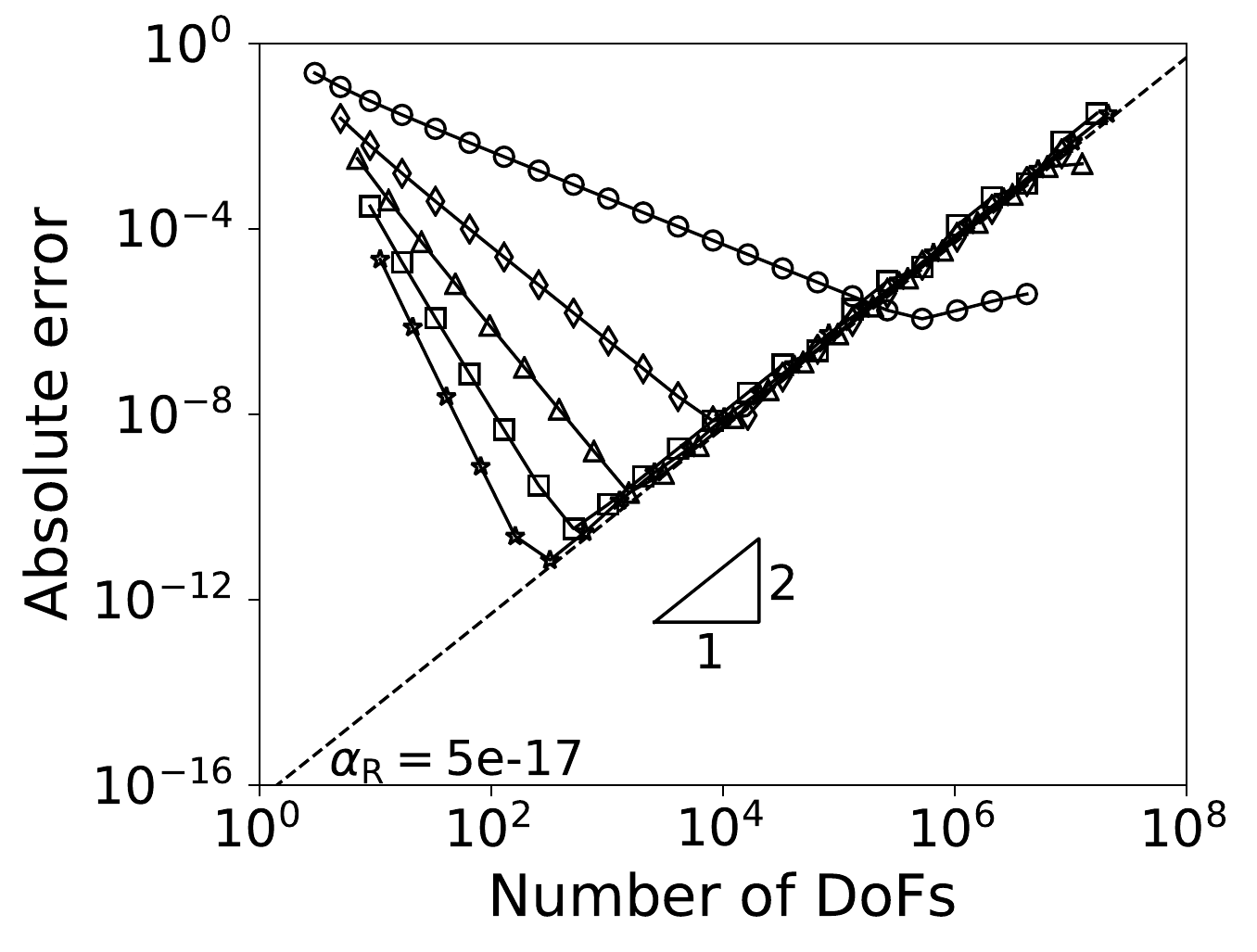}
        \caption{First derivative}
        \label{py_bench_Pois_SM_grad}
    \end{subfigure}
    \hspace{-0.2cm}
    \begin{subfigure}{5.5cm}
        \includegraphics[width=1.0\linewidth]{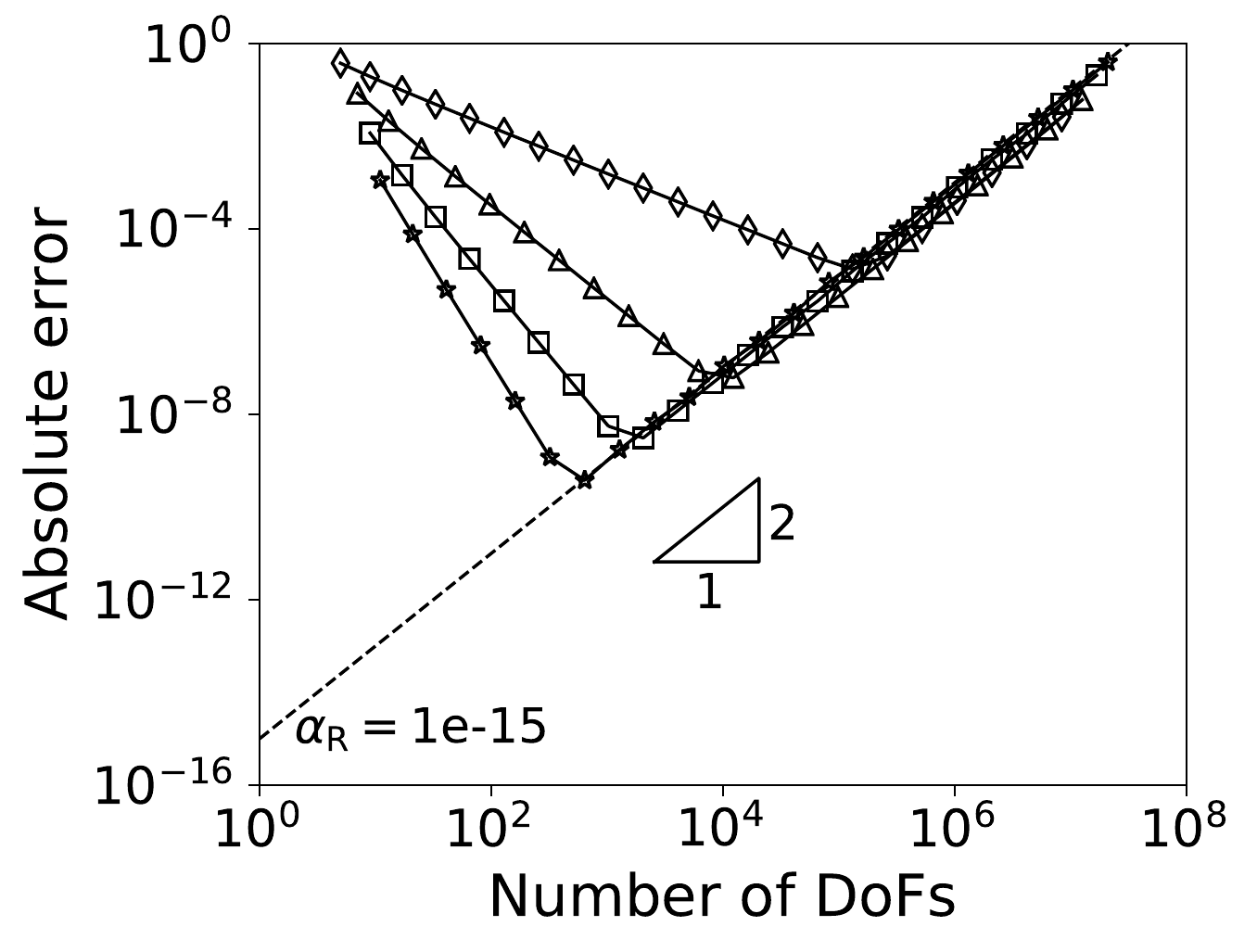}
        \caption{Second derivative}
        \label{py_bench_Pois_SM_2ndd}
    \end{subfigure}
\caption{Absolute errors for the benchmark Poisson equation using the standard FEM.}
\label{py_bench_Pois_SM}
\end{figure}

\begin{figure}[!ht]
    \begin{subfigure}{5.5cm}
        \includegraphics[width=1.0\linewidth]{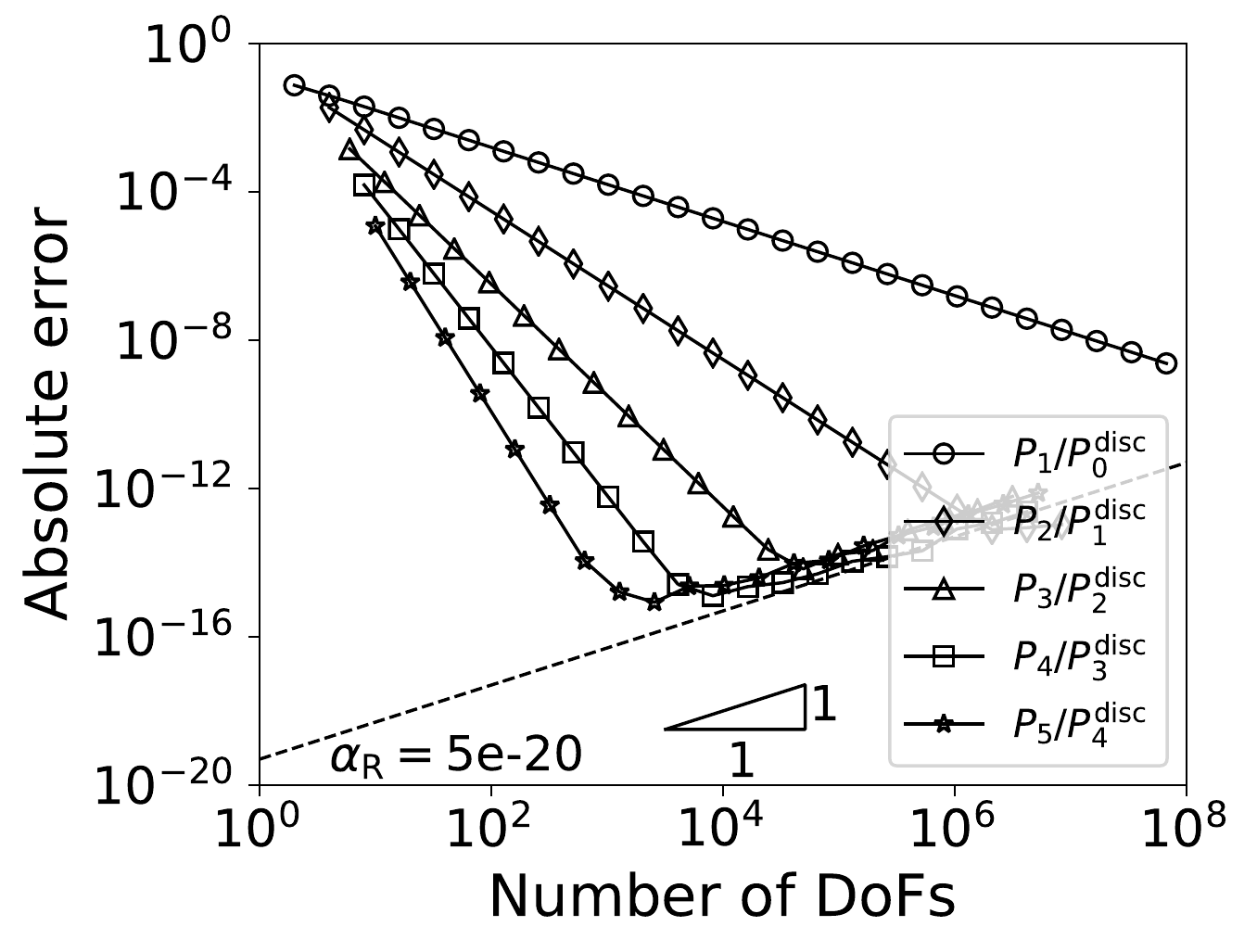}
        \caption{Solution}
        \label{py_bench_Pois_MM_solu}
    \end{subfigure}
    \hspace{-0.2cm}
    \begin{subfigure}{5.5cm}
        \includegraphics[width=1.0\linewidth]{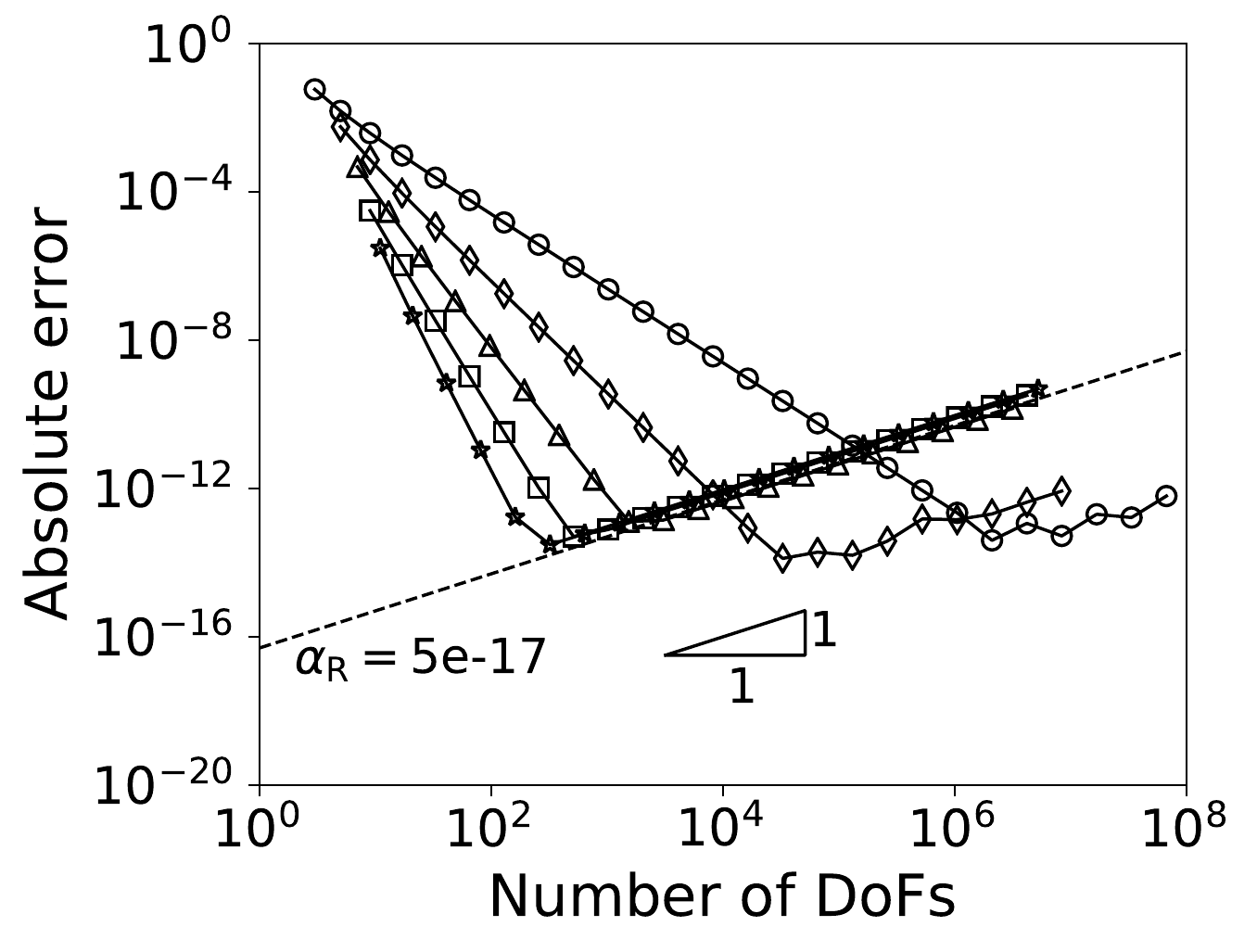}
        \caption{First derivative}
        \label{py_bench_Pois_MM_grad}
    \end{subfigure}
    \hspace{-0.2cm}
    \begin{subfigure}{5.5cm}
        \includegraphics[width=1.0\linewidth]{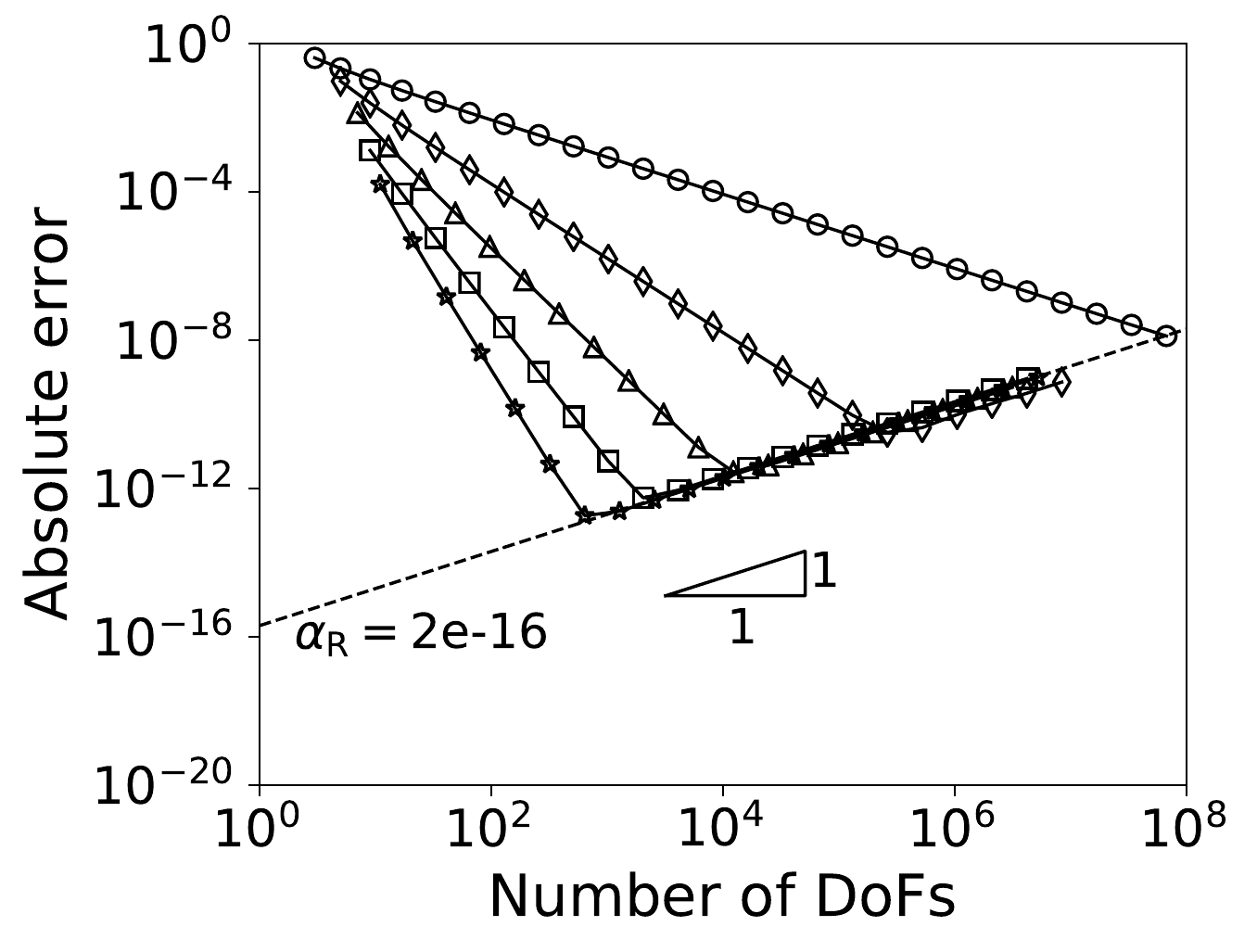}
        \caption{Second derivative}
        \label{py_bench_Pois_MM_2ndd}
    \end{subfigure}
\caption{Absolute errors for the benchmark Poisson equation using the mixed FEM.}
\label{py_bench_Pois_MM}
\end{figure}

For the benchmark Poisson equation, the discretization error $E_h$ for $u$, $u_x$ and $u_{xx}$ using both the standard FEM and the mixed FEM can be found in Fig.~\ref{py_bench_Pois_SM} and Fig.~\ref{py_bench_Pois_MM}, respectively. The offset $\alpha_{\rm R}$ and slope $\beta_{\rm R}$ are denoted in the figures, so are in the following figures of the same type.

Using both the standard FEM and the mixed FEM, for all the variables, the interesting point is that the values of $\alpha_{\rm R}$ and $\beta_{\rm R}$ for different element degrees tend to be the same. Notably, the value of the former is of order $10^{-16}$, which is as expected when using double precision.

\begin{figure}[!ht]
\hspace{0.0cm}
\begin{subfigure}[b]{0.35\textwidth}
\scalebox{0.85}{
\begin{tikzpicture}
\begin{axis}
[
    ymode=log,    
    ymin=1e-16,
    ymax=1e0,
    ytick={1e-20, 1e-16, 1e-12, 1e-8, 1e-4, 1e0},
    yticklabels={1e-20, 1e-16, 1e-12, 1e-8, 1e-4, 1e0},      
    legend style={nodes={scale=0.85}},
    legend cell align={left},
    height=5cm,
    width=6cm,
    ylabel={$E_{\rm min}$},
    ylabel style={at={(-0.01,0.5)}},    
    xtick={0,1,2,3,4},
    xticklabels={$1$,$2$, $3$, $4$, ${5}$},
    xlabel={$p$},
    xlabel style={at={(0.5,0.03)}},    
]
\addplot[black,mark=o,mark options={color=black,fill=black}] coordinates {(0,2.1e-25)};
\addplot[black,mark=diamond,mark options={color=black,fill=black}] coordinates {(0,2.1e-25)};
\addplot[black,mark=triangle,mark options={color=black,fill=black}] coordinates {(0,2.1e-25)};
\addplot[blue,mark=o,mark options={color=blue,fill=blue}] coordinates {(0,2.1e-10) (1,8.1e-11) (2,8.9e-12) (3,2.3e-12) (4,7.4e-13)};
\addplot[blue,mark=diamond,mark options={color=blue,fill=blue}] coordinates {(0,1.1e-6) (1,7.2e-9) (2,2.0e-10) (3,3.4e-11) (4,7.2e-12)};
\addplot[blue,mark=triangle,mark options={color=blue,fill=blue}] coordinates {(0) (1,1.4e-5) (2,6.2e-8) (3,3.1e-9) (4,3.8e-10)};

\addplot[red,mark=o,mark options={color=red,fill=red}] coordinates {(0,1.1e-9)}; 
\addplot[red,mark=o,mark options={color=red,fill=red}] coordinates {(1,7.7e-14) (2,9.3e-15) (3,2.6e-15) (4,7.4e-16)}; 
\draw [red,dashed] (axis cs:0,1.1e-9) -- (axis cs:1,7.7e-14);

\addplot[red,mark=diamond,mark options={color=red,fill=red}] coordinates {(0,4.0e-14) (1,1.3e-14) (2,1.3e-13) (3,5.1e-14) (4,3.1e-14)};

\addplot [red,mark=triangle,mark options={color=red,fill=red}] coordinates {(0,1.3e-8)};
\addplot[red,mark=triangle,mark options={color=red,fill=red}] coordinates {(1,3.5e-11) (2,2.4e-12) (3,5.1e-13) (4,2.4e-13)}; 
\draw [red,dashed] (axis cs:0,1.3e-8) -- (axis cs:1,3.5e-11);
\legend{$u$,$u_x$,$u_{xx}$};
\end{axis}
\end{tikzpicture}
}
\vspace{-0.2cm}
\caption{Poisson}
\label{E_min_benchmark_Poisson}
\end{subfigure}
\hspace{-0.7cm}
\begin{subfigure}[b]{0.35\textwidth}
\scalebox{0.85}{
\begin{tikzpicture} 
\begin{axis}
[
    ymode=log,    
    ymin=1e-16,
    ymax=1e0,
    ytick={1e-20, 1e-16, 1e-12, 1e-8, 1e-4, 1e0},
    yticklabels={1e-20, 1e-16, 1e-12, 1e-8, 1e-4, 1e0},      
    legend style={nodes={scale=0.85}},
    legend cell align={left},
    height=5cm,
    width=6cm,
    ylabel={$E_{\rm min}$},
    ylabel style={at={(-0.01,0.5)}},    
    xtick={0,1,2,3,4},
    xticklabels={$1$,$2$, $3$, $4$, ${5}$},
    xlabel={$p$},
    xlabel style={at={(0.5,0.03)}},   
]
\addplot[black,mark=o,mark options={color=black,fill=black}] coordinates {(0,2.1e-25)};
\addplot[black,mark=diamond,mark options={color=black,fill=black}] coordinates {(0,2.1e-25)};
\addplot[black,mark=triangle,mark options={color=black,fill=black}] coordinates {(0,2.1e-25)};

\addplot[blue,mark=o,mark options={color=blue,fill=blue}] coordinates {(0,1.3e-9) (1,1.0e-10) (2,5.5e-11) (3,3.8e-12) (4,2.9e-12)};
\addplot[blue,mark=diamond,mark options={color=blue,fill=blue}] coordinates {(0,7.7e-6) (1, 3.0e-8) (2, 1.6e-9) (3, 1.4e-10) (4, 2.1e-11)};
\addplot[blue,mark=triangle,mark options={color=blue,fill=blue}] coordinates {(0) (1,2.0e-4) (2,6.3e-7) (3, 1.6e-8) (4, 1.1e-9)};

\addplot[red,mark=o,mark options={color=red,fill=red}] coordinates {(0,8e-10)}; 
\addplot[red,mark=o,mark options={color=red,fill=red}] coordinates {(1,8.2e-14) (2,4.3e-14) (3,2.7e-15) (4,6.8e-16)}; 
\draw [red,dashed] (axis cs:0,8e-10) -- (axis cs:1,8.2e-14);

\addplot[red,mark=diamond,mark options={color=red,fill=red}] coordinates {(0, 3.2e-13) (1, 4.6e-14) (2, 2.5e-13) (3, 3.3e-14) (4,1.6e-14)};

\addplot [red,mark=triangle,mark options={color=red,fill=red}] coordinates {(0,8e-7)};
\addplot[red,mark=triangle,mark options={color=red,fill=red}] coordinates {(1, 1.3e-9) (2, 1.1e-10) (3, 2.0e-11) (4, 9.0e-12)}; 
\draw [red,dashed] (axis cs:0,8e-7) -- (axis cs:1,1.3e-9);
\end{axis}
\end{tikzpicture}
}
\vspace{-0.2cm}
\caption{Diffusion}
\label{E_min_benchmark_diffusion}
\end{subfigure}
\hspace{-0.7cm}
\begin{subfigure}[b]{0.35\textwidth}
\scalebox{0.85}{
\begin{tikzpicture} 
\begin{axis}
[
    ymode=log,    
    ymin=1e-16,
    ymax=1e0,
    ytick={1e-20, 1e-16, 1e-12, 1e-8, 1e-4, 1e0},
    yticklabels={1e-20, 1e-16, 1e-12, 1e-8, 1e-4, 1e0},      
    legend style={nodes={scale=0.85}},
    legend cell align={left},
    height=5cm,
    width=6cm,
    ylabel={$E_{\rm min}$},
    ylabel style={at={(-0.01,0.5)}},    
    xtick={0,1,2,3,4},
    xticklabels={$1$,$2$, $3$, $4$, ${5}$},
    xlabel={$p$},
    xlabel style={at={(0.5,0.03)}},  
]
\addplot[black,mark=o,mark options={color=black,fill=black}] coordinates {(0,2.1e-25)};
\addplot[black,mark=diamond,mark options={color=black,fill=black}] coordinates {(0,2.1e-25)};
\addplot[black,mark=triangle,mark options={color=black,fill=black}] coordinates {(0,2.1e-25)};

\addplot[blue,mark=o,mark options={color=blue,fill=blue}] coordinates {(0, 4.5e-9) (1, 1.2e-10) (2, 2.6e-11) (3, 2.5e-12) (4,6.2e-13)};
\addplot[blue,mark=diamond,mark options={color=blue,fill=blue}] coordinates {(0,2.8e-6) (1, 5.8e-9) (2, 3.1e-10) (3, 7.5e-12) (4, 4.8e-12)};
\addplot[blue,mark=triangle,mark options={color=blue,fill=blue}] coordinates {(0) (1,1.1e-5) (2, 4.0e-8) (3, 1.8e-9) (4, 1.9e-10)};

\addplot[red,mark=o,mark options={color=red,fill=red}] coordinates {(0,5e-9)}; 
\addplot[red,mark=o,mark options={color=red,fill=red}] coordinates {(1, 7.3e-12) (2, 3.3e-13) (3, 2.7e-14) (4, 2.3e-14)}; 
\draw [red,dashed] (axis cs:0,5e-9) -- (axis cs:1,7.3e-12);

\addplot[red,mark=diamond,mark options={color=red,fill=red}] coordinates {(0, 7.2e-14) (1, 1.6e-13) (2, 3.5e-14) (3, 3.1e-14) (4,2.9e-15)};

\addplot [red,mark=triangle,mark options={color=red,fill=red}] coordinates {(0,2e-8)};
\addplot[red,mark=triangle,mark options={color=red,fill=red}] coordinates {(1, 3.3e-11) (2, 2.5e-12) (3, 5.6e-13) (4, 1.6e-13)}; 
\draw [red,dashed] (axis cs:0,2e-8) -- (axis cs:1,3.3e-11);
\end{axis}
\end{tikzpicture}
}
\vspace{-0.2cm}
\caption{Helmholtz}
\label{E_min_benchmark_Helmholtz}
\end{subfigure}
\caption{$E_{\rm min}$ for the benchmark equations. The blue color denotes the standard FEM, and the red color denotes the mixed FEM.}
\label{E_min_benchmark}
\end{figure}
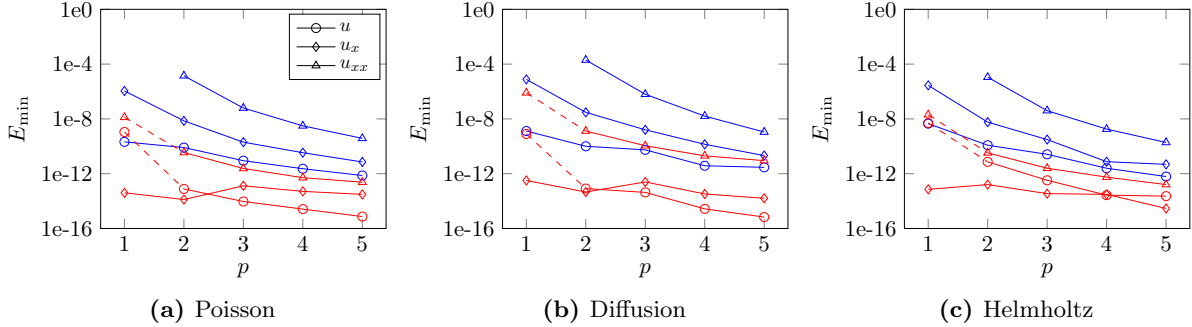

For the error of one particular variable using one particular FEM, since $E_{\rm T}$ decreases faster for larger $p$, smaller ${E}_{\text{min}}$ can be obtained using larger $p$. Since the slope $\beta_{\rm R}$ using the mixed FEM is half of that using the standard FEM\cite{kahan2013floating}, the mixed FEM gives smaller ${E}_{\text{min}}$ for each variable using the same $p$, see Fig.~\ref{E_min_benchmark_Poisson} for the statistics.

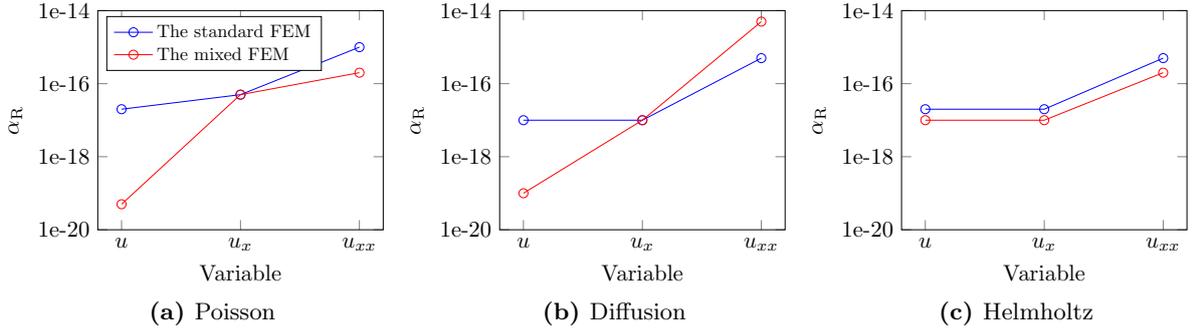
\begin{figure}[!ht]
\hspace{0.0cm}
\begin{subfigure}[b]{0.35\textwidth}
\scalebox{0.85}{
\begin{tikzpicture} 
\begin{axis}
[
    ymode=log,    
    ymin=1e-20,
    ymax=1e-14,
    ytick={1e-20, 1e-18, 1e-16, 1e-14},
    yticklabels={1e-20, 1e-18, 1e-16, 1e-14},      
    legend style={nodes={scale=0.8},at={(0.03,0.85)},anchor=west},
    legend cell align={left},
    height=5cm,
    width=6cm,
    ylabel={$\alpha_{\rm R}$},
    ylabel style={at={(-0.01,0.5)}},    
    xtick={0,1,2,3,4},
    xticklabels={$u$, $u_x$, $u_{xx}$},
    xlabel={Variable},
    xlabel style={at={(0.5,0.03)}},    
]
\addplot[blue,mark=o,mark options={color=blue,fill=blue}] coordinates {(0,2.0e-17) (1,5.0e-17) (2,1.0e-15)};
\addplot[red,mark=o,mark options={color=red,fill=red}] coordinates {(0,5.0e-20) (1,5.0e-17) (2,2.0e-16)};
\legend{The standard FEM, The mixed FEM};
\end{axis}
\end{tikzpicture}
}
\vspace{-0.2cm}
\caption{Poisson}
\label{alpha_R_benchmark_Poisson}
\end{subfigure}
\hspace{-0.7cm}
\begin{subfigure}[b]{0.35\textwidth}
\scalebox{0.85}{
\begin{tikzpicture} 
\begin{axis}
[
    ymode=log,    
    ymin=1e-20,
    ymax=1e-14,
    ytick={1e-20, 1e-18, 1e-16, 1e-14},
    yticklabels={1e-20, 1e-18, 1e-16, 1e-14},      
    legend style={nodes={scale=0.8},at={(0.03,0.85)},anchor=west},
    legend cell align={left},
    height=5cm,
    width=6cm,
    ylabel={$\alpha_{\rm R}$},
    ylabel style={at={(-0.01,0.5)}},    
    xtick={0,1,2,3,4},
    xticklabels={$u$, $u_x$, $u_{xx}$},
    xlabel={Variable},
    xlabel style={at={(0.5,0.03)}},    
]
\addplot[blue,mark=o,mark options={color=blue,fill=blue}] coordinates {(0,1.0e-17) (1,1.0e-17) (2,5.0e-16)};
\addplot[red,mark=o,mark options={color=red,fill=red}] coordinates {(0,1.0e-19) (1,1.0e-17) (2,5.0e-15)};
\end{axis}
\end{tikzpicture}
}
\vspace{-0.2cm}
\caption{Diffusion}
\label{alpha_R_benchmark_diffusion}
\end{subfigure}
\hspace{-0.7cm}
\begin{subfigure}[b]{0.35\textwidth}
\scalebox{0.85}{
\begin{tikzpicture} 
\begin{axis}
[
    ymode=log,    
    ymin=1e-20,
    ymax=1e-14,
    ytick={1e-20, 1e-18, 1e-16, 1e-14},
    yticklabels={1e-20, 1e-18, 1e-16, 1e-14},      
    legend style={nodes={scale=0.8},at={(0.03,0.85)},anchor=west},
    legend cell align={left},
    height=5cm,
    width=6cm,
    ylabel={$\alpha_{\rm R}$},
    ylabel style={at={(-0.01,0.5)}},    
    xtick={0,1,2,3,4},
    xticklabels={$u$, $u_x$, $u_{xx}$},
    xlabel={Variable},
    xlabel style={at={(0.5,0.03)}},    
]
\addplot[blue,mark=o,mark options={color=blue,fill=blue}] coordinates {(0,2.0e-17) (1,2.0e-17) (2,5.0e-16)};
\addplot[red,mark=o,mark options={color=red,fill=red}] coordinates {(0,1.0e-17) (1,1.0e-17) (2,2.0e-16)};
\end{axis}
\end{tikzpicture}
}
\vspace{-0.2cm}
\caption{Helmholtz}
\label{alpha_R_benchmark_Helmholtz}
\end{subfigure}
\caption{$\alpha_{\rm R}$ for the benchmark equations.}
\label{alpha_R_benchmark}
\end{figure}

It also shows that $\alpha_{\rm R}$ tends to increase slightly with increasing order of derivative, see Fig.~\ref{alpha_R_benchmark_Poisson}. Since $E_{\rm T}$ decreases slower after each differentiation, for the same element degree, using the standard FEM, ${E}_{\text{min}}$ tends to deteriorate with increasing order of derivative; using the mixed FEM, since the degree of the elements used for $u_{x}$ is one order higher than that used for $u$, ${E}_{\text{min}}$ for $u$ and $u_{x}$ tend to be of the same order, but ${E}_{\text{min}}$ for $u_{xx}$ is still larger than that for $u_{x}$, see Fig.~\ref{E_min_benchmark_Poisson} for the statistics. 

\subsubsection{Benchmark diffusion and Helmholtz equations}      \label{section_equation_type}

For the benchmark diffusion and Helmholtz equations, the discretization errors are shown in \ref{discretization_error_bench_diff}.
The slopes $\beta_{\rm R}$ remain the same with that of the Poisson equation.
The offsets $\alpha_{\rm R}$ and  $E_{\rm min}$ also follow the same trend, see the rest of Fig.~\ref{alpha_R_benchmark} and Fig.~\ref{E_min_benchmark}, respectively.

Summarizing this section, $\alpha_{\rm R}$ varies not only with the variable, but also with the equation and FEM method; $\beta_{\rm R}$ is relatively fixed, which is 2 using the standard FEM and 1 using the mixed FEM. In what follows, we will take $\beta_{\rm R}$ as constant if not stated otherwise.

\subsection{Sensitivity analysis}       \label{section_sensitivity}

We focus on the benchmark Poisson equation, for which $P_2$ elements are used for the standard FEM, and $P_4/P_3^{\rm disc}$ elements are used for the mixed FEM.

\subsubsection{Solution strategy}		\label{section_solver}

In this section, we investigate the influence of the solution strategy on the accuracy of the numerical solution. In particular, we compare the outcome when applying the direct solver UMFPACK with that of using the iterative Conjugate Gradient (CG) method \cite{ginsburg1963cg}, which can be applied since the system matrix $A$ in Eq. (\ref{matrix equation std FEM}) is symmetric and positive definite. The tolerance of the CG solver is set to be the product of a parameter, denoted by $tol_{prm}$, and the $L_2$ norm of the discrete right-hand side $\|F\|_2$.
When the $L_2$ norm of the residual, i.e. $\|F-Au\|_2$ in Eq. (\ref{matrix equation std FEM}), is smaller than the tolerance, the iteration is stopped.
For the mixed FEM, we additionally investigate the impact of using a segregated solution approach based on the Schur complement instead of a fully coupled approach.

\paragraph{The standard FEM}

\begin{figure}[!ht]
    \begin{subfigure}{5.5cm}
        \includegraphics[width=1.0\linewidth]{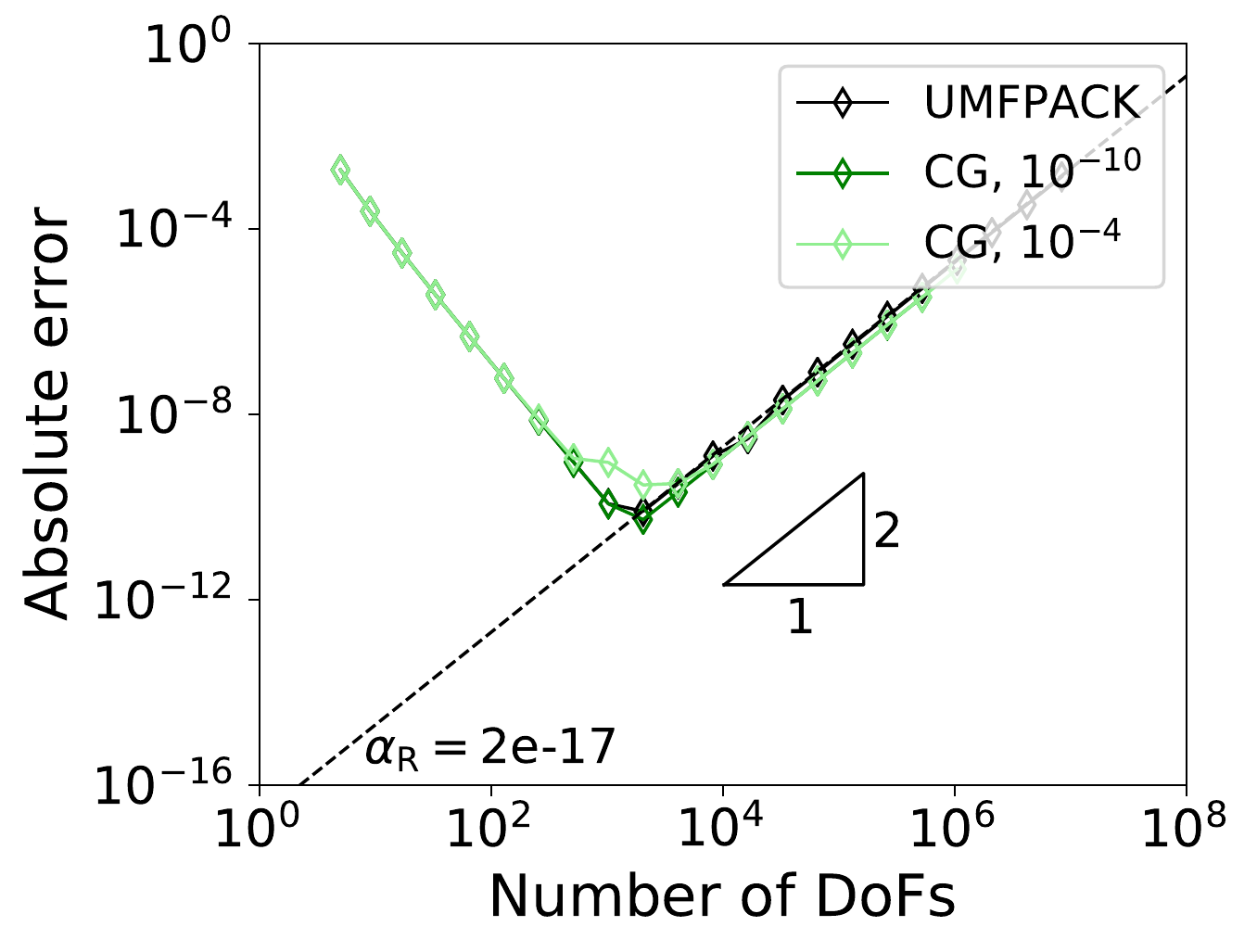}
        \caption{Solution}
        \label{py_bench_Pois_SM_error_solution_strategy_solu}
    \end{subfigure}
    \hspace{-0.2cm}
    \begin{subfigure}{5.5cm}
        \includegraphics[width=1.0\linewidth]{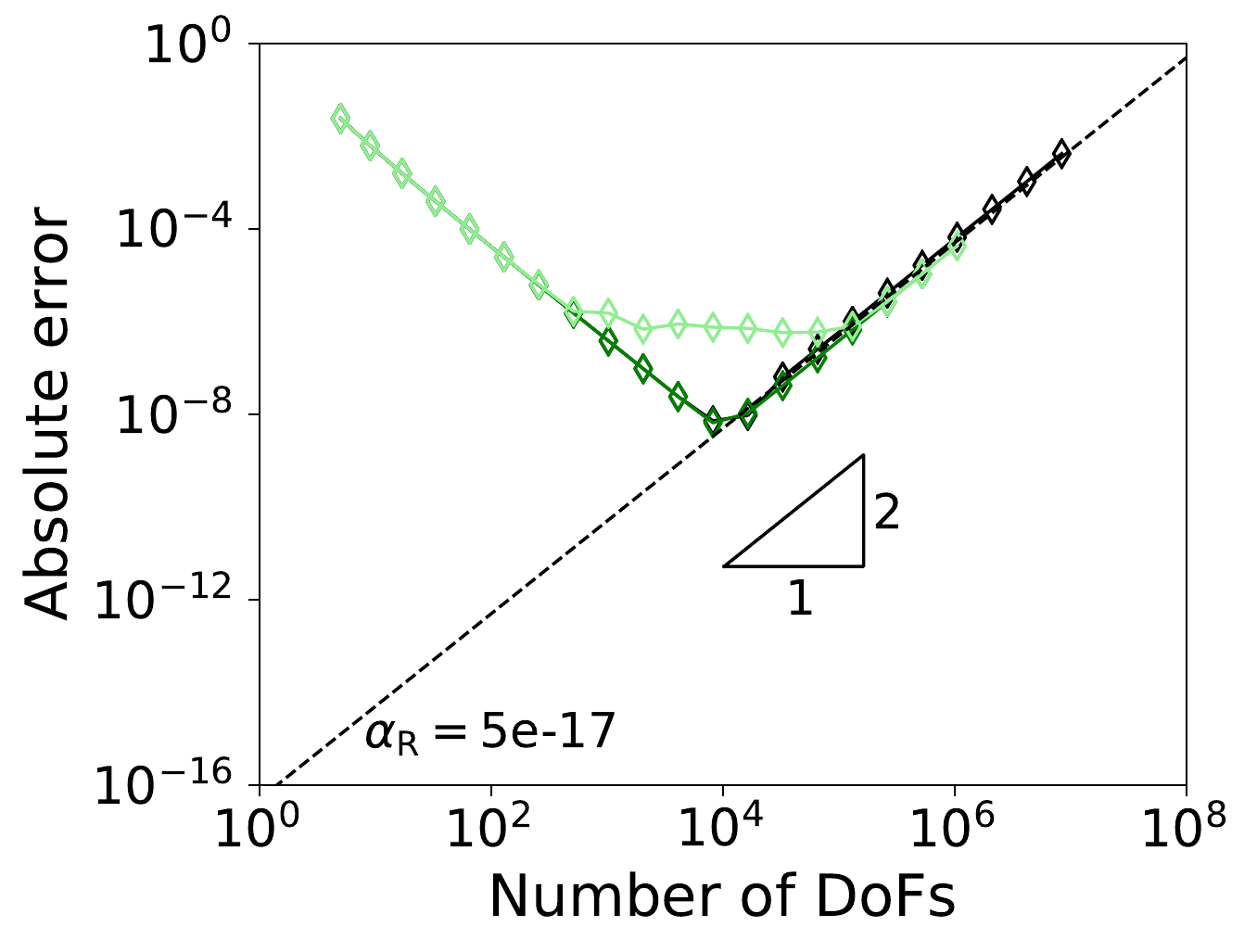}
        \caption{First derivative}
        \label{py_bench_Pois_SM_error_solution_strategy_grad}
    \end{subfigure}
    \hspace{-0.2cm}
    \begin{subfigure}{5.5cm}
        \includegraphics[width=1.0\linewidth]{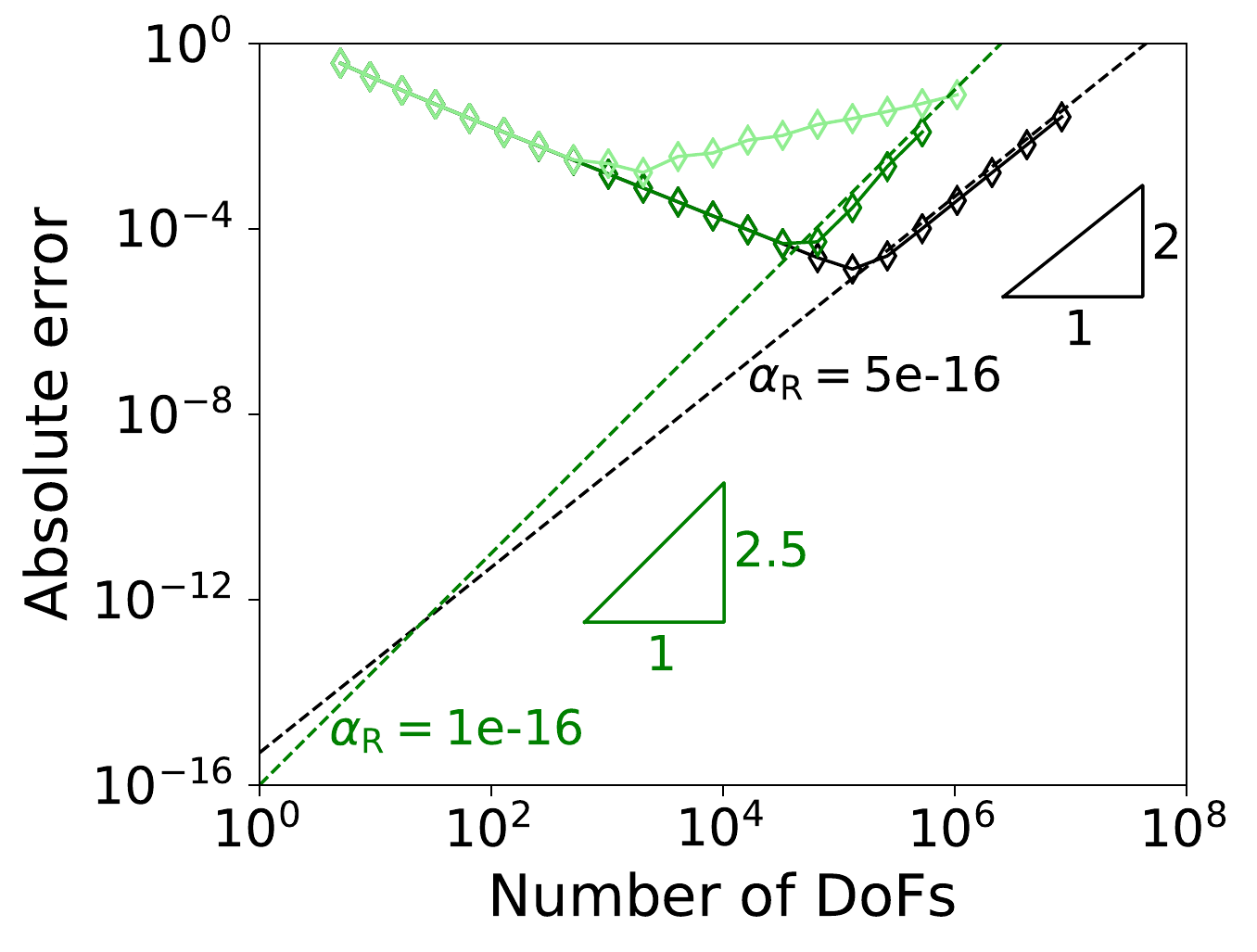}
        \caption{Second derivative}
        \label{py_bench_Pois_SM_error_solution_strategy_2ndd}
    \end{subfigure}
\caption{Comparison of the errors using the CG solver and the UMFPACK solver.}
\label{py_bench_Pois_SM_error_solution_strategy}
\end{figure}

The CG solver is stopped once $\|F-Au\|_2 \leqslant tol_{prm}\|F\|_2$, with $tol_{prm}=$ $10^{-10}$ and $10^{-4}$, respectively. 
The absolute errors for $u$, $u_{x}$ and $u_{xx}$ using the CG solver are shown in Fig.~\ref{py_bench_Pois_SM_error_solution_strategy}, in comparison with that using the direct solver UMFPACK. 

When $tol_{prm}$ is adequately small, i.e. $tol_{prm}=10^{-10}$, the round-off error for the solution and the first derivative using the CG solver is the same with that using the UMFPACK solver; the round-off error for the second derivative using the CG solver increases faster than that using the UMFPACK solver.
When $tol_{prm}$ is too large, i.e. $tol_{prm}=10^{-4}$, the error contribution due to the iterative solver dominates both truncation and round-off errors. 

\paragraph{The mixed FEM}

Since the resulting matrix Eq.~(\ref{matrix equation mix FEM}) is indefinite, a widely used alternative is to decouple the fully coupled monolithic approach
\begin{subequations}
 \begin{align}
  B^{\top} M^{-1} B U &= B^{\top} M^{-1} G - H, 	\label{schur_complement_solution} \\
  MV&=G-BU						\label{schur_complement_gradient}
\end{align}						\label{schur_complement_solu_grad}%
\end{subequations}
and solve both equations in segregated manner, i.e. Eq. (\ref{schur_complement_solution}) is solved in the first place to obtain $U$, and then it is substituted into Eq. (\ref{schur_complement_gradient}) to obtain $V$.

Eq. (\ref{schur_complement_solution}) involves the term $M^{-1} G$ in the right-hand side, which is computed by solving the auxiliary linear system $MY=G$ by using either the UMFPACK or the CG solver. The same options are available for solving Eq. (\ref{schur_complement_gradient}). 

The difficulty in solving Eq. (\ref{schur_complement_solution}) lies in not assembling the Schur complement matrix explicitly since it comprises $M^{-1}$. 
The CG solver only makes use of matrix-vector products of the form $(B^{\top}M^{-1}B)W$, which can be computed by the following three-step algorithm: $X=BW$, $MY=X$ and $Z=B^{\top}Y$. As before, the linear system $MY=X$ can be solved by the UMFPACK or the CG solver.
 
\begin{figure}[!ht]
    \begin{subfigure}{5.5cm}
        \includegraphics[width=1.0\linewidth]{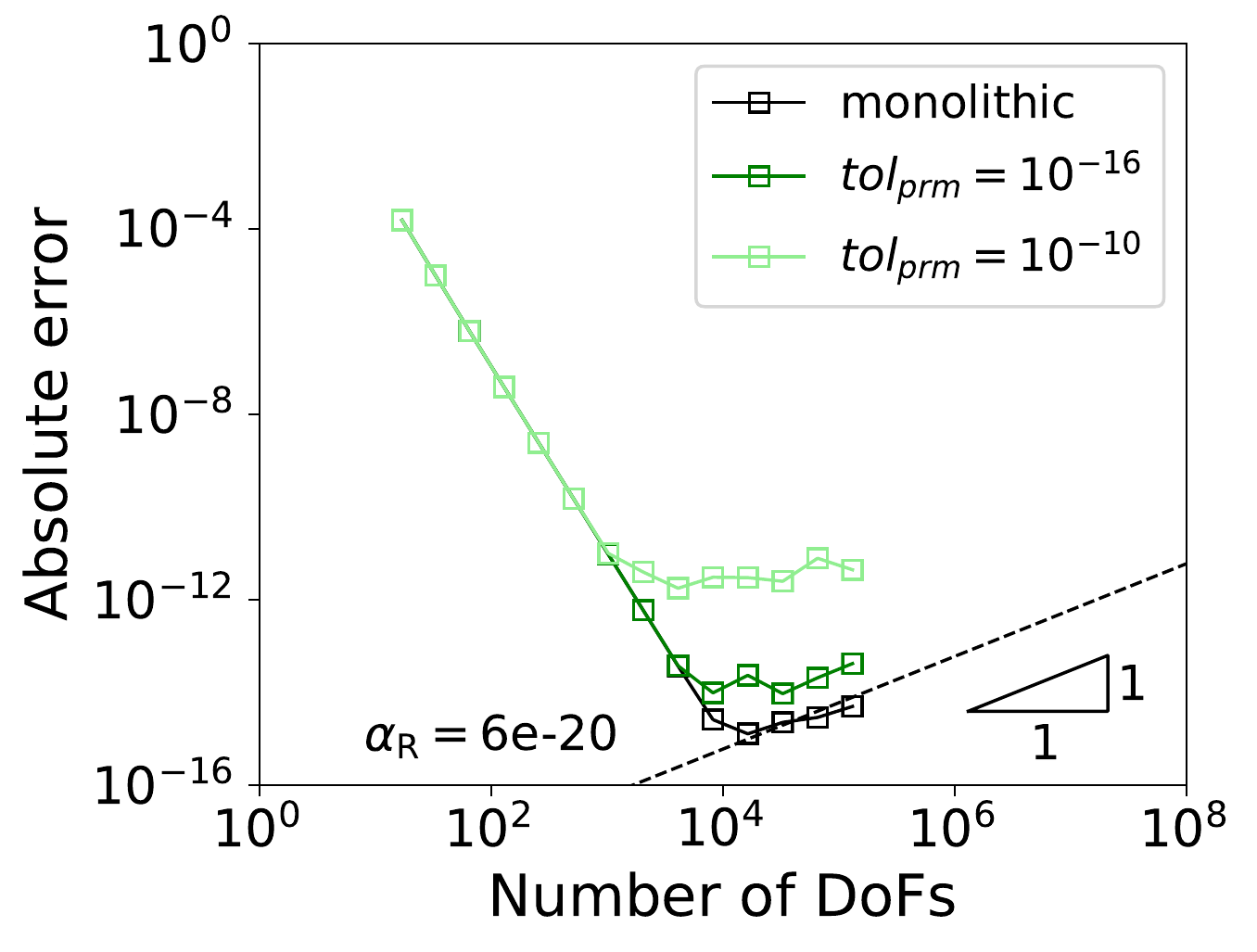}
        \caption{Solution}
        \label{Fig:py_bench_Pois_MM_error_solution_strategy_schur_variant_other_UMF_solu}
    \end{subfigure}
    \hspace{-0.2cm}
    \begin{subfigure}{5.5cm}
        \includegraphics[width=1.0\linewidth]{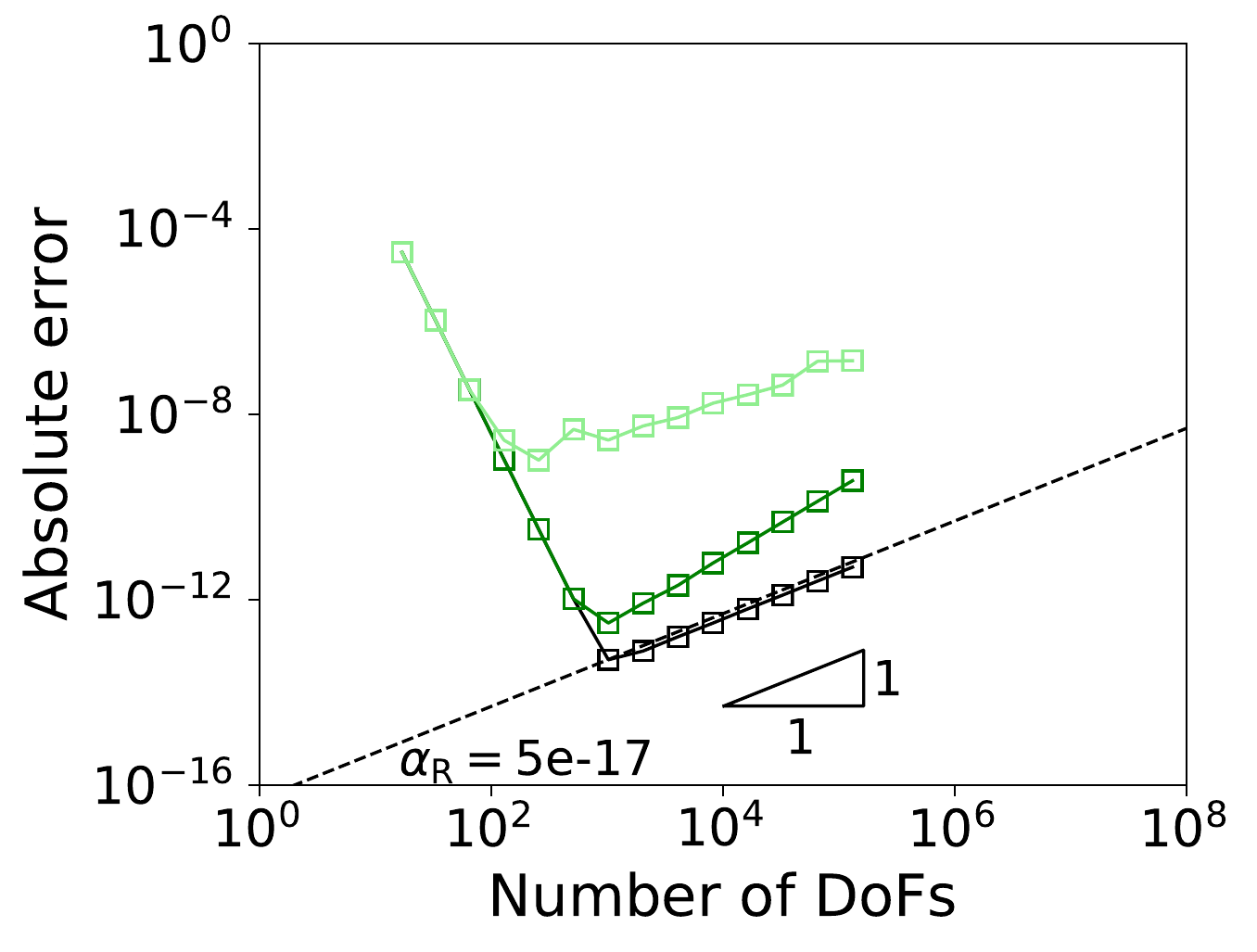}
        \caption{First derivative}
        \label{Fig:py_bench_Pois_MM_error_solution_strategy_schur_variant_other_UMF_grad}
    \end{subfigure}
    \hspace{-0.2cm}
    \begin{subfigure}{5.5cm}
        \includegraphics[width=1.0\linewidth]{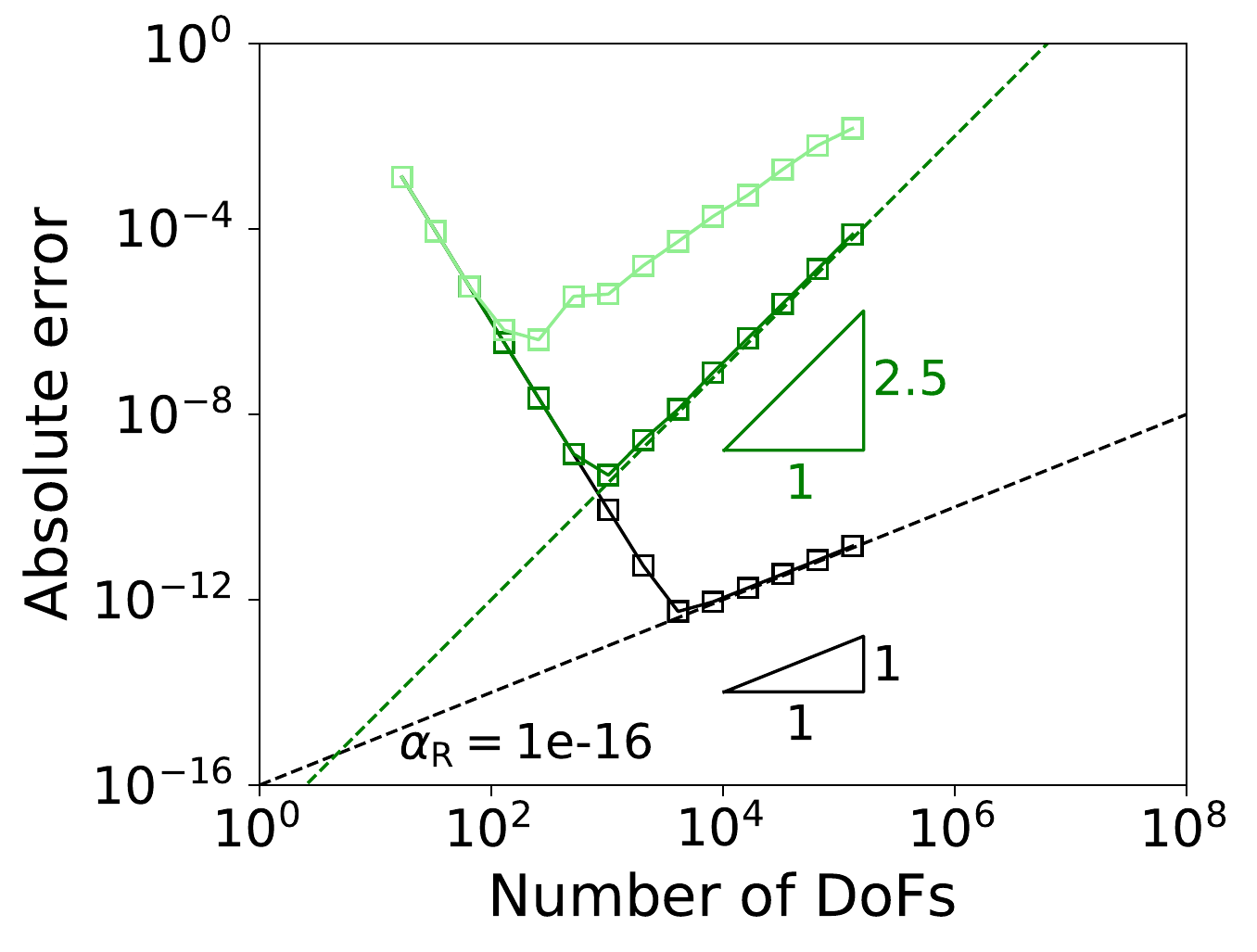}
        \caption{Second derivative}
        \label{Fig:py_bench_Pois_MM_error_solution_strategy_schur_variant_other_UMF_2ndd}
    \end{subfigure}
\caption{Influence of the CG solver on the accuracy when the left-hand side is the Schur complement using the mixed FEM.}
\label{Fig:py_bench_Pois_MM_error_solution_strategy_schur_variant_other_UMF}
\end{figure}

We first investigate the influence of $tol_{prm}$ of the CG solver on the accuracy of the solutions when the left-hand side is $B^{\top}M^{-1}B$. In this case, the UMFPACK solver is used to solve the matrix equations when the left-hand side is $M$.  
For $tol_{prm}$ being $10^{-16}$ and $10^{-10}$, the results are shown in Fig. \ref{Fig:py_bench_Pois_MM_error_solution_strategy_schur_variant_other_UMF}, in comparison with that obtained from solving the monolithic Eq. (\ref{matrix equation mix FEM}) directly using the UMFPACK solver.
It shows that, for the problem at hand, the monolithic solution approach yields by far the most accurate solution and derivative values.
Remarkably, the round-off error for $v_{x}$ increases fastest using the Schur complement approach even though $tol_{prm}$ is sufficiently small, i.e. $tol_{prm}=10^{-16}$, which makes the highest attainable accuracy much lower.
When $tol_{prm}$ is less strict, i.e. $tol_{prm}=10^{-10}$, the iteration error dominates the total error instead of the round-off error.

\begin{figure}[!ht]
    \begin{subfigure}{5.5cm}
        \includegraphics[width=1.0\linewidth]{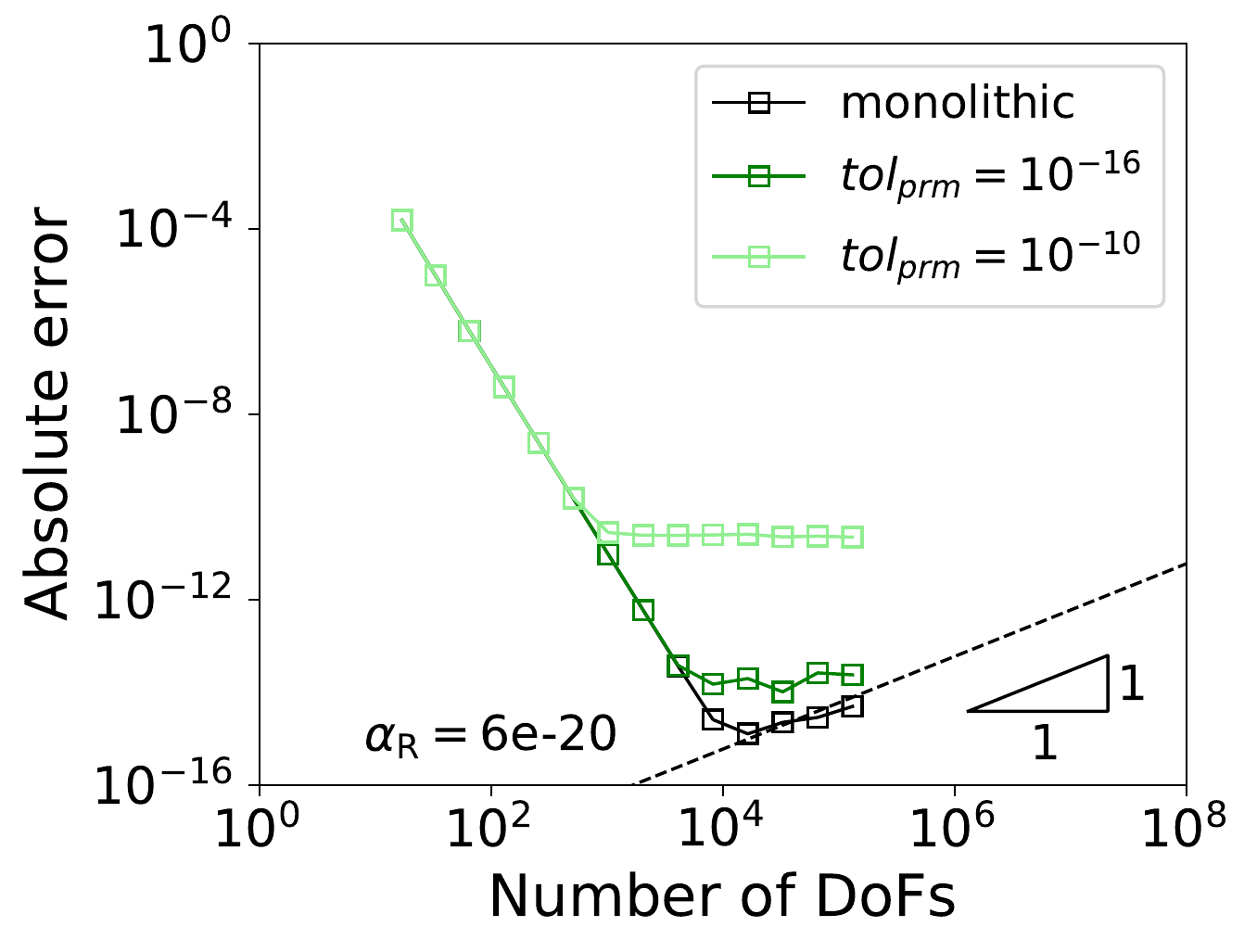}
        \caption{Solution}
        \label{py_bench_Pois_MM_error_solution_strategy_schur_1em16_M_variant_solu}
    \end{subfigure}
    \hspace{-0.2cm}
    \begin{subfigure}{5.5cm}
        \includegraphics[width=1.0\linewidth]{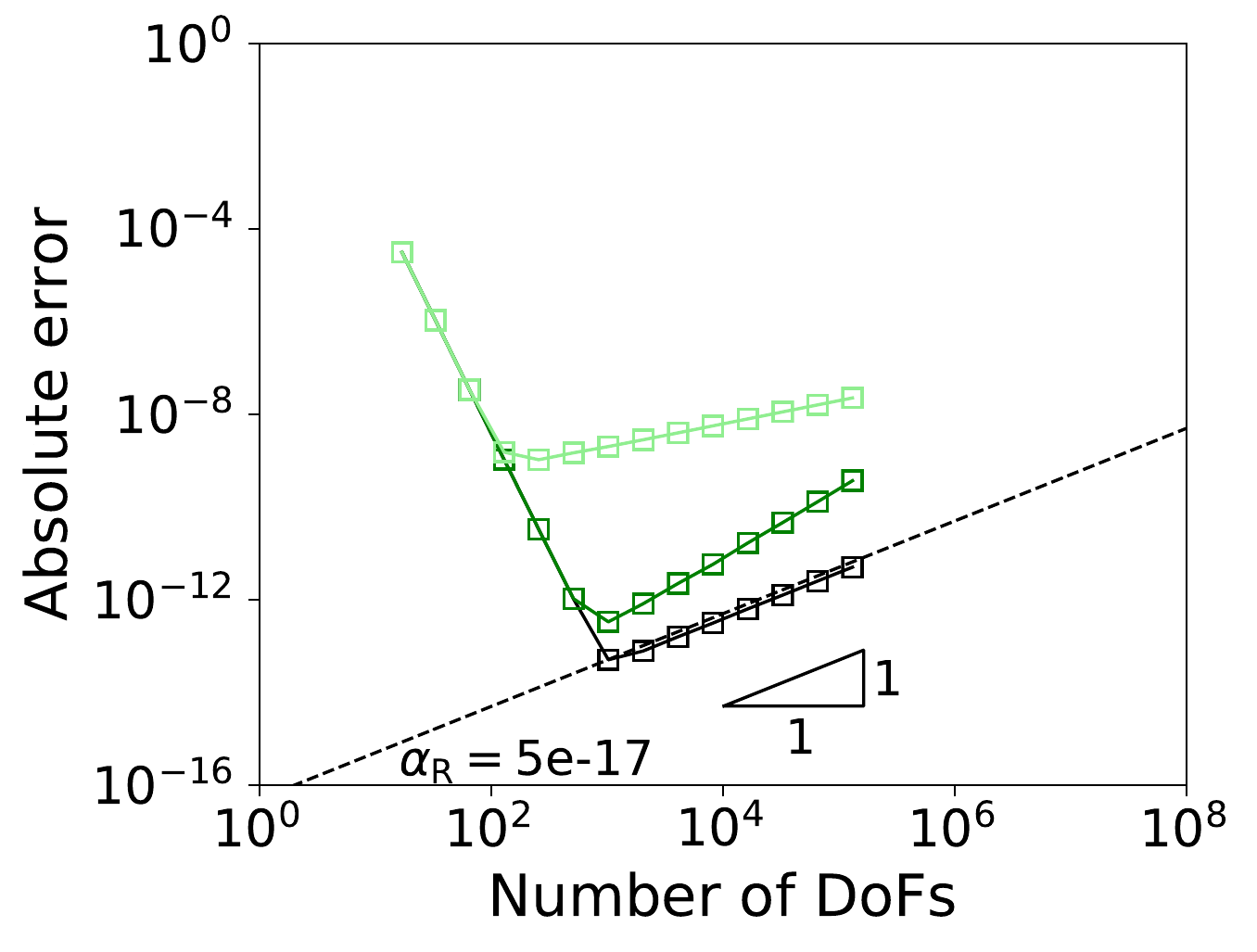}
        \caption{First derivative}
        \label{py_bench_Pois_MM_error_solution_strategy_schur_1em16_M_variant_grad}
    \end{subfigure}
    \hspace{-0.2cm}
    \begin{subfigure}{5.5cm}
        \includegraphics[width=1.0\linewidth]{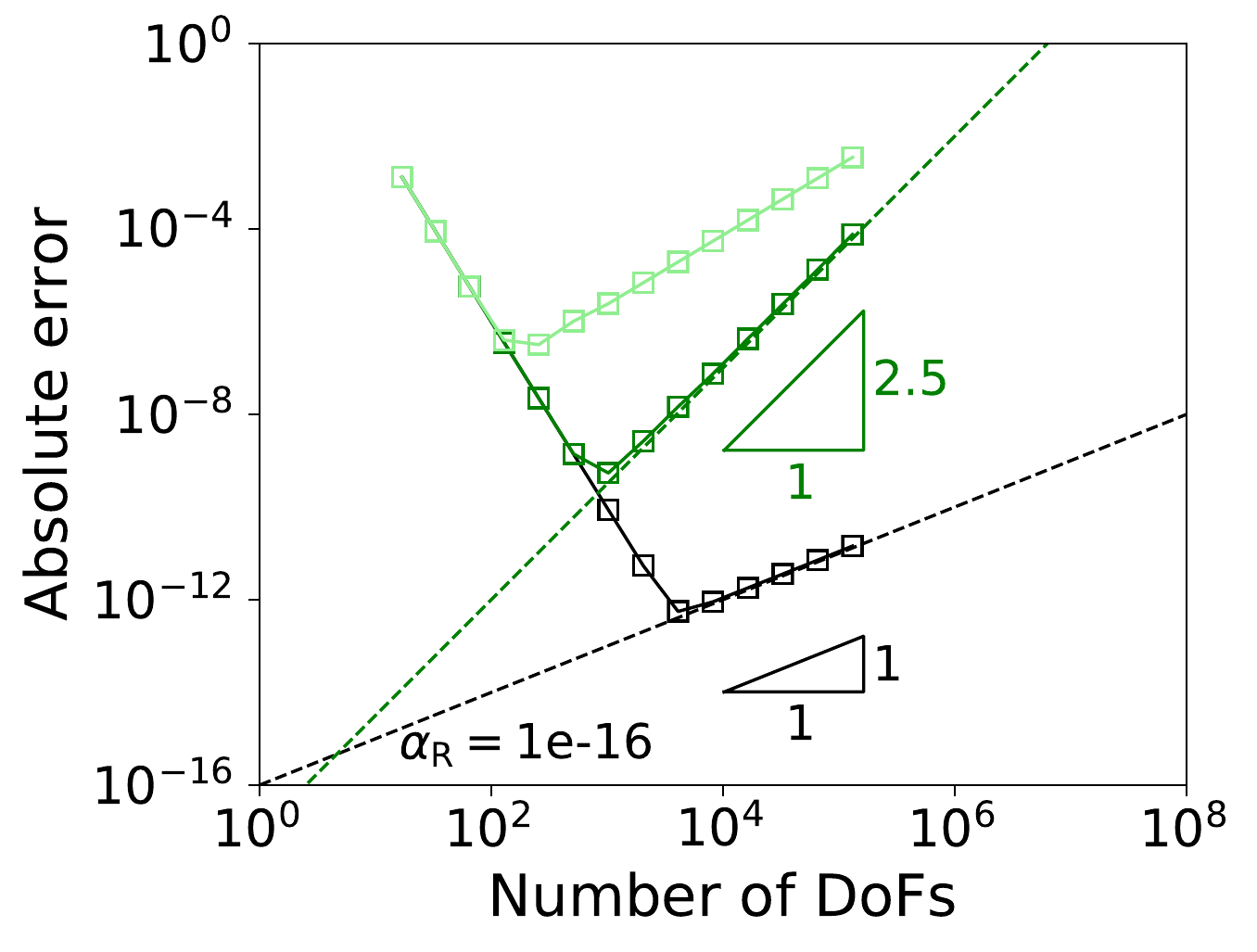}
        \caption{Second derivative}
        \label{py_bench_Pois_MM_error_solution_strategy_schur_1em16_M_variant_2ndd}
    \end{subfigure}
\caption{Influence of the CG solver on the accuracy when the left-hand side is $M$ using the mixed FEM.}		%
\label{py_bench_Pois_MM_error_solution_strategy_schur_1em16_M_variant}
\end{figure}

Next, we investigate the influence of $tol_{prm}$ of the CG solver when the left-hand side is $M$. In this case, the CG solver with $tol_{prm}$ being $10^{-16}$ is used to solve the matrix equation with the left-hand side being $B^{\top}M^{-1}B$. For $tol_{prm}$ being $10^{-16}$ and $10^{-10}$, the results are shown in Fig.~\ref{py_bench_Pois_MM_error_solution_strategy_schur_1em16_M_variant}, in comparison with that obtained from solving the monolithic Eq. (\ref{matrix equation mix FEM}) directly using the UMFPACK solver.
It also shows that, when the tolerance is less strict, i.e. $tol_{prm}=10^{-10}$, the iteration error dominates the total error before the round-off error.

In summary, for the standard FEM, the CG solver gives the same accuracy for $u$ and $u_{x}$ as the UMFPACK solver when $tol_{prm}$ is strict enough, while the UMFPACK solver is recommended for computing $u_{xx}$; for the mixed FEM, the accuracy for all the three variables is the highest when using the UMFPACK solver to solve the monolithic Eq. (\ref{matrix equation mix FEM}) directly. Moreover, the application of the CG solver on both the standard and mixed FEM methods shows that less strict values for $tol_{prm}$ introduce iteration errors.

\subsubsection{Order of magnitude}	    \label{section_scaling}

\begin{figure}[!ht]
\centering
      \includegraphics[width=0.4\linewidth]{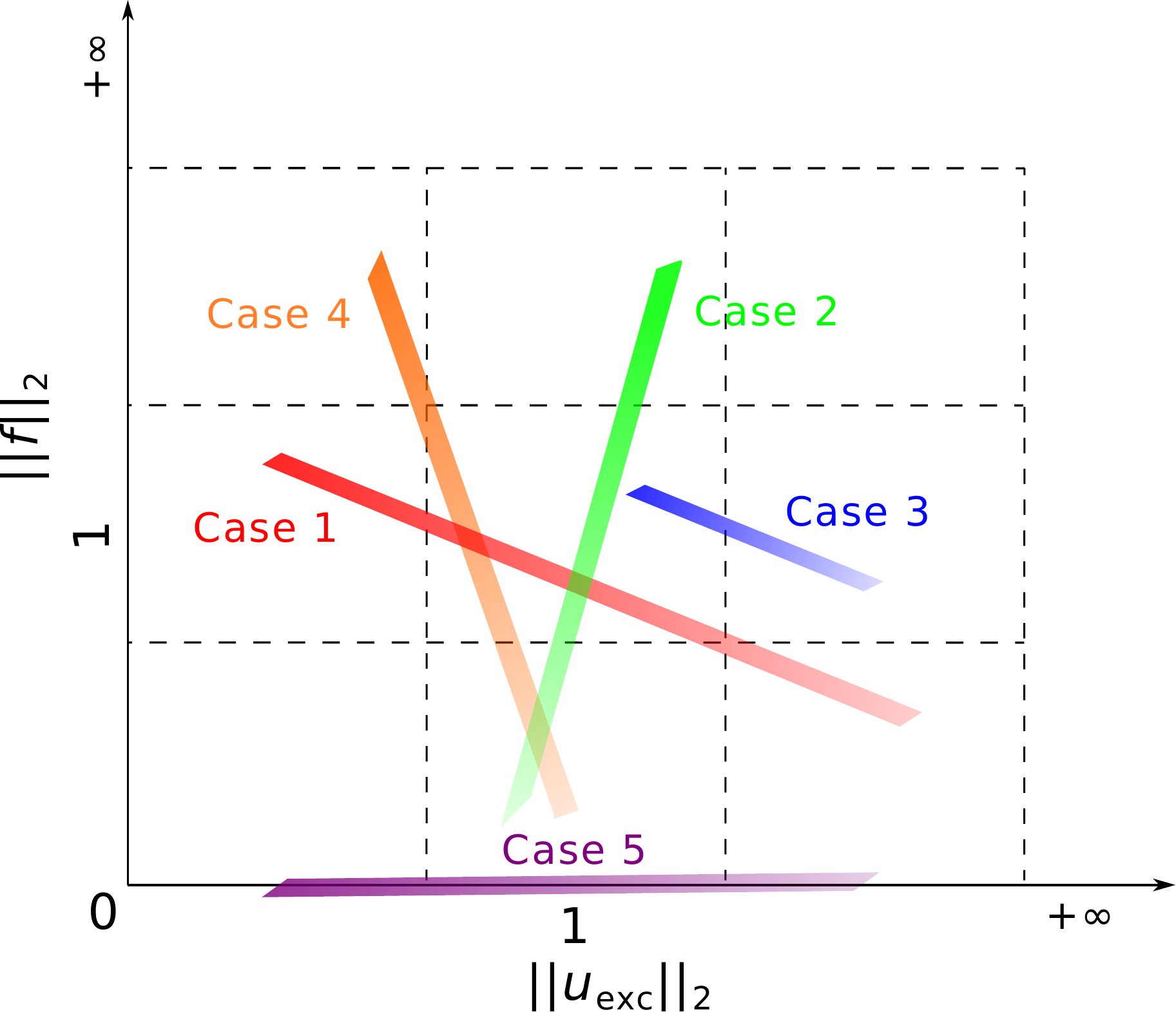}
\caption{Distribution of $\|u_{\rm exc}\|_2$ and $\|f\|_2$ for the test cases with the settings from Table \ref{scaling_cases_Poisson}. The color density increases with the value of the coefficient $c_i$.}
\label{L2 norm u f}
\end{figure}

In this section, we investigate the influence of the order of magnitude of the solution and the right-hand side on the offset $\alpha_{\rm R}$ of the round-off error and propose different scaling schemes to mitigate this influence factor.
To cover a wide range of scenarios, we choose the right-hand sides shown in the second column of Table \ref{scaling_cases_Poisson}. The corresponding boundary conditions and analytical solutions are given in the remaining columns of the table. 
Each case contains a coefficient $c_i,~i=1,2, \ldots , 5$, which is varied over several orders of magnitude so that the $L_2$ norm of the exact solution, denoted by $\|u_{\rm exc}\|_2$, and the $L_2$ norm of the right-hand side, denoted by $\|f\|_2$, extend over a wide range of magnitudes. Fig.~\ref{L2 norm u f} gives an overview of the distribution of $\|u_{\rm exc}\|_2$ and $\|f\|_2$ for different cases, and the more detailed information can be found in Fig.~\ref{L2_norms_cases_1_to_5}.


\begin{table}[!ht]
\centering
\caption [w]{Setting of the Poisson equation with different right-hand sides.} 
\label{scaling_cases_Poisson}
 \begin{tabular}{c c c c c} \hline      
\multirow{2}{*}{Case} & \multirow{2}{*}{$f(x)$}  & \multicolumn{2}{c}{Boundary conditions} & \multirow{2}{*}{$u_{\text{exc}}(x)$} \\
\cline{3-4}
& & $u(0)$ & $u(1)$ & \\ \hline
{1} & {$\sin (2 \pi c_1x)$} & {0}& ${(2 \pi c_1)}^{-2} \sin (2 \pi c_1)$ & ${(2 \pi c_1)}^{-2} \sin (2 \pi c_1x)$\\ \hline
2 & $\makecell{-e^{-{c_2}{(x-1/2)^2}} \cdot \\ \left({4{c_2}^2(x-1/2)^2 -2c_2} \right)}$ & $e^{-c_2/4}$ & $e^{-c_2/4}$ & $e^{-{c_2}{{(x-1/2)^2}}}$ \\ \hline
3 & $\sin (2 \pi c_3 x) +1$ & $0$ & ${(2 \pi c_3)}^{-2}\sin (2 \pi c_3)-\frac{1}{2}$ & ${(2 \pi c_3)}^{-2}\sin (2 \pi c_3 x)-\frac{x^2}{2}$ \\ \hline
4 & $2 \pi c_4 \sin (2 \pi c_4 x)$ & $0$ & ${(2 \pi c_4)}^{-1} \sin (2 \pi c_4)$ & ${(2 \pi c_4)}^{-1} \sin (2 \pi c_4x)$ \\ \hline
5 & $0$ & $0$ & ${c_5}^{-1}$ & ${c_5}^{-1} x$ \\ \hline
\end{tabular}
\end{table}

For case 1, the results are given below, and that of other cases can be found in \ref{results_order_of_magnitude_other_cases}, which shows qualitatively the same behavior.

\paragraph{The standard FEM}		\label{scaling_std_FEM}

\begin{figure}[!ht]
    \begin{subfigure}{5.5cm}
        \includegraphics[width=1.0\linewidth]{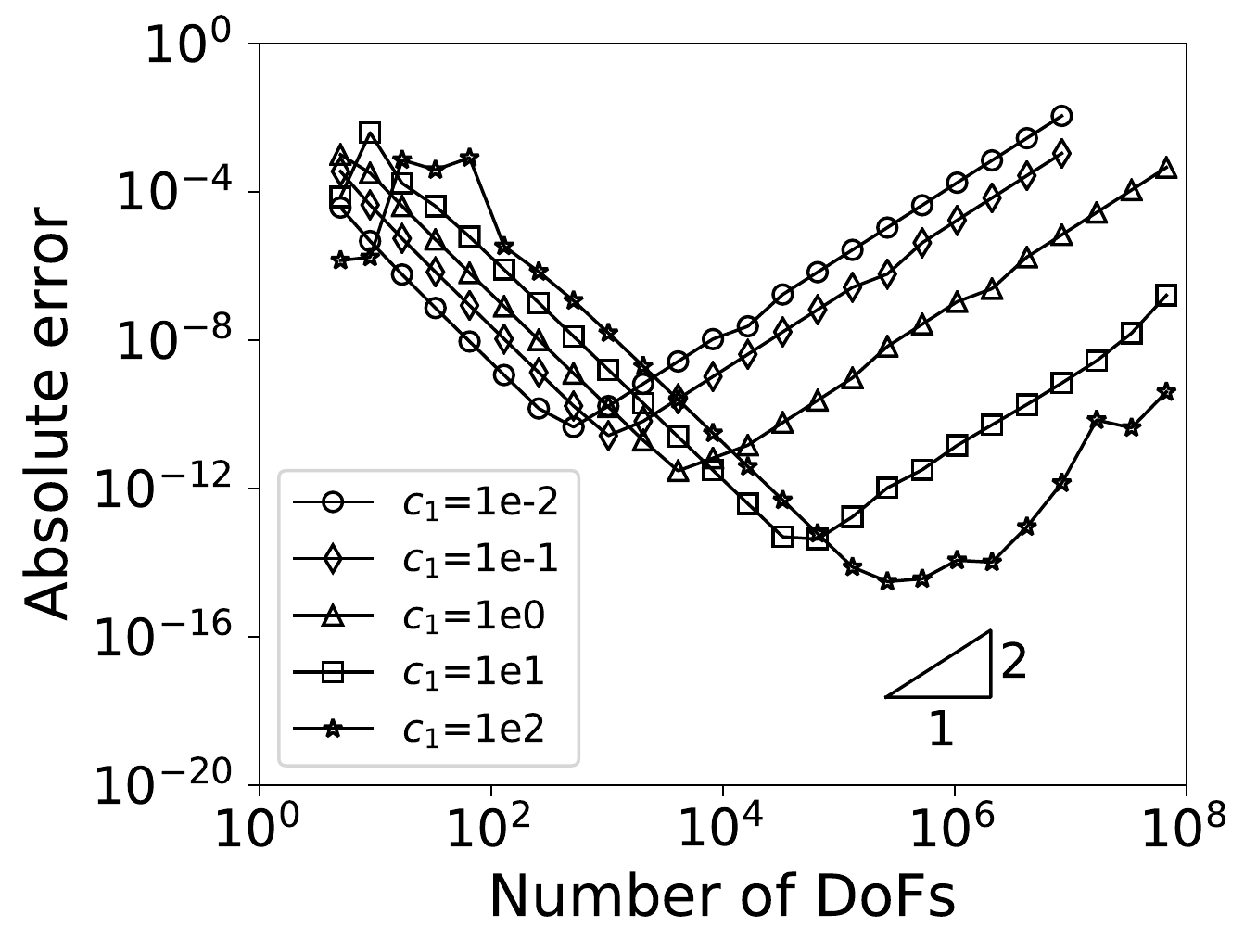}
        \caption{Solution}
        \label{py_L2_Pois1_SM_scaling_no_solu}
    \end{subfigure}
    \hspace{-0.2cm}
    \begin{subfigure}{5.5cm}
        \includegraphics[width=1.0\linewidth]{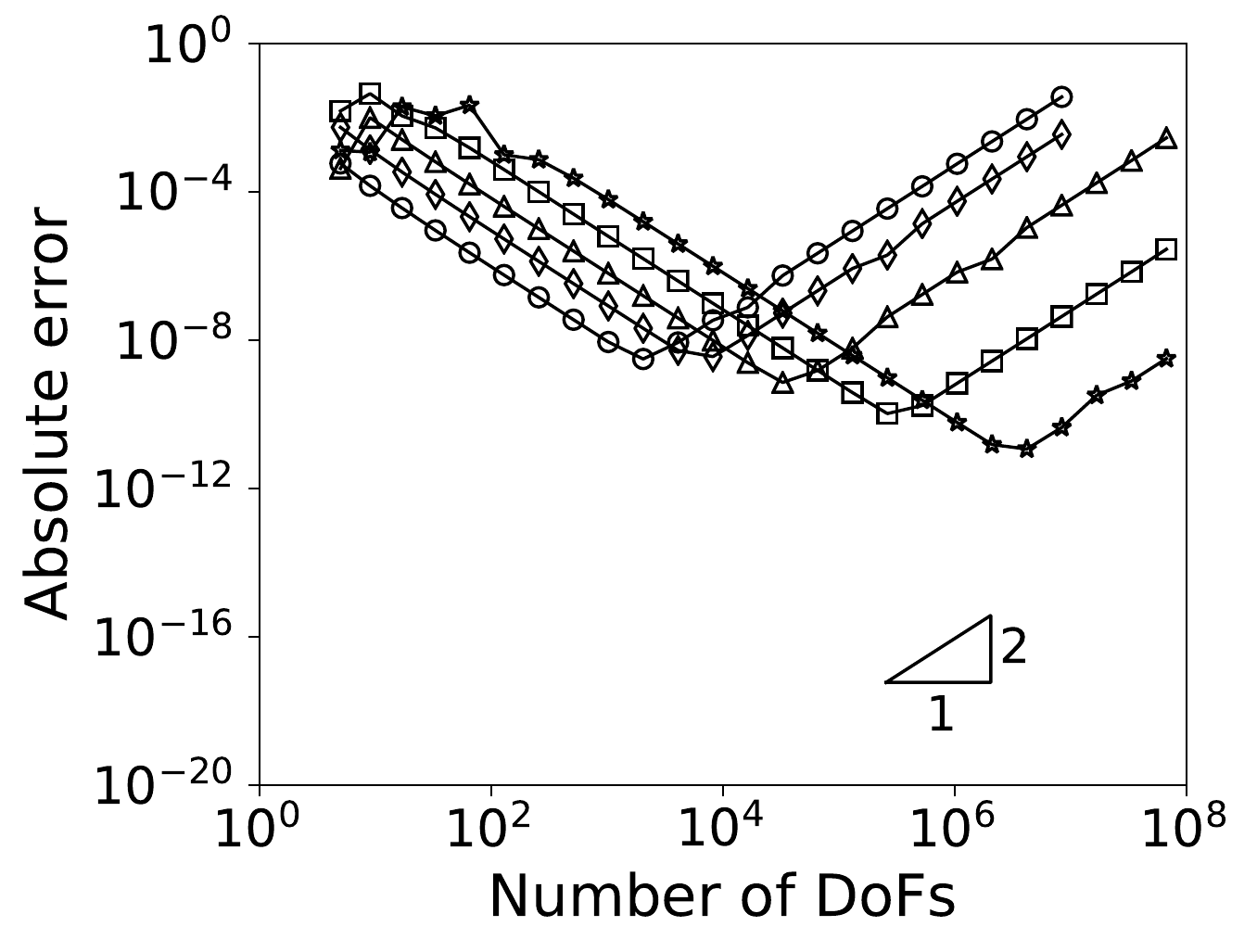}
        \caption{First derivative}
        \label{py_L2_Pois1_SM_scaling_no_grad}
    \end{subfigure}
    \hspace{-0.2cm}
    \begin{subfigure}{5.5cm}
        \includegraphics[width=1.0\linewidth]{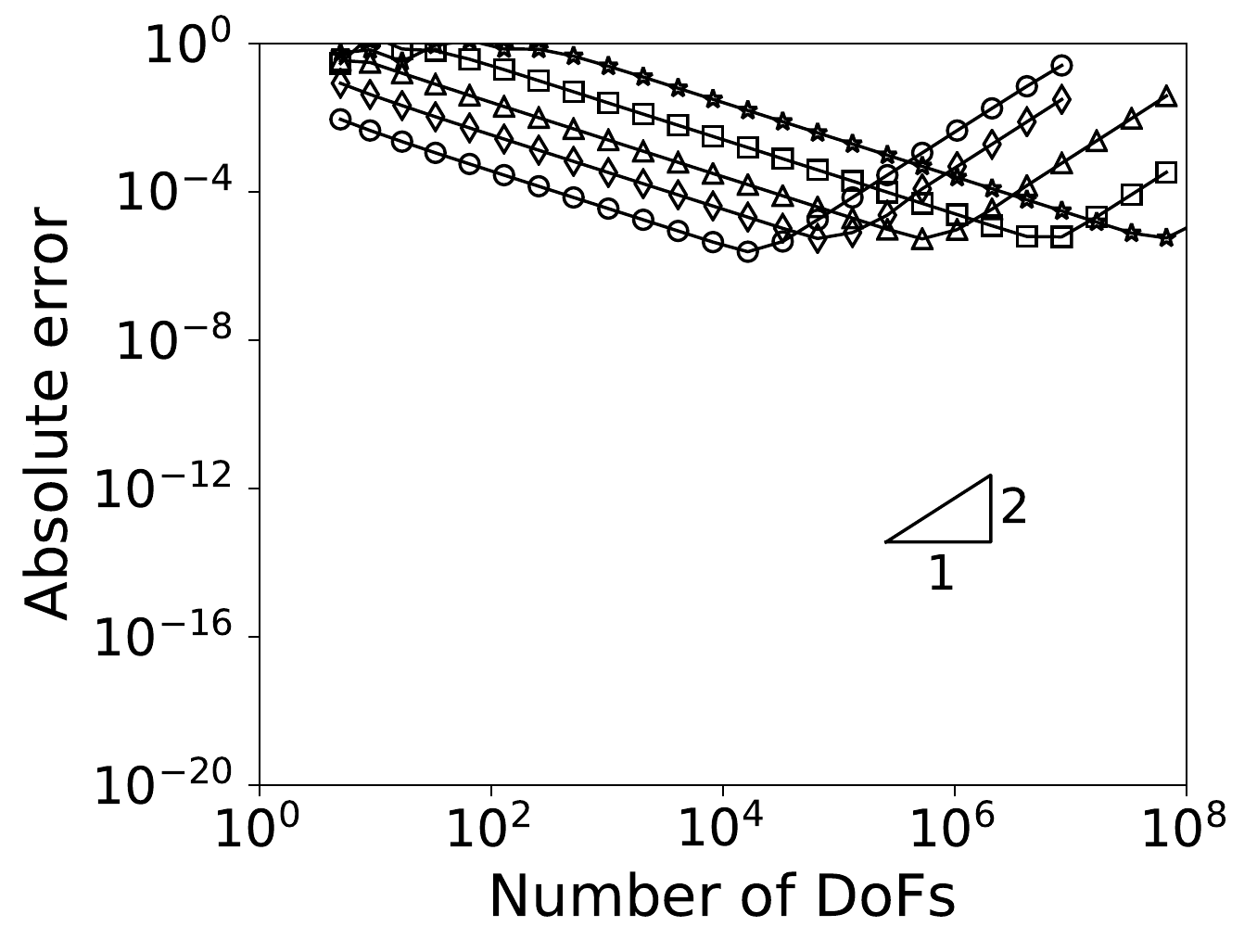}
        \caption{Second derivative}
        \label{py_L2_Pois1_SM_scaling_no_2ndd}
    \end{subfigure}
\caption{Absolute errors for different $c_1$ using the standard FEM without scaling the right-hand side.}   
\label{Pois_SM_rhs_solu orig}
\end{figure}

\begin{figure}[!ht]
    \begin{subfigure}{5.5cm}
        \includegraphics[width=1.0\linewidth]{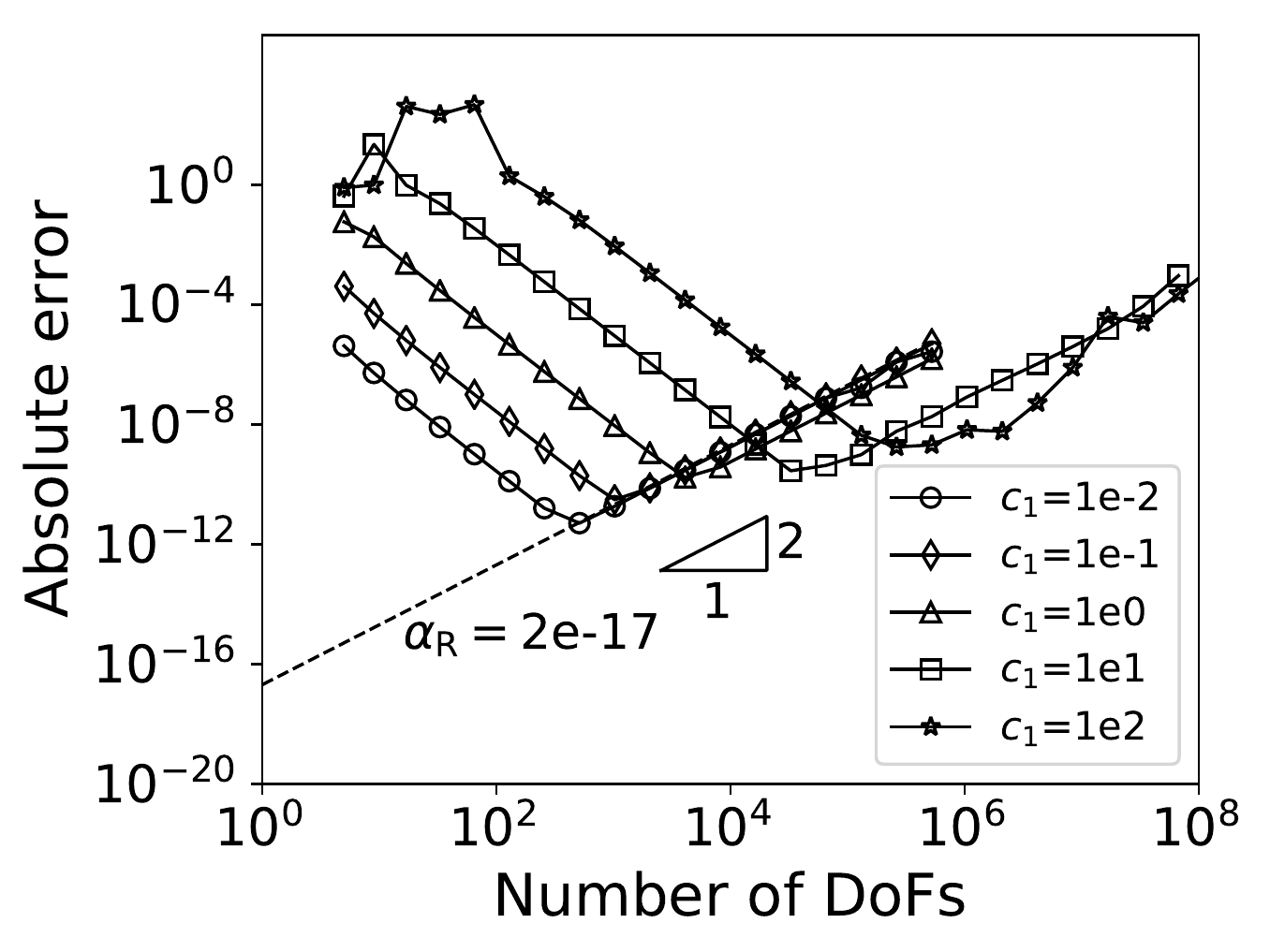}
        \caption{Solution}
        \label{py_L2_Pois1_SM_scaling_S_solu}
    \end{subfigure}
    \hspace{-0.2cm}
    \begin{subfigure}{5.5cm}
        \includegraphics[width=1.0\linewidth]{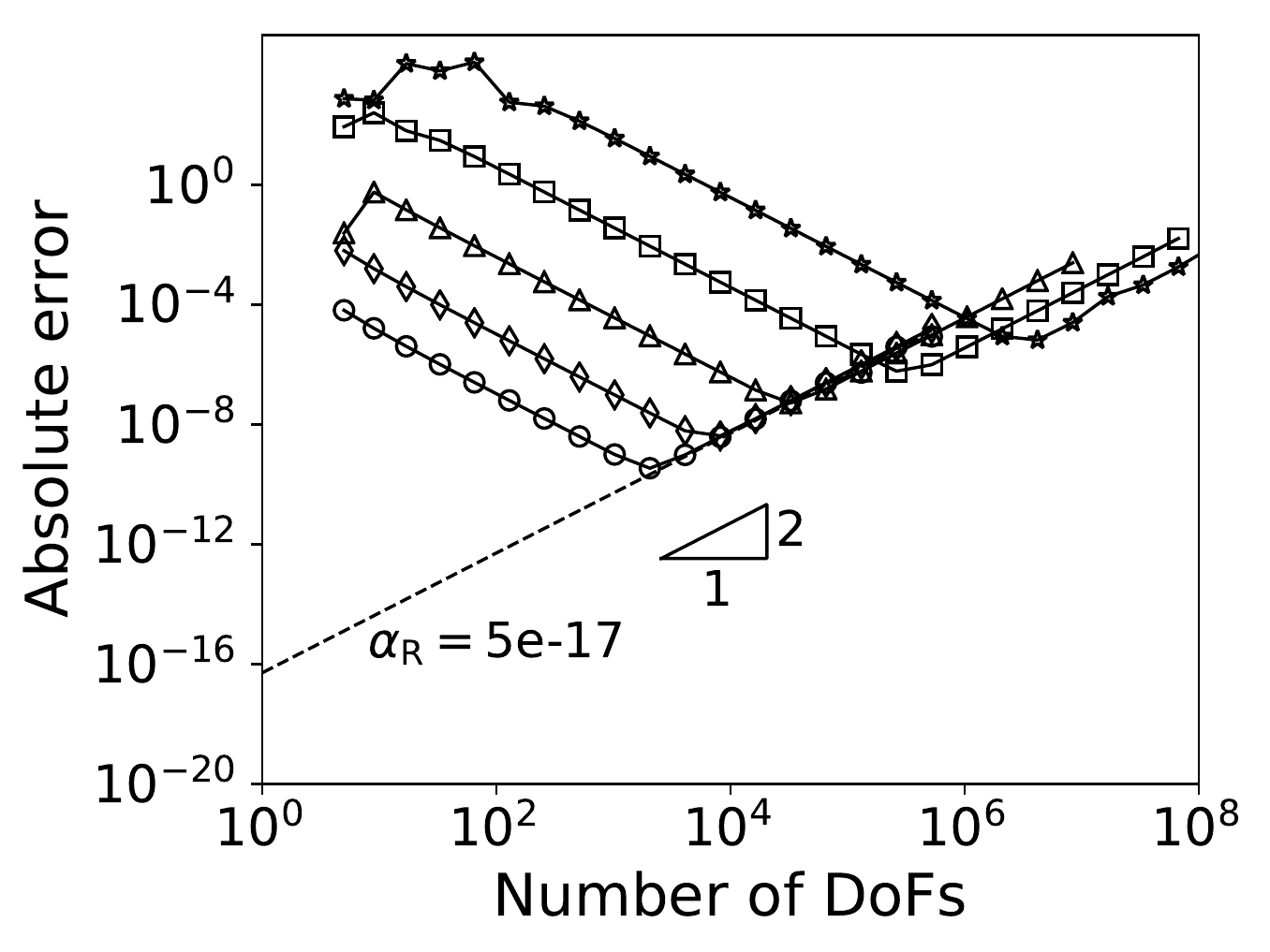}
        \caption{First derivative}
        \label{py_L2_Pois1_SM_scaling_S_grad}
    \end{subfigure}
    \hspace{-0.2cm}
    \begin{subfigure}{5.5cm}
        \includegraphics[width=1.0\linewidth]{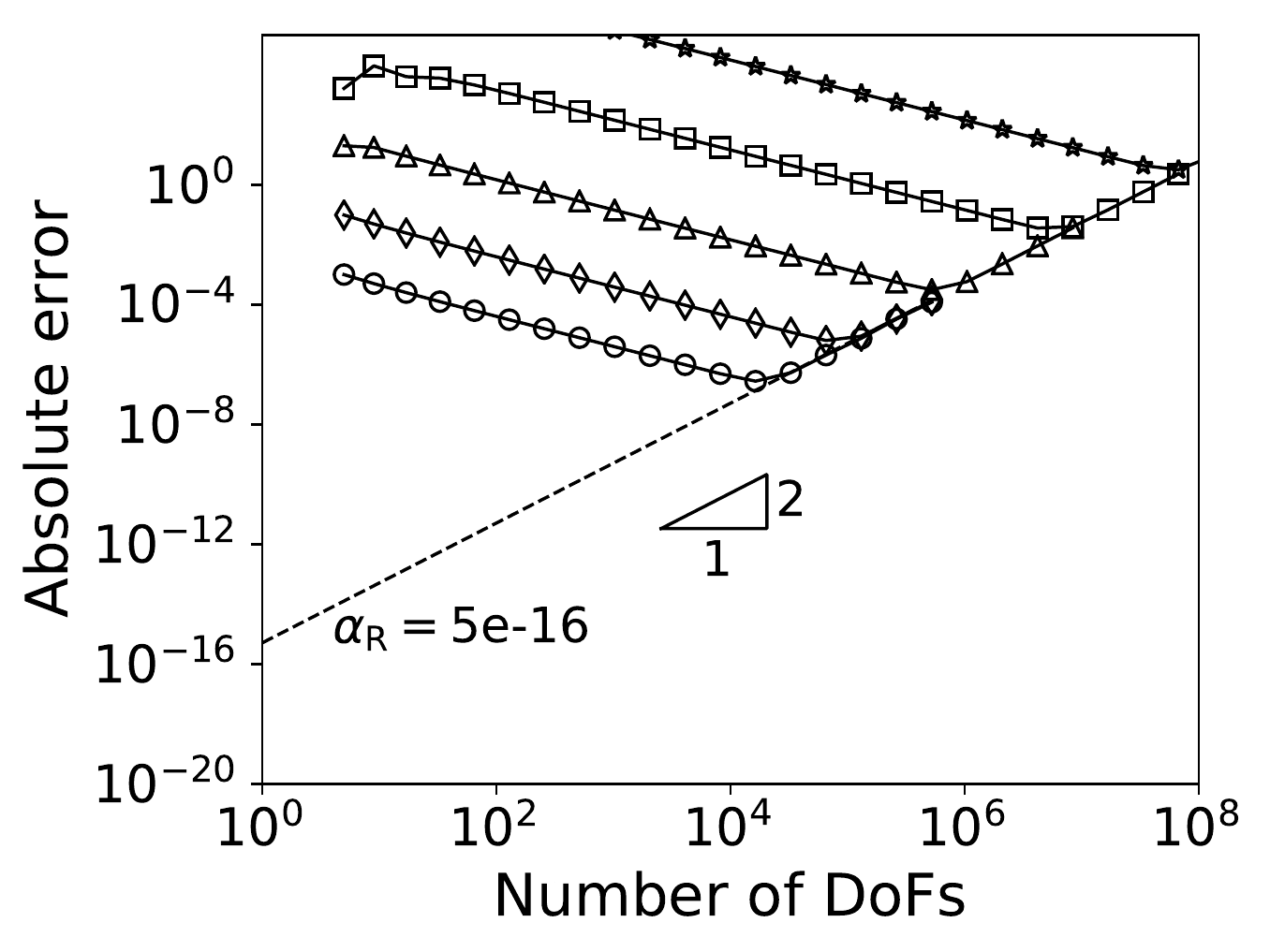}
        \caption{Second derivative}
        \label{py_L2_Pois1_SM_scaling_S_2ndd}
    \end{subfigure}
\caption{Absolute errors for different $c_1$ using the standard FEM with scheme $S$.}    
\label{py_L2_Pois1_SM_scaling_S}
\end{figure}

The absolute errors for $u$, $u_{x}$ and $u_{xx}$ for different values of $c_1$ using the standard FEM are depicted in Fig.~\ref{Pois_SM_rhs_solu orig}. It shows that, for all the three variables, the offsets $\alpha_{\rm R}$ increase with increasing $\|u\|_2$ (decreasing $c_1$), which makes it impossible to determine the break-even point between truncation and round-off error in a generic, that is, problem independent way. 

This is because the number of accurate significant digits that the double-precision floating-point format can hold is 17 at most, and hence, more significant digits in the fractional part will be rounded with increasing $\|u\|_2$. To eliminate this influence factor, we scale the $L_2$ norm of $u$ to 1, which is achieved by dividing the right-hand side $F$ of the linear system of equations (\ref{matrix equation std FEM}) by $\|u\|_{2}$.
The scaling scheme can be found in the second row of Table \ref{scaling schemes std and mix FEM}, which is denoted as $S$.
Note that, the scaling factor is approximated from the numerical solution through an a posteriori algorithm presented in Section \ref{section_algorithm}.

Using scheme $S$ for Case 1, the absolute errors are depicted in Fig.~\ref{py_L2_Pois1_SM_scaling_S}. It shows that $\alpha_{\rm R}$ for different $c_1$ converge to common values, which are $2\times10^{-17}$, $5\times10^{-17}$ and $5\times10^{-16}$ for $u$, $u_{x}$ and $u_{xx}$, respectively. These values also apply to Cases 2$\sim$5 when using scheme $S$. 


%

\begin{table}[!ht]
\centering
\caption [sss] {Scaling schemes.}
\label{scaling schemes std and mix FEM} 
\begin{tabular}{c c c c}
\hline  
{Scheme}& Left-hand side & Solution & Right-hand side \\	\hline
$S$ & {$A$} & $\frac{1}{\|u\|_{2}} U$ & $\frac{1}{\|u\|_{2}} F$ \\	\hline
$M_1$ & {$\left[ \begin{array}{cc} M & \frac{\|u\|_{2}}{\|v\|_{2}} B  \\ B^T & 0 \end{array}\right]$ } & $\left[ \begin{array}{cc} \frac{1}{\|v\|_{2}} {V} \\ \frac{1}{\|u\|_{2}} {U} \end{array}\right]$ & $\left[ \begin{array}{cc} \frac{1}{\|v\|_{2}} G \\ {\frac{1}{\|v\|_{2}} H} \end{array}\right]$\\	\hline
$M_2$ & {$\left[ \begin{array}{cc} M & B  \\ B^T & 0 \end{array}\right]$ } & $\frac{1}{\|u\|_{2}} \left[ \begin{array}{cc} {V} \\ {U} \end{array}\right]$ & $\frac{1}{\|u\|_{2}} \left[ \begin{array}{cc}  G \\ H \end{array}\right]$ \\	\hline
\end{tabular}
\end{table}

\paragraph{The mixed FEM}			\label{scaling_mix_FEM}

\begin{figure}[!ht]
    \begin{subfigure}{5.5cm}
        \includegraphics[width=1.0\linewidth]{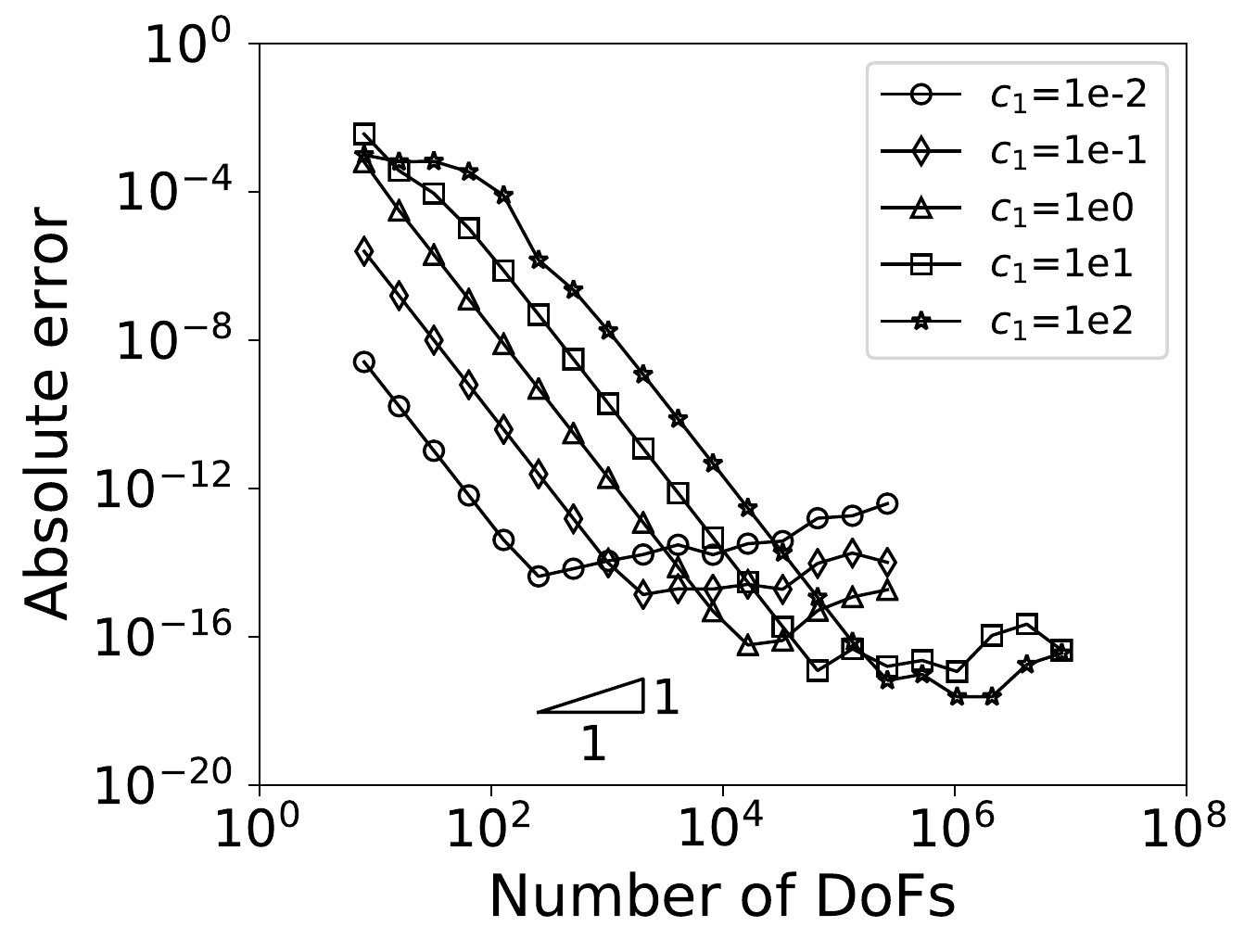}
        \caption{Solution}
        \label{py_L2_Pois1_MM_scaling_no_solu}
    \end{subfigure}
    \hspace{-0.2cm}
    \begin{subfigure}{5.5cm}
        \includegraphics[width=1.0\linewidth]{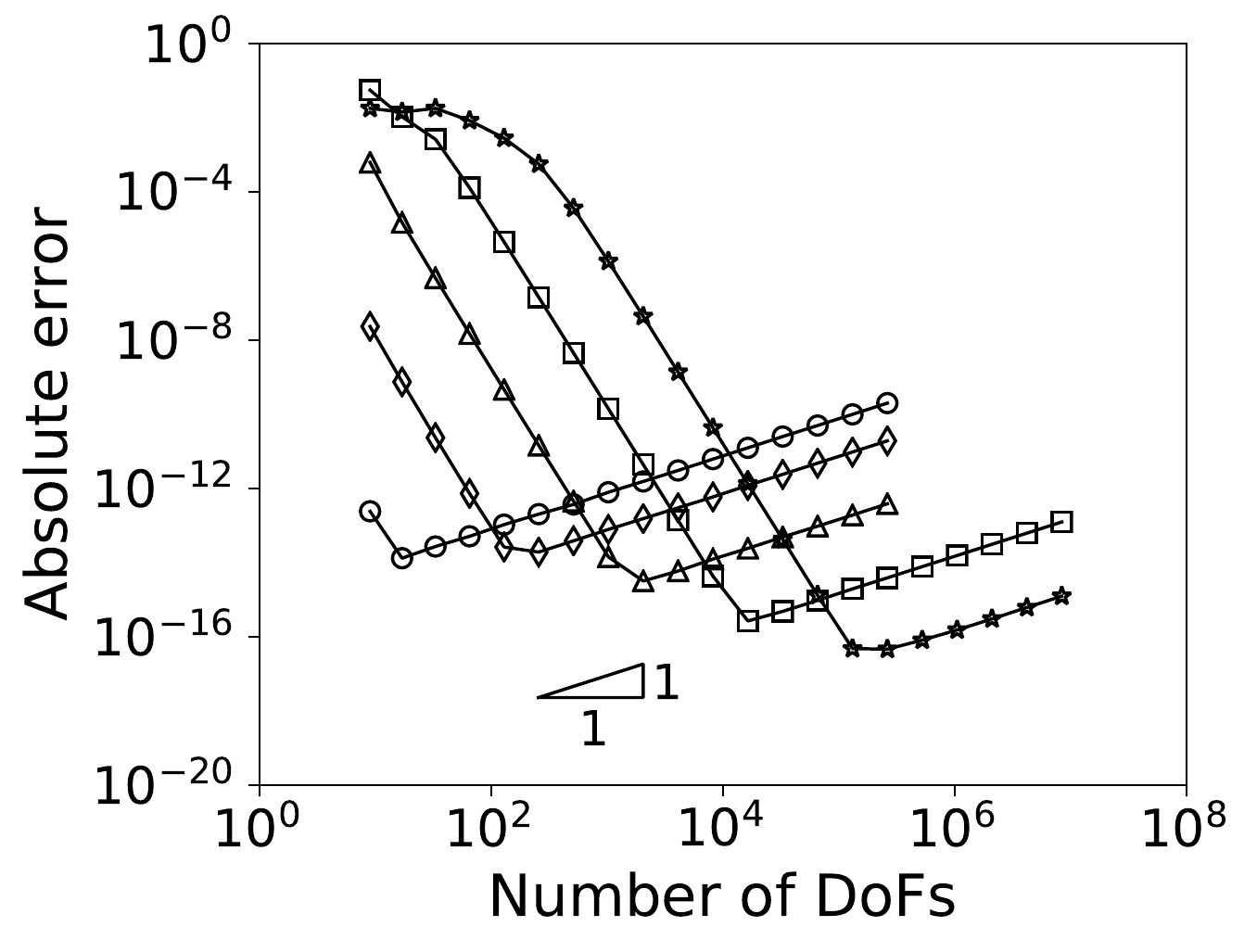}
        \caption{First derivative}
        \label{py_L2_Pois1_MM_scaling_no_grad}
    \end{subfigure}
    \hspace{-0.2cm}
    \begin{subfigure}{5.5cm}
        \includegraphics[width=1.0\linewidth]{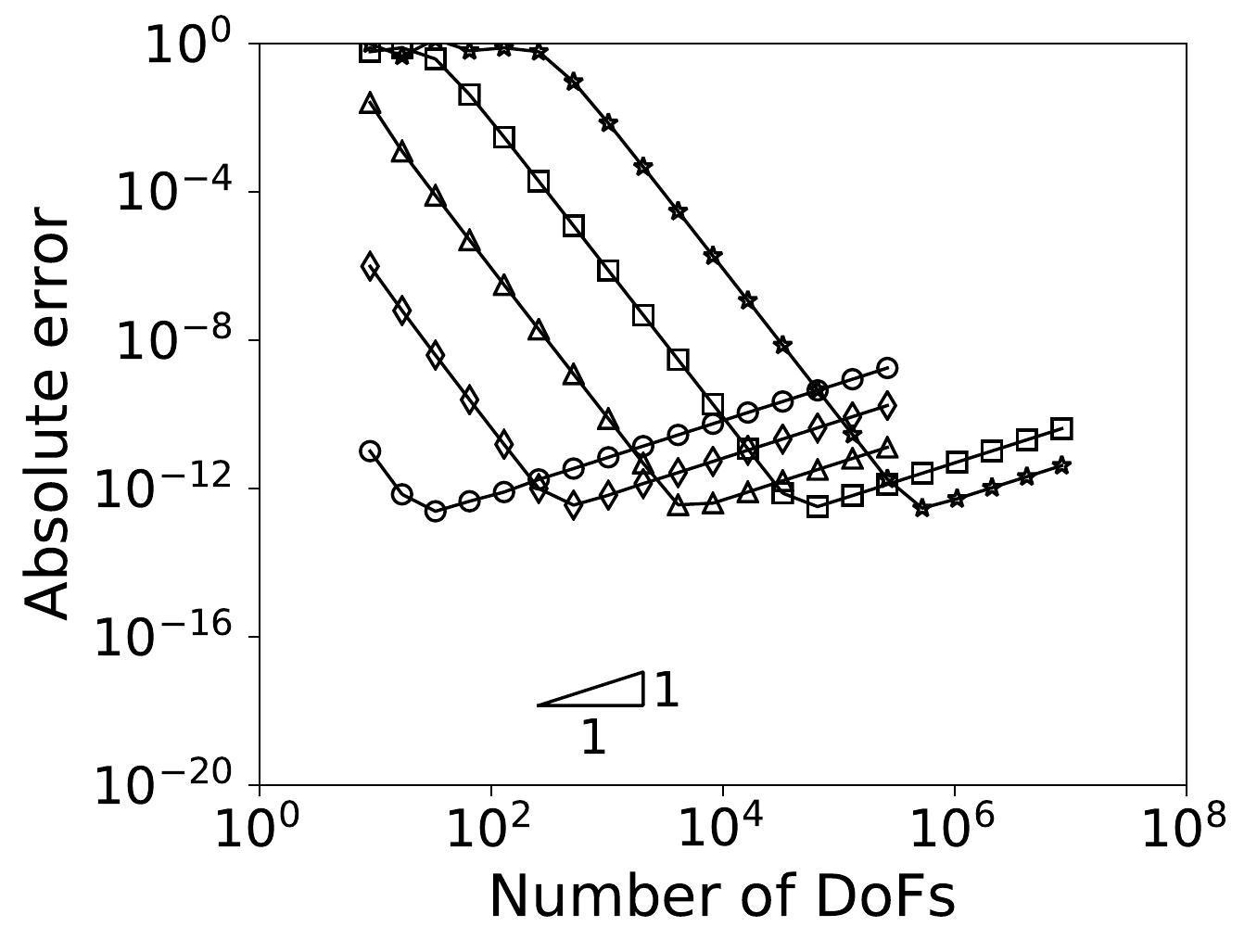}
        \caption{Second derivative}
        \label{py_L2_Pois1_MM_scaling_no_2ndd}
    \end{subfigure}
\caption{Absolute errors for different $c_1$ using the mixed FEM without scaling the right-hand side.}    
\label{Pois_pLov2pi2sin_MM_rhs orig}
\end{figure}

The outcome of the numerical experiments performed with the mixed-FEM formulation Eq. (\ref{matrix equation mix FEM}) are presented in Fig.~\ref{Pois_pLov2pi2sin_MM_rhs orig}.
Like with the standard FEM, the offsets $\alpha_{\rm R}$ for $u$ and $u_x$ increase whenever $\|u\|_2$ and $\|u_x\|_2$ are increased.

\begin{figure}[!ht]
    \begin{subfigure}{5.5cm}
        \includegraphics[width=1.0\linewidth]{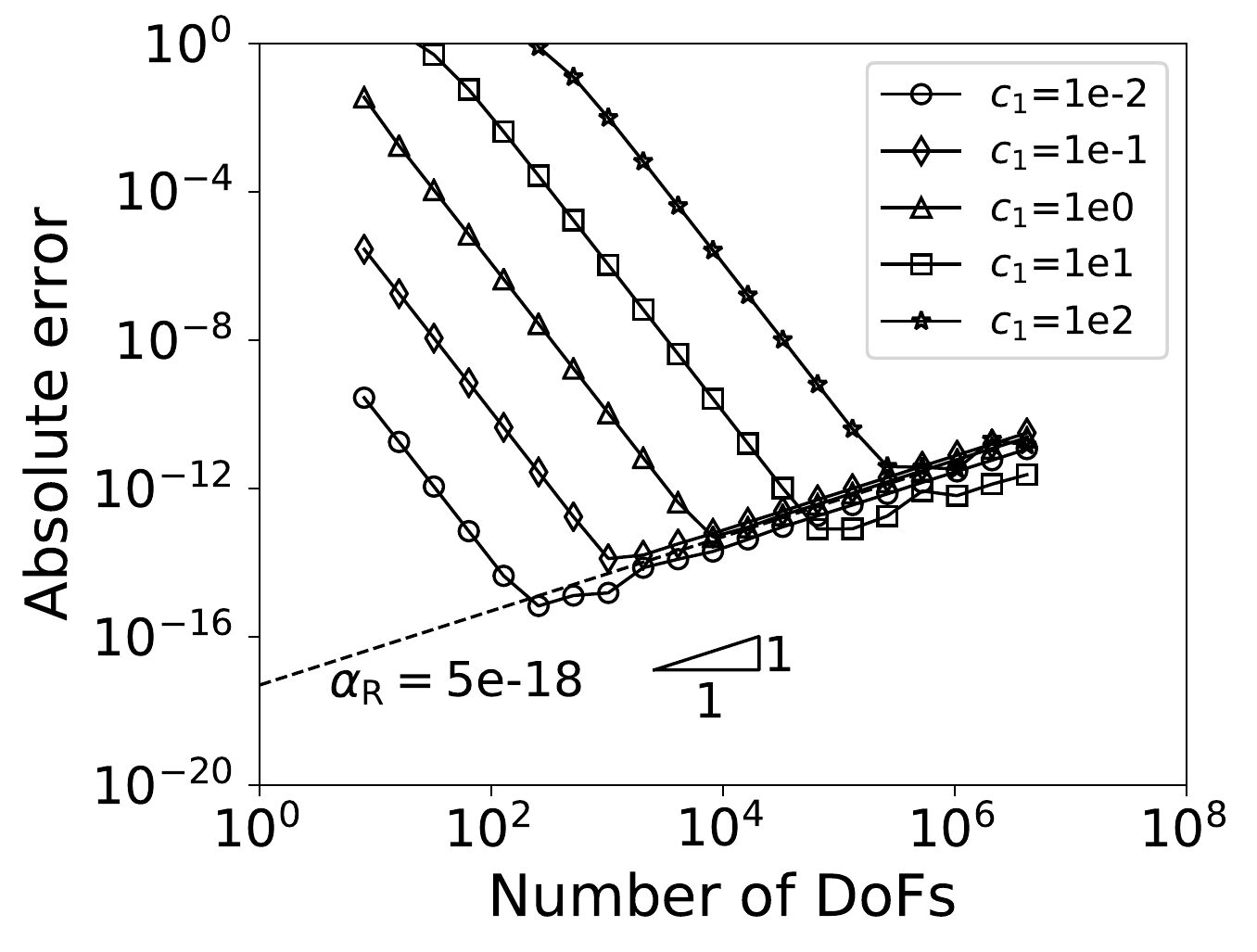}
        \caption{Solution}
        \label{py_L2_Pois1_MM_scaling_M1_solu}
    \end{subfigure}
    \hspace{-0.2cm}
    \begin{subfigure}{5.5cm}
        \includegraphics[width=1.0\linewidth]{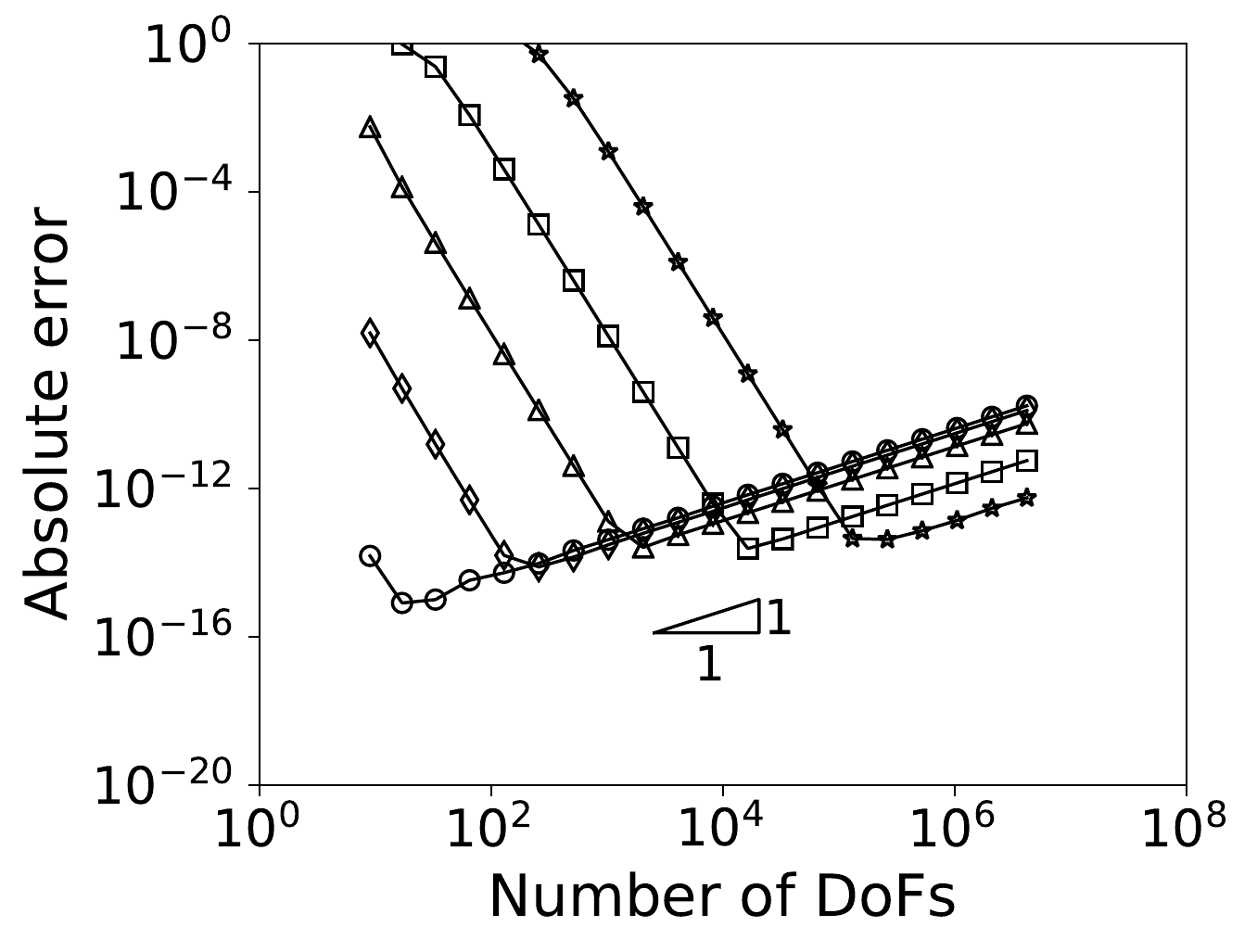}
        \caption{First derivative}
        \label{py_L2_Pois1_MM_scaling_M1_grad}
    \end{subfigure}
    \hspace{-0.2cm}
    \begin{subfigure}{5.5cm}
        \includegraphics[width=1.0\linewidth]{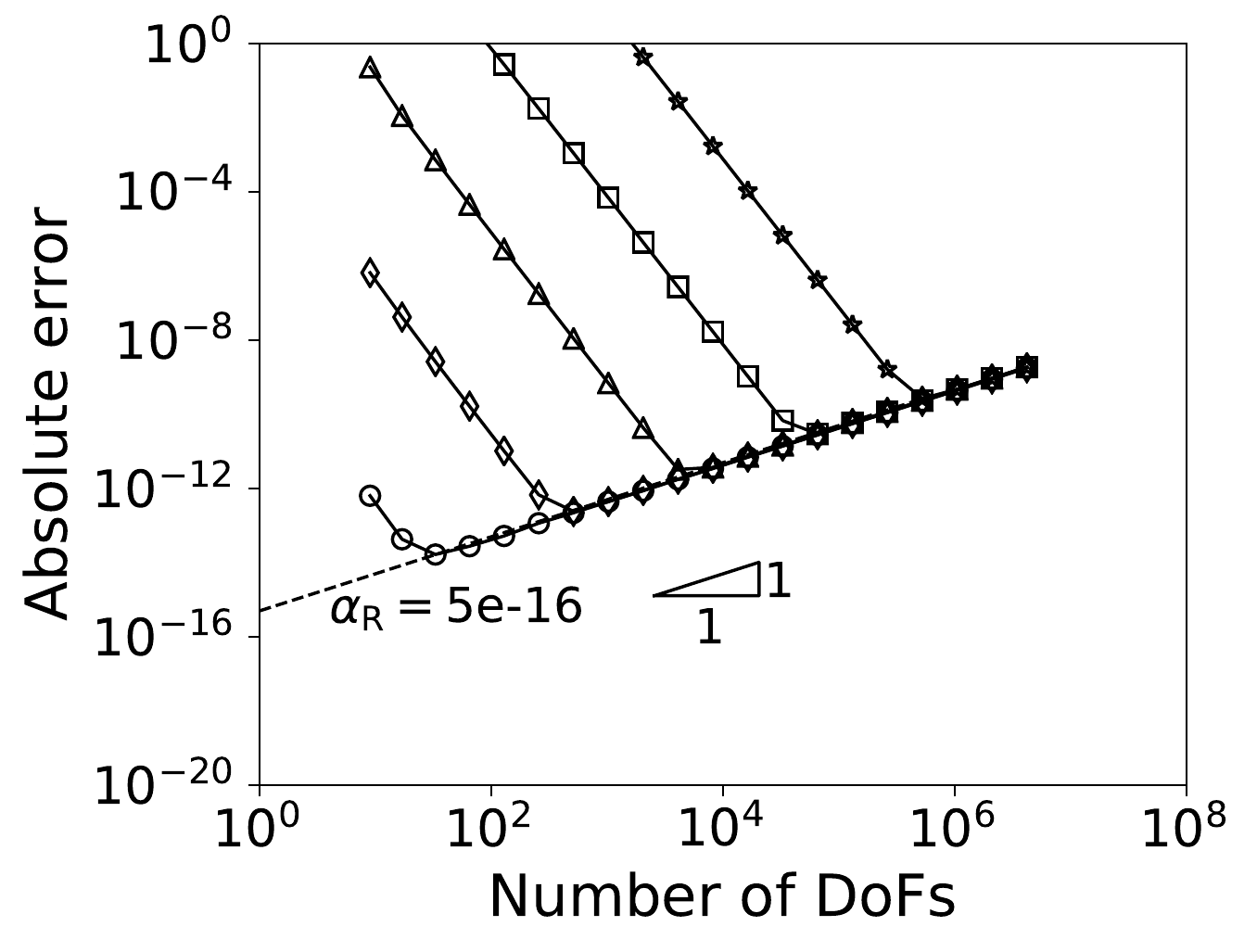}
        \caption{Second derivative}
        \label{py_L2_Pois1_MM_scaling_M1_2ndd}
    \end{subfigure}
\caption{Absolute errors for different $c_1$ using the mixed FEM with scheme $M_1$.}    
\label{Pois_pLov2pi2sin_MM_rhs M1}
\end{figure}

Instinctively, to mitigate the influence of the magnitude of the solution $u$ and the first derivative $v$ on the offset $\alpha_{\rm R}$, one would scale the $L_2$ norm of $u$ and $v$ to 1. This can be achieved by dividing the right-hand sides $G$ and $H$ by the $L_2$ norm of the first derivative $\|v\|_{2}$ and multiplying the discrete first derivative operator $B$ by $\frac{\|u\|_{2}}{\|v\|_{2}}$, see scheme $M_1$ shown in Table \ref{scaling schemes std and mix FEM}.
Using this scheme, the absolute errors of $u$, $v$ and $v_{x}$ are shown in Fig.~\ref{Pois_pLov2pi2sin_MM_rhs M1}.
As expected, the offsets $\alpha_{\rm R}$ for $u$ and $v_{x}$ converge, but that for $v$ only converge when $c_1<1$.
For $c_1>1$, no convergence of $\alpha_{\rm R}$ is seen for $v$. It further indicates that we need smaller scaling factor for $v$ when $c_1>1$.

\vspace{0.0cm}

\begin{figure}[!ht]
    \begin{subfigure}{5.5cm}
        \includegraphics[width=1.0\linewidth]{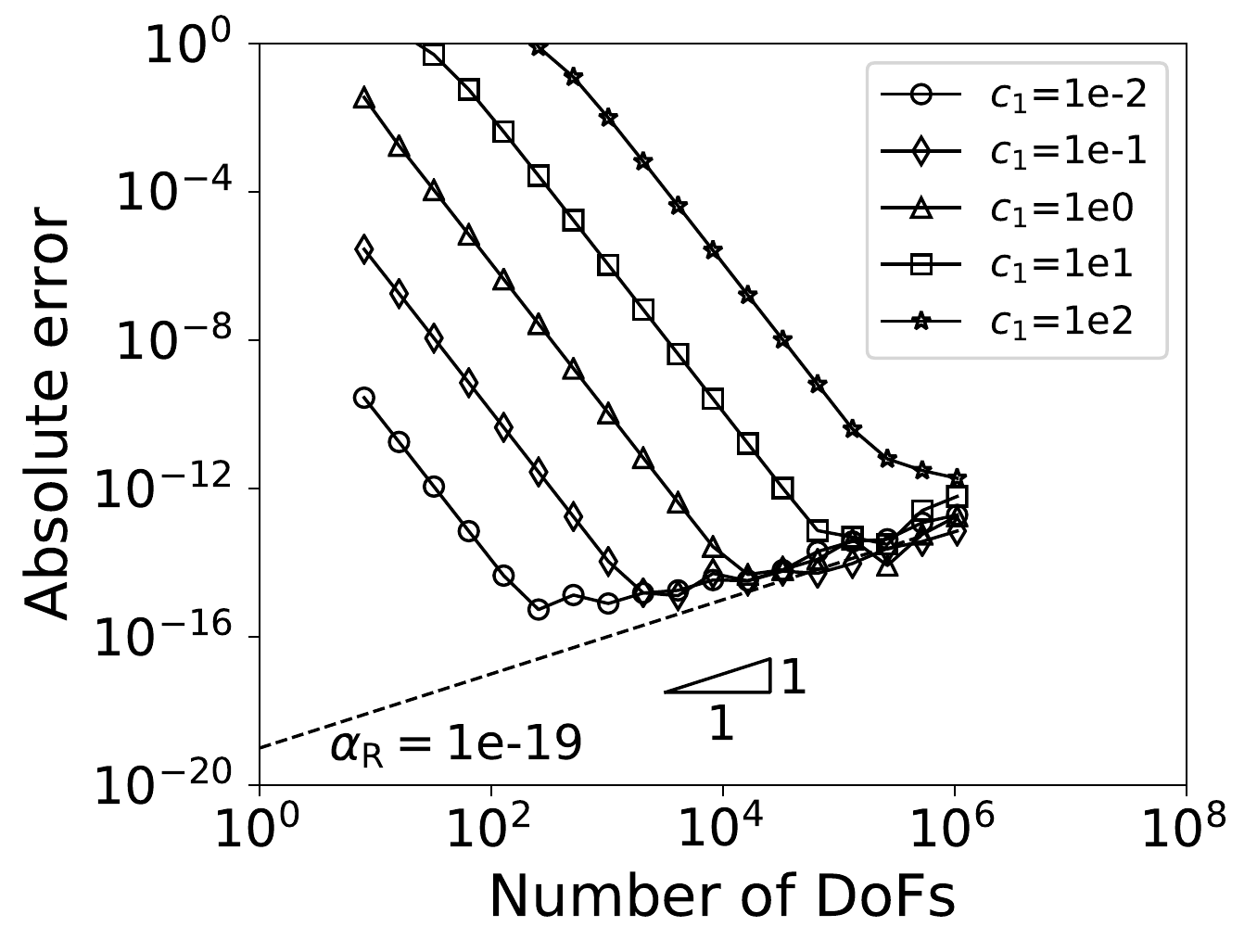}
        \caption{Solution}
        \label{py_L2_Pois1_MM_scaling_M2_solu}
    \end{subfigure}
    \hspace{-0.2cm}
    \begin{subfigure}{5.5cm}
        \includegraphics[width=1.0\linewidth]{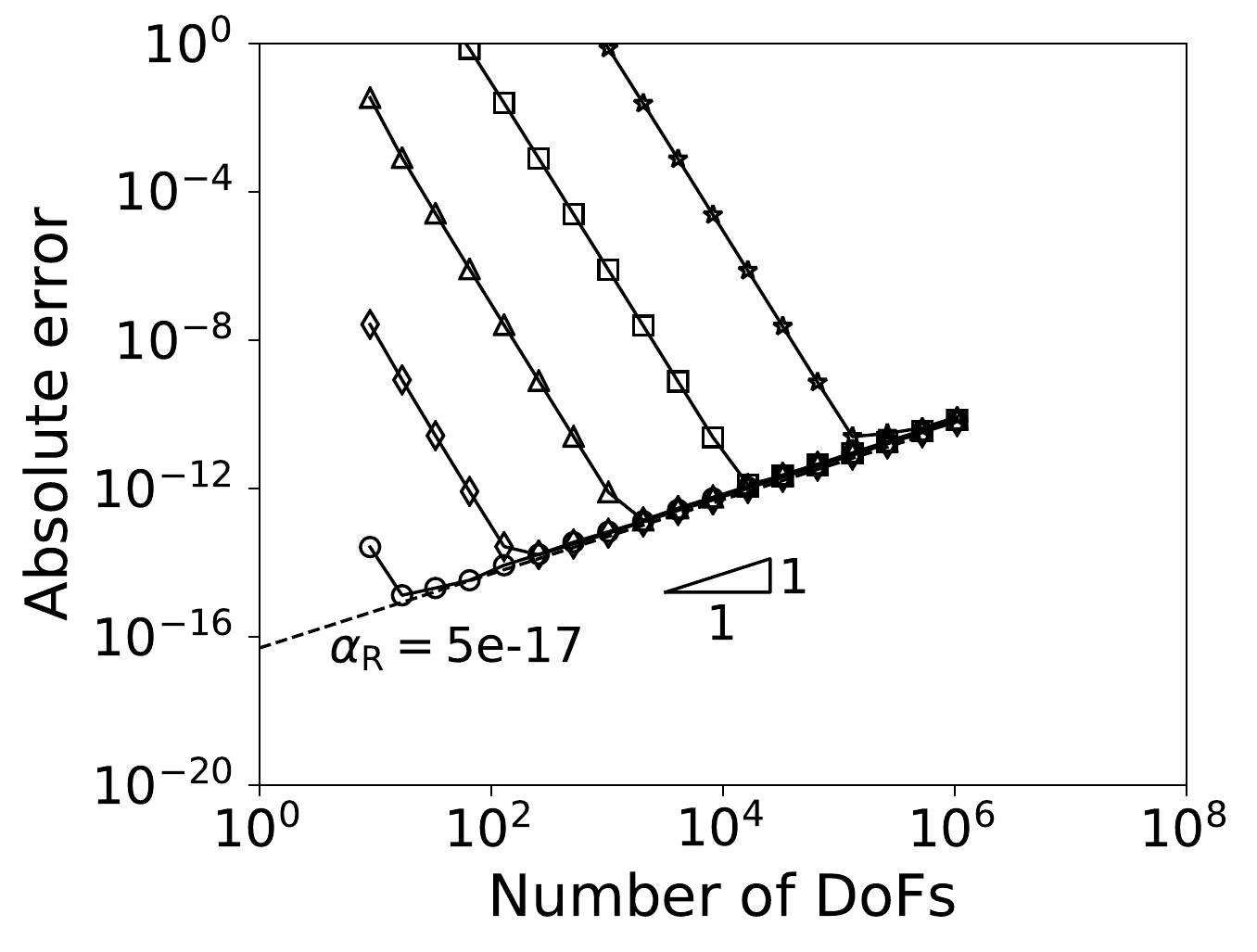}
        \caption{First derivative}
        \label{py_L2_Pois1_MM_scaling_M2_grad}
    \end{subfigure}
    \hspace{-0.2cm}
    \begin{subfigure}{5.5cm}
        \includegraphics[width=1.0\linewidth]{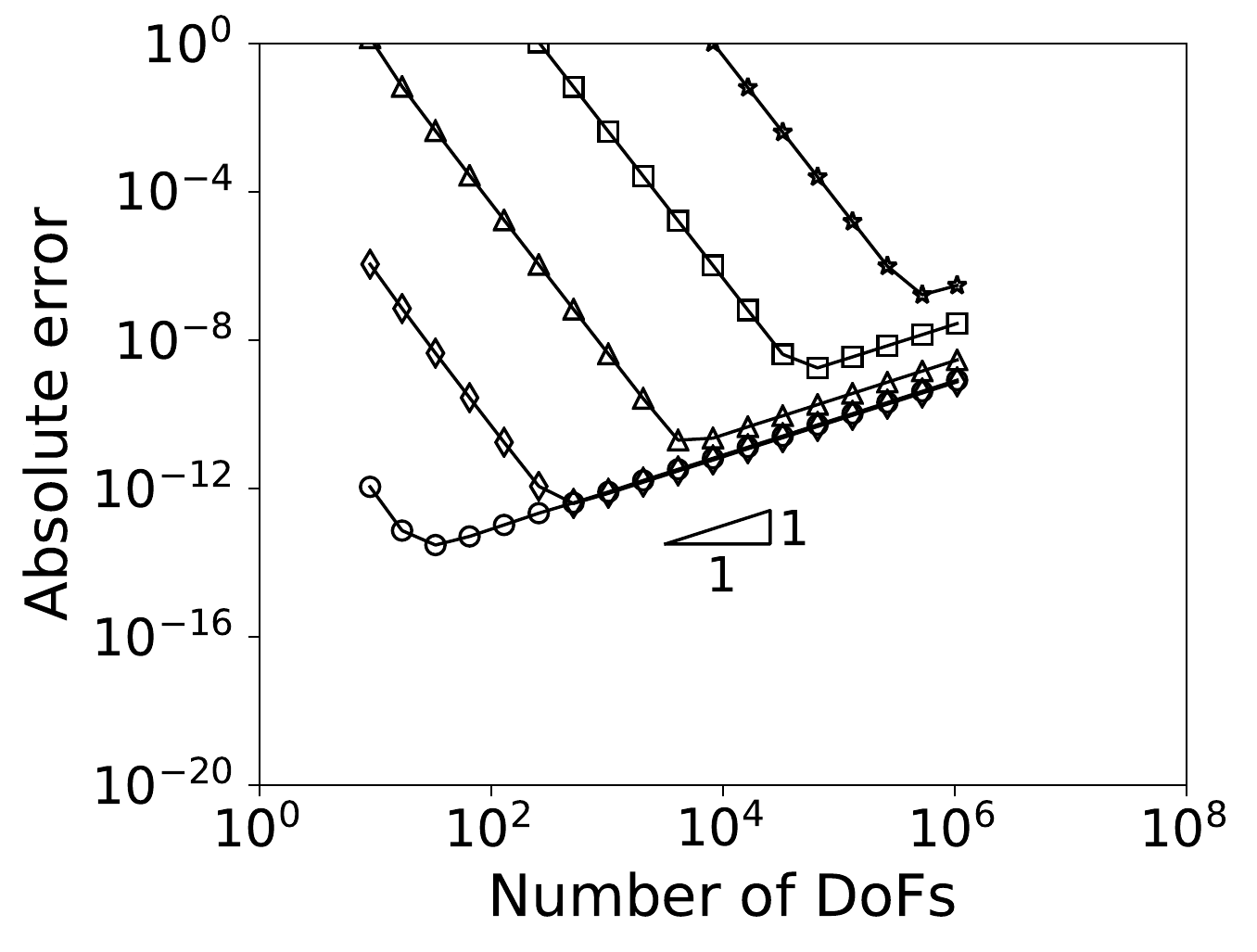}
        \caption{Second derivative}
        \label{py_L2_Pois1_MM_scaling_M2_2ndd}
    \end{subfigure}
\caption{Absolute errors for different $c_1$ using the mixed FEM with scheme $M_2$.}
\label{Pois_pLov2pi2sin_MM_rhs M2}
\end{figure}

\newpage
Given that ${\|u\|_{2}}$ is of the same order with ${\|v\|_{2}}$ when $c_1<1$, while it is smaller than ${\|v\|_{2}}$ when $c_1>1$, we scale both $u$ and $v$ by $\|u\|_{2}$. This scaling scheme, which divides both the right-hand sides $G$ and $H$ by $\|u\|_{2}$, is denoted as scheme $M_2$ as shown in Table \ref{scaling schemes std and mix FEM}. The absolute errors obtained by using this scheme are shown in Fig.~\ref{Pois_pLov2pi2sin_MM_rhs M2}, where the offsets $\alpha_{\rm R}$ for both $u$ and $v$ converge. However, not for $v_{x}$.

Therefore, scheme $M_2$ is preferable if $u$ and $v$ are of primary interest, and scheme $M_1$ is more suitable when $v_{x}$ is of interest. 
If all three quantities need to be computed with required accuracy, both schemes $M_1$ and $M_2$ need to be applied side by side. 
The generalized values of $\alpha_{\rm R}$ using the mixed FEM are $1\times10^{-19}$, $5\times10^{-17}$ and $5\times10^{-16}$ for $u$, $v$ and $v_{x}$, respectively. 
These two scaling schemes also work for Cases 2$\sim$5, but the resulting $\alpha_{\rm R}$ are slightly different. After being amended by Cases 2$\sim$5, $\alpha_{\rm R}$ become $1\times10^{-18}$, $1\times10^{-16}$ and $5\times10^{-16}$.

Generalizing $\alpha_{\rm R}$ using both the standard FEM and the mixed FEM for the Poisson equations, we obtain $2.0\times10^{-17}$, $5.0\times10^{-17}$ and $5.0\times10^{-16}$ for $u$, $u_{x}$ and $u_{xx}$, respectively. Using the above scaling schemes for the benchmark diffusion and Helmholtz equations, we obtain values $2.0\times10^{-17}$, $2.0\times10^{-17}$ and $1.0\times10^{-15}$.
Generalizing these two sets of values, we obtain $\alpha_{\rm R}$ shown in Table \ref{value_alpha_R_generalized}.

\begin{table}[!ht]
\captionof{table}{Generalized values of $\alpha_{\rm R}$ for Eq. (\ref{1D_general_Helmholtz_equation}).}
\centering
\begin{tabular}{c c c c}
\hline
 & $u$ & $u_x$ & $u_{xx}$ \\ \hline 
$\alpha_{\rm R}$ & 2e-17 & 5e-17 & 1e-15 \\ \hline
\end{tabular}
\label{value_alpha_R_generalized}
\end{table}


Summarizing this section, to mitigate the influence of the order of magnitude of the different variables on $\alpha_{\rm R}$, we have proposed and validated three different scaling schemes $S$, $M_1$ and $M_2$, resulting in the common values for $\alpha_{\rm R}$.
This is an essential prerequisite for our a posteriori refinement strategy to be robust and generally applicable.

\subsubsection{Boundary conditions}	\label{section_BC}

In this section, two aspects of the influence of the boundary conditions on the round-off error are investigated: first the method of implementing the Dirichlet boundary conditions, and secondly types of boundary conditions. 

\begin{figure}[!ht]
    \begin{subfigure}{5.5cm}
        \includegraphics[width=1.0\linewidth]{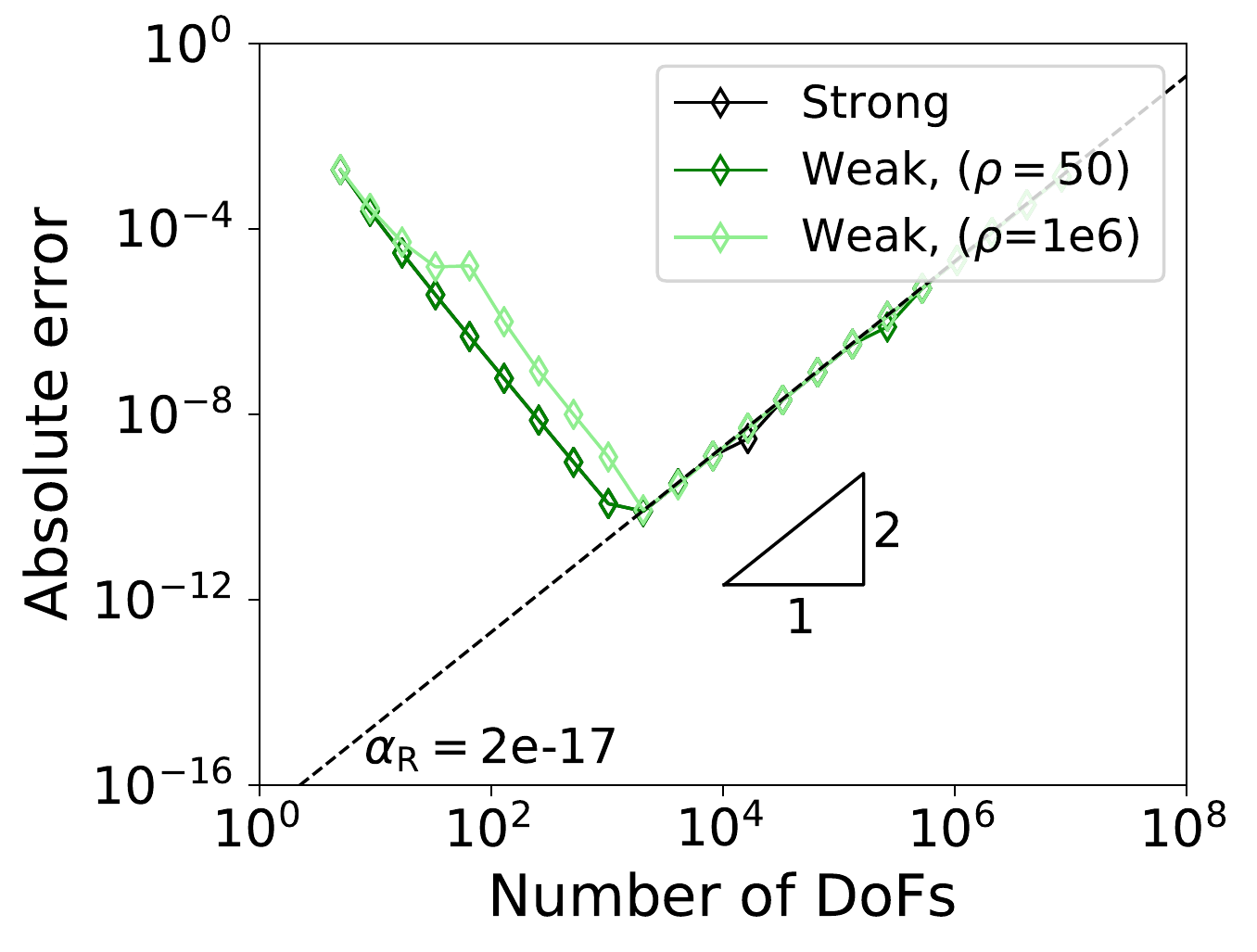}
        \caption{Solution}
        \label{py_bench_Pois_SM_error_boundary_weak_solu}
    \end{subfigure}
    \hspace{-0.2cm}
    \begin{subfigure}{5.5cm}
        \includegraphics[width=1.0\linewidth]{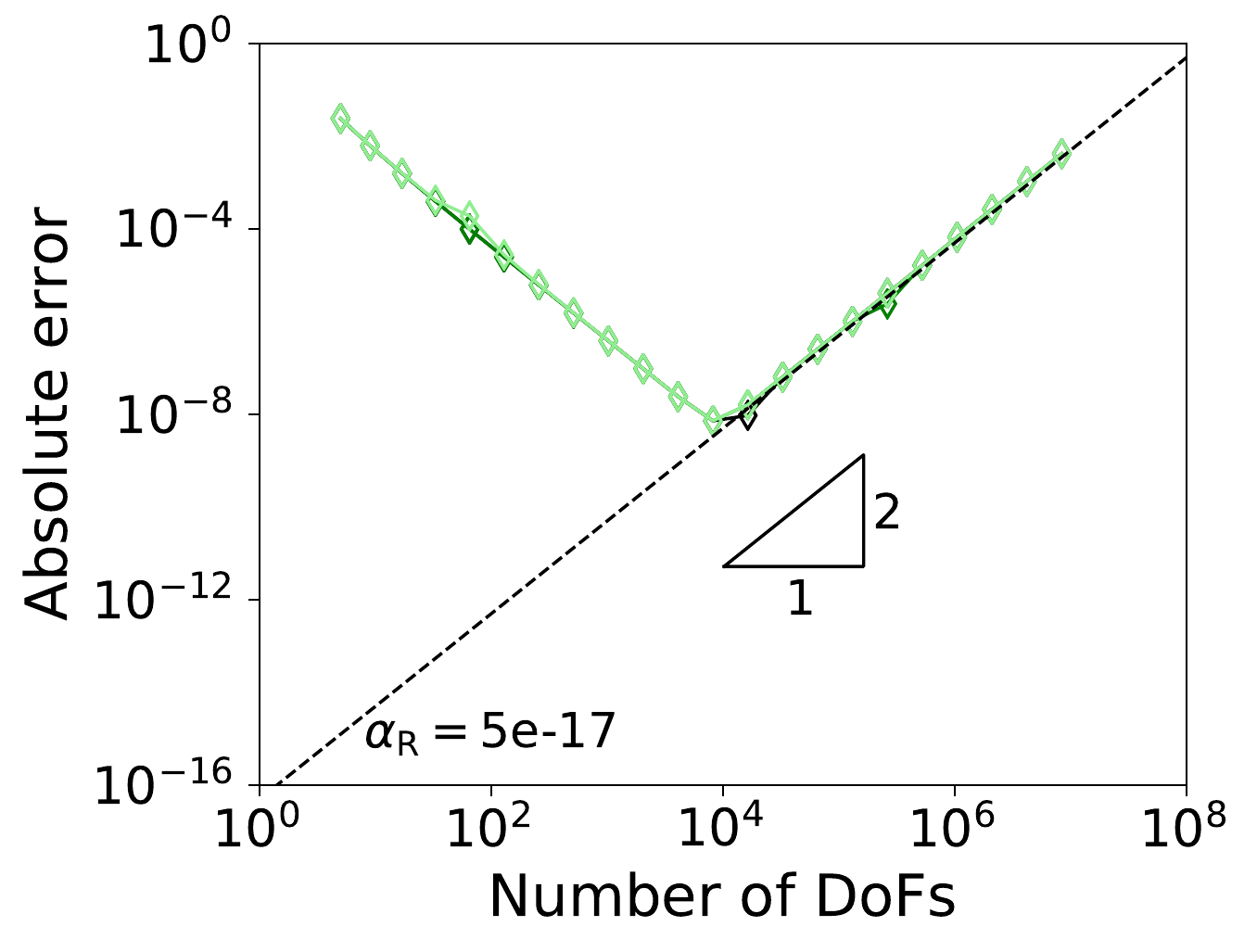}
        \caption{First derivative}
        \label{py_bench_Pois_SM_error_boundary_weak_grad}
    \end{subfigure}
    \hspace{-0.2cm}
    \begin{subfigure}{5.5cm}
        \includegraphics[width=1.0\linewidth]{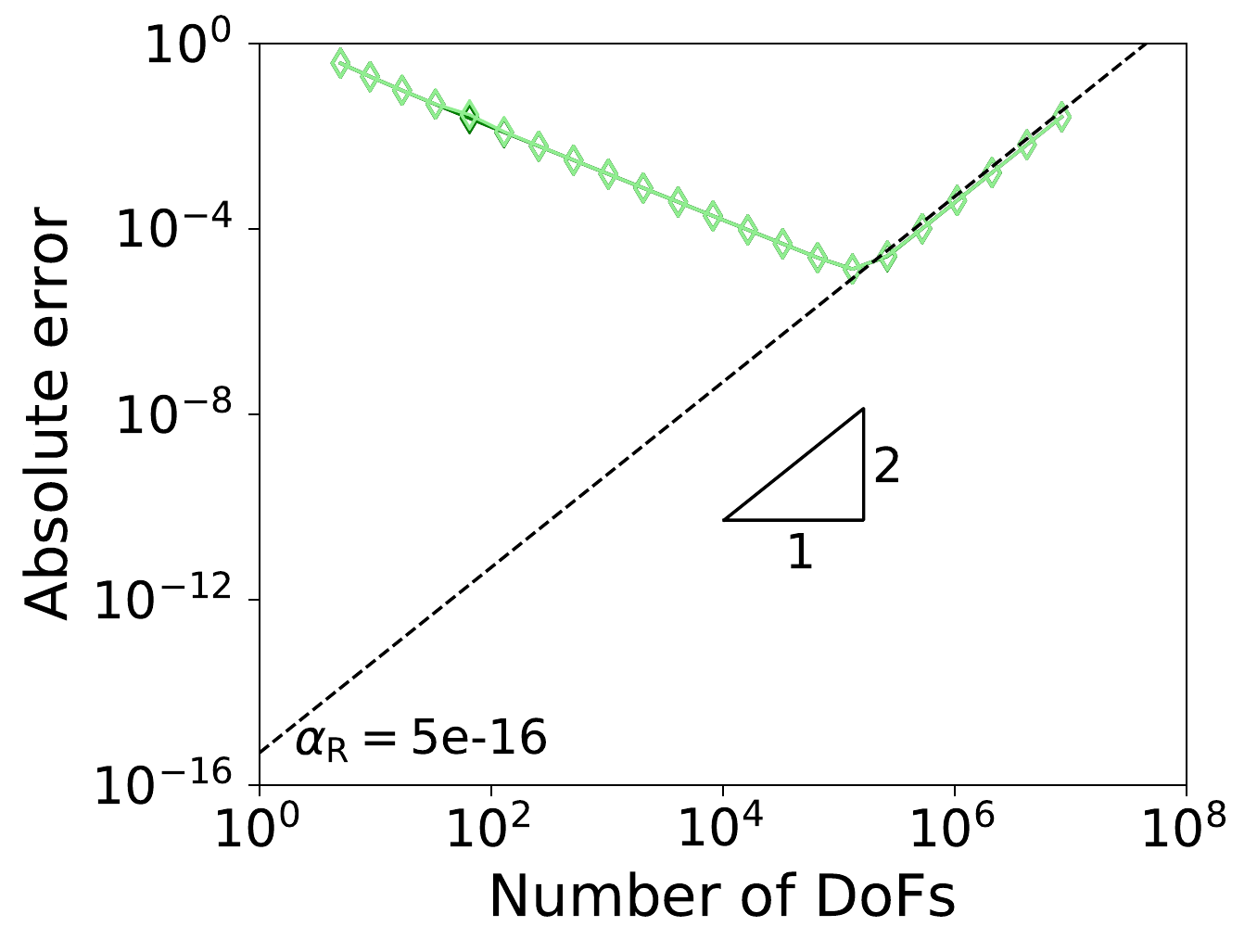}
        \caption{Second derivative}
        \label{py_bench_Pois_SM_error_boundary_weak_2ndd}
    \end{subfigure}
\caption{Comparison of the errors for imposing the Dirichlet boundary condition strongly and weakly.}
\label{py_bench_Pois_SM_error_boundary_weak}
\end{figure}

For the first aspect, using Weak form 2 for $\rho=50$ and $10^6$, the discretization errors are depicted in Fig.~\ref{py_bench_Pois_SM_error_boundary_weak}, in comparison with that using Weak form 1. As can be seen, both weak and strong imposition of the Dirichlet boundary condition yield the same trend line for the round-off error for the solution and its derivatives, and the magnitude of the penalty parameter in the weak imposition makes no difference. In addition, small penalty parameters might lead to larger truncation errors for $u$, but the difference diminishes when the penalty parameter is large enough.

To construct the problem for the second aspect, the Dirichlet boundary condition at the left boundary ($x=0$) is kept while the Dirichlet boundary condition at the right boundary ($x=1$) has been replaced by the Neumann boundary condition $u_x (1) = -e^{-1/4}$, leading to the same solution and derivative profiles.

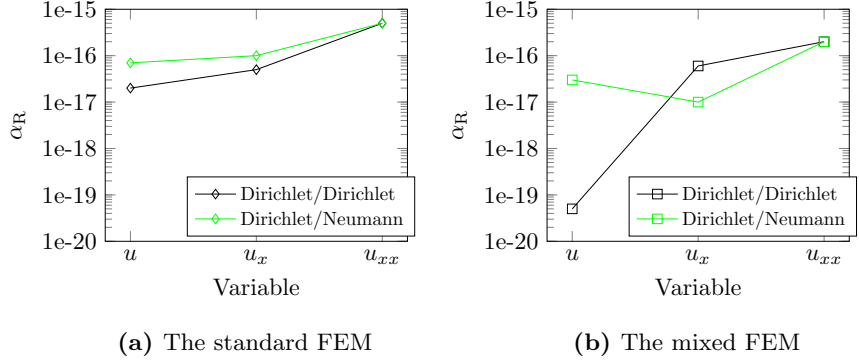
\begin{figure}[!ht]
\hspace{2.2cm}
\begin{subfigure}[b]{0.4\textwidth}
\scalebox{0.9}{
\begin{tikzpicture} 
\begin{axis}
[
    ymode=log,    
    ymin=1e-20,
    ymax=1e-15,
    ytick={1e-20,1e-19,1e-18,1e-17,1e-16,1e-15},
    yticklabels={1e-20,1e-19,1e-18,1e-17,1e-16,1e-15},      
    legend style={nodes={scale=0.8},at={(0.27,0.15)},anchor=west},
    legend cell align={left},
    height=5cm,
    width=6cm,
    ylabel={$\alpha_{\rm R}$},
    ylabel style={at={(-0.01,0.5)}},    
    xtick={0,1,2,3,4},
    xticklabels={$u$,$u_x$, $u_{xx}$, $4$, ${5}$},
    xlabel={Variable},
    xlabel style={at={(0.5,0.03)}},    
]
\addplot[black,mark=diamond,mark options={color=black,fill=black}] coordinates {(0,2.0e-17) (1,5.0e-17) (2,5.0e-16)};
\addplot[green,mark=diamond,mark options={color=green,fill=green}] coordinates {(0,7.0e-17) (1,1.0e-16) (2,5.0e-16)};
\legend{Dirichlet/Dirichlet, Dirichlet/Neumann};
\end{axis}
\end{tikzpicture}
}
\caption{The standard FEM}
\label{boundary_type_benchmark_Poisson_std}
\end{subfigure}
\hspace{-1.0cm}
\begin{subfigure}[b]{0.4\textwidth}
\scalebox{0.9}{
\begin{tikzpicture} 
\begin{axis}
[
    ymode=log,    
    ymin=1e-20,
    ymax=1e-15,
    ytick={1e-20,1e-19,1e-18,1e-17,1e-16,1e-15},
    yticklabels={1e-20,1e-19,1e-18,1e-17,1e-16,1e-15},        
    legend style={nodes={scale=0.8},at={(0.27,0.15)},anchor=west},
    legend cell align={left},
    height=5cm,
    width=6cm,
    ylabel={$\alpha_{\rm R}$},
    ylabel style={at={(-0.01,0.5)}},    
    xtick={0,1,2,3,4},
    xticklabels={$u$,$u_x$, $u_{xx}$, $4$, ${5}$},
    xlabel={Variable},
    xlabel style={at={(0.5,0.03)}},    
]
\addplot[black,mark=square,mark options={color=black,fill=black}] coordinates {(0, 5.0e-20) (1, 6.0e-17) (2, 2.0e-16)};
\addplot[green,mark=square,mark options={color=green,fill=green}] coordinates {(0, 3.0e-17) (1, 1.0e-17) (2, 2.0e-16)};
\legend{Dirichlet/Dirichlet, Dirichlet/Neumann};
\end{axis}
\end{tikzpicture}
}
\caption{The mixed FEM}
\label{boundary_type_benchmark_Poisson_mix}
\end{subfigure}
\caption{Comparison of the errors for imposing Dirichlet/Dirichlet and Dirichlet/Neumann boundary conditions.}
\label{boundary_type_benchmark_Poisson}
\end{figure}

\paragraph{The standard FEM}
Using the standard FEM, the offsets $\alpha_{\rm R}$ for the two types of boundary conditions are depicted in Fig.~\ref{boundary_type_benchmark_Poisson_std}. 
For the Dirichlet/Neumann boundary condition, the offsets $\alpha_{\rm R}$ for $u$ and $u_x$ are slightly larger than that for the Dirichlet/Dirichlet boundary condition by a factor of 3.5 and 2, respectively. The offsets $\alpha_{\rm R}$ for $u_{xx}$ are identical for the two types of boundary conditions.

\paragraph{The mixed FEM}
Using the mixed FEM, the offsets $\alpha_{\rm R}$ for the two types of boundary conditions are depicted in Fig.~\ref{boundary_type_benchmark_Poisson_mix}.
As can be seen, the type of boundary conditions plays a more important role for $\alpha_{\rm R}$ for the solution than $\alpha_{\rm R}$ for other variables.

In summary, $\alpha_{\rm R}$ are relatively independent of the variations in the type of boundary conditions and the method Dirichlet boundary conditions are implemented, which is an important prerequisite for our a posteriori refinement strategy to be applicable for a wide range of problems.

To conclude the sections on sensitivity analysis, the factors that cannot be mitigated are the tolerances for the iterative linear solver, that can be mitigated are the order of magnitude, and that are relatively irrelevant are the boundary conditions.


\section{A posteriori algorithm for finding the optimal number of degrees of freedom}		\label{section_algorithm}

Based on the validation experiments from the previous section, we introduce a novel a posteriori algorithm for determining $E_{\rm min}$ for the solution and its first and second derivative without performing the brute-force mesh refinement. 
Table \ref{settings_algorithm} gives the default settings and the required custom input of the algorithm.

\begin{table}[!ht]
\captionof{table}{Settings of the algorithm.}
\label{settings_algorithm}
  \centering
  \begin{tabular}{lL{5cm}L{6cm}}
    \toprule
    Item & Default & Custom  \\
    \midrule
    Problem & - & \tabitem the differential equation to be solved \\
     &  & \tabitem its associated boundary conditions \\\hline
    Grid & \tabitem initial number of vertices: 2 & - \\
     & \tabitem the vertices are equidistant &  \\\hline
    FEM & \tabitem the maximum $N_h$, denoted by $N_{\rm max}$, : $10^8$ & \tabitem standard or mixed formulation \\
    & \tabitem Dirichlet boundary conditions are imposed strongly & \tabitem an ordered array of element degrees $\{p_{\min}, \ldots, p_{\rm max}\}$ \\\hline
    Computer precision & IEEE-754 double precision & - \\\hline
    Solver & UMFPACK & - \\\hline
    $var$ & - & \tabitem chosen from $\{u,~u_x,~u_{xx}\}$ \\     
    & & \tabitem error tolerance $tol_{var}$ \\     
    \bottomrule
  \end{tabular}
\end{table}

Furthermore, we use the following coefficients in the algorithm:		
\begin{itemize}
  \renewcommand\labelitemi{--}
  \item a minimal number of $h$-refinements before `\textit{NORMALIZATION}' and carrying out `\textit{PREDICTION}', denoted by $REF_{\rm min}$, with the following default values:
  \begin{equation}
  \begin{aligned}
      REF_{\rm min} &=
      \begin{cases*}
	9-p & for p $<$ 6, \\
	4 & otherwise.
      \end{cases*}
  \end{aligned}
  \end{equation} 
  We choose this parameter mainly because the error might increase, or decrease faster than the theoretical order of convergence for coarse refinements, especially for lower-order elements.
  \item a stopping criterion $c_s$ for seeking the scaling factor $\|var_{\rm exc}\|_{2}$ in Table~\ref{scaling schemes std and mix FEM}, its value is 0.001 by default. We choose this parameter because the analytical solution does not exist for most practical problems.
  \item a relaxation coefficient $c_r$ for seeking the theoretical order of convergence, with the following default values: 
    \begin{equation}
    \begin{aligned}
	c_r &=
	\begin{cases*}
	  0.9 & for p $<$ 4, \\
	  0.7 & for 4 $\leqslant$ p $<$ 10, \\
	  0.5 & otherwise.
	\end{cases*}
    \end{aligned}
    \end{equation}
  \item the offset $\alpha _{\rm R}$, see Table \ref{value_alpha_R_generalized} for the default values.
\end{itemize}

The procedure of our algorithm consists of four steps, which are explained below:

\paragraph{Step-1} `\textit{INPUT}'. In this step, the custom input has to be provided.
\paragraph{Step-2} `\textit{NORMALIZATION}'. The function of this step is to find the scaling factor to normalize problems of different orders of magnitude for the variable. The specific procedure can be found in Algorithm \ref{algo_scaling_factor}, where elements of degree $p_{\rm min}$ are used. 

\vspace{0.2cm}
\begin{algorithm}[H]
\caption{NORMALIZATION}
\label{algo_scaling_factor}
\While{$N_h<N_{\rm max}$}
{
    \eIf{$\left|\frac{\|var_{h}\|_{2} - \|var_{2h}\|_{2}}{\|var_{h}\|_{2}} \right| < c_s$}
    {
        $\|var_{\rm exc}\|_{2}$ $\gets$ $\|var_{h}\|_{2}$\;
        break\;
    }
    {
        $h$ $\gets$ $h/2$\;
        calculate $\|var_h\|_{2}$ using Eq. (\ref{formula_abs_error_analytical}) without scaling\;    
    }
}
\end{algorithm}
                                                                   
\paragraph{Step-3} `\textit{PREDICTION}'. This step finds $E_{\rm min}$ for each $var$ and $p$ of interest, as illustrated in Fig.~\ref{sketch_discretization_error_one_p}.
The procedure for carrying out this step can be found in Algorithm \ref{block_PREDICTION}.

\vspace{0.2cm}
\begin{algorithm}[H]
\caption{PREDICTION}			
\label{block_PREDICTION}
    \While{$\widetilde {E_{h}}>E_{\rm R}$ \textbf{\textup{and}} $N_h<N_{\rm max}$}
    {
        $\widetilde{Q}$ $\gets$ $\log _2 \left( {\widetilde {E_{2h}}}/{\widetilde {E_{h}}} \right)$\;
        \eIf
        {
            $\widetilde{Q} \geqslant \beta_{\rm T} \times c_r$
        }
        {
            $N_{\rm c} \gets N_h$\;
            $E_{\rm c} \gets \widetilde {E_{h}}$\;
            $\alpha_{\rm T}$ $\gets$ ${E_{\rm c}}/{N_{\rm c}}^{- \beta_{\rm T}}$\;
            $N_{\rm opt} \gets \left( \frac{\alpha_{\rm T} \beta_{\rm T}}{\alpha _{\rm R} \beta_{\rm R}} \right)^{\frac{1}{\beta_{\rm R} + \beta_{\rm T}}}$\;
            $E_{\rm min} \gets \alpha_{\rm T} {N_{\rm opt}}^{- {\beta _{\rm T}}} + \alpha_{\rm R} {N_{\rm opt}}^{{\beta _{\rm R}}}$\;

        }
        {
            $h$ $\gets$ $h/2$\;
            calculate $\widetilde {E_{h}}$ using Eq.~(\ref{formula_abs_error_numerical}) with proper scaling schemes\;
        }
	}    
\end{algorithm}

\paragraph{Step-4} `\textit{OUTPUT}'. In this step, we output $E_{\rm min}$ obtained from $Step$-3.

\section{Validation}		\label{section_validation}

In what follows, we validate the strategy discussed in Section \ref{section_behaviour_discretization_error_and_prediction} by using the following Helmholtz problem:
\begin{equation}
  \left((0.01+x)(1.01-x) u_x \right)_x -(0.01i) u(x) = 1.0,\qquad x \in I = (0,1),	\label{1D_Helmholtz_equation_application}
\end{equation}
with homogeneous Dirichlet and Neumann boundary conditions imposed as follows: $u(0)=0$ and $u_x(1)=0$.

Both the standard FEM and the mixed FEM are investigated, and the element degree $p$ has a range of $\{1, 2, \ldots, 5\}$. Variables $u$, $u_x$ and $u_{xx}$ are all investigated, for which $tol_{var}$ is set to be $10^{-9}$. 

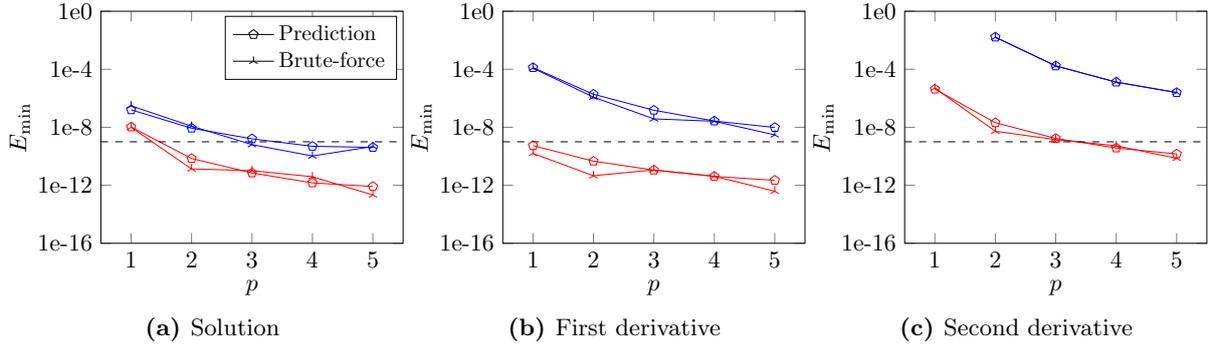
\begin{figure}[!ht]
\hspace{0.0cm}
\begin{subfigure}[b]{0.35\textwidth}
\scalebox{0.9}{
\begin{tikzpicture} 
\begin{axis}
[
    xmin=-0.5,
    xmax=4.5,
    ymode=log,
    ymin=1e-16,
    ymax=1,
    ytick={1e-16, 1e-12, 1e-8, 1e-4, 1e0},
    yticklabels={1e-16, 1e-12, 1e-8, 1e-4, 1e0},    
    legend style={nodes={scale=0.9}},
    legend cell align={left},    
    height=5cm,
    width=6cm,
    ylabel={$E_{\min}$}, 
    ylabel style={at={(0.01,0.5)}},  
    xtick={0,1,2,3,4},
    xticklabels={$1$, $2$, $3$, $4$, ${5}$},
    xlabel={$p$},
    xlabel style={at={(0.5,0.03)}}, 
]
\addplot[black,mark=pentagon,mark options={color=black,fill=black}] coordinates {(0,2.1e-25)};
\addplot[black,mark=Mercedes star,mark options={color=black,fill=black}] coordinates {(0,2.1e-25)};

\addplot[blue,mark=pentagon,mark options={color=blue,fill=blue}] coordinates {(0,1.6e-7) (1,8.7e-9) (2,1.6e-9) (3,4.9e-10) (4,4.1e-10)}; 
\addplot[blue,mark=Mercedes star,mark options={color=blue,fill=blue}] coordinates {(0,3e-7) (1,1.2e-8) (2,6.4e-10) (3,1.05e-10) (4,4.7e-10)}; 

\addplot[red,mark=pentagon,mark options={color=red,fill=red}] coordinates {(0,1.06e-8) (1,7.05e-11) (2,7.13e-12) (3,1.47e-12) (4,8.27e-13)}; 
\addplot[red,mark=Mercedes star,mark options={color=red,fill=red}] coordinates {(0,1.0e-8) (1,1.35e-11) (2,1.02e-11) (3,3.9e-12) (4,2.16e-13)}; 
\addplot[color=black, dashed] coordinates {(-0.5, 1e-9) (4.5, 1e-9)};
\legend{Prediction,Brute-force};
\end{axis}
\end{tikzpicture}
}
\vspace{-0.2cm}
\caption{Solution}
\end{subfigure}
\hspace{-0.7cm}
\begin{subfigure}[b]{0.35\textwidth}
\scalebox{0.9}{
\begin{tikzpicture} 
\begin{axis}
[
    xmin=-0.5,
    xmax=4.5,
    ymode=log,
    ymin=1e-16,
    ymax=1,
    ytick={1e-16, 1e-12, 1e-8, 1e-4, 1e0},
    yticklabels={1e-16, 1e-12, 1e-8, 1e-4, 1e0},
    legend style={nodes={scale=0.8}},
    legend cell align={left},    
    height=5cm,
    width=6cm,
    ylabel={$E_{\min}$},
    ylabel style={at={(0.01,0.5)}},  
    xtick={0,1,2,3,4},
    xticklabels={$1$, $2$, $3$, $4$, ${5}$},
    xlabel={$p$},
    xlabel style={at={(0.5,0.03)}},    
]
\addplot[blue,mark=pentagon,mark options={color=blue,fill=blue}] coordinates {(0,1.3e-4) (1,1.9e-6) (2,1.5e-7) (3,2.7e-8) (4,9.7e-9)}; 
\addplot[blue,mark=Mercedes star,mark options={color=blue,fill=blue}] coordinates {(0,1.2e-4) (1,1.2e-6) (2,3.8e-8) (3,2.6e-8) (4,3.0e-9)}; 

\addplot[red,mark=pentagon,mark options={color=red,fill=red}] coordinates {(0,5.47e-10) (1,4.62e-11) (2,1.12e-11) (3,4.11e-12) (4,2.21e-12)}; 
\addplot[red,mark=Mercedes star,mark options={color=red,fill=red}] coordinates {(0,1.57e-10) (1,4.56e-12) (2,1.15e-11) (3,4.2e-12) (4,3.91e-13)}; 
\addplot[color=black, dashed] coordinates {(-0.5, 1e-9) (4.5, 1e-9)};
\end{axis}
\end{tikzpicture}
}
\vspace{-0.2cm}
\caption{First derivative}
\end{subfigure}
\hspace{-0.7cm}
\begin{subfigure}[b]{0.35\textwidth}
\scalebox{0.9}{
\begin{tikzpicture} 
\begin{axis}
[
    xmin=-0.5,
    xmax=4.5,
    ymode=log,
    ymin=1e-16,
    ymax=1,
    ytick={1e-16, 1e-12, 1e-8, 1e-4, 1e0},
    yticklabels={1e-16, 1e-12, 1e-8, 1e-4, 1e0},
    legend style={nodes={scale=0.8}},
    legend cell align={left},    
    height=5cm,
    width=6cm,
    ylabel={$E_{\min}$},
    ylabel style={at={(0.01,0.5)}},  
    xtick={0,1,2,3,4},
    xticklabels={$1$, $2$, $3$, $4$, ${5}$},
    xlabel={$p$},
    xlabel style={at={(0.5,0.03)}}, 
]
\addplot[blue,mark=pentagon,mark options={color=blue,fill=blue}] coordinates {(0) (1,1.66e-2) (2,1.71e-4) (3,1.3e-5) (4,2.46e-6)}; 
\addplot[blue,mark=Mercedes star,mark options={color=blue,fill=blue}] coordinates {(1,1.6e-2) (2,1.7e-4) (3,1.3e-5) (4,2.5e-6)}; 

\addplot[red,mark=pentagon,mark options={color=red,fill=red}] coordinates {(0,4.25e-6) (1,2.07e-8) (2,1.70e-9) (3,3.75e-10) (4,1.45e-10)}; 
\addplot[red,mark=Mercedes star,mark options={color=red,fill=red}] coordinates {(0,5e-6) (1,5.36e-9) (2,1.41e-9) (3,5.33e-10) (4,7.29e-11)}; 
\addplot[color=black, dashed] coordinates {(-0.5, 1e-9) (4.5, 1e-9)};
\end{axis}
\end{tikzpicture}
}
\vspace{-0.2cm}
\caption{Second derivative}
\end{subfigure}
\caption{Comparison of $E_{\rm min}$ for Eq. (\ref{1D_Helmholtz_equation_application}) using the algorithm and the brute-force refinement. The blue color denotes the standard FEM, and the red color denotes the mixed FEM.}
\label{E_min_application}
\end{figure}

Using the prediction approach and the brute-force approach, $E_{\rm min}$ are compared in Fig.~\ref{E_min_application}. As can be seen, $E_{\rm min}$ can be predicted correctly.

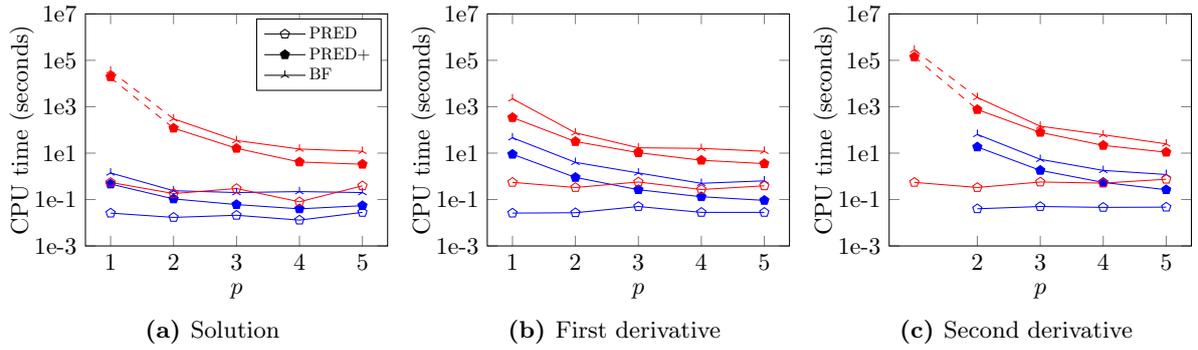
\begin{figure}[!ht]
\hspace{0.0cm}
\begin{subfigure}[b]{0.35\textwidth}
\scalebox{0.9}{
\begin{tikzpicture} 
\begin{axis}
[
    ymode=log,
    ymin=1e-3,
    ymax=1e7,
    ytick={1e-3, 1e-1, 1e1, 1e3, 1e5, 1e7},
    yticklabels={1e-3, 1e-1, 1e1, 1e3, 1e5, 1e7},
    legend style={nodes={scale=0.7}},
    legend cell align={left},
    height=5cm,
    width=6cm,
    ylabel={CPU time (seconds)},
    ylabel style={below},    
    xtick={0,1,2,3,4},
    xticklabels={$1$, $2$, $3$, $4$, ${5}$},
    xlabel={$p$},
    xlabel style={at={(0.5,0.03)}},   
]
\addplot[black,mark=pentagon,mark options={color=black,fill=black}] coordinates {(0,2.1e-25)}; 
\addplot[black,mark=pentagon*,mark options={color=black,fill=black}] coordinates {(0,2.1e-25)};
\addplot[black,mark=Mercedes star,mark options={color=black,fill=black}] coordinates {(0,2.1e-25)}; 

\addplot[blue,mark=pentagon,mark options={color=blue,fill=blue}] coordinates {(0,0.026) (1,0.017) (2,0.021) (3,0.013) (4,0.028)}; 
\addplot[blue,mark=pentagon*,mark options={color=blue,fill=blue}] coordinates {(0,0.026+0.439343) (1,0.017+0.08932) (2,0.021+0.0397) (3,0.013+0.026735) (4,0.028+0.02639)};
\addplot[blue,mark=Mercedes star,mark options={color=blue,fill=blue}] coordinates {(0,1.4) (1,0.24) (2,0.2) (3,0.22) (4,0.2)}; 

\addplot[red,mark=pentagon,mark options={color=red,fill=red}] coordinates {(0,0.546) (1,0.181) (2,0.298) (3,0.08) (4,0.39)}; 
\addplot[red,mark=pentagon*,mark options={color=red,fill=red}] coordinates {(0,0.546+19600)};
\addplot[red,mark=pentagon*,mark options={color=red,fill=red}] coordinates {(1,0.181+119) (2,0.298+15.8) (3,0.08+4.07) (4,0.39+2.92)};
\draw [red,dashed] (axis cs:0,0.546+17600) -- (axis cs:1,0.181+119);
\addplot[red,mark=Mercedes star,mark options={color=red,fill=red}] coordinates {(0,35000)};
\addplot[red,mark=Mercedes star,mark options={color=red,fill=red}] coordinates {(1,300) (2,35) (3,15) (4,12)};
\draw [red,dashed] (axis cs:0,35000) -- (axis cs:1,300);

\legend{PRED, PRED+, BF};
\end{axis}
\end{tikzpicture}
}
\vspace{-0.2cm}
\caption{Solution}
\end{subfigure}
\hspace{-0.7cm}
\begin{subfigure}[b]{0.35\textwidth}
\scalebox{0.9}{
\begin{tikzpicture} 
\begin{axis}
[
    ymode=log,
    ymin=1e-3,
    ymax=1e7,
    ytick={1e-3, 1e-1, 1e1, 1e3, 1e5, 1e7},
    yticklabels={1e-3, 1e-1, 1e1, 1e3, 1e5, 1e7},
    legend style={nodes={scale=0.7}},
    legend cell align={left},
    height=5cm,
    width=6cm,
    ylabel={CPU time (seconds)},
    ylabel style={below},    
    xtick={0,1,2,3,4},
    xticklabels={$1$, $2$, $3$, $4$, ${5}$},
    xlabel={$p$},
    xlabel style={at={(0.5,0.03)}},   
]
\addplot[blue,mark=pentagon,mark options={color=blue,fill=blue}] coordinates {(0,0.026) (1,0.027) (2,0.050) (3,0.028) (4,0.028)}; 
\addplot[blue,mark=pentagon*,mark options={color=blue,fill=blue}] coordinates {(0,0.026+8.9) (1,0.027+0.87) (2,0.05+0.215) (3,0.028+0.104) (4,0.028+0.064)};
\addplot[blue,mark=Mercedes star,mark options={color=blue,fill=blue}] coordinates {(0,46) (1,4) (2,1.4) (3,0.5) (4,0.64)}; 

\addplot[red,mark=pentagon,mark options={color=red,fill=red}] coordinates {(0,0.55) (1,0.33) (2,0.57) (3,0.27) (4,0.39)}; 
\addplot[red,mark=pentagon*,mark options={color=red,fill=red}] coordinates {(0,0.55+340) (1,0.33+31.3) (2,0.57+9.88) (3,0.27+4.66) (4,0.39+3.13)};
\addplot[red,mark=Mercedes star,mark options={color=red,fill=red}] coordinates {(0,2200) (1,75) (2,17) (3,16) (4,12)}; 
\end{axis}
\end{tikzpicture}
}
\vspace{-0.2cm}
\caption{First derivative}
\end{subfigure}
\hspace{-0.7cm}
\begin{subfigure}[b]{0.35\textwidth}
\scalebox{0.9}{
\begin{tikzpicture} 
\begin{axis}
[
    ymode=log,
    ymin=1e-3,
    ymax=1e7,
    ytick={1e-3, 1e-1, 1e1, 1e3, 1e5, 1e7},
    yticklabels={1e-3, 1e-1, 1e1, 1e3, 1e5, 1e7},
    legend style={nodes={scale=0.7}},
    legend cell align={left},
    height=5cm,
    width=6cm,
    ylabel={CPU time (seconds)},
    ylabel style={below},     
    xtick={1,2,3,4},
    xticklabels={$2$, $3$, $4$, ${5}$},
    xlabel={$p$},
    xlabel style={at={(0.5,0.03)}},     
]
\addplot[blue,mark=pentagon,mark options={color=blue,fill=blue}] coordinates {(1,0.04) (2,0.05) (3,0.046) (4,0.047)}; 
\addplot[blue,mark=pentagon*,mark options={color=blue,fill=blue}] coordinates {(1,0.04+18.4184) (2,0.05+1.76016) (3,0.046+0.500414) (4,0.047+0.217396)}; 
\addplot[blue,mark=Mercedes star,mark options={color=blue,fill=blue}] coordinates {(1,64) (2,5.4) (3,1.8) (4,1.2)}; 

\addplot[red,mark=pentagon,mark options={color=red,fill=red}] coordinates {(0,0.55) (1,0.33) (2,0.57) (3,0.51) (4,0.75)}; 
\addplot[red,mark=pentagon*,mark options={color=red,fill=red}] coordinates {(0,0.55+140000)};
\addplot[red,mark=pentagon*,mark options={color=red,fill=red}] coordinates {(1,0.33+755) (2,0.57+78) (3,0.51+21) (4,0.75+10.2)};
\draw [red,dashed] (axis cs:0,0.55+140000) -- (axis cs:1,0.33+755);
\addplot[red,mark=Mercedes star,mark options={color=red,fill=red}] coordinates {(0,280000)};
\addplot[red,mark=Mercedes star,mark options={color=red,fill=red}] coordinates {(1,2500) (2,145) (3,62) (4,25)};
\draw [red,dashed] (axis cs:0,280000) -- (axis cs:1,2500);
\end{axis}
\end{tikzpicture}
}
\vspace{-0.2cm}
\caption{Second derivative}
\end{subfigure}
\caption{Comparison of the CPU time to obtain $E_{\rm min}$ for Eq. (\ref{1D_Helmholtz_equation_application}) using the algorithm and the brute-force refinement. The blue color denotes the standard FEM, and the red color denotes the mixed FEM.}
\label{CPU_algorithm_brute_force}
\end{figure}

\begin{figure}[!ht]
\centering
\scalebox{0.9}{
\begin{tikzpicture}
\begin{axis}
[
    ymin=20,
    ymax=100,
    legend style={nodes={scale=0.8},at={(0.65,0.2)},anchor=west},
    legend cell align={left},
    height=5cm,
    width=6cm,
    ylabel={CPU time saved},
    ylabel style={at={(0.04,0.5)},font=\small},
    y tick label style={font=\footnotesize},
    yticklabel=\pgfmathprintnumber{\tick}\,$\%$,
    xtick={0,1,2,3,4},
    xticklabels={$1$, $2$, $3$, $4$, ${5}$},
    xlabel={$p$},
    xlabel style={at={(0.5,0.03)}}, 
]
\addplot[black,mark=o,mark options={color=black,fill=black}] coordinates {(0,2.1e-25)};
\addplot[black,mark=diamond,mark options={color=black,fill=black}] coordinates {(0,2.1e-25)};
\addplot[black,mark=triangle,mark options={color=black,fill=black}] coordinates {(0,2.1e-25)};

\addplot[blue,mark=o,mark options={color=blue,fill=blue}] coordinates {(0,66.8) (1,55.8) (2,69.5) (3,82.0) (4,73.0)};
\addplot[blue,mark=diamond,mark options={color=blue,fill=blue}] coordinates {(0,80.6) (1,77.6) (2,81.1) (3,73.6) (4,85.6)};
\addplot[blue,mark=triangle,mark options={color=blue,fill=blue}] coordinates {(0) (1,71.2) (2,66.5) (3,69.7) (4,78)};

\addplot[red,mark=o,mark options={color=red,fill=red}] coordinates {(0,44) (1,60.3) (2,54) (3,72.3) (4,72.4)};
\addplot[red,mark=diamond,mark options={color=red,fill=red}] coordinates {(0,84.5) (1, 57.8) (2, 38.5) (3, 69.3) (4, 70.7)};
\addplot[red,mark=triangle,mark options={color=red,fill=red}] coordinates {(0,50) (1,69.8) (2,45.8) (3, 65.3) (4, 56.1)};
\legend{$u$,$u_x$,$u_{xx}$};
\end{axis}
\end{tikzpicture}
}
\caption{Percentage of CPU time saved using the algorithm. The blue color denotes the standard FEM, and the red color denotes the mixed FEM.}
\label{Percentage_CPU_Saved}
\end{figure}
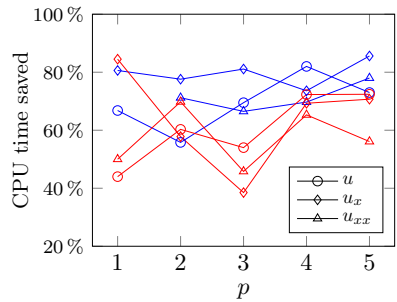

The CPU time required by the prediction approach (PRED) and the brute-force approach (BF) is shown in Fig.~\ref{CPU_algorithm_brute_force}. Next to time PRED, and the computation time for the optimal grid (PRED+) using the prediction approach is also given.
As can be seen, both time BF and time PRED+ decrease with increasing element degree. Time PRED+ is much smaller compared to time BF, see Fig.~\ref{Percentage_CPU_Saved} for the percentage of the CPU time saved by PRED+, which shows a saving of the CPU time basically more than 60\% and 40\% for the standard FEM and the mixed FEM, respectively. Last but not least, time PRED is negligible compared to time PRED+.

\vspace{-0.5cm}

\newpage
Furthermore, the dashed line indicating the desired error tolerance in Fig.~\ref{E_min_application} cannot be reached using the standard FEM, whereas it can be reached using the mixed FEM with $P_4/P_3^{\rm disc}$ or betters. When using $P_4/P_3^{\rm disc}$, $N_{\rm opt}$ for $u$, $u_x$ and $u_{xx}$ are predicted to be 6042, 9812 and 123486, respectively.

\section{Conclusions}		\label{paragraph on conclusion}

A novel approach is presented to predict the highest attainable accuracy for second-order ordinary differential equations using the finite element methods.
In contrast to the brute-force approach, which uses successive $h$-refinements, this approach uses only a few coarse grid refinements. 
This approach is viable for the solution and its first and second derivative, for both the standard FEM and the mixed FEM, and different element degrees.
The algorithm for implementing the approach shows that the highest attainable accuracy can be accurately predicted and the CPU time is significantly reduced.
To compute the solution of the highest attainable accuracy using our approach, the CPU time can be saved more than 60\% for the standard FEM and 40\% for the mixed FEM.

Future research will focus on the validation of the approach for 2D second-order problems, where the influence of the linear system solver, local mesh refinement and boundary conditions might be significantly different from 1D problems. 

\appendix

\section{Derivation of the weak form}		\label{weak form appendix}

\subsection{The standard FEM}		\label{derivation_weak_form_SM}

Multiply Eq. (\ref{1D_general_Helmholtz_equation}) by a test function $\eta \in H ^1 (I)$, and integrate it over $I$ yields
\begin{equation}
(\eta, \, \left(D u_x \right)_x + ru) = (\eta, \, f). \label{1D_general_inte}
\end{equation}

By applying Gauss's theorem for the first term of the left-hand side of Eq. (\ref{1D_general_inte}), we obtain
\begin{equation}
 -({\eta} _x, \, D {u} _{ x }) + (\eta, \, ru) = (\eta, \, f) - \left( \eta, \, D u_x n \right)_{ {\Gamma_N}}.		\label{1D_general_gauss}
\end{equation}

Therefore, without considering the boundary conditions, the weak form reads
\begin{equation}
\centering
\boxed{ 
\begin{aligned}
&\text{Find $u \in H ^1 (I)$ such that:} \\
&-({ \eta} _{ x }, \,D { u} _{ x }  ) + (\eta, \, ru) = (\eta, \, f ) - (\eta, \, D u_x n )_{\Gamma _N} \qquad \forall \eta \in H ^1 (I),\\
&\text{where } {n} \text{ is 1 at $x=1$, and -1 at $x=0$.}
\end{aligned}	\label{1D_general_SM_weak_form_no_boundary} 
}
\end{equation}

Imposing the original Dirichlet boundary conditions on $u$ and the corresponding homogeneous Dirichlet boundary conditions on $\eta$ in Eq. (\ref{1D_general_SM_weak_form_no_boundary}), which is called the strong imposition of the Dirichlet boundary conditions, the weak form can be found in Eq. (\ref{1D_general_SM_weak_form_Diri_strong}).
Instead of imposing the Dirichlet boundary conditions directly on the variables $u$ and $\eta$ in Eq. (\ref{1D_general_SM_weak_form_no_boundary}), by adding auxiliary terms \cite{bazilevs2007weak}, which is called the weak imposition of the Dirichlet boundary conditions, we obtain the weak form Eq. (\ref{1D_general_SM_weak_form_Diri_weak}).

\subsection{The mixed FEM}		\label{derivation_weak_form_MM}

To obtain the weak form of Eq. (\ref{1D_general_MM_2in1}), Eq. (\ref{Gene_MM_strong1}) is multiplied by a test function of $v$, i.e. $w \in H _{N0}^{1}(I)$, and integrated over $I$, yielding
\begin{subequations}
\begin{align}
  ( v + u _x, w) = 0,	\label{Gene_MM_weak1_inte}
\end{align}
and Eq. (\ref{Gene_MM_strong2}) is multiplied by a test function of $u$, i.e. $q \in L^2 (I)$, and integrated over $I$, yielding 
\begin{align}
  -( q, \, D_x v) - ( q , \, D v_x) + (q, \, ru) = (q, \, f ). \label{1D_Poisson_classical_Dint2}
\end{align}
By applying Gauss's theorem to Eq. (\ref{Gene_MM_weak1_inte}) and imposing the natural boundary condition $u(x)=g(x)$ on $\Gamma_D$, we obtain
\begin{align}
 (w, \, v) - (w_x, \,  u ) = -(w, \, g n)_{\Gamma_D}.
\end{align}
\end{subequations}
The resulting weak form can be found in Eq. (\ref{1D_General_MM_weak_2in1}).

\newpage

\section{Numerical results of the benchmark diffusion and Helmholtz equations}         \label{discretization_error_bench_diff_Helm}

\subsection{The diffusion equation}         \label{discretization_error_bench_diff}

\begin{figure}[!ht]
    \begin{subfigure}{5.5cm}
        \includegraphics[width=1.0\linewidth]{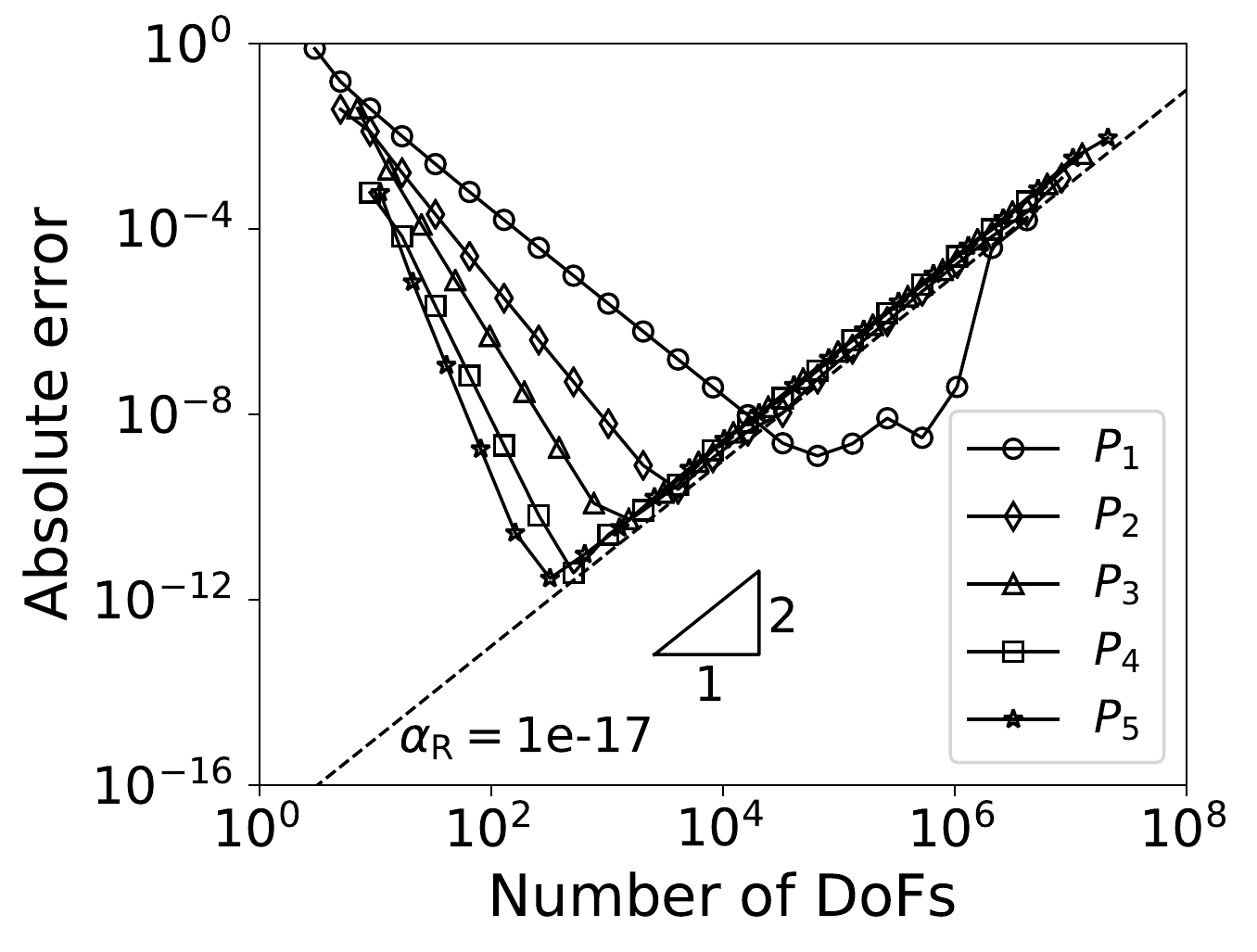}
        \caption{Solution}
        \label{py_bench_diff_SM_solu}
    \end{subfigure}
    \hspace{-0.2cm}
    \begin{subfigure}{5.5cm}
        \includegraphics[width=1.0\linewidth]{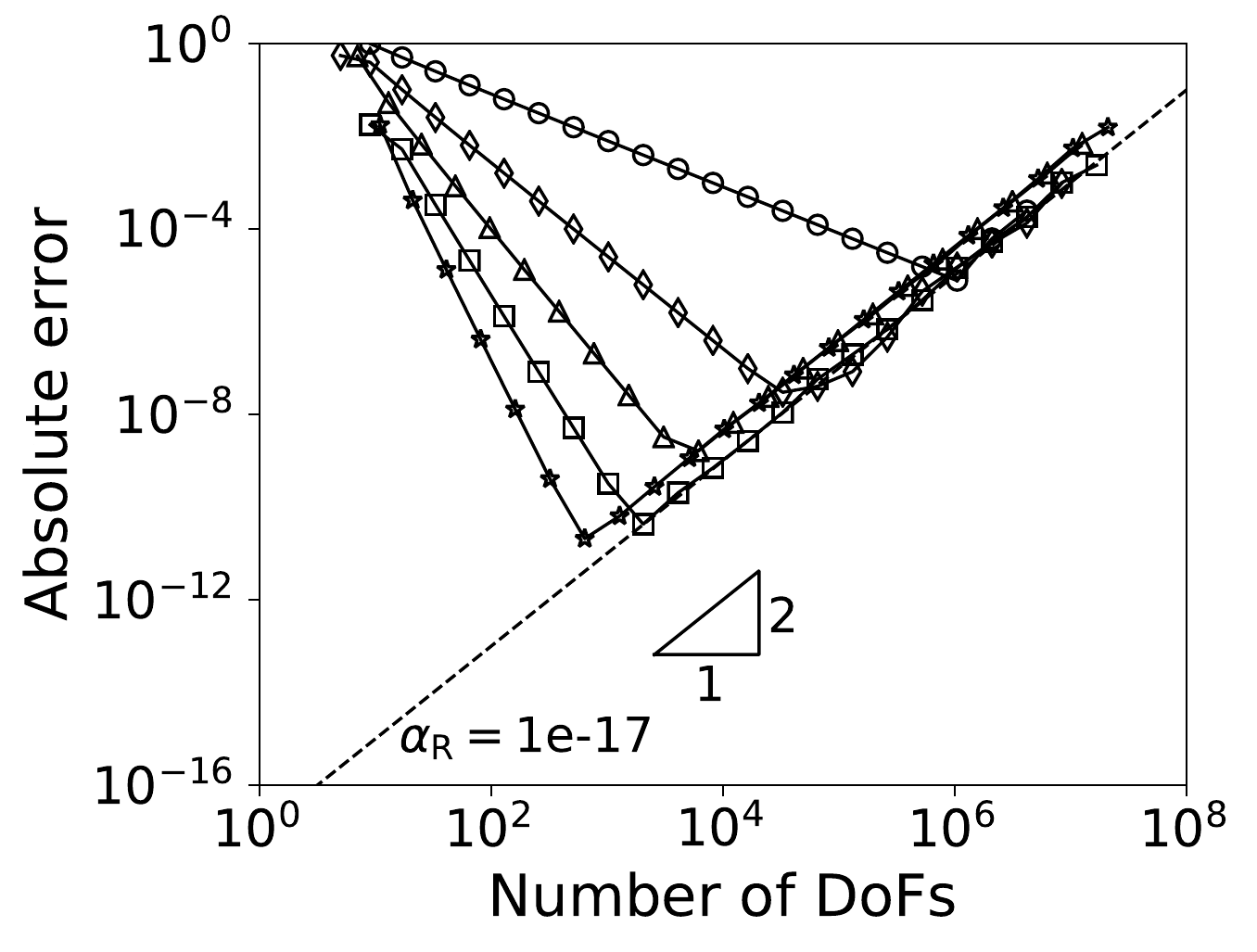}
        \caption{First derivative}
        \label{py_bench_diff_SM_grad}
    \end{subfigure}
    \hspace{-0.2cm}
    \begin{subfigure}{5.5cm}
        \includegraphics[width=1.0\linewidth]{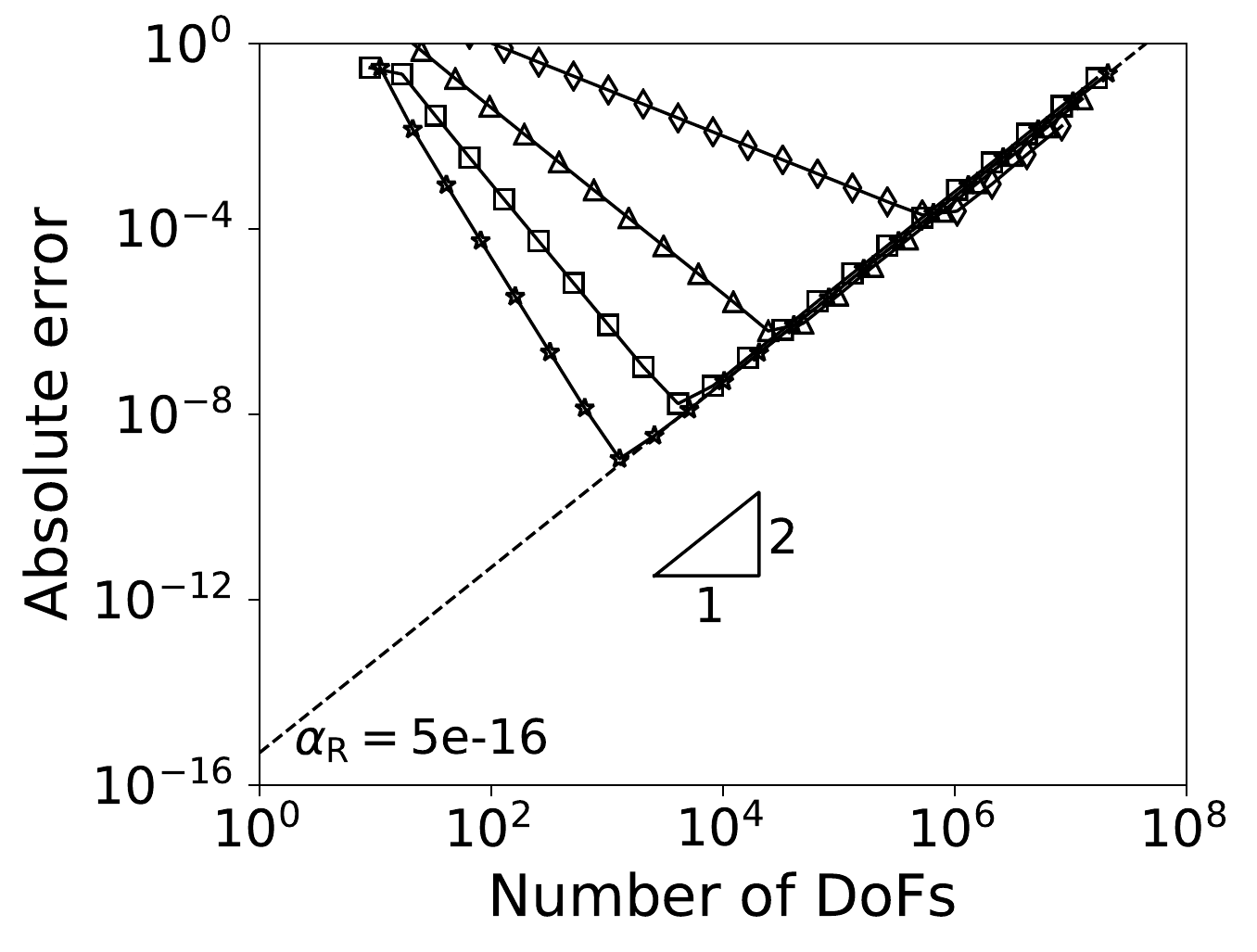}
        \caption{Second derivative}
        \label{py_bench_diff_SM_2ndd}
    \end{subfigure}
\caption{Absolute errors for the benchmark diffusion equation using the standard FEM.}
\label{py_bench_diff_SM}
\end{figure}

\begin{figure}[!ht]
    \begin{subfigure}{5.5cm}
        \includegraphics[width=1.0\linewidth]{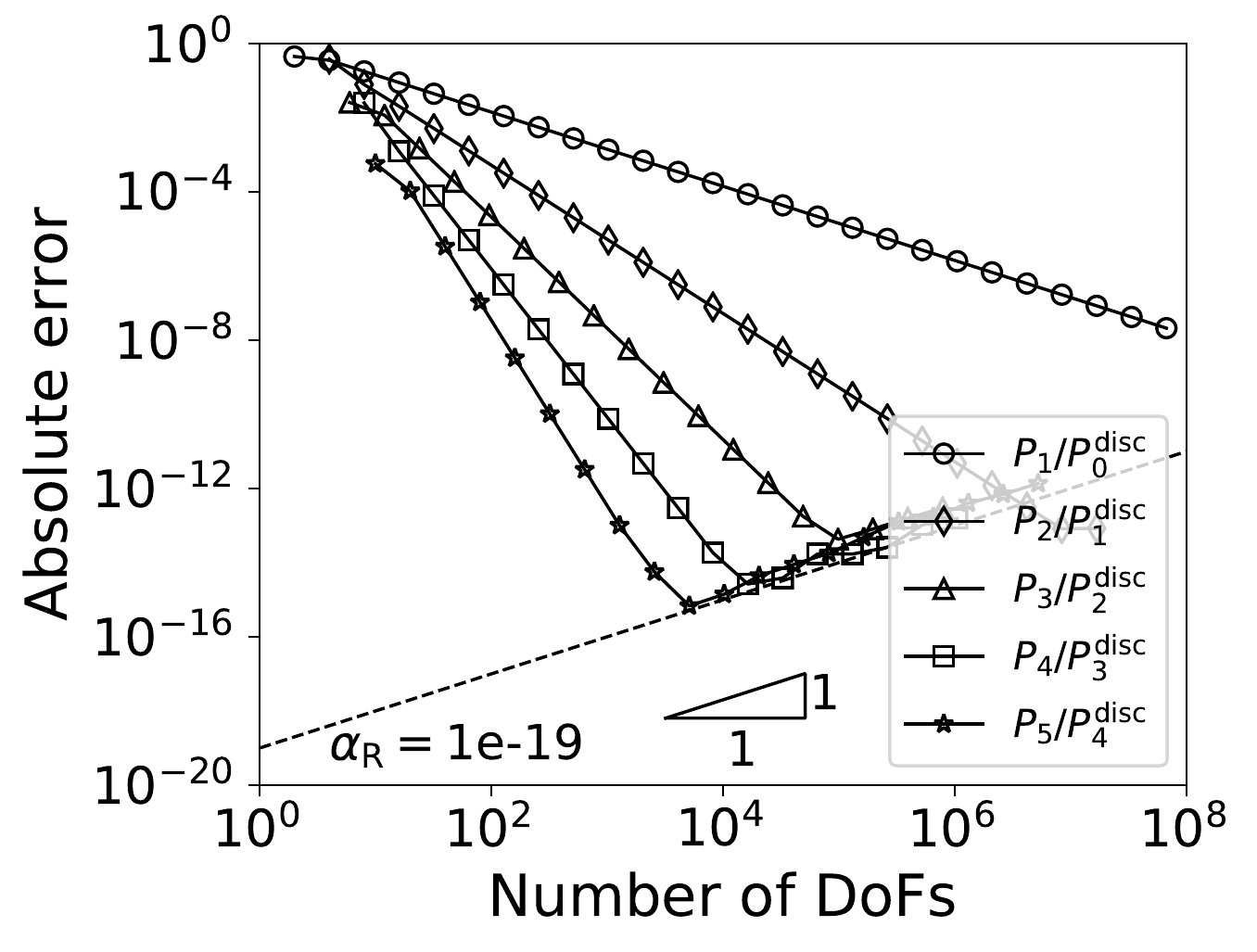}
        \caption{Solution}
        \label{py_bench_diff_MM_solu}
    \end{subfigure}
    \hspace{-0.2cm}
    \begin{subfigure}{5.5cm}
        \includegraphics[width=1.0\linewidth]{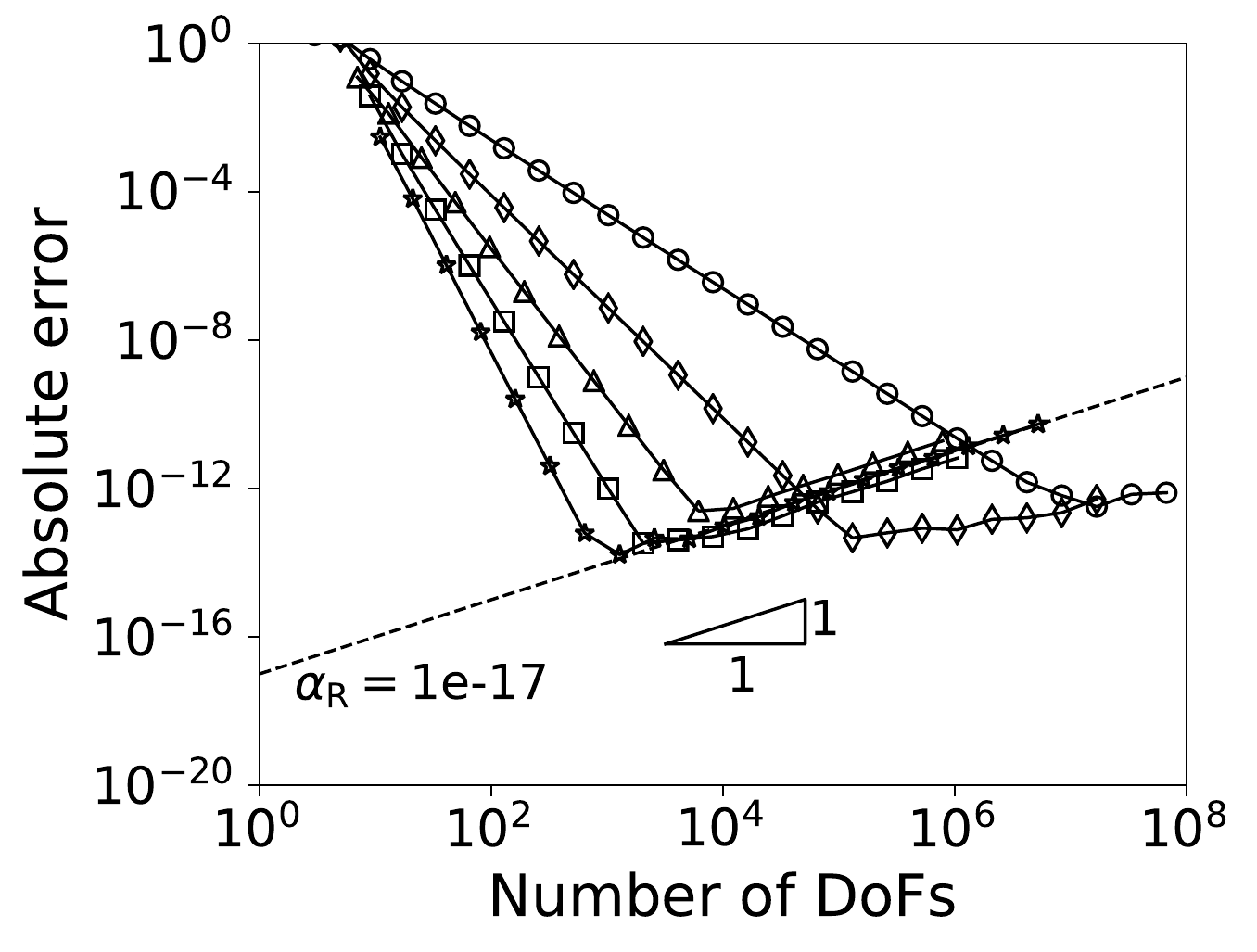}
        \caption{First derivative}
        \label{py_bench_diff_MM_grad}
    \end{subfigure}
    \hspace{-0.2cm}
    \begin{subfigure}{5.5cm}
        \includegraphics[width=1.0\linewidth]{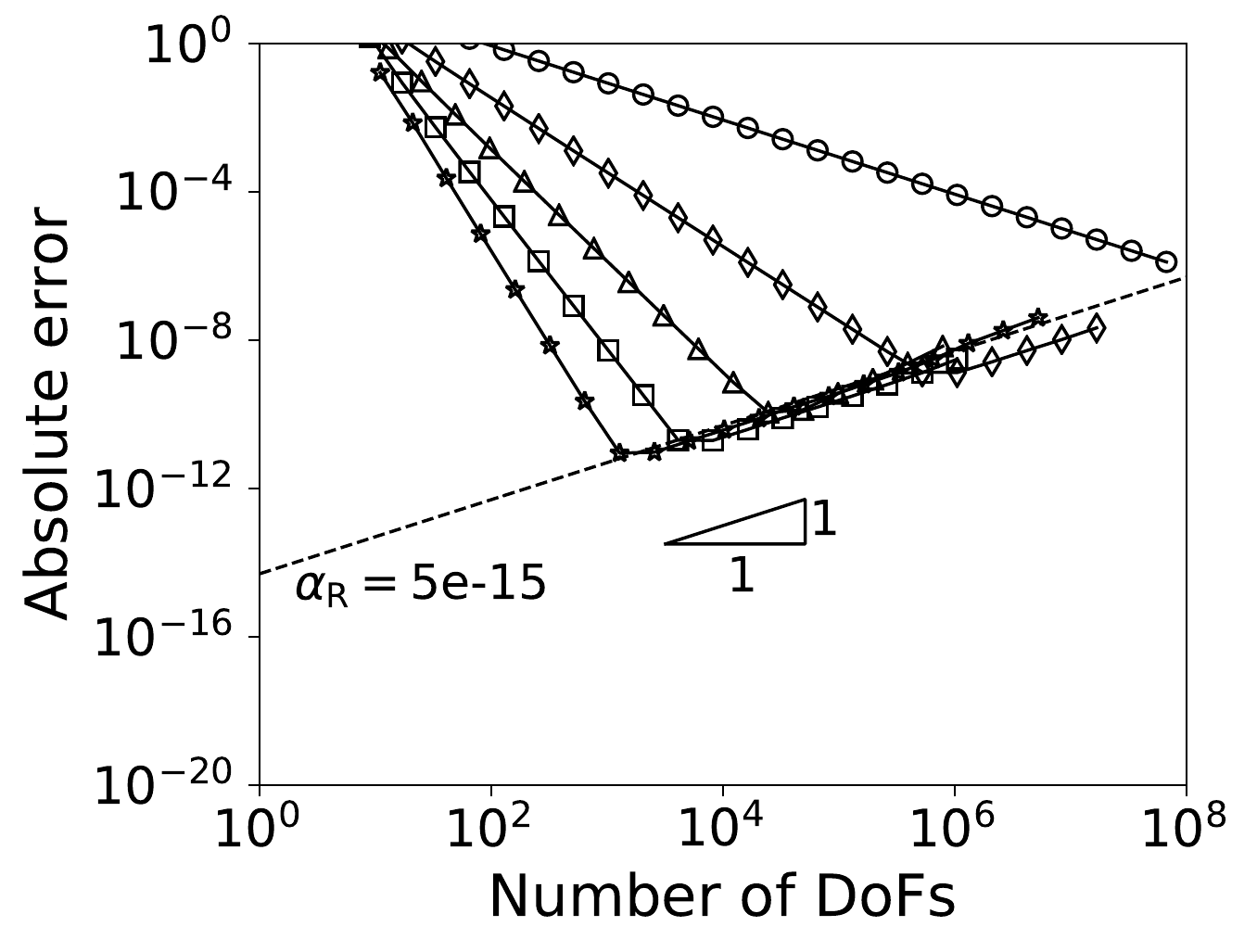}
        \caption{Second derivative}
        \label{py_bench_diff_MM_2ndd}
    \end{subfigure}
\caption{Absolute errors for the benchmark diffusion equation using the mixed FEM}
\label{py_bench_diff_MM}    
\end{figure}

For the benchmark diffusion equation, using both the standard FEM and the mixed FEM, the absolute errors for all \emph{three} variables are shown in Fig.~\ref{py_bench_diff_SM} and  Fig.~\ref{py_bench_diff_MM}, respectively. 

\newpage

\subsection{The Helmholtz equation}             \label{discretization_error_bench_Helm}

\begin{figure}[!ht]
    \begin{subfigure}{5.5cm}
        \includegraphics[width=1.0\linewidth]{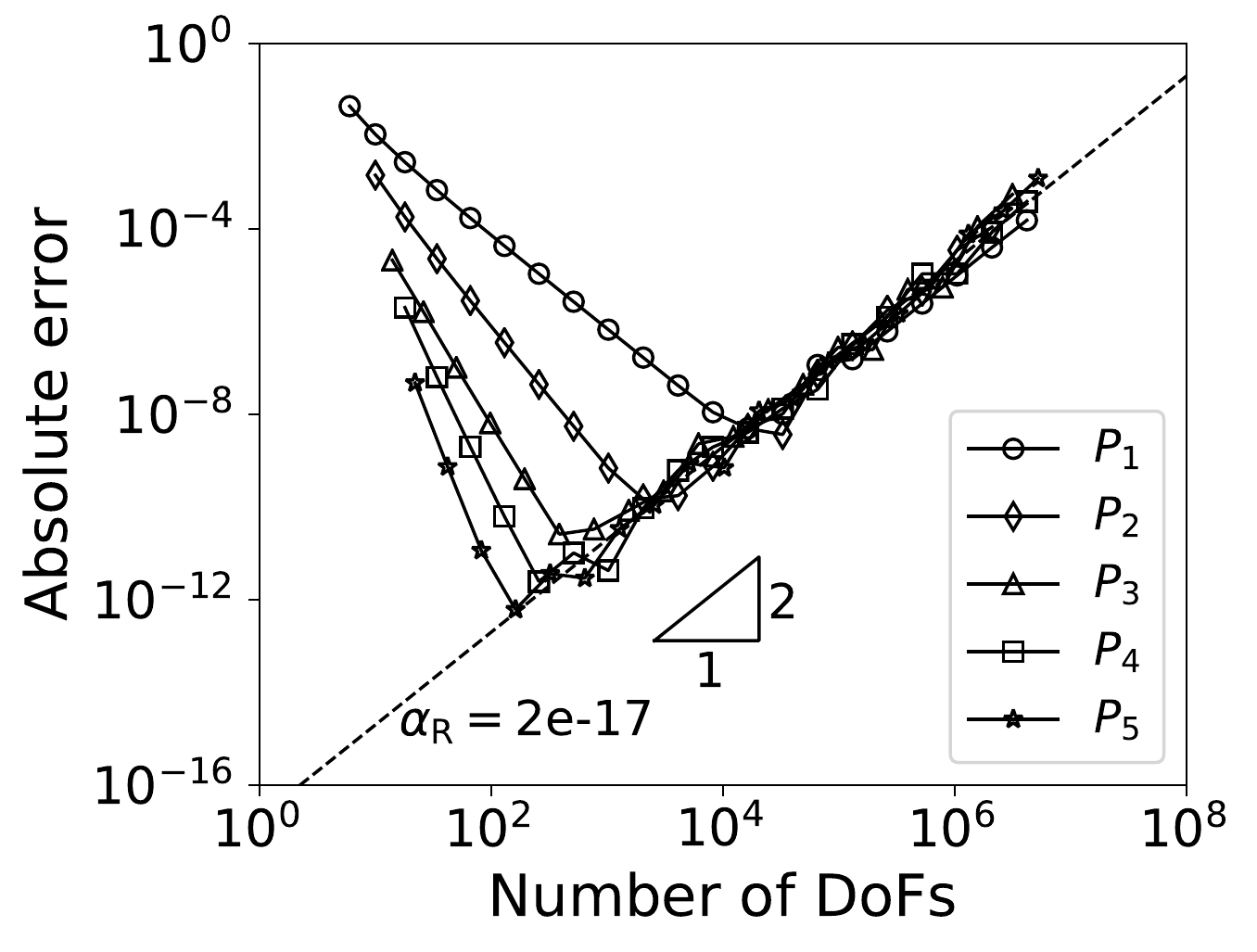}
        \caption{Solution}
        \label{py_bench_Helm_SM_solu}
    \end{subfigure}
    \hspace{-0.2cm}
    \begin{subfigure}{5.5cm}
        \includegraphics[width=1.0\linewidth]{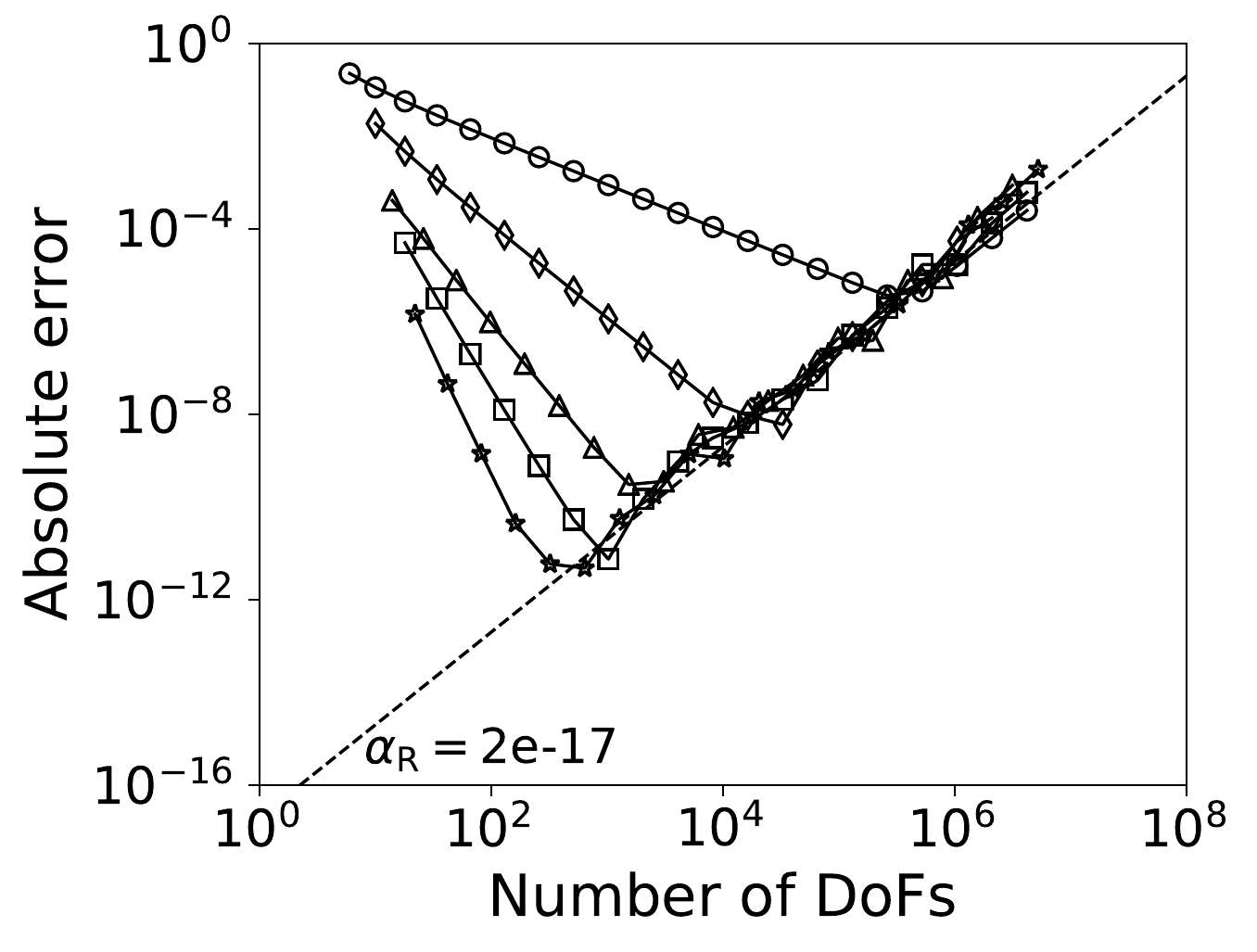}
        \caption{First derivative}
        \label{py_bench_Helm_SM_grad}
    \end{subfigure}
    \hspace{-0.2cm}
    \begin{subfigure}{5.5cm}
        \includegraphics[width=1.0\linewidth]{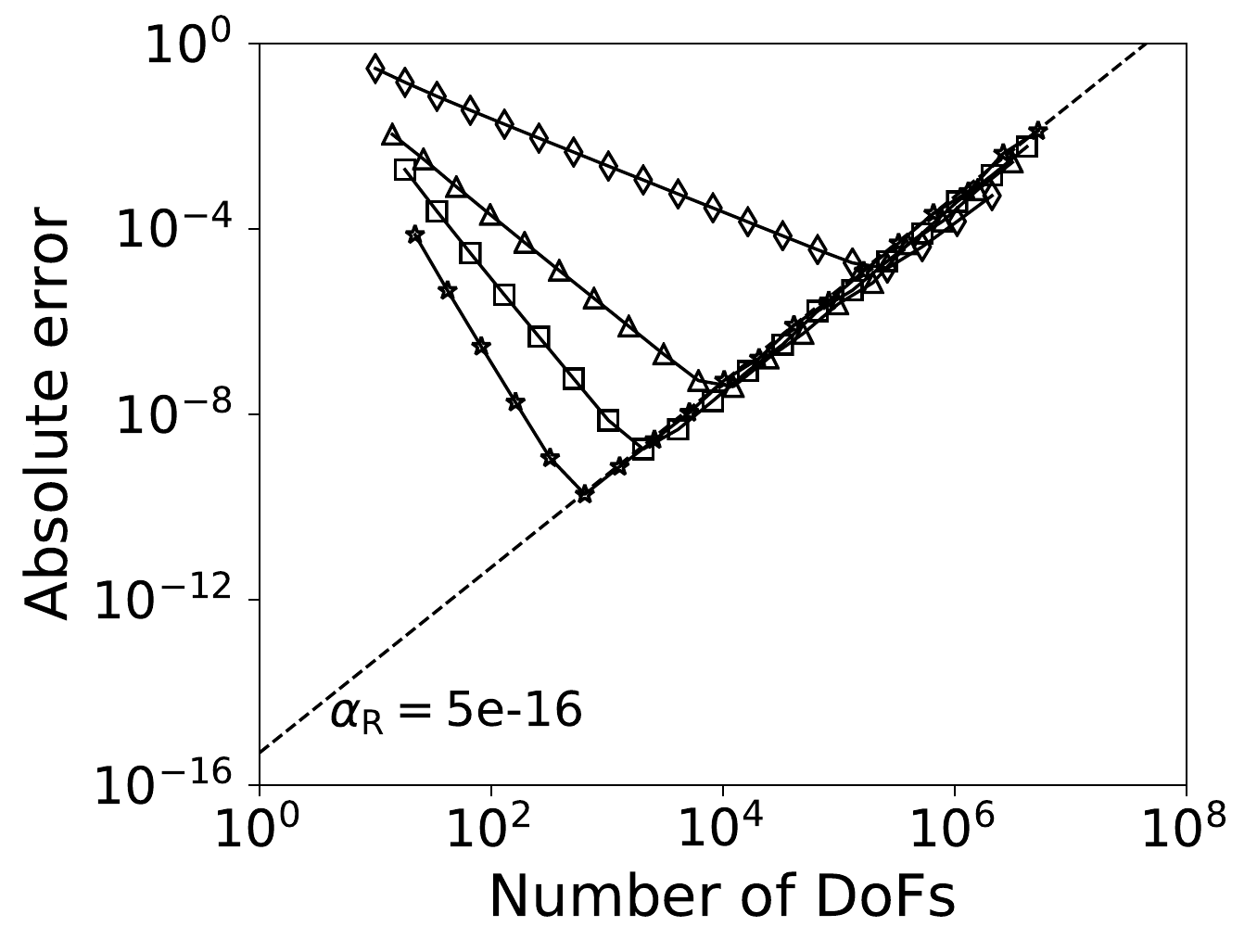}
        \caption{Second derivative}
        \label{py_bench_Helm_SM_2ndd}
    \end{subfigure}
\caption{Absolute errors for the benchmark Helmholtz equation using the standard FEM.}
\label{py_bench_Helm_SM}
\end{figure}

\begin{figure}[!ht]
    \begin{subfigure}{5.5cm}
        \includegraphics[width=1.0\linewidth]{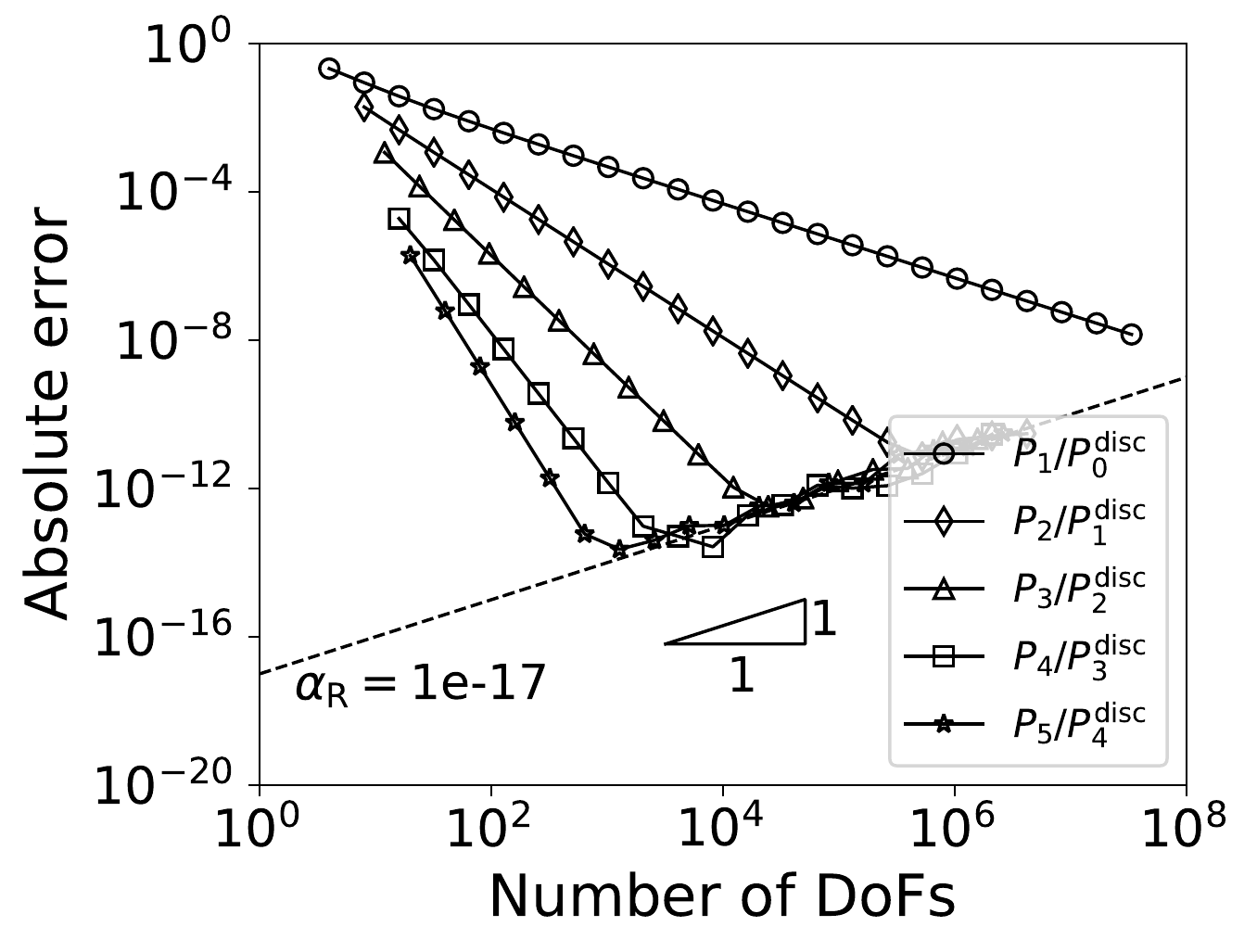}
        \caption{Solution}
        \label{py_bench_Helm_MM_solu}
    \end{subfigure}
    \hspace{-0.2cm}
    \begin{subfigure}{5.5cm}
        \includegraphics[width=1.0\linewidth]{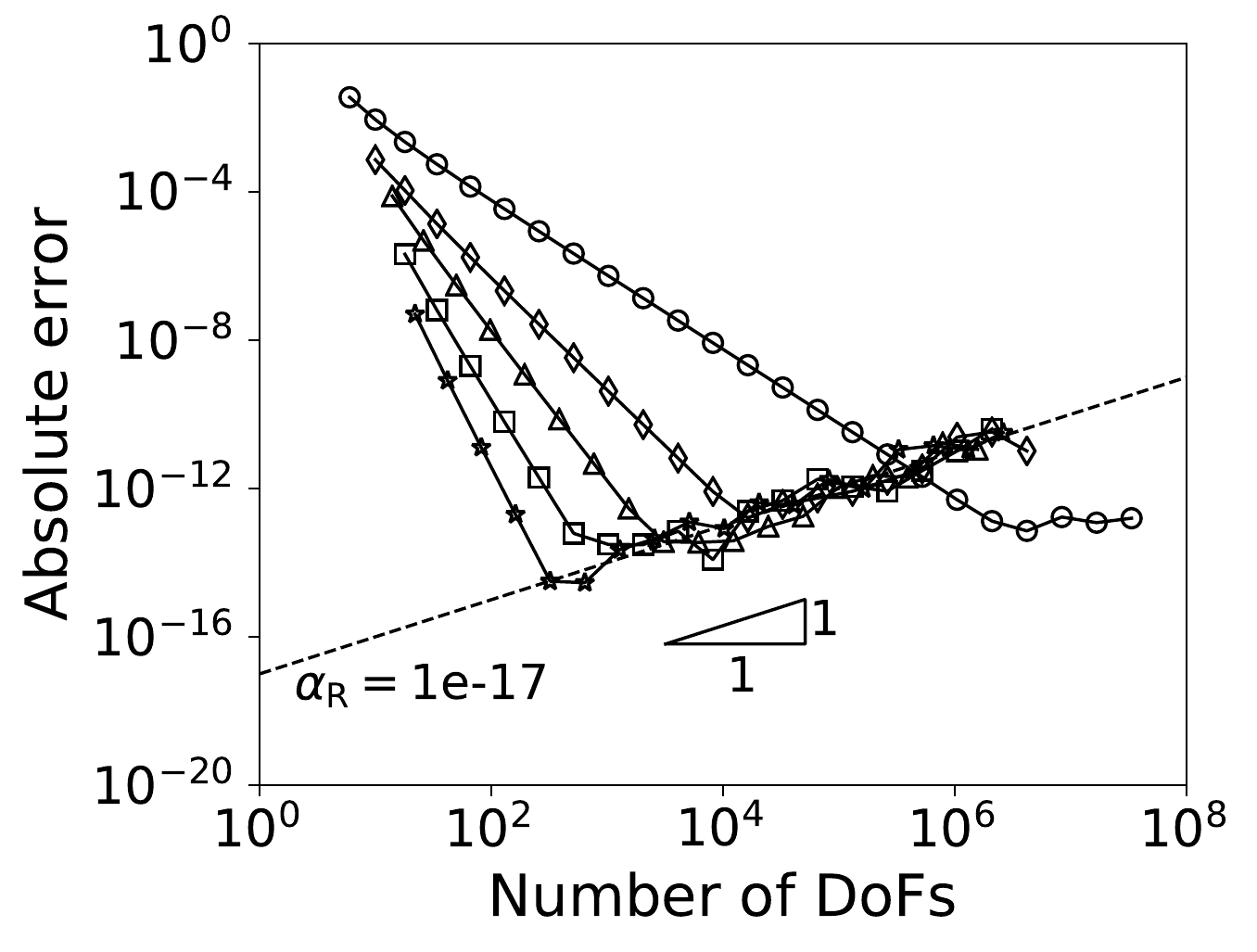}
        \caption{First derivative}
        \label{py_bench_Helm_MM_grad}
    \end{subfigure}
    \hspace{-0.2cm}
    \begin{subfigure}{5.5cm}
        \includegraphics[width=1.0\linewidth]{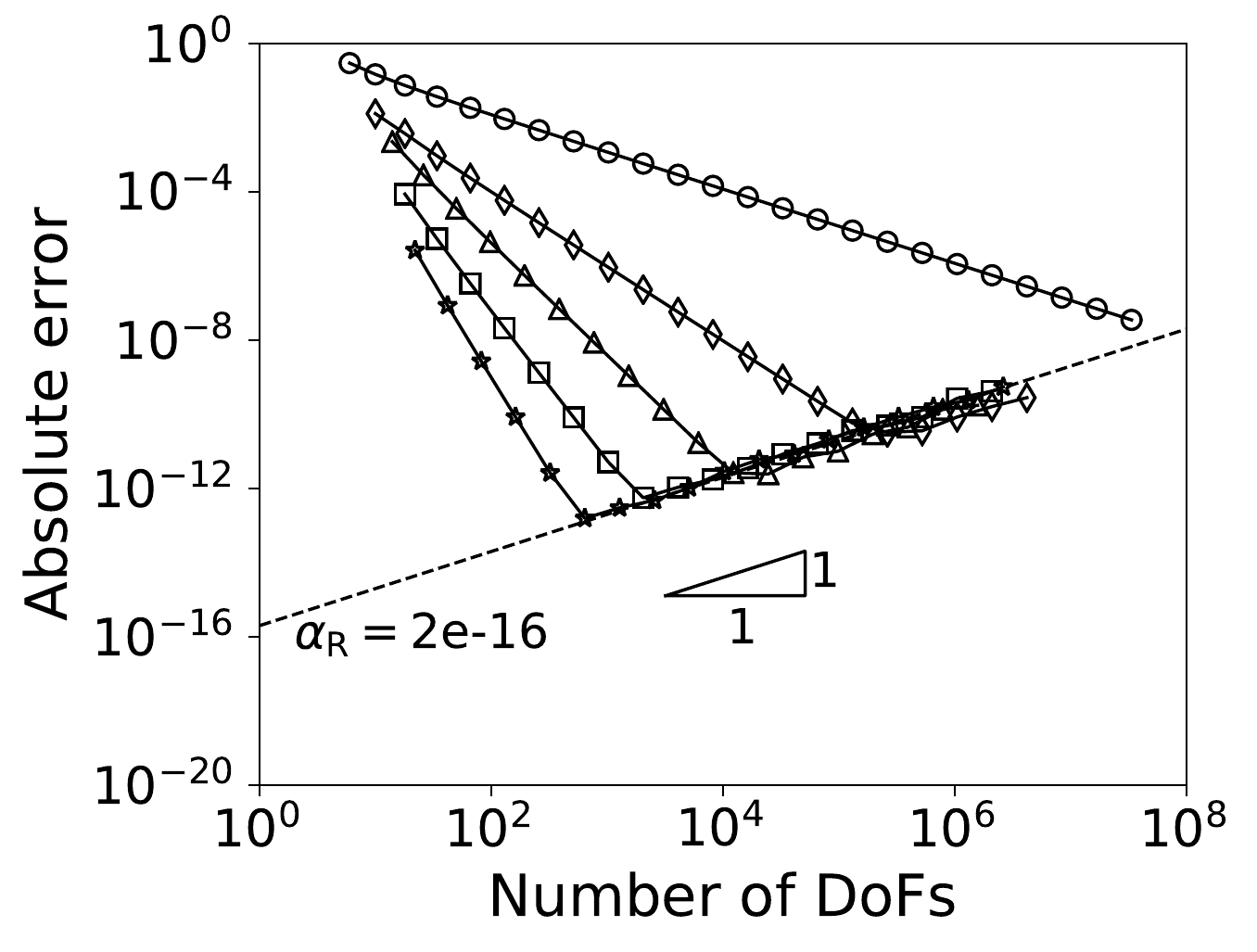}
        \caption{Second derivative}
        \label{py_bench_Helm_MM_2ndd}
    \end{subfigure}
\caption{Absolute errors for the benchmark Helmholtz equation using the mixed FEM.}
\label{py_bench_Helm_MM}
\end{figure}

For the benchmark Helmholtz equation, using both the standard FEM and the mixed FEM, the absolute errors for all \emph{three} variables are shown in Fig.~\ref{py_bench_Helm_SM} and  Fig.~\ref{py_bench_Helm_MM}, respectively. 

\newpage
\section{\texorpdfstring{$L_2$}{} norms and absolute errors for different cases in Table \ref{scaling_cases_Poisson}}   \label{results_order_of_magnitude_other_cases}

\subsection{\texorpdfstring{$L_2$}{} norms}     \label{L2_norm_scaling_cases}

\begin{figure}[!ht]
\hspace{0.0cm}
\begin{subfigure}[b]{0.3\textwidth}
\scalebox{0.9}{
\begin{tikzpicture} 
\begin{axis}
[
    ymode=log,
    ymin=1e-8,
    ymax=1e4,
    ytick={1e-8, 1e-4, 1e0, 1e4},
    yticklabels={1e-8, 1e-4, 1e0, 1e4},
    legend style={at={(0.4,0.4)},nodes={scale=1}},
    legend cell align={left},
    height=5cm,
    width=6cm,
    xlabel={$c_1$},
    xlabel style={below},      
    ylabel={$L_2$ norm},
    ylabel style={below},    
    xtick=data,
    xticklabels={$10^{-2}$, $10^{-1}$, $10^{-0}$, $10^{1}$, $10^{2}$}    
]
\addplot[mark=o] coordinates {(0,9.2e0) (1,8.8e-1) (2,1.8e-2) (3,1.8e-4) (4,1.8e-6)}; 
\addplot[mark=diamond] coordinates {(0,1.6e1) (1,1.5e0) (2,1.1e-1) (3,1.1e-2) (4,1.1e-3)}; 
\addplot[mark=triangle] coordinates {(0,3.6e-2) (1,3.5e-1) (2,7.1e-1) (3,7.1e-1) (4,7.1e-1)}; 
\legend{$u$,$u_x$,$u_{xx}$};
\end{axis}
\end{tikzpicture}
}
\caption{Case 1}
\label{L2_norm_case_1}
\end{subfigure}
\hspace{0.15cm}
\begin{subfigure}[b]{0.3\textwidth}
\scalebox{0.9}{
\begin{tikzpicture} 
\begin{axis}
[
    ymode=log,
    ymin=1e-8,
    ymax=1e4,
    ytick={1e-8, 1e-4, 1e0, 1e4},
    yticklabels={1e-8, 1e-4, 1e0, 1e4},    
    height=5cm,
    width=6cm,
    xlabel={$c_2$},
    xlabel style={below},      
    ylabel={$L_2$ norm},
    ylabel style={below},    
    xtick=data,
    xticklabels={$10^{-4}$, $10^{-2}$, $10^{-0}$, $10^{2}$, $10^{4}$}    
]
\addplot[mark=o] coordinates {(0, 1e0) (1, 1e0) (2, 0.92) (3, 0.35) (4, 0.11)}; 
\addplot[mark=diamond] coordinates {(0, 5.8e-5) (1, 5.8e-3) (2, 5.0e-1) (3, 3.5) (4, 11)}; 
\addplot[mark=triangle] coordinates {(0, 2.0e-4) (1, 2.0e-2) (2, 1.6) (3, 61) (4, 1.9e3)}; 
\end{axis}
\end{tikzpicture}
}
\caption{Case 2}
\label{L2_norm_case_2}
\end{subfigure}
\hspace{0.15cm}
\begin{subfigure}[b]{0.3\textwidth}
\scalebox{0.9}{
\begin{tikzpicture} 
\begin{axis}
[
    ymode=log,
    ymin=1e-8,
    ymax=1e4,
    ytick={1e-8, 1e-4, 1e0, 1e4},
    yticklabels={1e-8, 1e-4, 1e0, 1e4},    
    height=5cm,
    width=6cm,
    xlabel={$c_3$},
    xlabel style={below},      
    ylabel={$L_2$ norm},
    ylabel style={below},    
    xtick=data,
    xticklabels={$10^{-4}$, $10^{-2}$, $10^{-0}$, $10^{2}$, $10^{4}$}    
]
\addplot[mark=o] coordinates {(0, 920) (1, 9) (2, 0.23) (3, 0.22) (4, 0.22)}; 
\addplot[mark=diamond] coordinates {(0, 1.6e3) (1, 15.4) (2, 0.59) (3, 0.58) (4, 0.58)}; 
\addplot[mark=triangle] coordinates {(0, 1.0) (1, 1.0) (2, 1.2) (3, 1.2) (4, 1.2)}; 
\end{axis}
\end{tikzpicture}
}
\caption{Case 3}
\label{L2_norm_case_3}
\end{subfigure}

\hspace{2.5cm}
\begin{subfigure}[b]{0.3\textwidth}
\scalebox{0.9}{
\begin{tikzpicture} 
\begin{axis}
[
    ymode=log,
    ymin=1e-8,
    ymax=1e4,
    ytick={1e-8, 1e-4, 1e0, 1e4},
    yticklabels={1e-8, 1e-4, 1e0, 1e4},    
    height=5cm,
    width=6cm,
    xlabel={$c_4$},
    xlabel style={below},      
    ylabel={$L_2$ norm},
    ylabel style={below},    
    xtick=data,
    xticklabels={$10^{-2}$, $10^{-1}$, $10^{-0}$, $10^{1}$, $10^{2}$}    
]
\addplot[mark=o] coordinates {(0, 0.58) (1, 0.8) (2, 1.1e-1) (3, 1.1e-2) (4, 1.1e-3)}; 
\addplot[mark=diamond] coordinates {(0, 1.0) (1, 0.94) (2, 0.71) (3, 0.71) (4, 0.71)}; 
\addplot[mark=triangle] coordinates {(0, 2.3e-3) (1, 0.22) (2, 4.4) (3, 4.4e1) (4, 4.4e2)}; 
\end{axis}
\end{tikzpicture}
}
\caption{Case 4}
\label{L2_norm_case_4}
\end{subfigure}
\hspace{0.15cm}
\begin{subfigure}[b]{0.3\textwidth}
\scalebox{0.9}{
\begin{tikzpicture} 
\begin{axis}
[
    ymode=log,
    ymin=1e-8,
    ymax=1e4,
    ytick={1e-8, 1e-4, 1e0, 1e4},
    yticklabels={1e-8, 1e-4, 1e0, 1e4},    
    height=5cm,
    width=6cm,
    xlabel={$c_5$},
    xlabel style={below},      
    ylabel={$L_2$ norm},
    ylabel style={below},    
    xtick=data,
    xticklabels={$10^{-4}$, $10^{-2}$, $10^{0}$, $10^{2}$, $10^{4}$} 
]
\addplot[mark=o] coordinates {(0, 5800) (1, 58) (2, 0.58) (3, 5.8e-3) (4, 5.8e-5)}; 
\addplot[mark=diamond] coordinates {(0, 1e4) (1, 1e2) (2, 1) (3, 1e-2) (4, 1e-4)}; 
\addplot[mark=triangle] coordinates {(0, 0.0) (1, 0.0) (2, 0.0) (3, 0.0) (4, 0.0)}; 
\end{axis}
\end{tikzpicture}
}
\caption{Case 5}
\label{L2_norm_case_5}
\end{subfigure}
\caption{$L_2$ norms of $u$, $u_{x}$ and $u_{xx}$ for different cases in Table \ref{scaling_cases_Poisson}.}
\label{L2_norms_cases_1_to_5}
\end{figure}
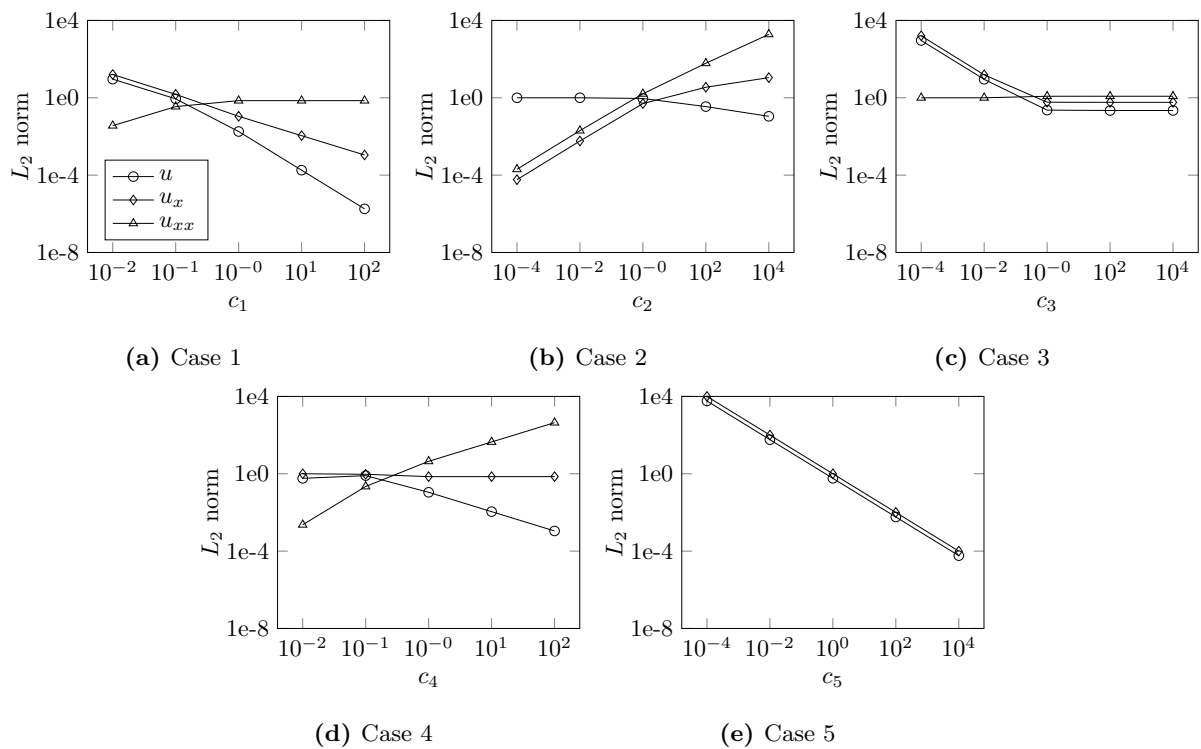

The $L_2$ norms of $u$, $u_x$ and/or $u_{xx}(f)$ for different cases are shown in Fig.~\ref{L2_norms_cases_1_to_5}.

\newpage

\subsection{Absolute errors}

\subsubsection{The standard FEM}

\begin{figure}[!ht]
    \begin{subfigure}{5.5cm}
        \includegraphics[width=1.0\linewidth]{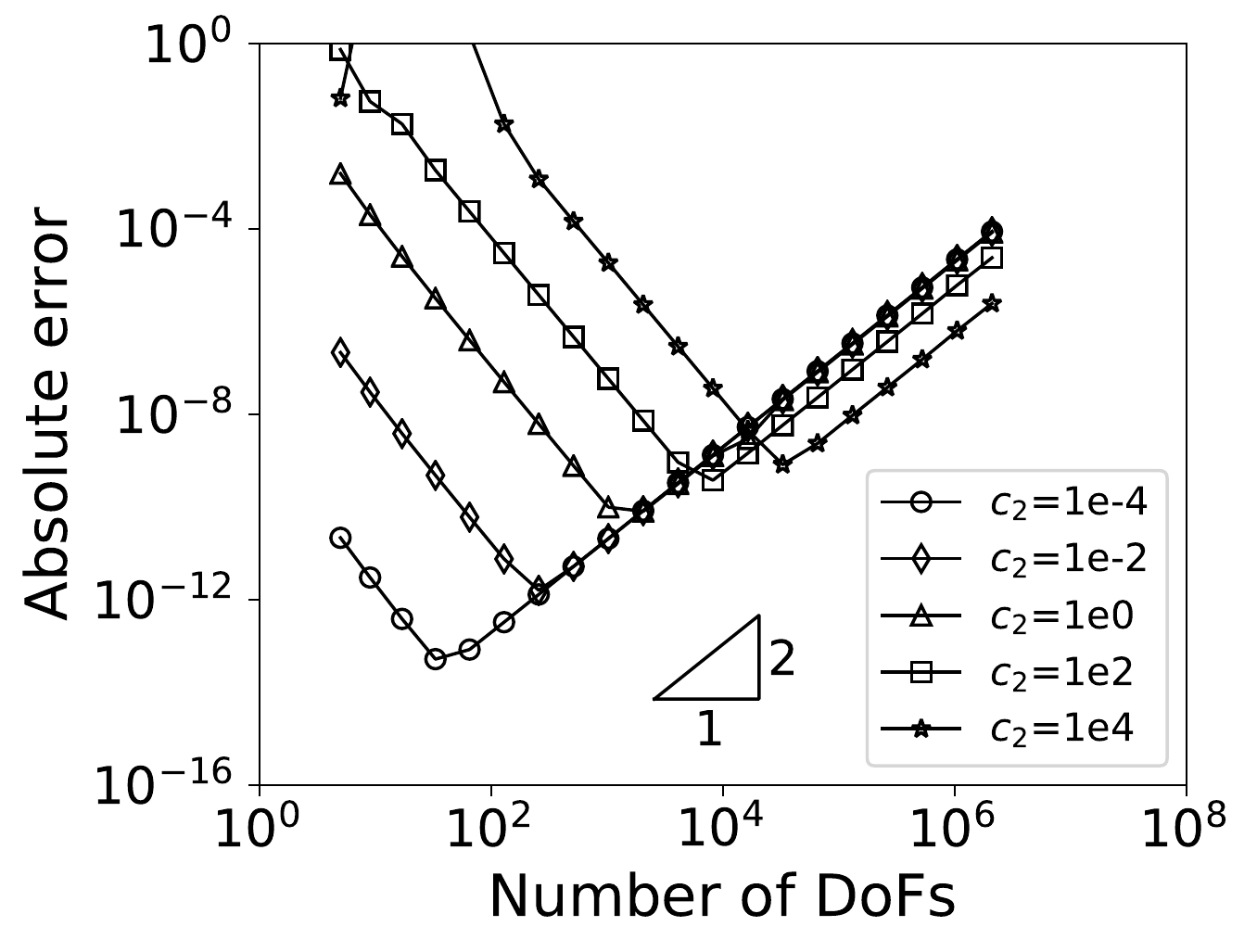}
        \caption{Solution}
        \label{py_L2_Pois2_SM_scaling_no_solu}
    \end{subfigure}
    \hspace{-0.2cm}
    \begin{subfigure}{5.5cm}
        \includegraphics[width=1.0\linewidth]{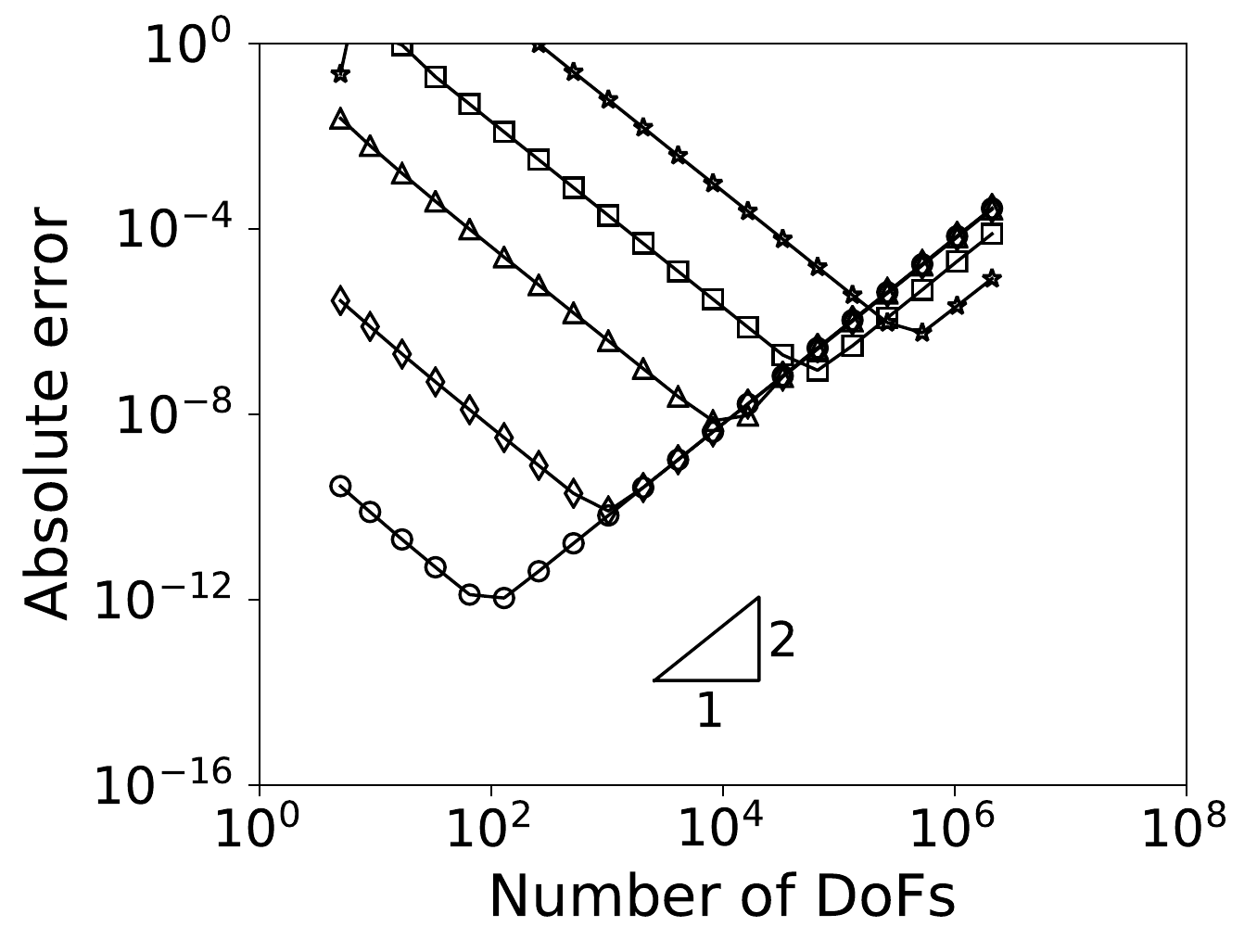}
        \caption{First derivative}
        \label{py_L2_Pois2_SM_scaling_no_grad}
    \end{subfigure}
    \hspace{-0.2cm}
    \begin{subfigure}{5.5cm}
        \includegraphics[width=1.0\linewidth]{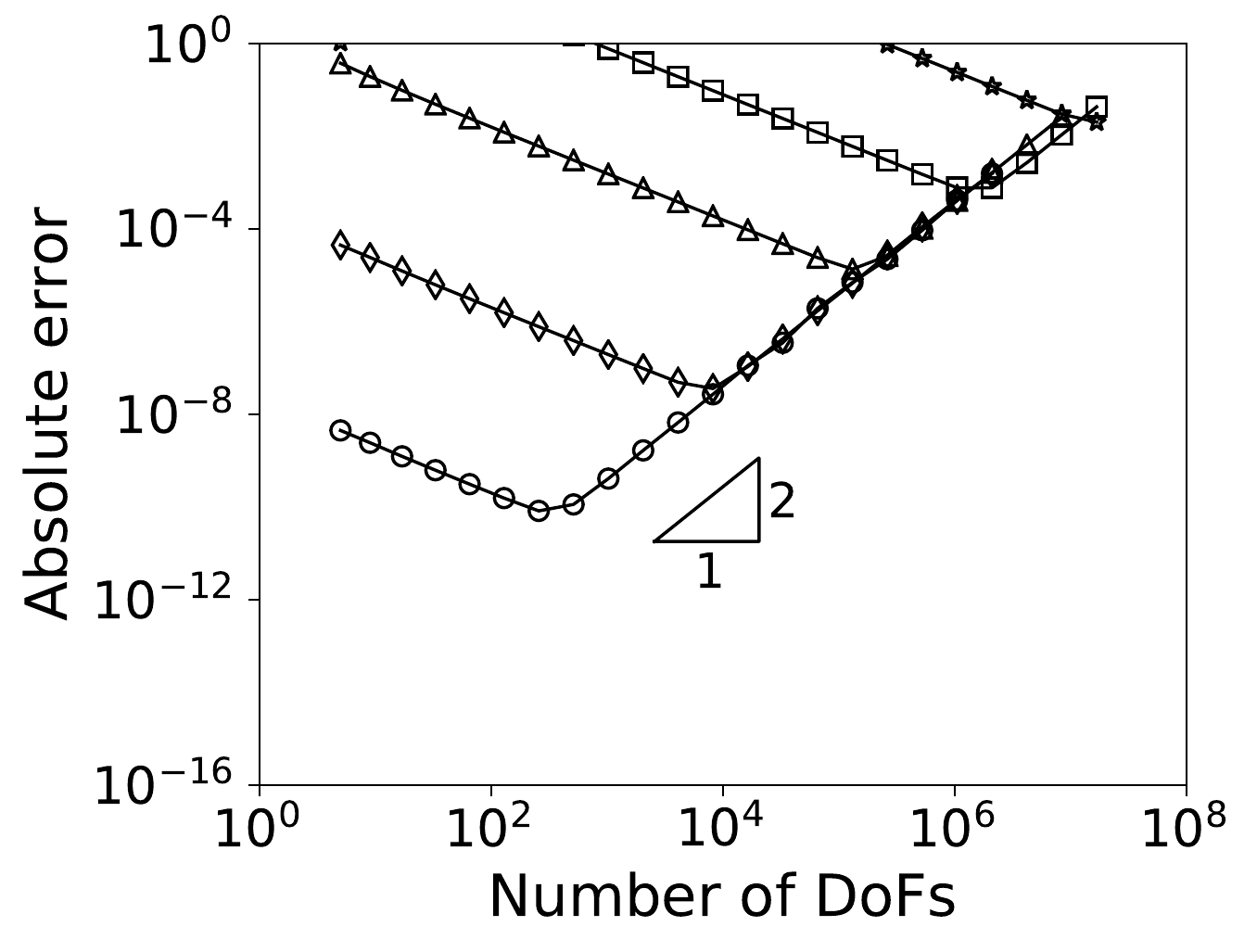}
        \caption{Second derivative}
        \label{py_L2_Pois2_SM_scaling_no_2ndd}
    \end{subfigure}
\caption{Absolute errors for Case 2 in Table \ref{scaling_cases_Poisson} using the standard FEM without scaling the right-hand side.}
\label{py_L2_Pois2_SM_scaling_no}
\end{figure}

\begin{figure}[!ht]
    \begin{subfigure}{5.5cm}
        \includegraphics[width=1.0\linewidth]{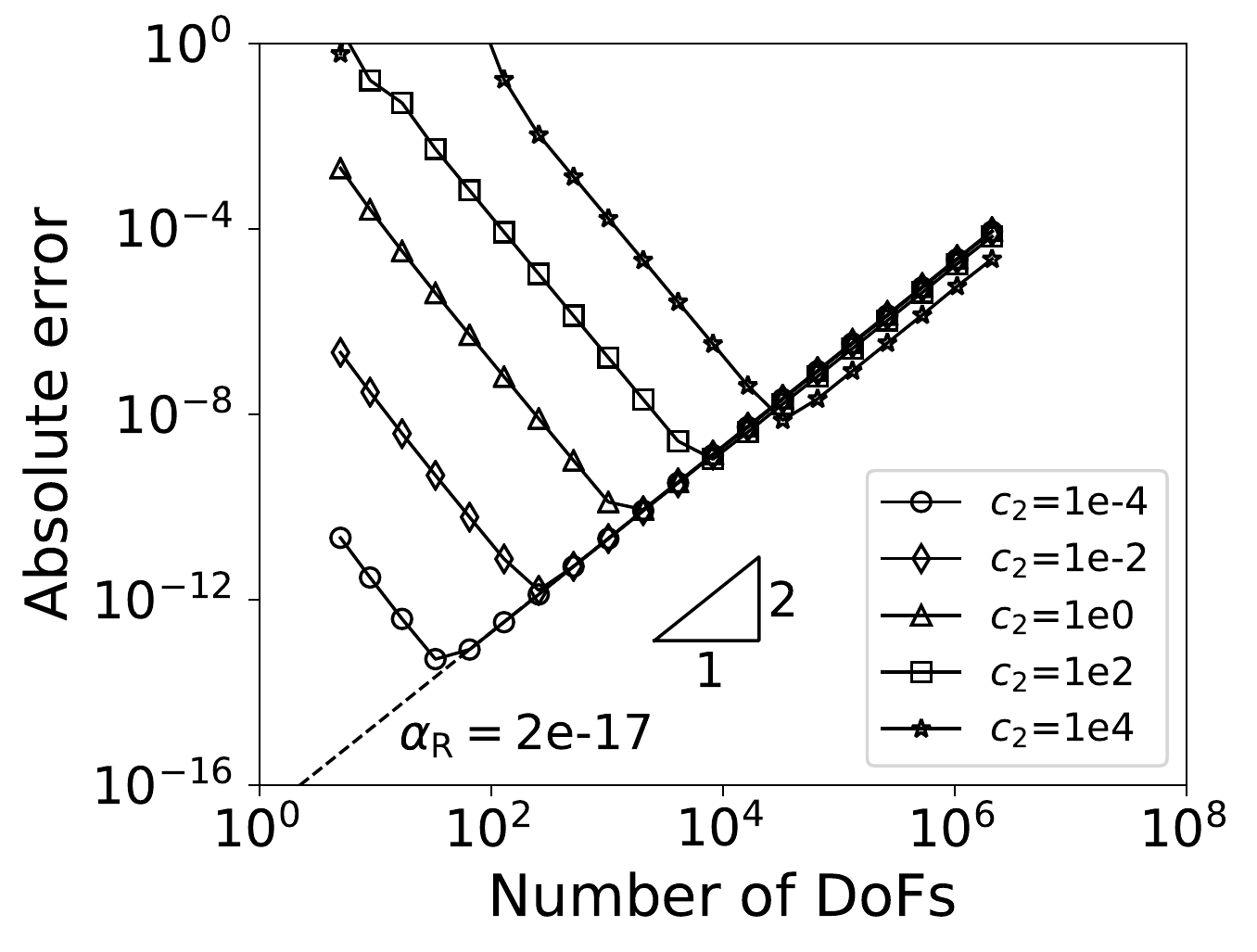}
        \caption{Solution}
        \label{py_L2_Pois2_SM_scaling_S_solu}
    \end{subfigure}
    \hspace{-0.2cm}
    \begin{subfigure}{5.5cm}
        \includegraphics[width=1.0\linewidth]{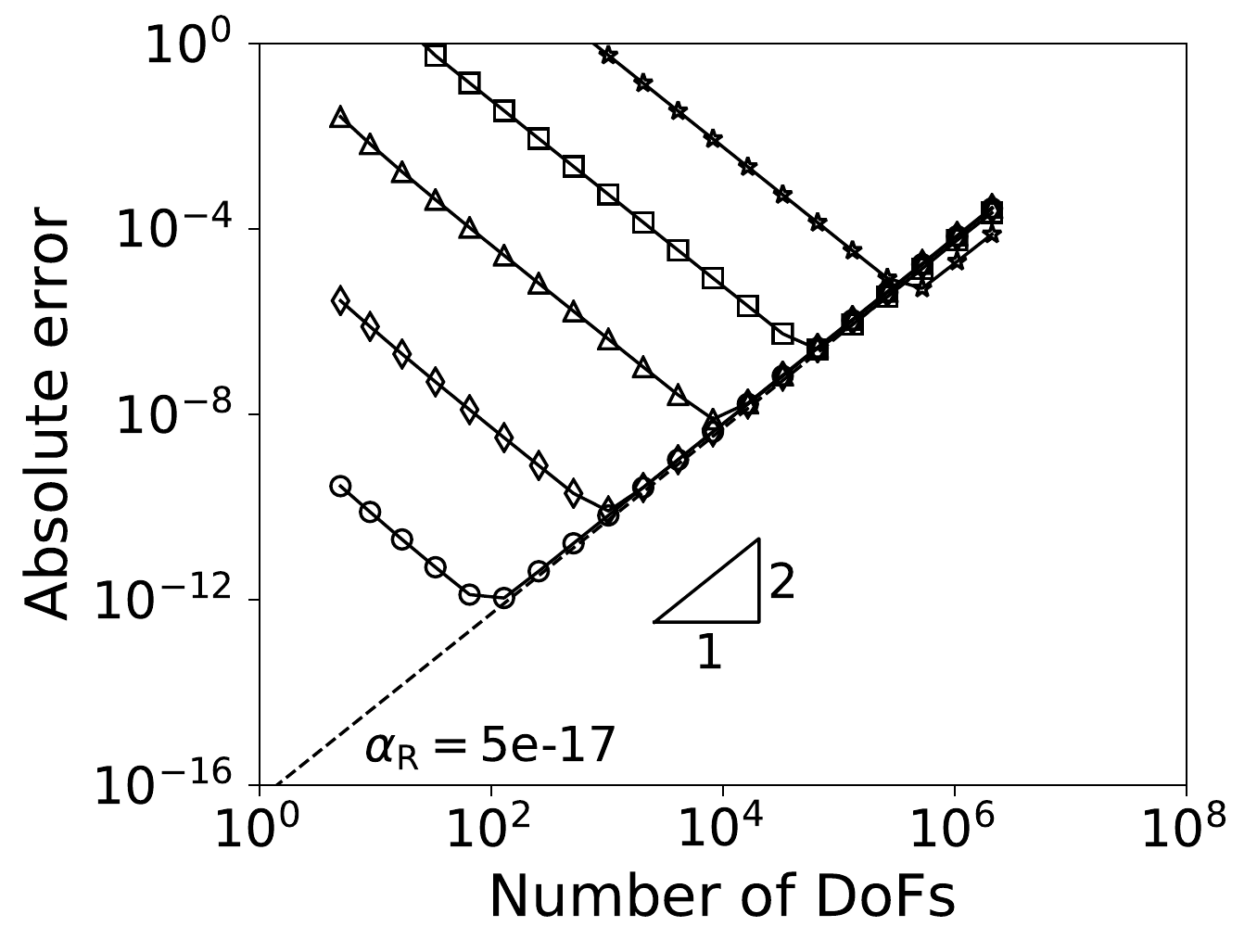}
        \caption{First derivative}
        \label{py_L2_Pois2_SM_scaling_S_grad}
    \end{subfigure}
    \hspace{-0.2cm}
    \begin{subfigure}{5.5cm}
        \includegraphics[width=1.0\linewidth]{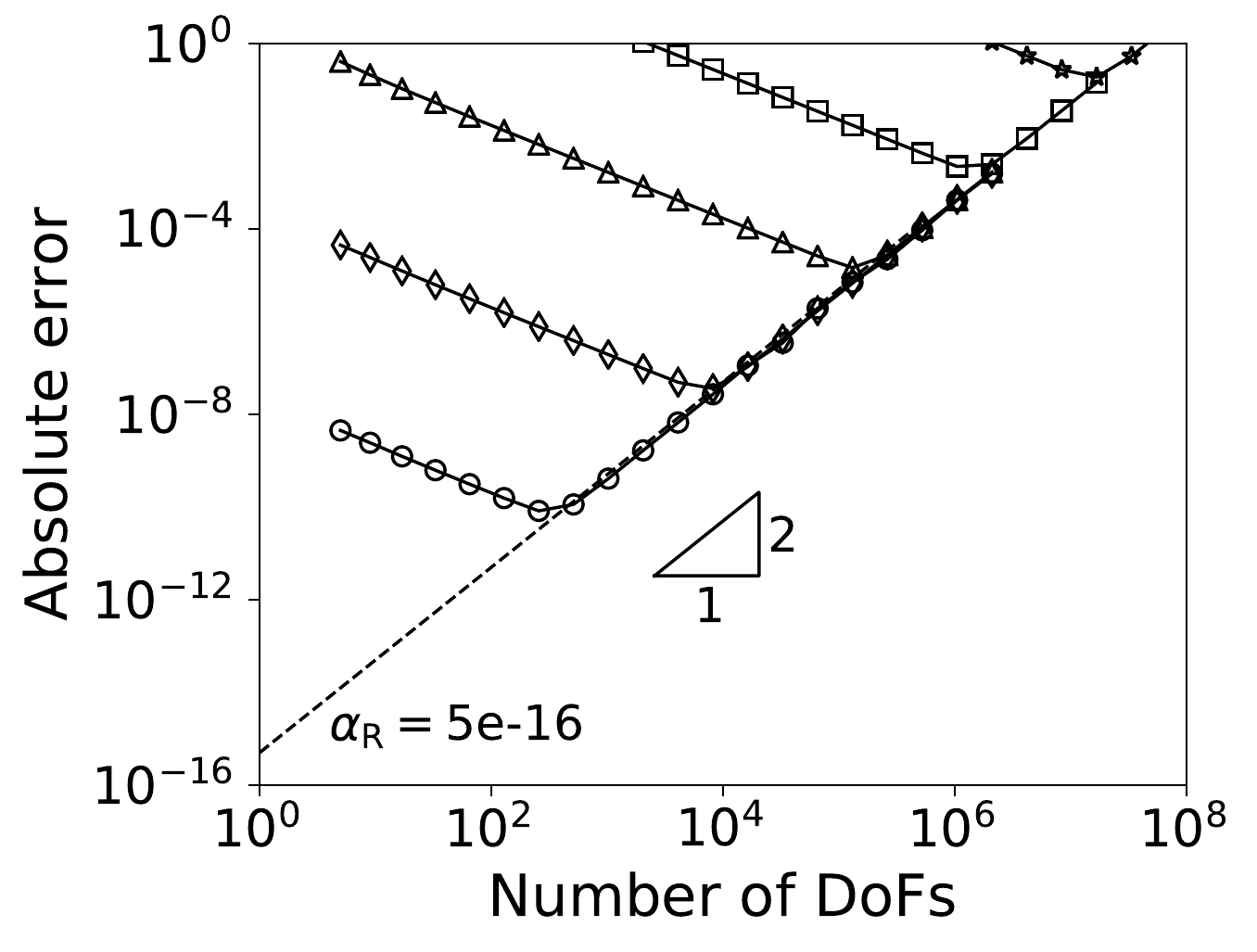}
        \caption{Second derivative}
        \label{py_L2_Pois2_SM_scaling_S_2ndd}
    \end{subfigure}
\caption{Absolute errors for Case 2 in Table \ref{scaling_cases_Poisson} using scheme $S$.}
\label{py_L2_Pois2_SM_scaling_S}
\end{figure}

\paragraph{Case 2}
For Case 2, using the standard FEM without scaling the right-hand side and scheme S, the absolute errors are shown in Figs. \ref{py_L2_Pois2_SM_scaling_no}--\ref{py_L2_Pois2_SM_scaling_S}.

\begin{figure}[!ht]
    \begin{subfigure}{5.5cm}
        \includegraphics[width=1.0\linewidth]{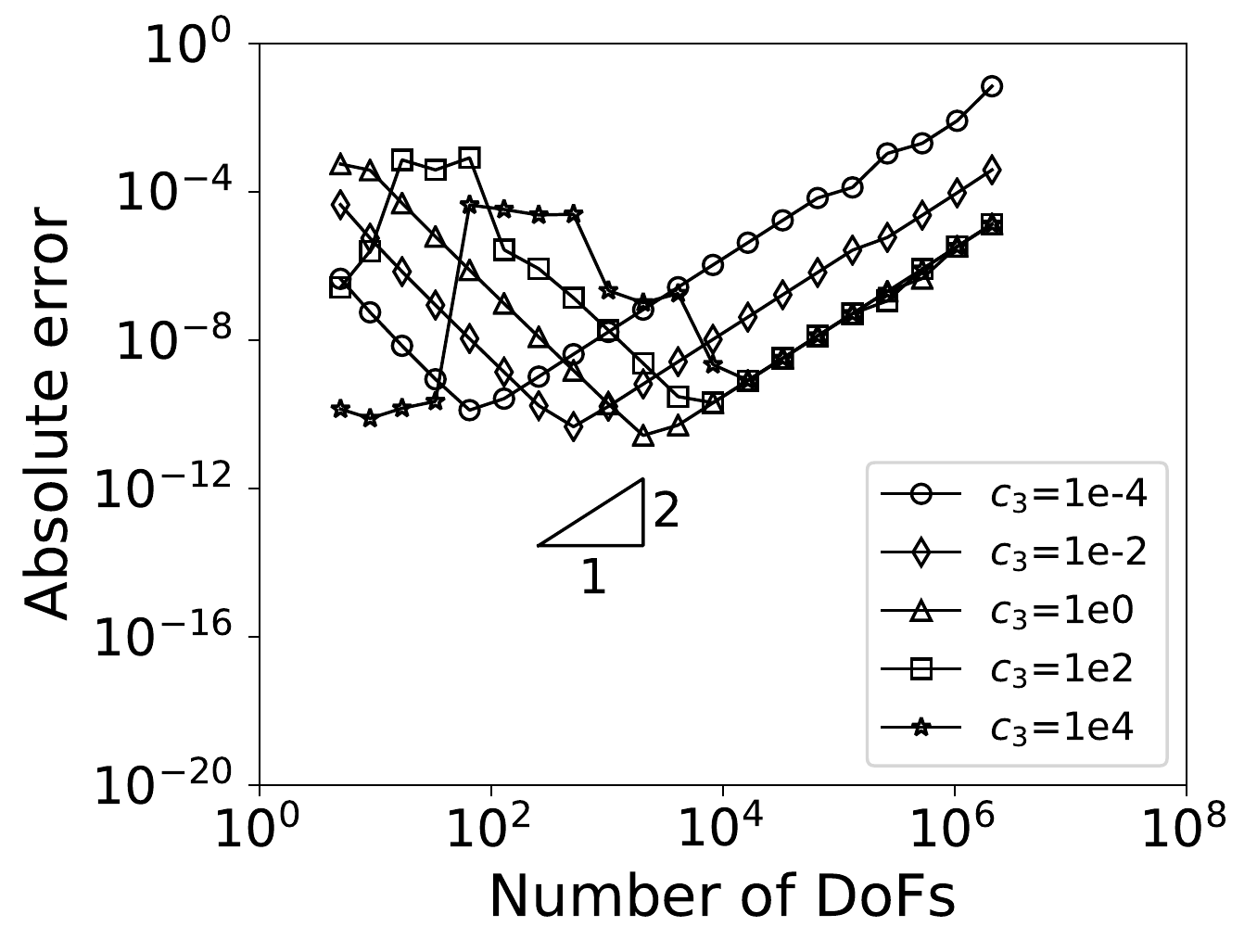}
        \caption{Solution}
        \label{py_L2_Pois3_SM_scaling_no_solu}
    \end{subfigure}
    \hspace{-0.2cm}
    \begin{subfigure}{5.5cm}
        \includegraphics[width=1.0\linewidth]{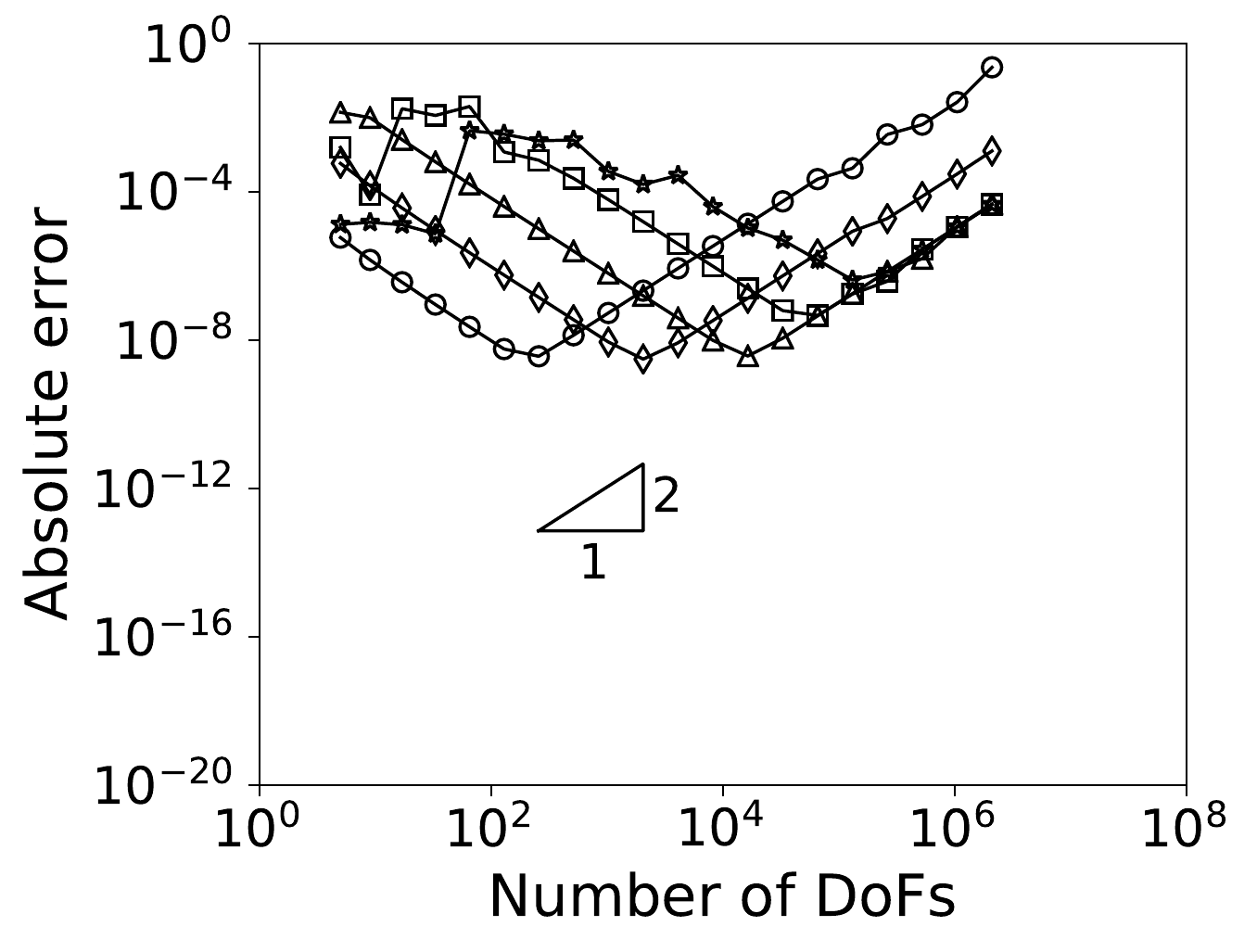}
        \caption{First derivative}
        \label{py_L2_Pois3_SM_scaling_no_grad}
    \end{subfigure}
    \hspace{-0.2cm}
    \begin{subfigure}{5.5cm}
        \includegraphics[width=1.0\linewidth]{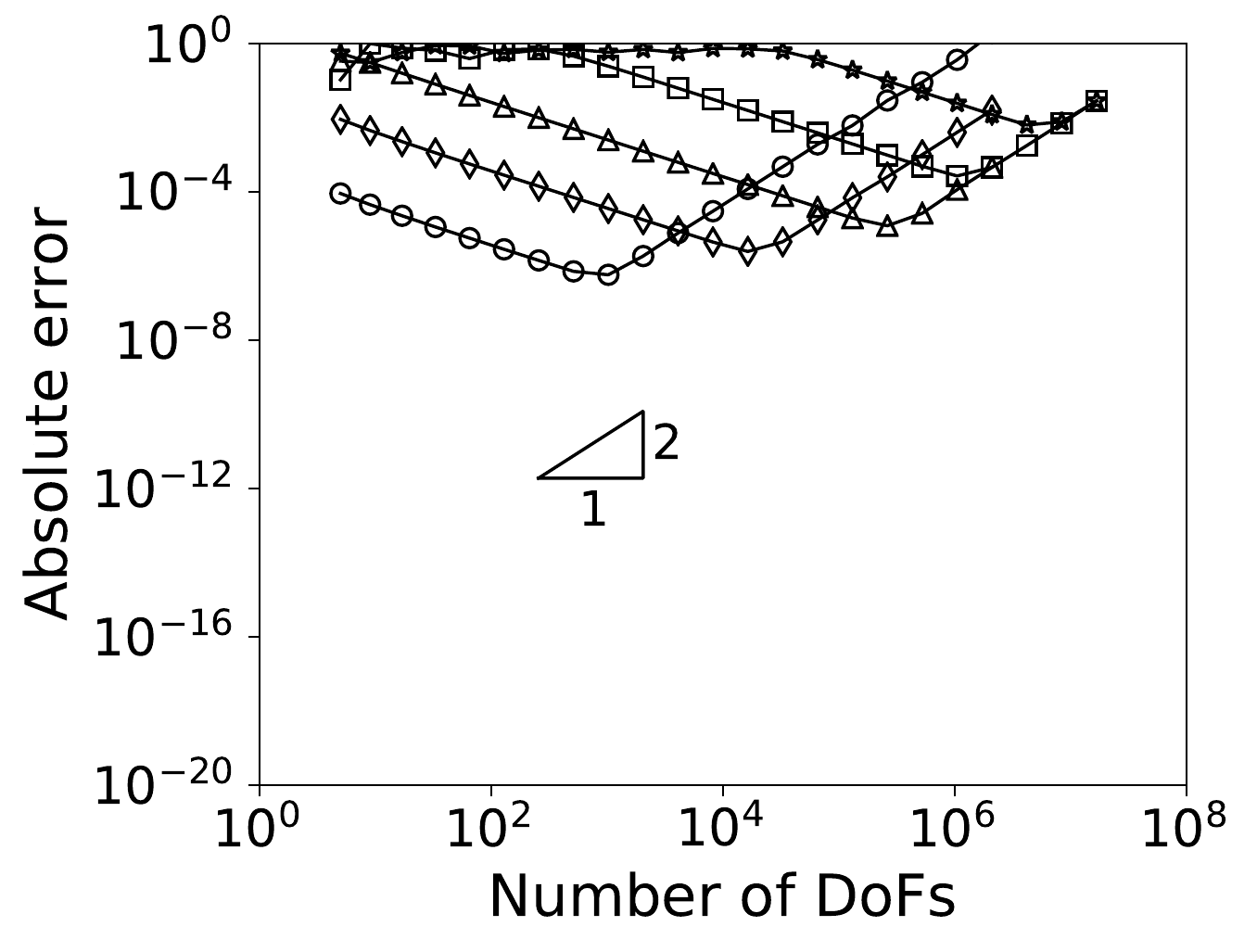}
        \caption{Second derivative}
        \label{py_L2_Pois3_SM_scaling_no_2ndd}
    \end{subfigure}
\caption{Absolute errors for Case 3 in Table \ref{scaling_cases_Poisson} using the standard FEM without scaling the right-hand side.}
\label{py_L2_Pois3_SM_scaling_no}
\end{figure}

\begin{figure}[!ht]
    \begin{subfigure}{5.5cm}
        \includegraphics[width=1.0\linewidth]{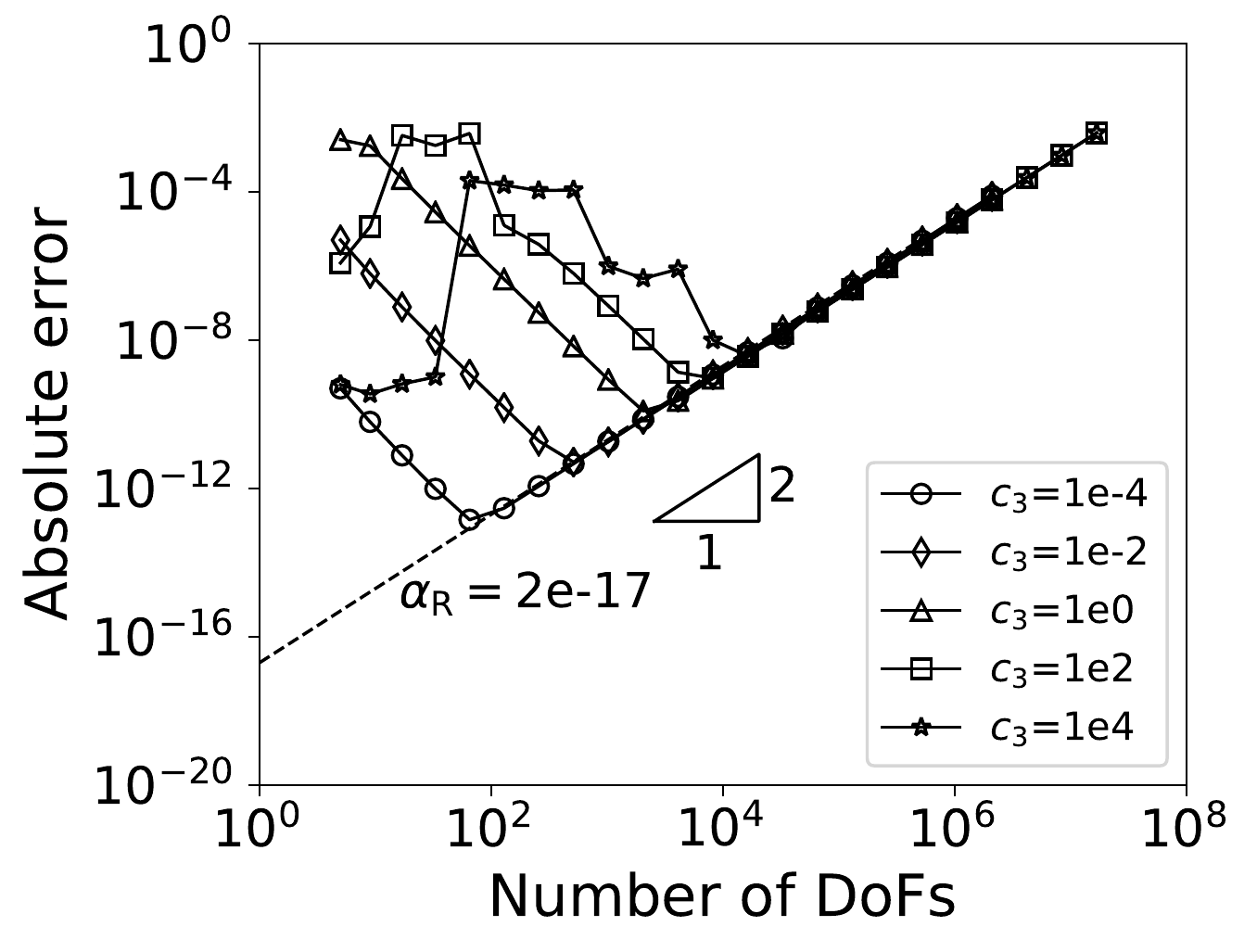}
        \caption{Solution}
        \label{py_L2_Pois3_SM_scaling_S_solu}
    \end{subfigure}
    \hspace{-0.2cm}
    \begin{subfigure}{5.5cm}
        \includegraphics[width=1.0\linewidth]{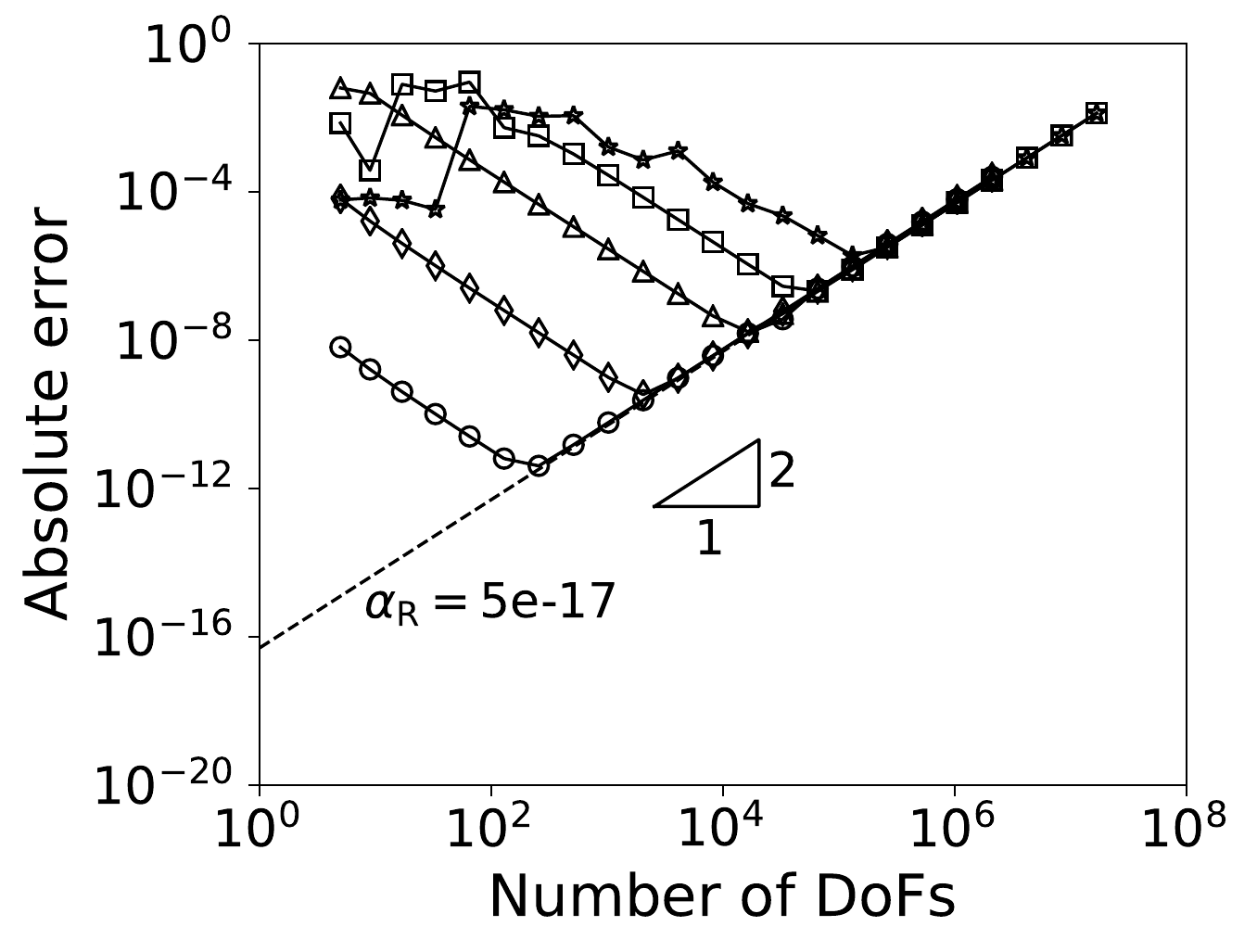}
        \caption{First derivative}
        \label{py_L2_Pois3_SM_scaling_S_grad}
    \end{subfigure}
    \hspace{-0.2cm}
    \begin{subfigure}{5.5cm}
        \includegraphics[width=1.0\linewidth]{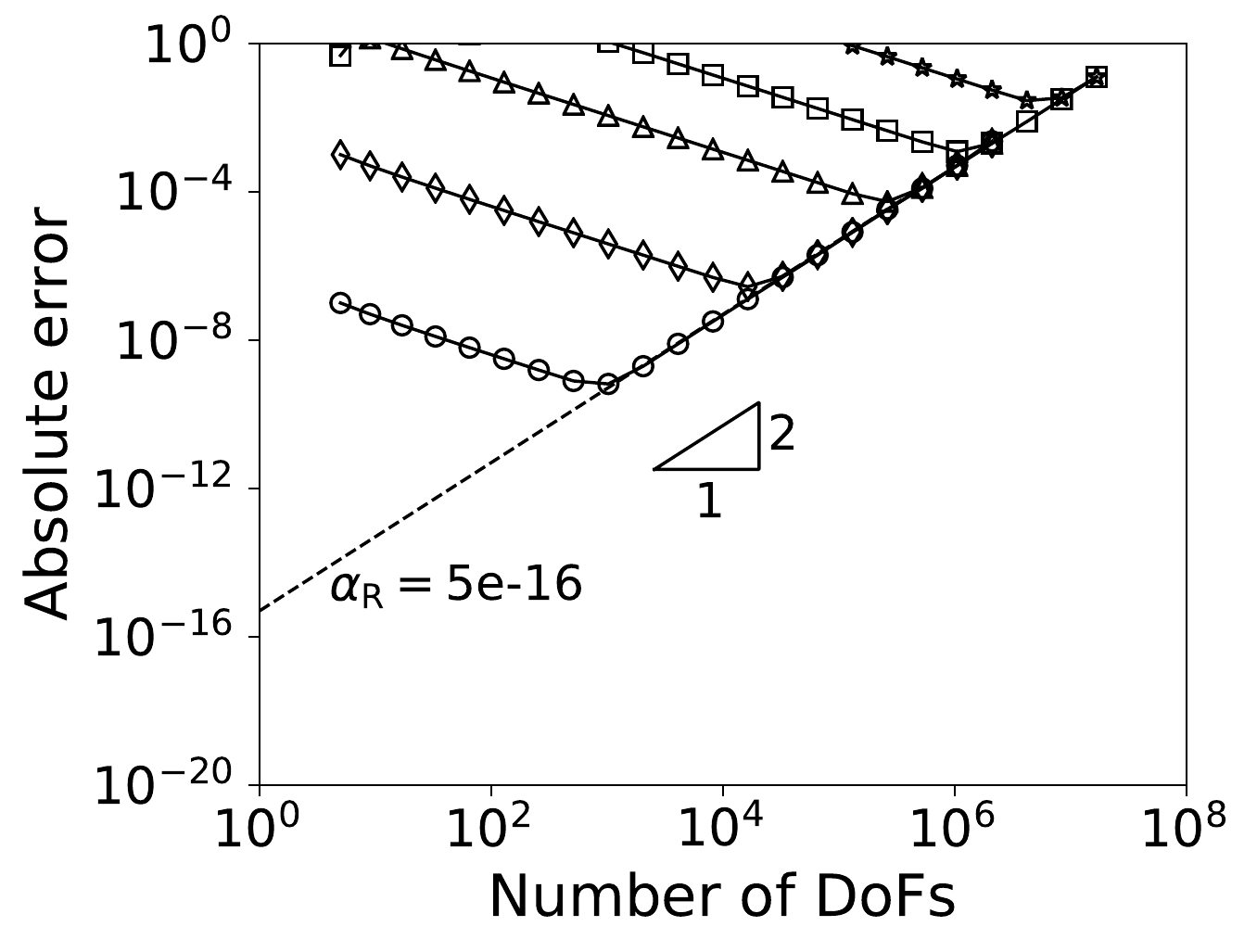}
        \caption{Second derivative}
        \label{py_L2_Pois3_SM_scaling_S_2ndd}
    \end{subfigure}
\caption{Absolute errors of Case 3 in Table \ref{scaling_cases_Poisson} using scheme $S$.}
\label{py_L2_Pois3_SM_scaling_S}
\end{figure}

\paragraph{Case 3}
For Case 3, using the standard FEM without scaling the right-hand side and scheme S, the absolute errors are shown in Figs. \ref{py_L2_Pois3_SM_scaling_no}--\ref{py_L2_Pois3_SM_scaling_S}.

\begin{figure}[!ht]
    \begin{subfigure}{5.5cm}
        \includegraphics[width=1.0\linewidth]{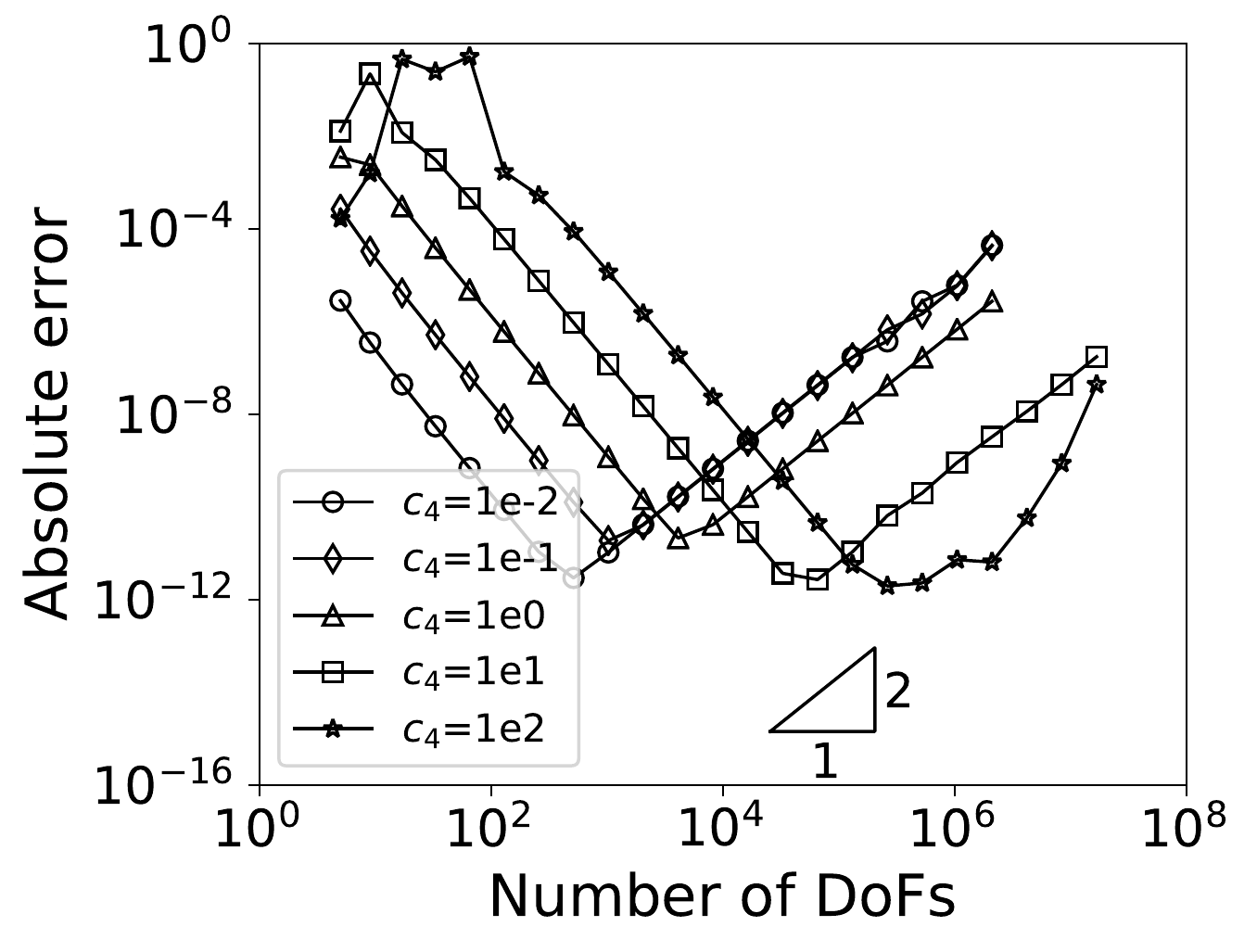}
        \caption{Solution}
        \label{py_L2_Pois4_SM_scaling_no_solu}
    \end{subfigure}
    \hspace{-0.2cm}
    \begin{subfigure}{5.5cm}
        \includegraphics[width=1.0\linewidth]{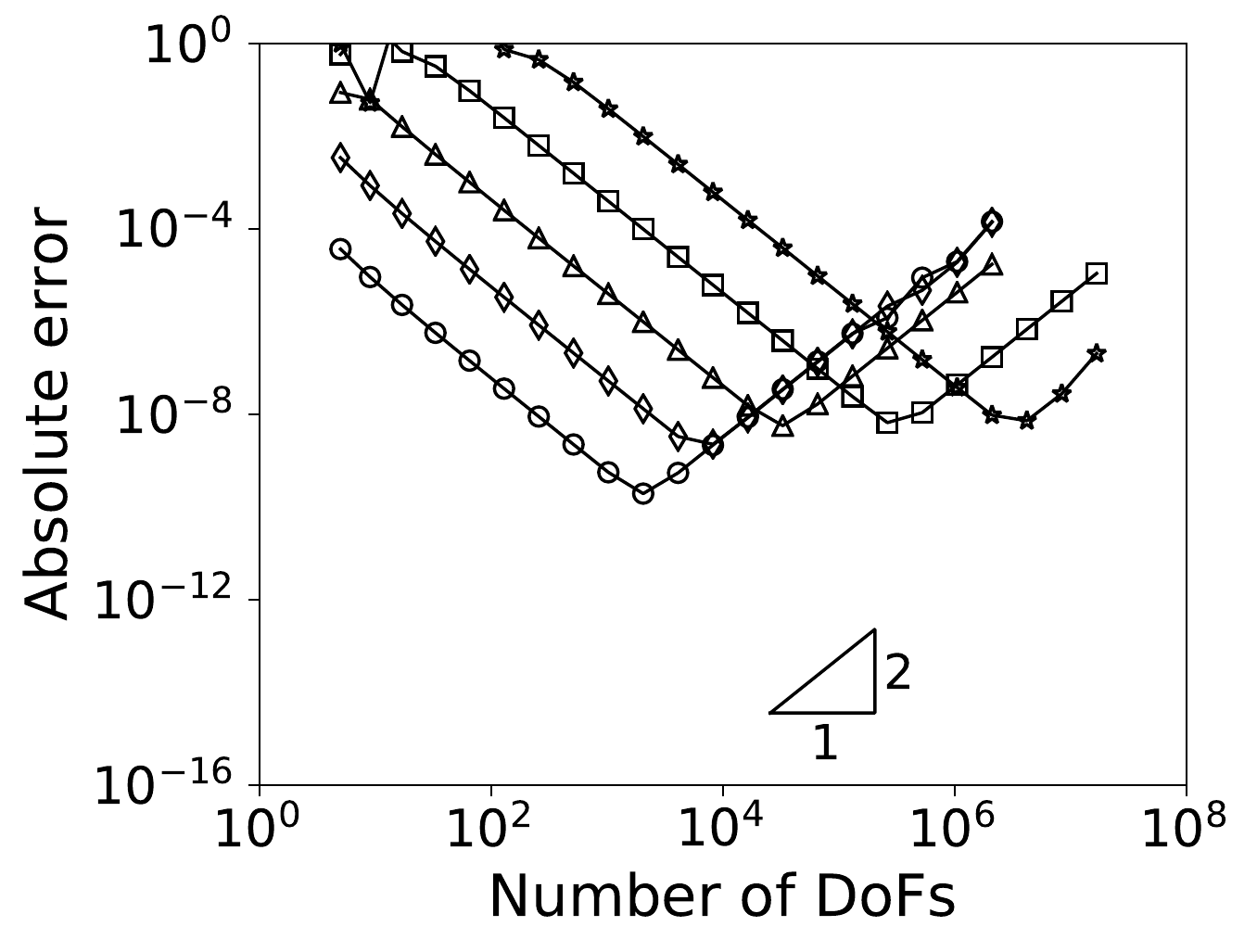}
        \caption{First derivative}
        \label{py_L2_Pois4_SM_scaling_no_grad}
    \end{subfigure}
    \hspace{-0.2cm}
    \begin{subfigure}{5.5cm}
        \includegraphics[width=1.0\linewidth]{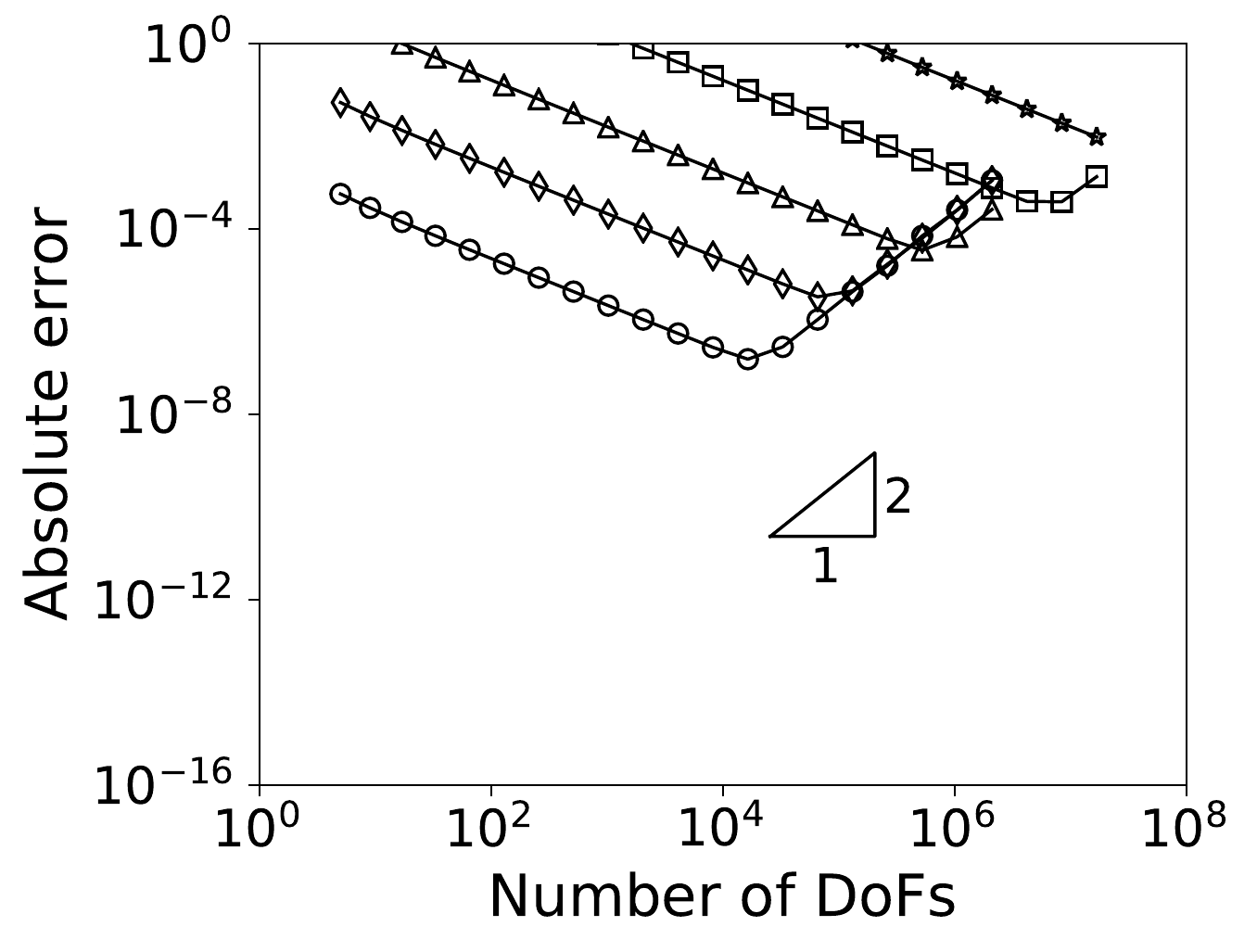}
        \caption{Second derivative}
        \label{py_L2_Pois4_SM_scaling_no_2ndd}
    \end{subfigure}
\caption{Absolute errors for Case 4 in Table \ref{scaling_cases_Poisson} using the standard FEM without scaling the right-hand side.}
\label{py_L2_Pois4_SM_scaling_no}
\end{figure}

\begin{figure}[!ht]
    \begin{subfigure}{5.5cm}
        \includegraphics[width=1.0\linewidth]{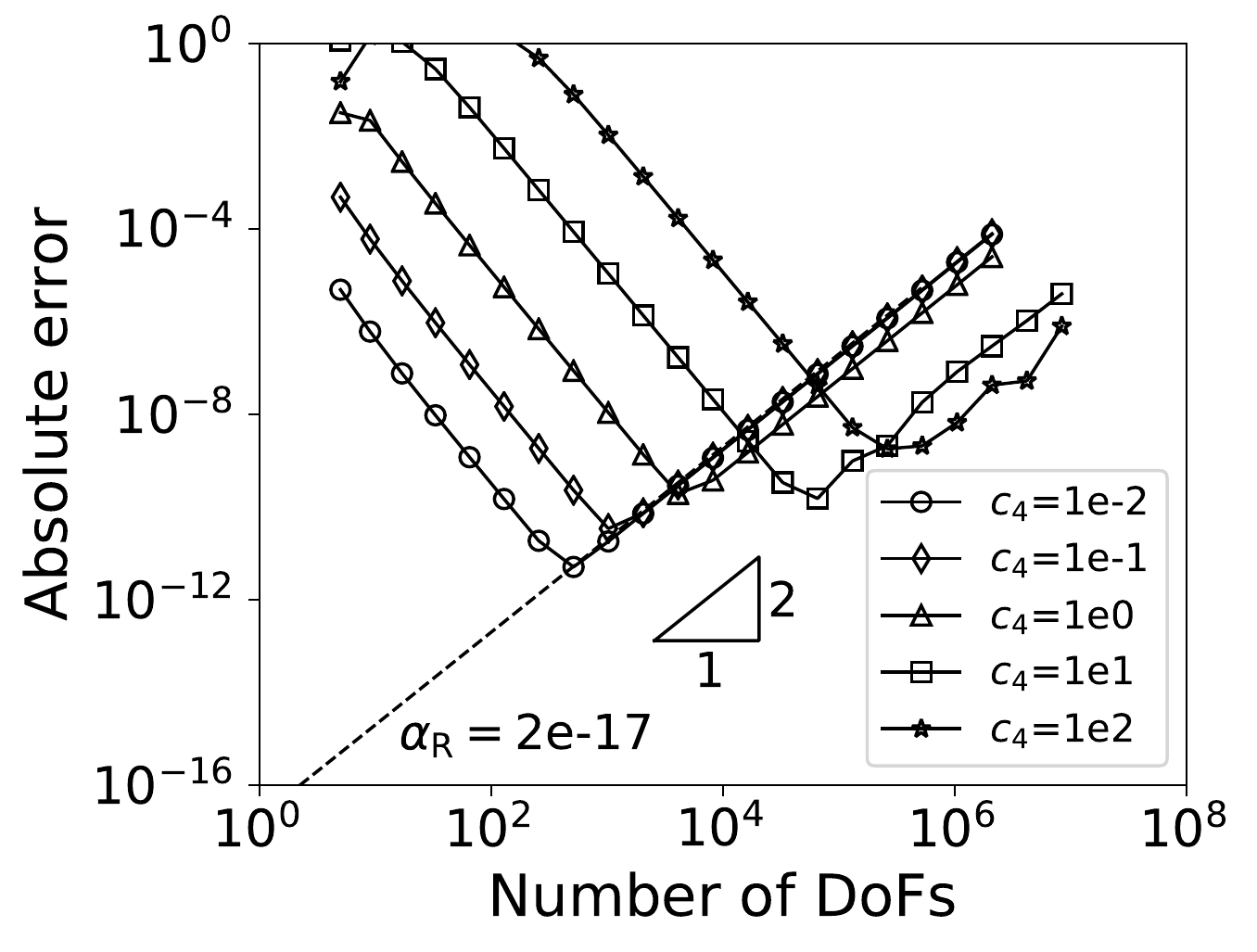}
        \caption{Solution}
        \label{py_L2_Pois4_SM_scaling_S_solu}
    \end{subfigure}
    \hspace{-0.2cm}
    \begin{subfigure}{5.5cm}
        \includegraphics[width=1.0\linewidth]{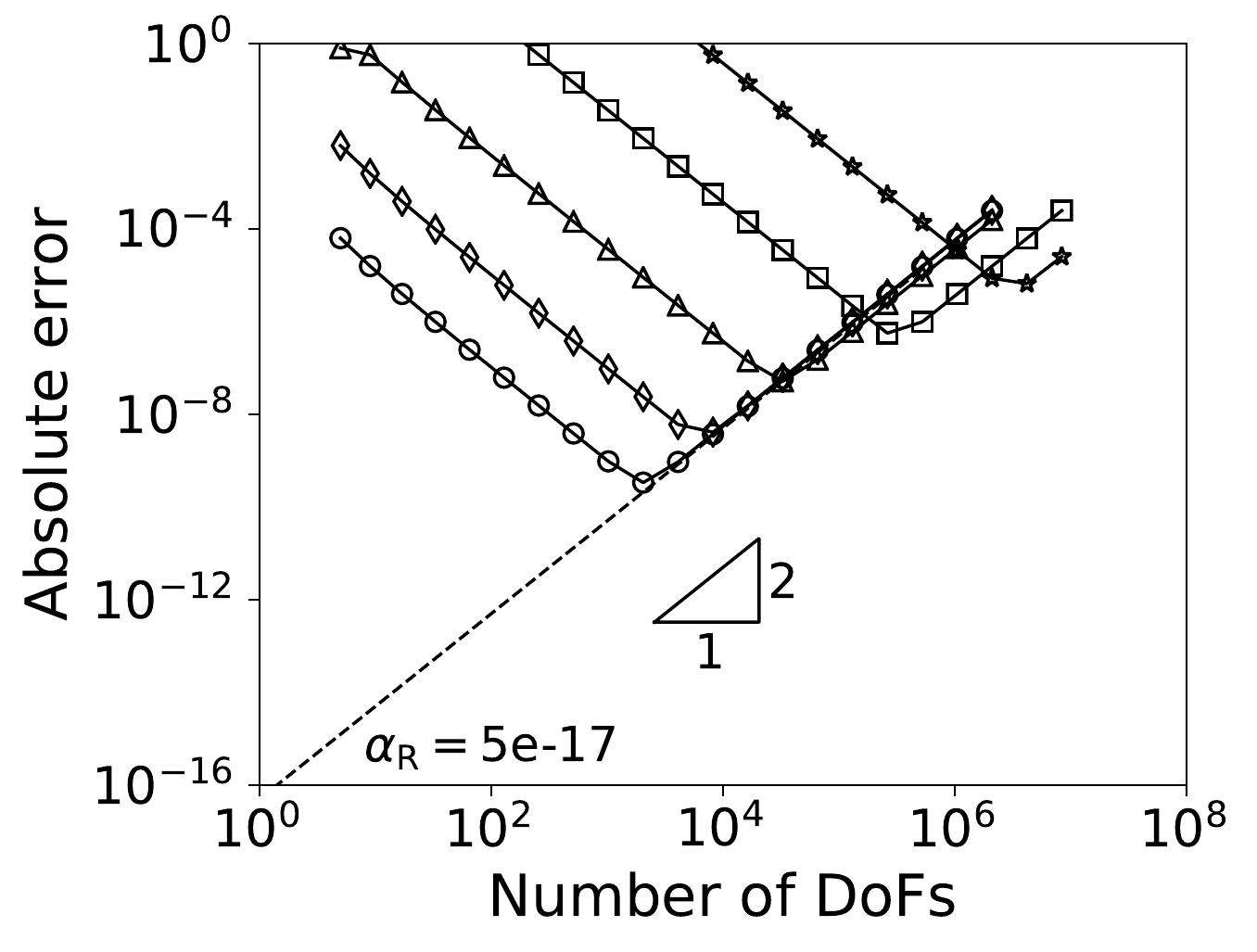}
        \caption{First derivative}
        \label{py_L2_Pois4_SM_scaling_S_grad}
    \end{subfigure}
    \hspace{-0.2cm}
    \begin{subfigure}{5.5cm}
        \includegraphics[width=1.0\linewidth]{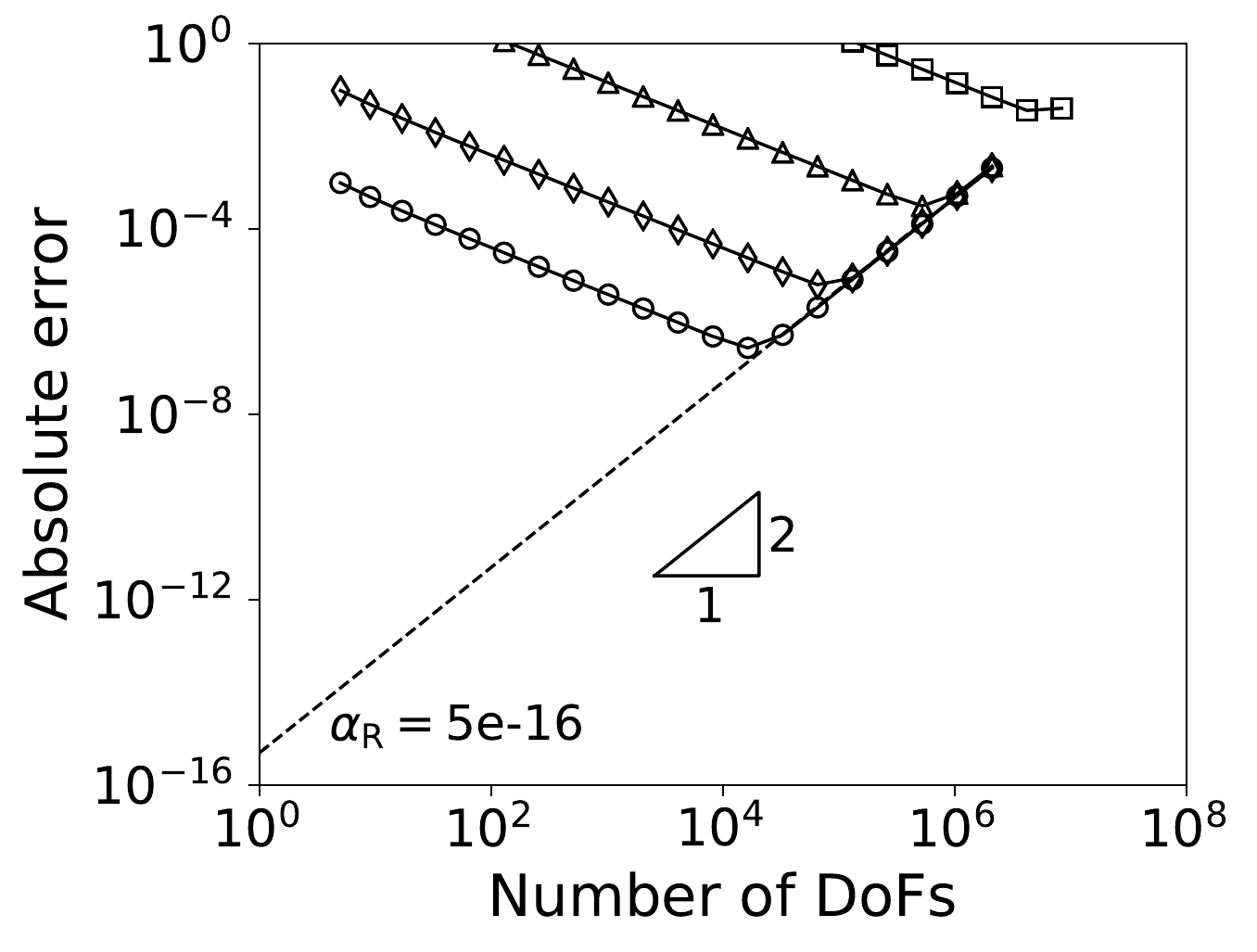}
        \caption{Second derivative}
        \label{py_L2_Pois4_SM_scaling_S_2ndd}
    \end{subfigure}
\caption{Absolute errors of Case 4 in Table \ref{scaling_cases_Poisson} using scheme $S$.}
\label{py_L2_Pois4_SM_scaling_S}
\end{figure}

\newpage
\paragraph{Case 4}
For Case 4, using the standard FEM without scaling the right-hand side and scheme S, the absolute errors are shown in Figs. \ref{py_L2_Pois4_SM_scaling_no}--\ref{py_L2_Pois4_SM_scaling_S}.

\begin{figure}[!ht]
    \begin{subfigure}{5.5cm}
        \includegraphics[width=1.0\linewidth]{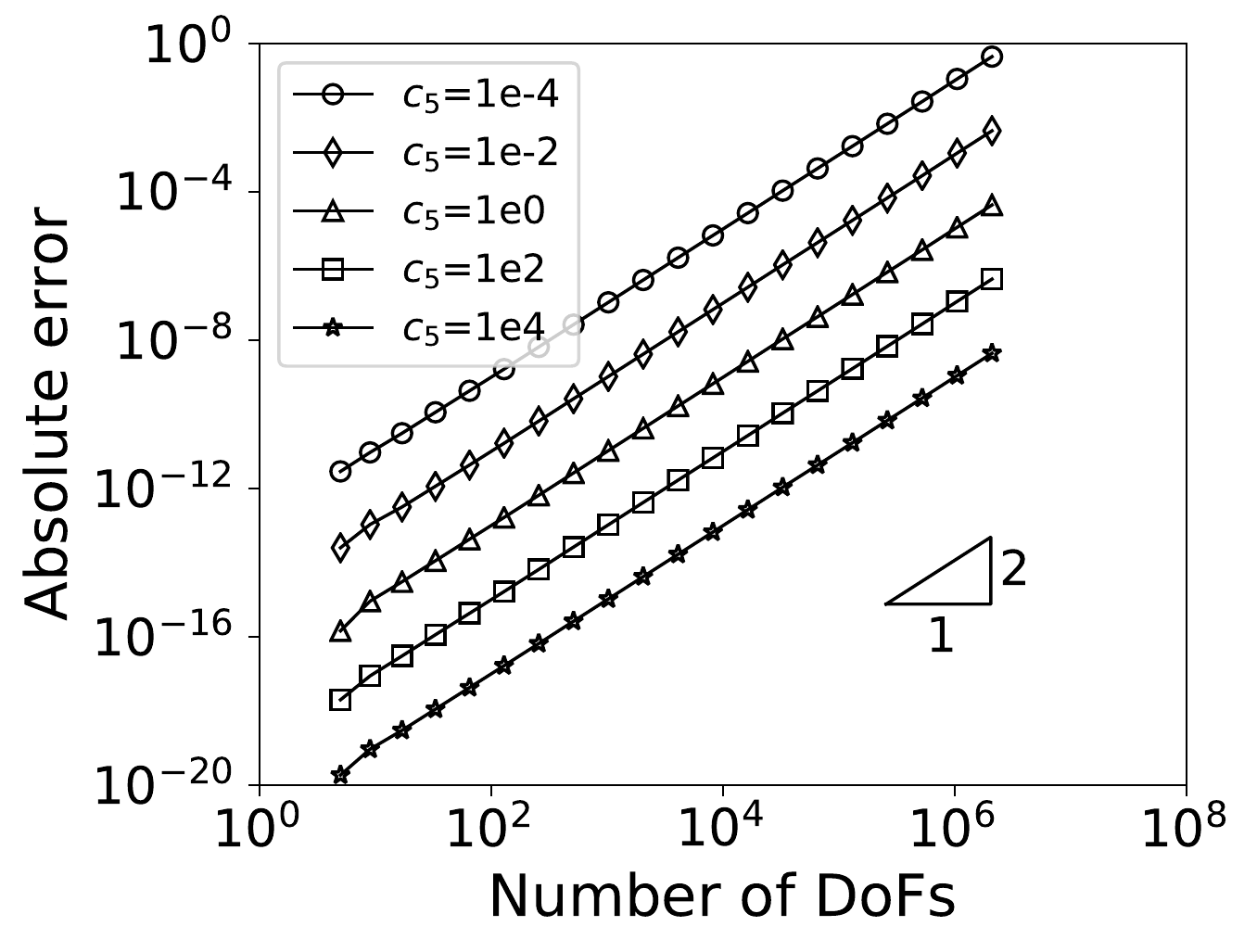}
        \caption{Solution}
        \label{py_L2_Pois5_SM_scaling_no_solu}
    \end{subfigure}
    \hspace{-0.2cm}
    \begin{subfigure}{5.5cm}
        \includegraphics[width=1.0\linewidth]{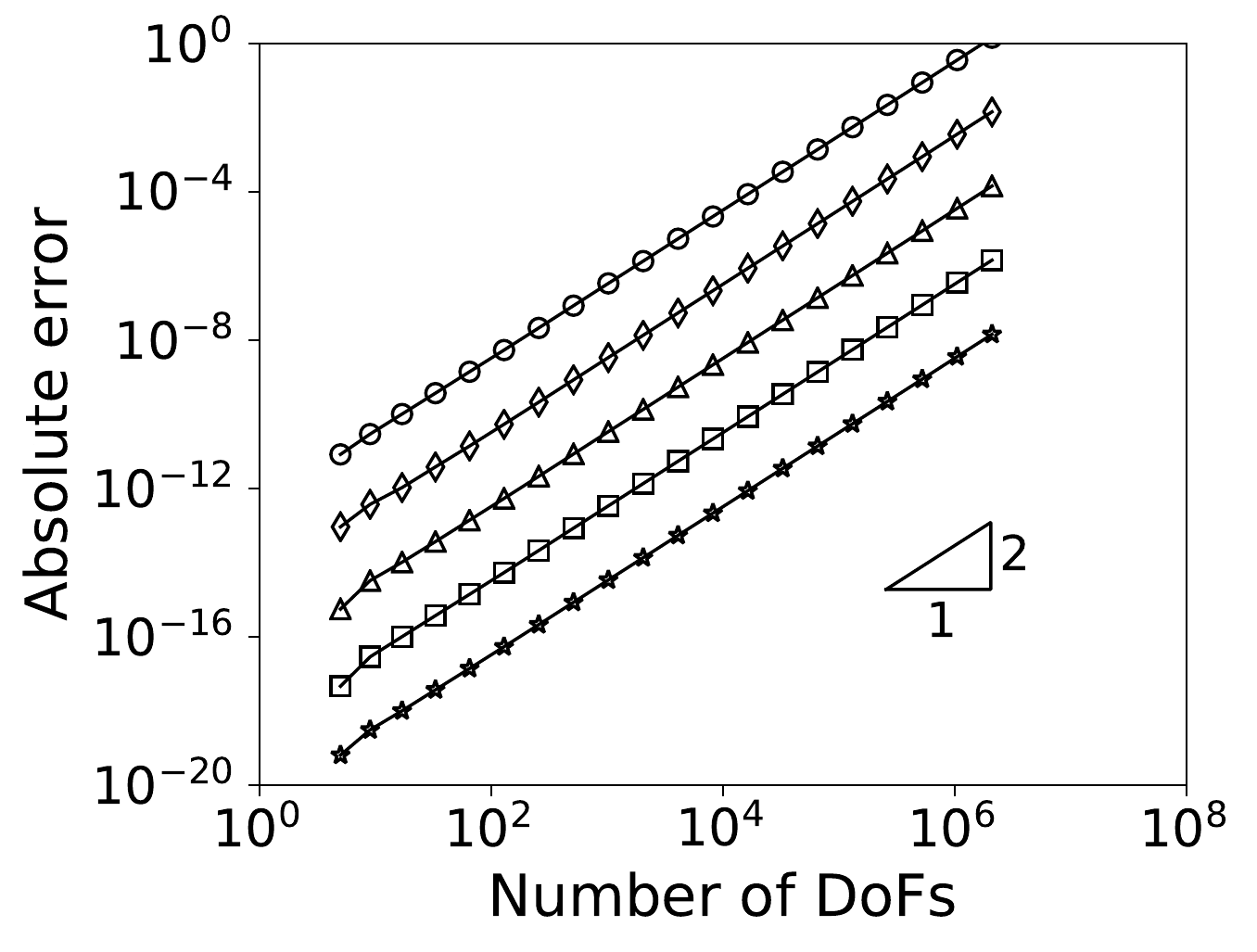}
        \caption{First derivative}
        \label{py_L2_Pois5_SM_scaling_no_grad}
    \end{subfigure}
    \hspace{-0.2cm}
    \begin{subfigure}{5.5cm}
        \includegraphics[width=1.0\linewidth]{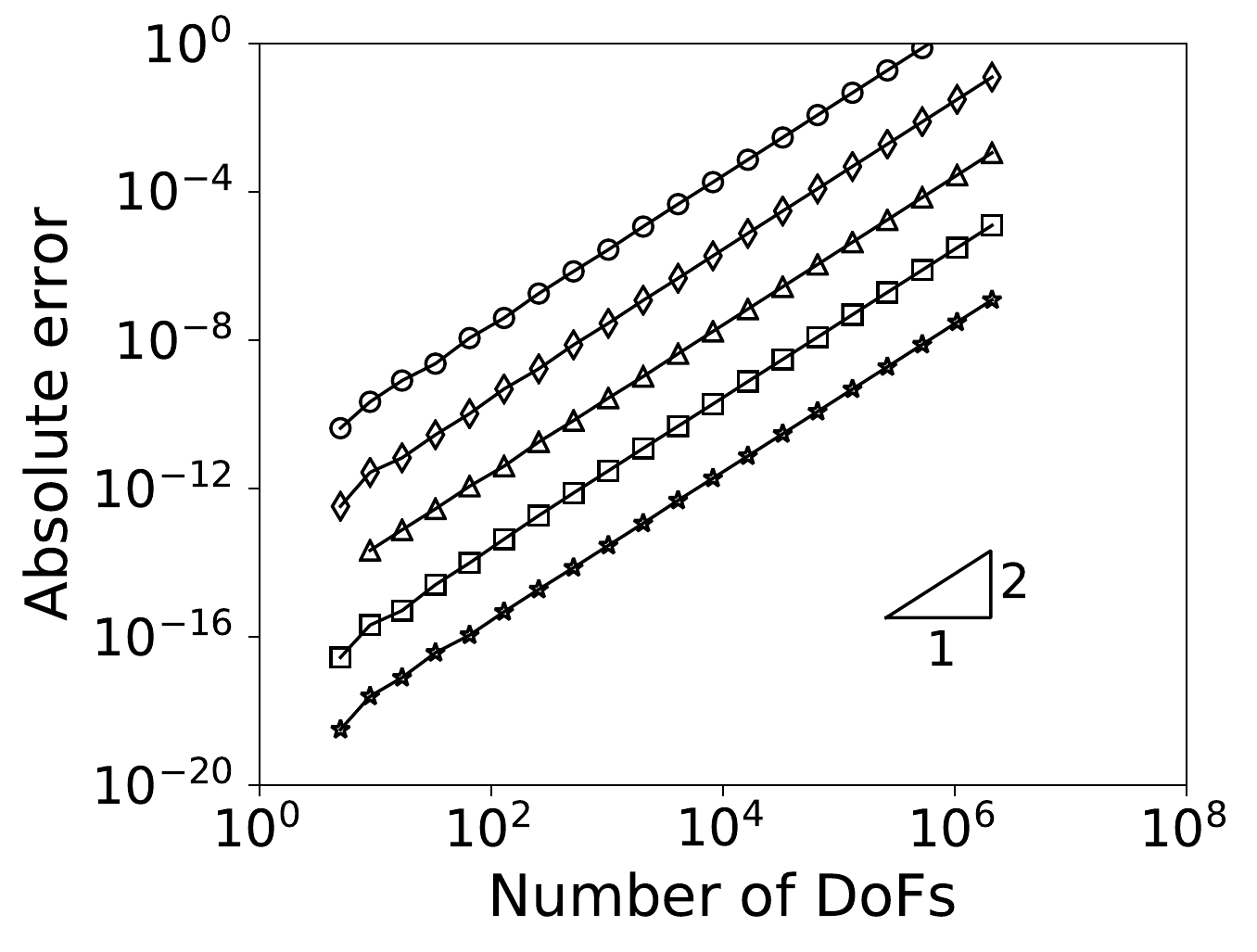}
        \caption{Second derivative}
        \label{py_L2_Pois5_SM_scaling_no_2ndd}
    \end{subfigure}
\caption{Absolute errors for Case 5 in Table \ref{scaling_cases_Poisson} using the standard FEM without scaling the right-hand side.}
\label{py_L2_Pois5_SM_scaling_no}
\end{figure}

\begin{figure}[!ht]
    \begin{subfigure}{5.5cm}
        \includegraphics[width=1.0\linewidth]{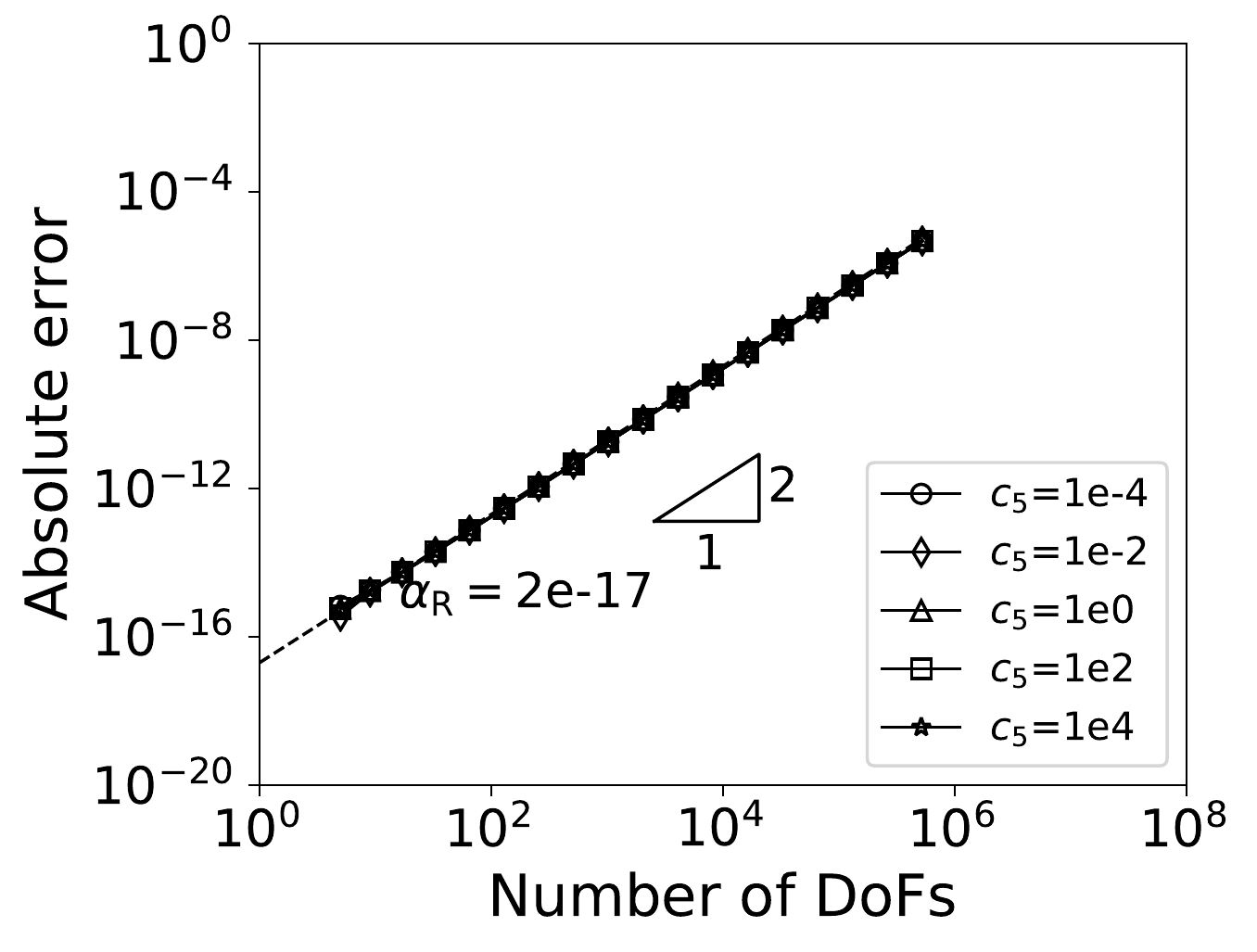}
        \caption{Solution}
        \label{py_L2_Pois5_SM_scaling_S_solu}
    \end{subfigure}
    \hspace{-0.2cm}
    \begin{subfigure}{5.5cm}
        \includegraphics[width=1.0\linewidth]{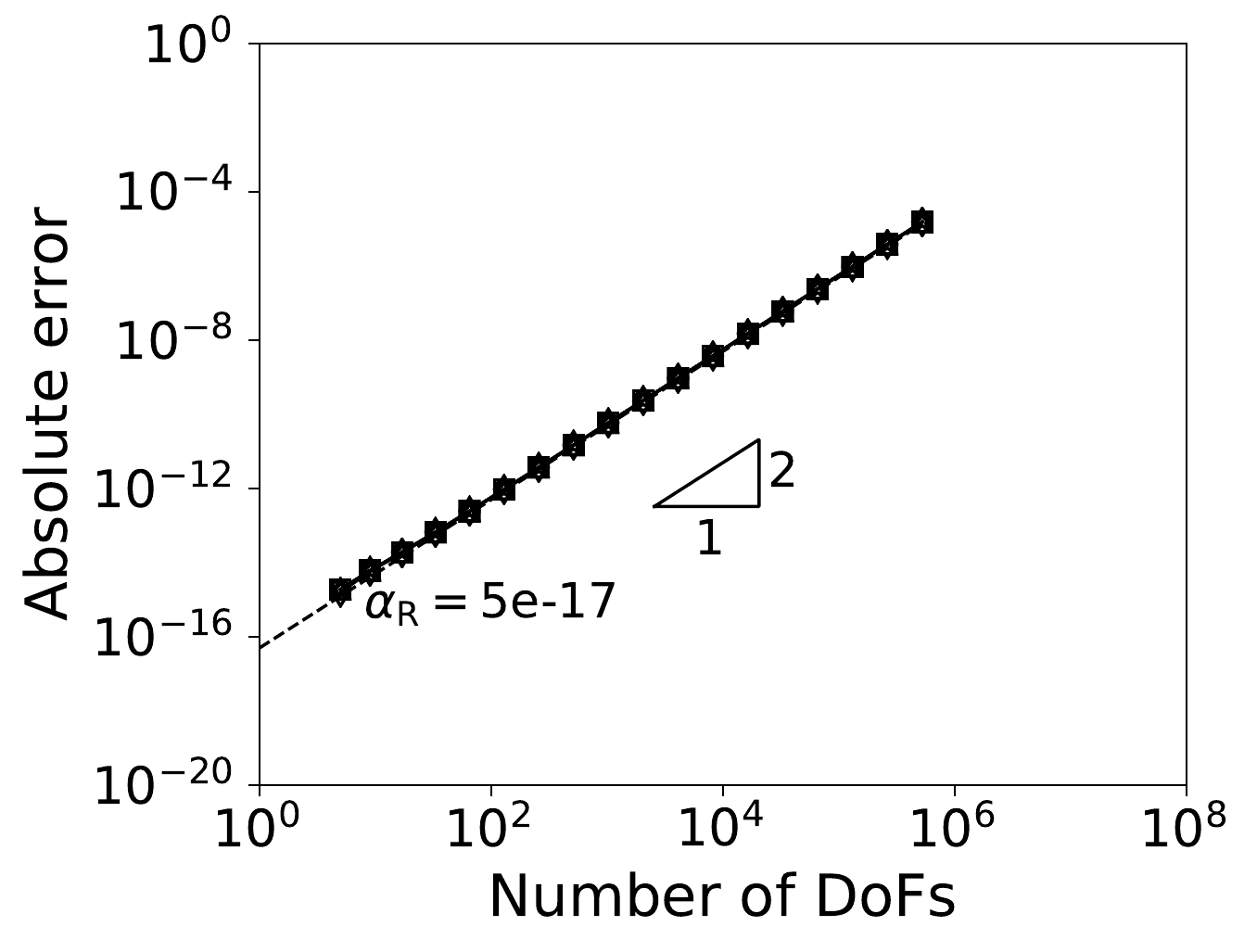}
        \caption{First derivative}
        \label{py_L2_Pois5_SM_scaling_S_grad}
    \end{subfigure}
    \hspace{-0.2cm}
    \begin{subfigure}{5.5cm}
        \includegraphics[width=1.0\linewidth]{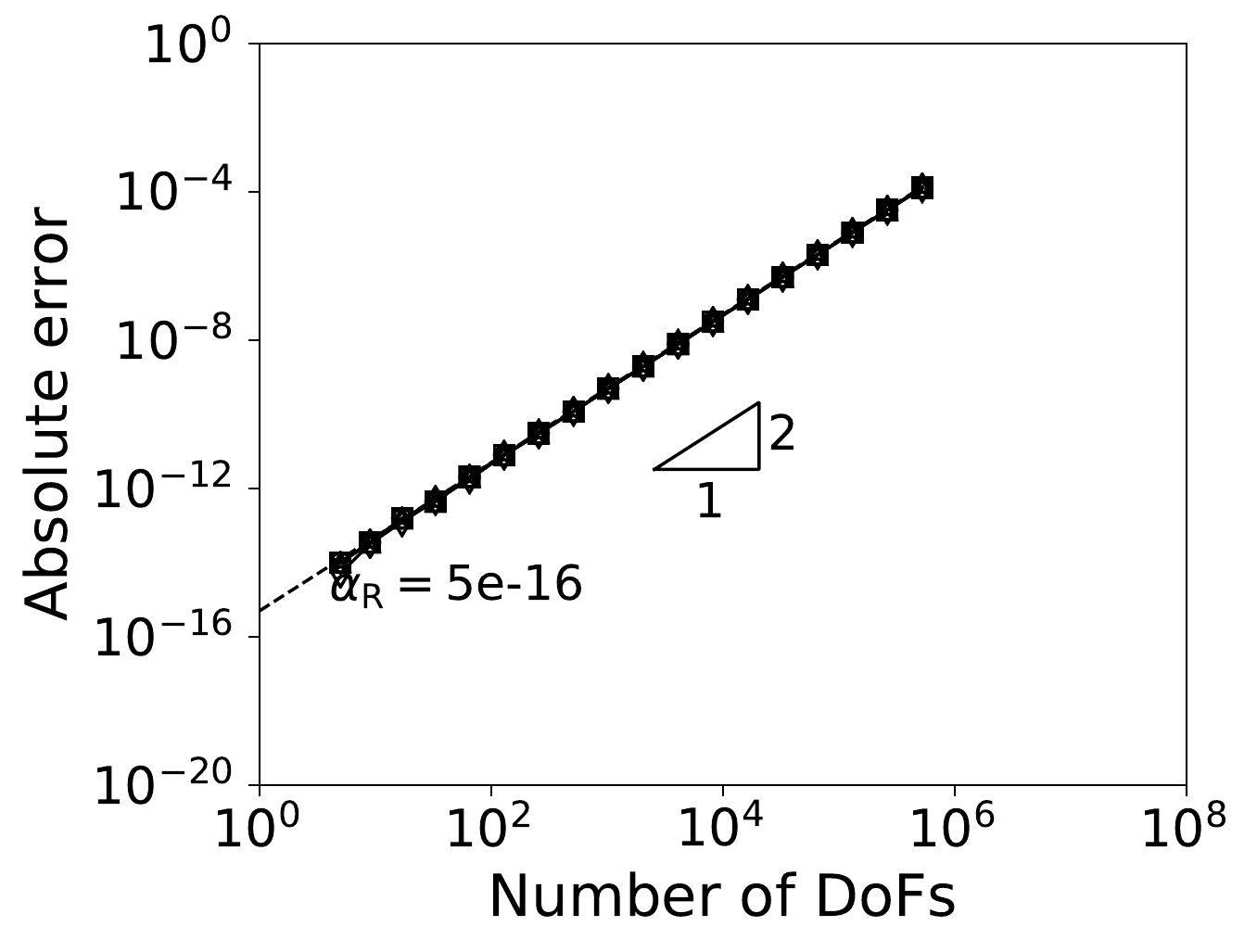}
        \caption{Second derivative}
        \label{py_L2_Pois5_SM_scaling_S_2ndd}
    \end{subfigure}
\caption{Absolute errors of Case 5 in Table \ref{scaling_cases_Poisson} using scheme $S$.}
\label{py_L2_Pois5_SM_scaling_S}
\end{figure}

\paragraph{Case 5}
For Case 5, using the standard FEM without scaling the right-hand side and scheme S, the absolute errors are shown in Figs. \ref{py_L2_Pois5_SM_scaling_no}--\ref{py_L2_Pois5_SM_scaling_S}.

\newpage

\subsubsection{The mixed FEM}

\begin{figure}[!ht]
    \begin{subfigure}{5.5cm}
        \includegraphics[width=1.0\linewidth]{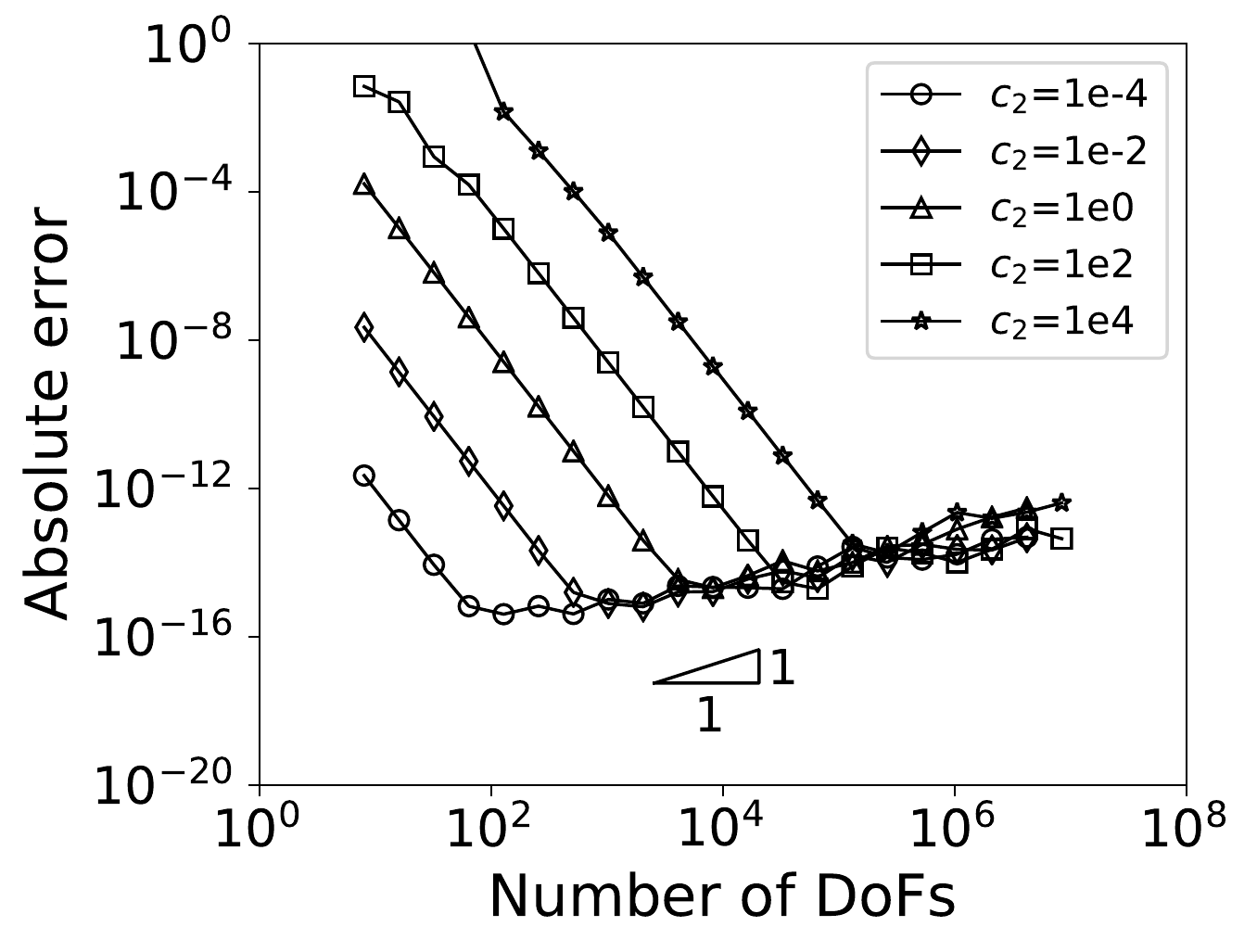}
        \caption{Solution}
        \label{py_L2_Pois2_MM_scaling_no_solu}
    \end{subfigure}
    \hspace{-0.2cm}
    \begin{subfigure}{5.5cm}
        \includegraphics[width=1.0\linewidth]{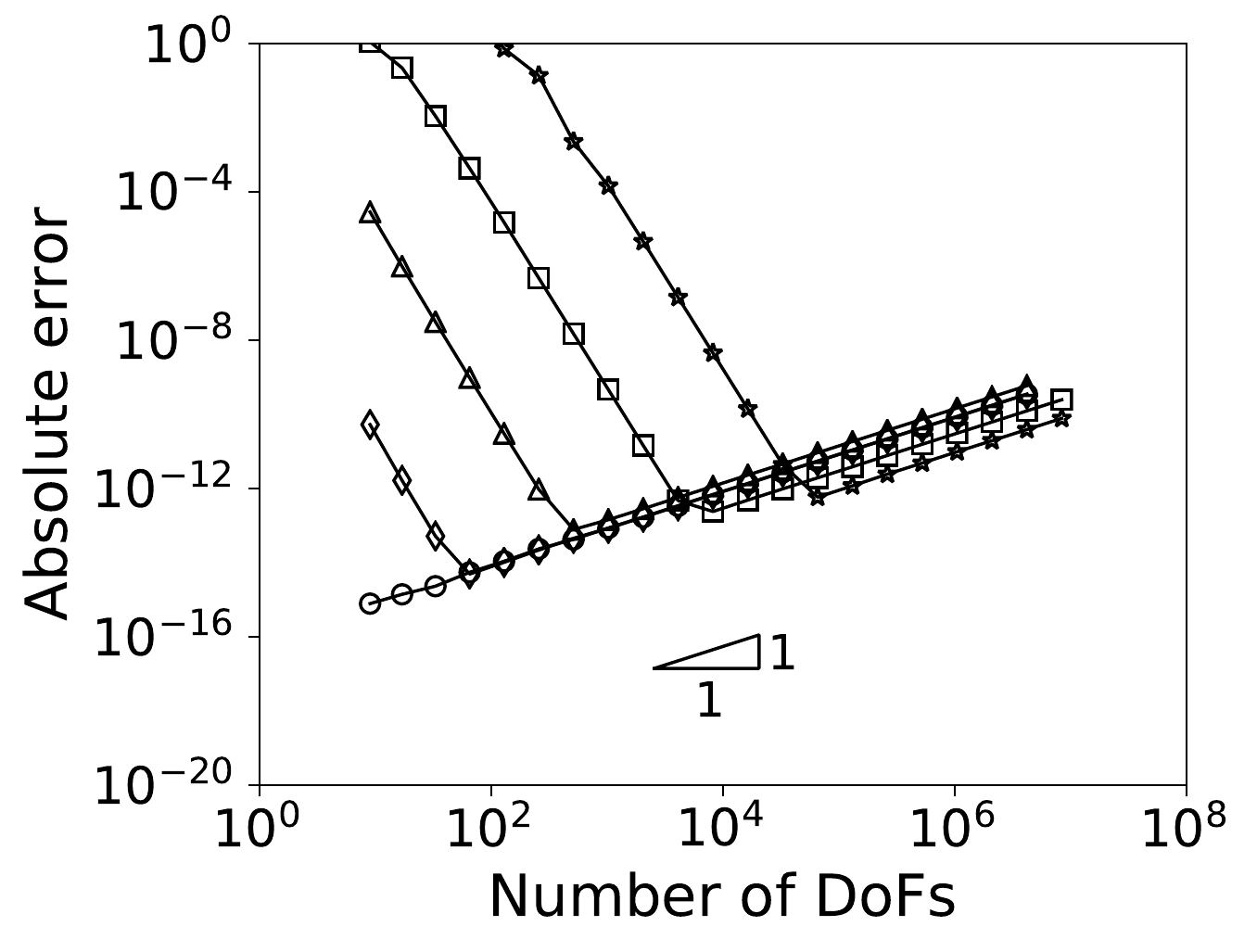}
        \caption{First derivative}
        \label{py_L2_Pois2_MM_scaling_no_grad}
    \end{subfigure}
    \hspace{-0.2cm}
    \begin{subfigure}{5.5cm}
        \includegraphics[width=1.0\linewidth]{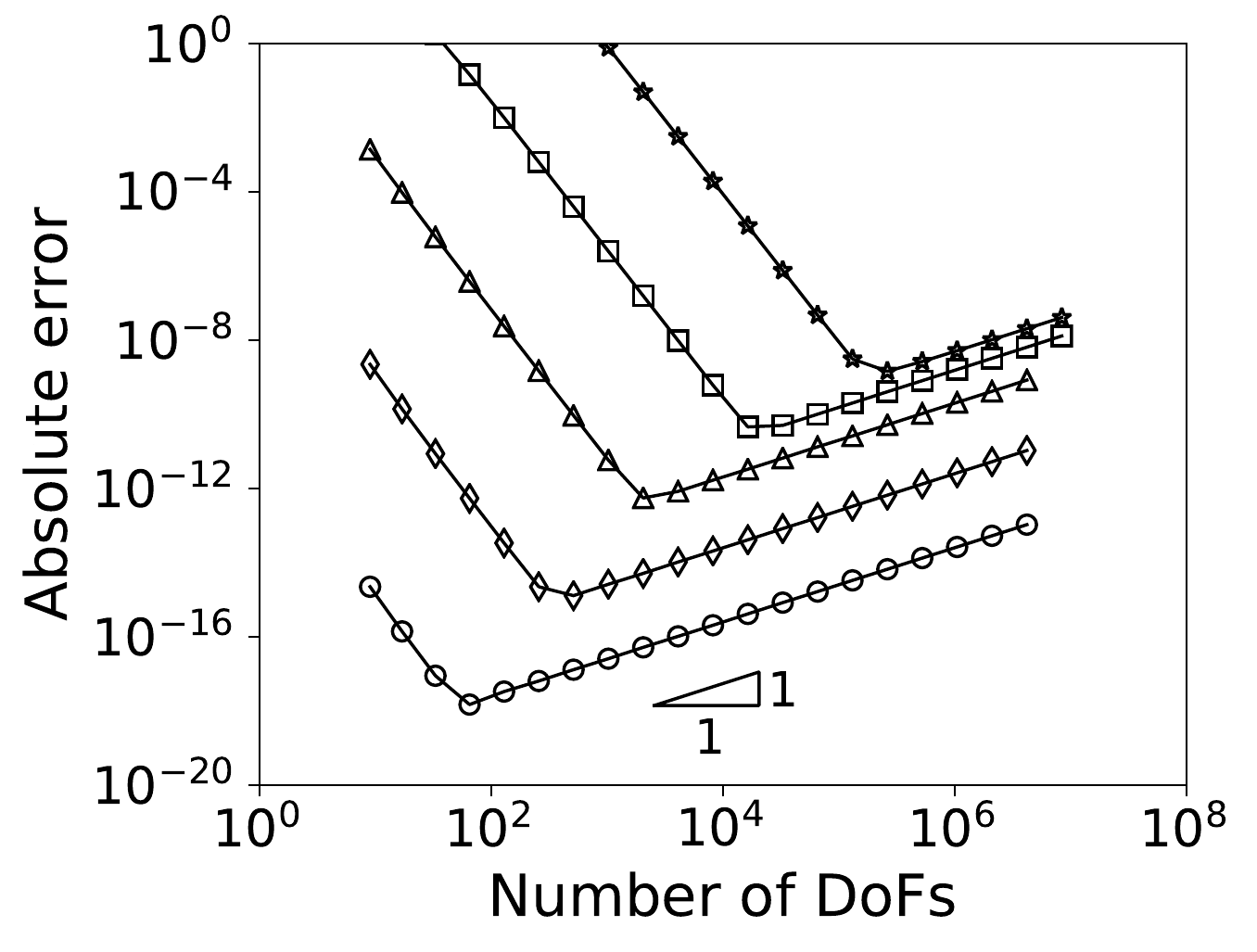}
        \caption{Second derivative}
        \label{py_L2_Pois2_MM_scaling_no_2ndd}
    \end{subfigure}
\caption{Absolute errors for Case 2 in Table \ref{scaling_cases_Poisson} using the mixed FEM without scaling the right-hand side.}
\label{py_L2_Pois2_MM_scaling_no}
\end{figure}

\begin{figure}[!ht]
    \begin{subfigure}{5.5cm}
        \includegraphics[width=1.0\linewidth]{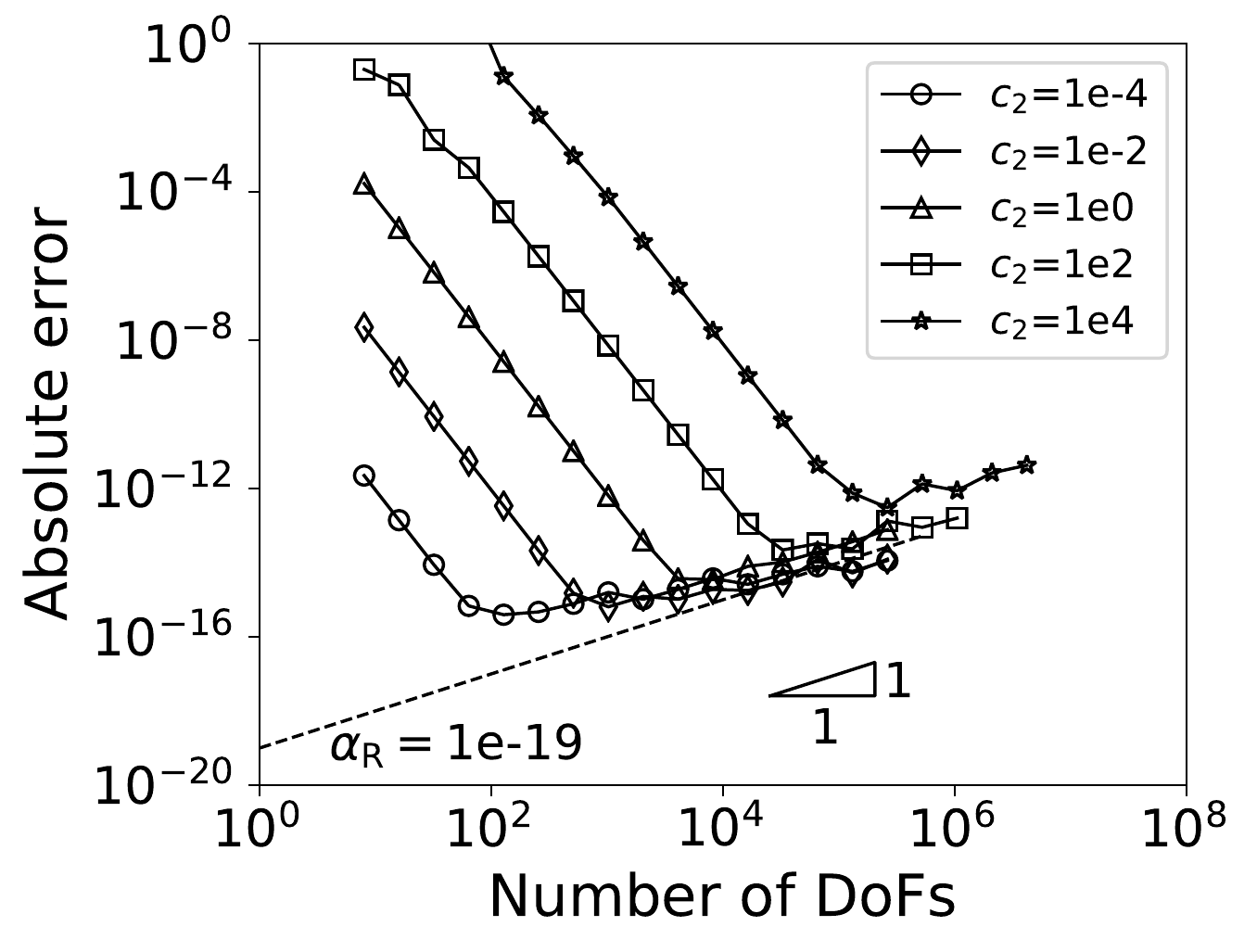}
        \caption{Solution}
        \label{py_L2_Pois2_MM_scaling_M1_solu}
    \end{subfigure}
    \hspace{-0.2cm}
    \begin{subfigure}{5.5cm}
        \includegraphics[width=1.0\linewidth]{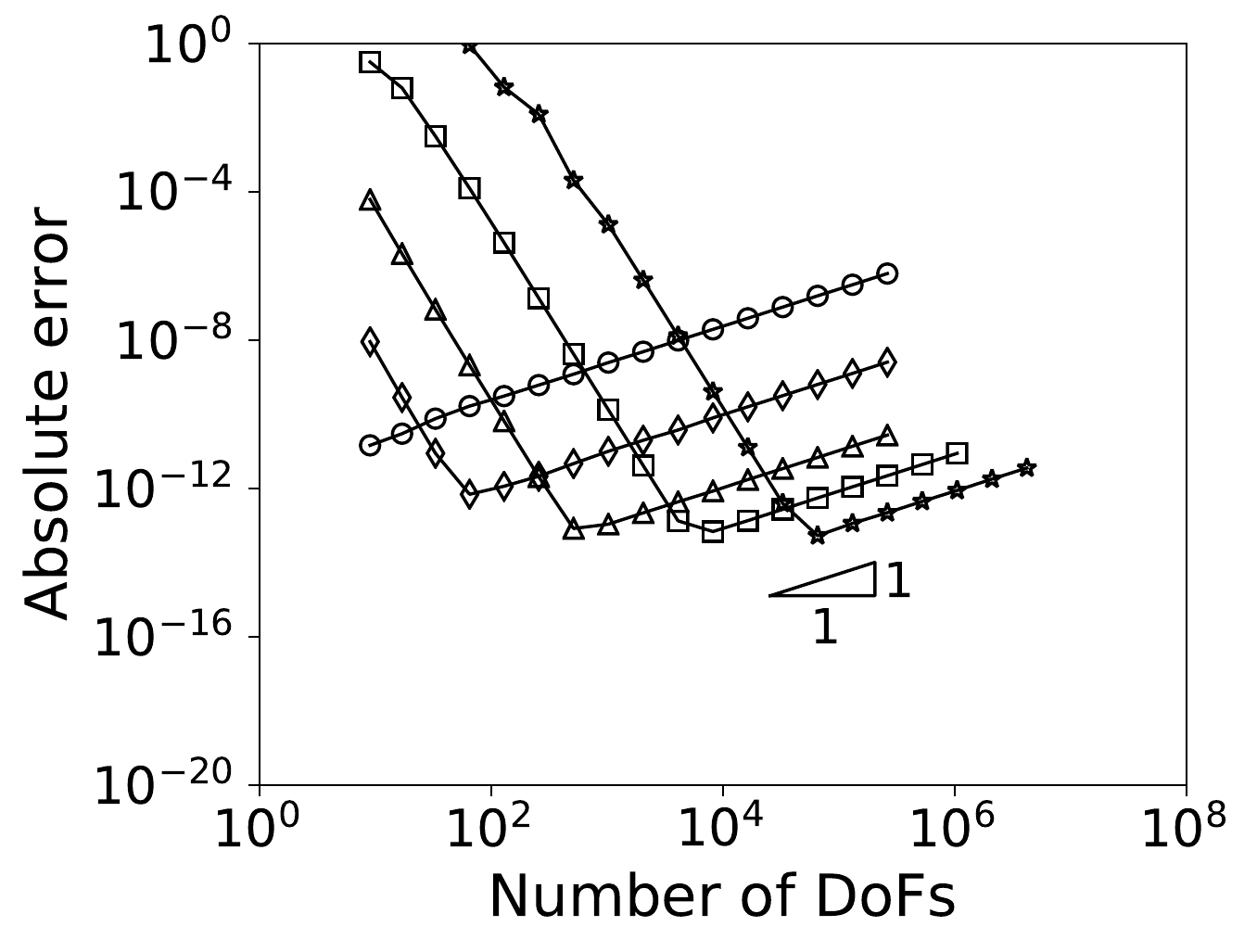}
        \caption{First derivative}
        \label{py_L2_Pois2_MM_scaling_M1_grad}
    \end{subfigure}
    \hspace{-0.2cm}
    \begin{subfigure}{5.5cm}
        \includegraphics[width=1.0\linewidth]{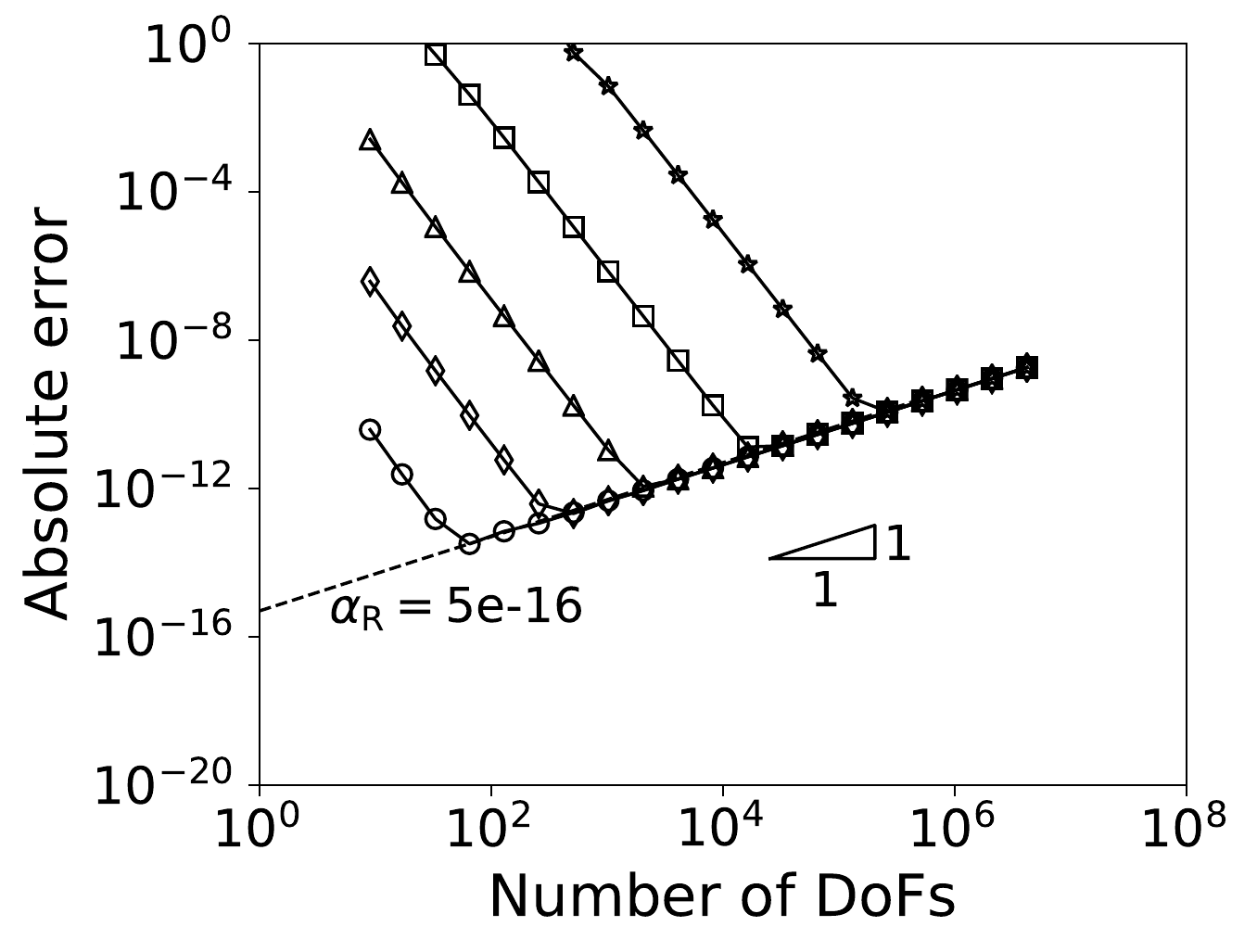}
        \caption{Second derivative}
        \label{py_L2_Pois2_MM_scaling_M1_2ndd}
    \end{subfigure}
\caption{Absolute errors for Case 2 in Table \ref{scaling_cases_Poisson} using scheme $M_1$.}
\label{py_L2_Pois2_MM_scaling_M1}
\end{figure}

\begin{figure}[!ht]
    \begin{subfigure}{5.5cm}
        \includegraphics[width=1.0\linewidth]{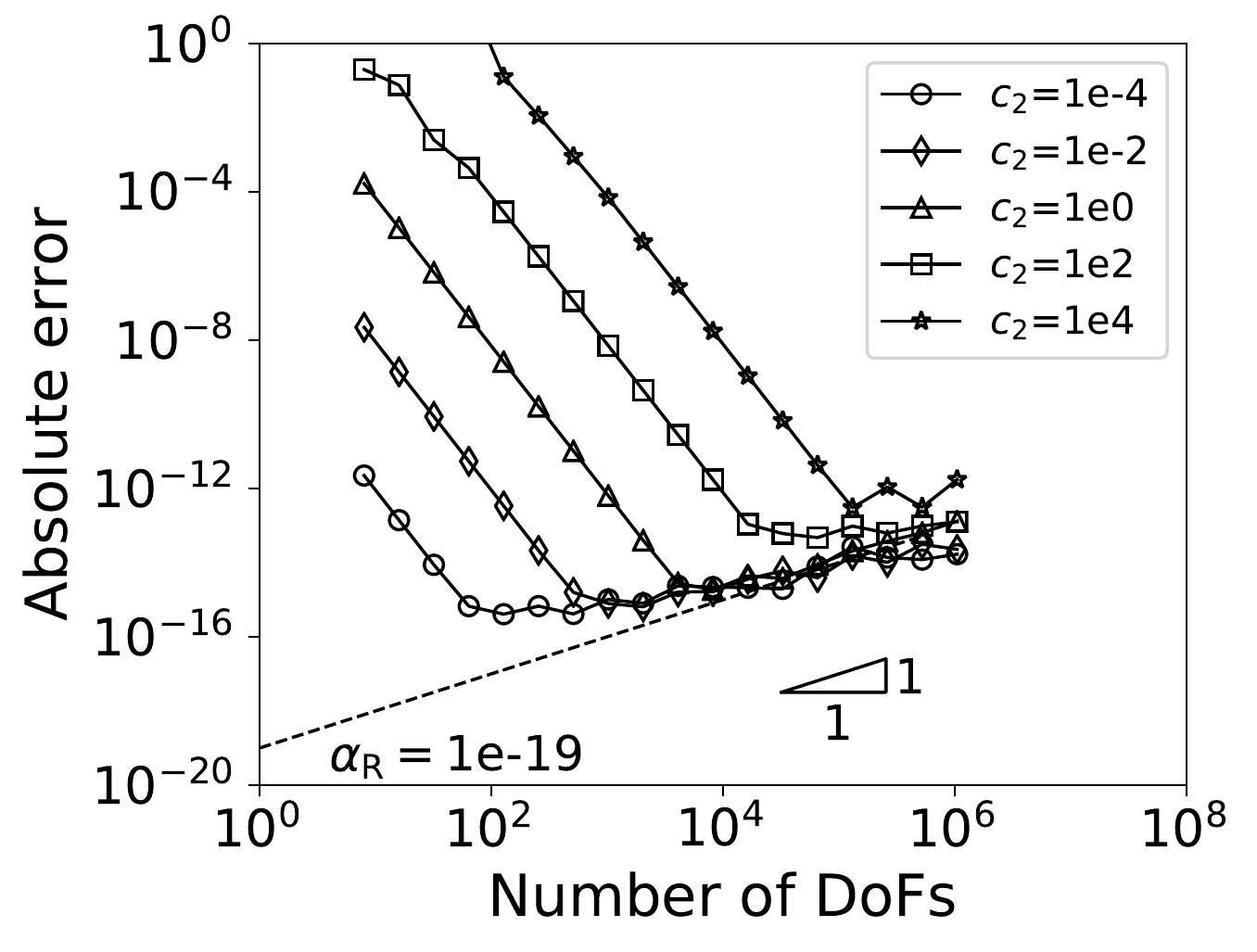}
        \caption{Solution}
        \label{py_L2_Pois2_MM_scaling_M2_solu}
    \end{subfigure}
    \hspace{-0.2cm}
    \begin{subfigure}{5.5cm}
        \includegraphics[width=1.0\linewidth]{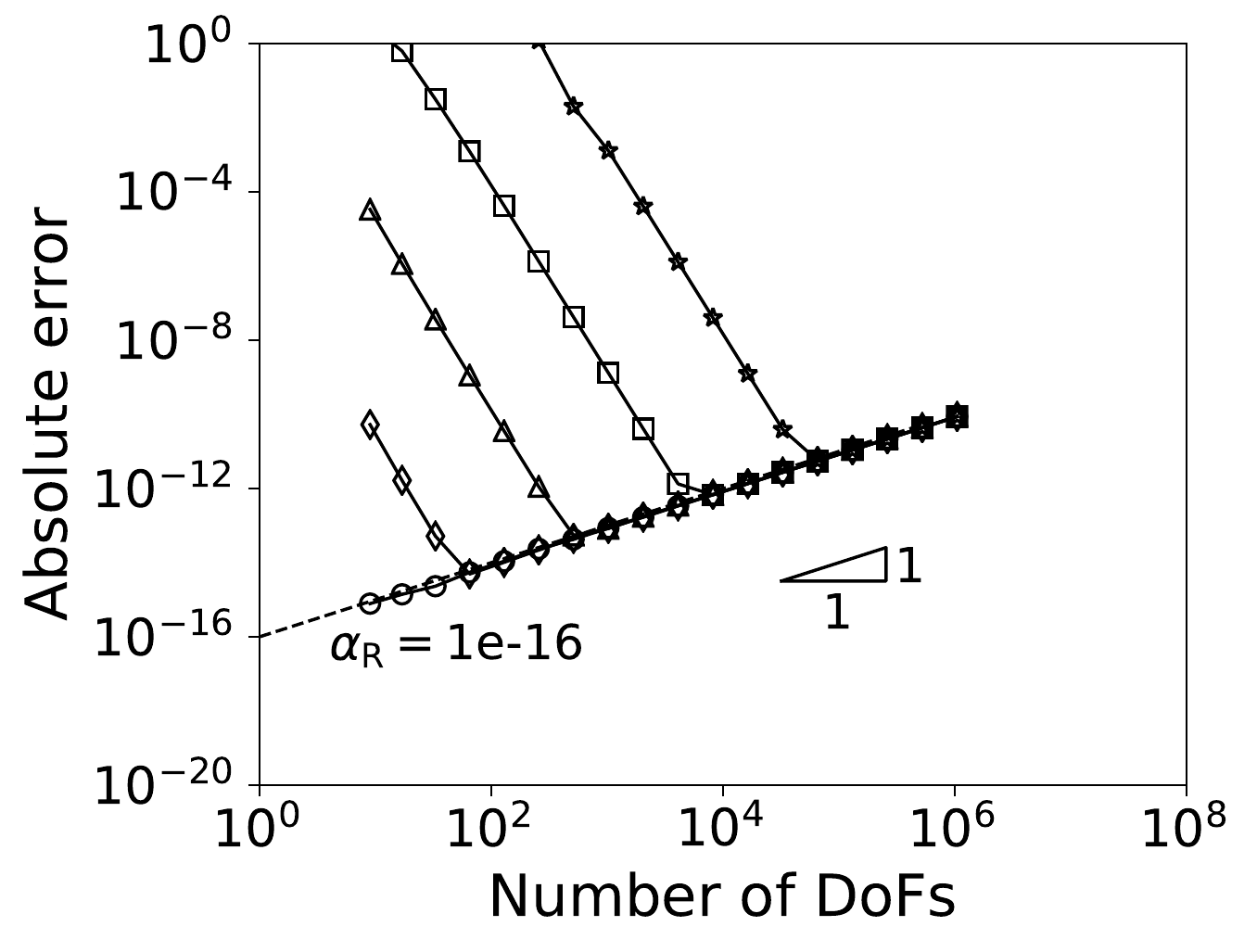}
        \caption{First derivative}
        \label{py_L2_Pois2_MM_scaling_M2_grad}
    \end{subfigure}
    \hspace{-0.2cm}
    \begin{subfigure}{5.5cm}
        \includegraphics[width=1.0\linewidth]{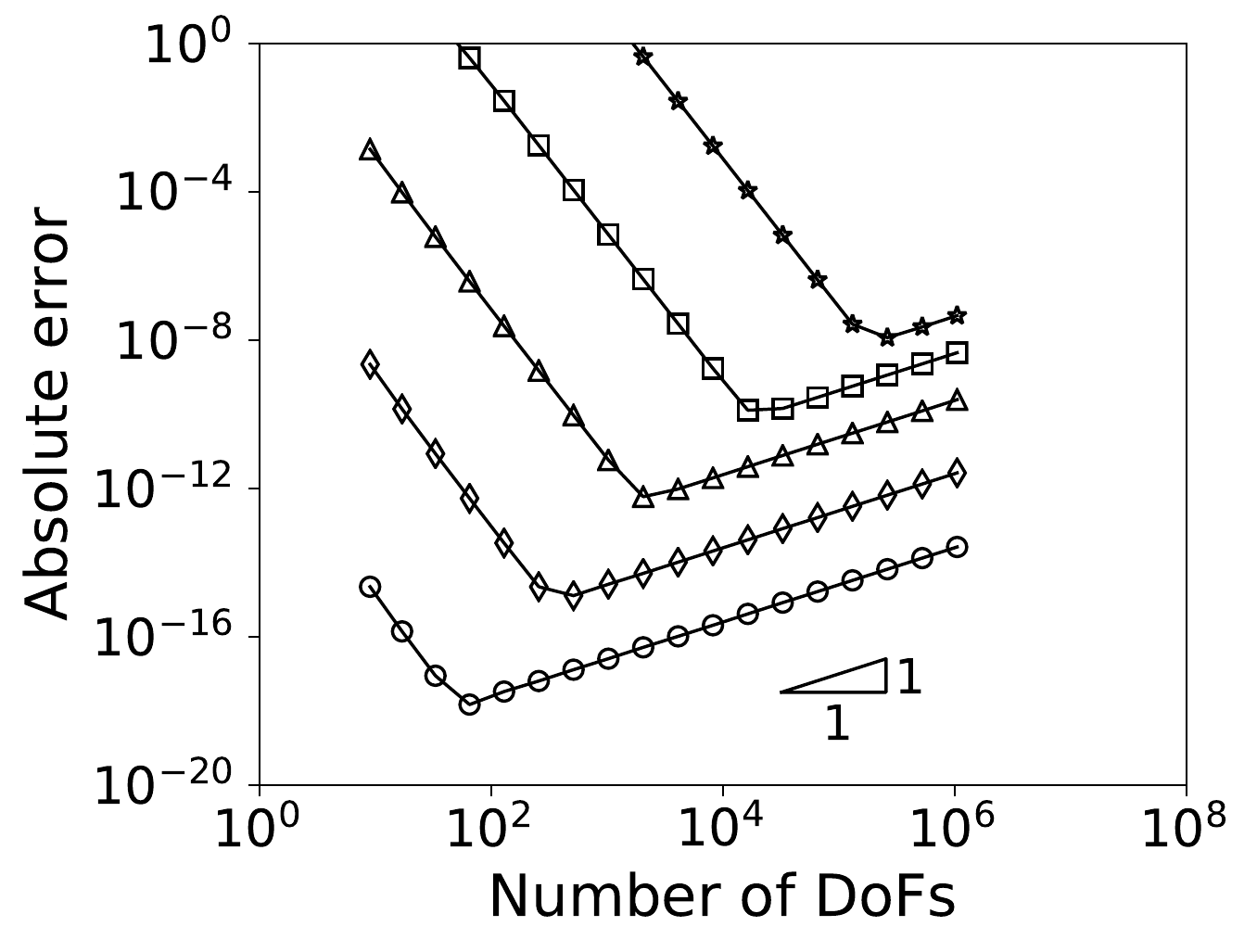}
        \caption{Second derivative}
        \label{py_L2_Pois2_MM_scaling_M2_2ndd}
    \end{subfigure}
\caption{Absolute errors for Case 2 in Table \ref{scaling_cases_Poisson} using scheme $M_2$.}
\label{py_L2_Pois2_MM_scaling_M2}
\end{figure}

\paragraph{Case 2}
For Case 2, using the mixed FEM without scaling the right-hand side and schemes ${\rm M}_1$ and ${\rm M}_2$, the absolute errors are shown in Figs. \ref{py_L2_Pois2_MM_scaling_no}--\ref{py_L2_Pois2_MM_scaling_M2}.

\begin{figure}[!ht]
    \begin{subfigure}{5.5cm}
        \includegraphics[width=1.0\linewidth]{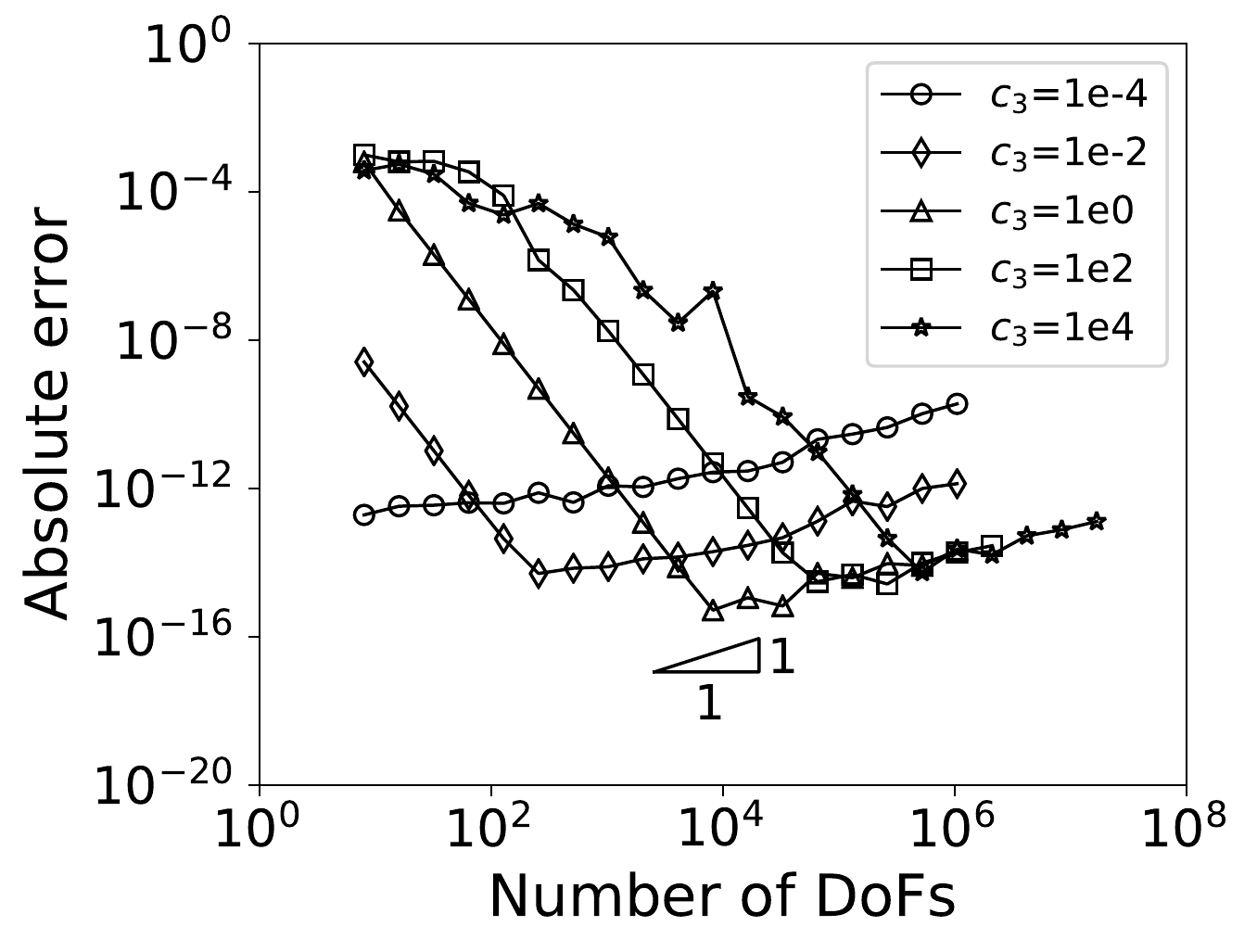}
        \caption{Solution}
        \label{py_L2_Pois3_MM_scaling_no_solu}
    \end{subfigure}
    \hspace{-0.2cm}
    \begin{subfigure}{5.5cm}
        \includegraphics[width=1.0\linewidth]{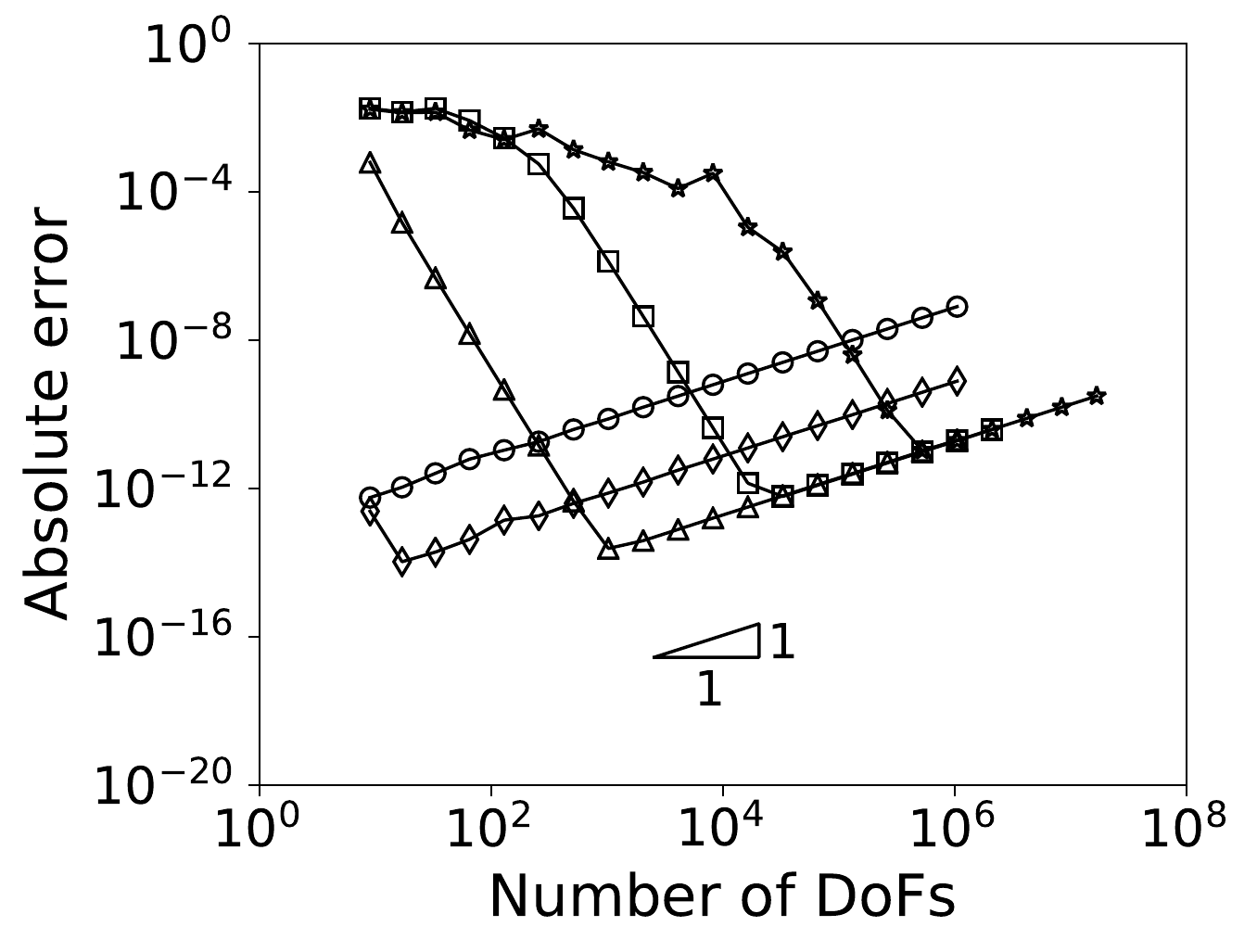}
        \caption{First derivative}
        \label{py_L2_Pois3_MM_scaling_no_grad}
    \end{subfigure}
    \hspace{-0.2cm}
    \begin{subfigure}{5.5cm}
        \includegraphics[width=1.0\linewidth]{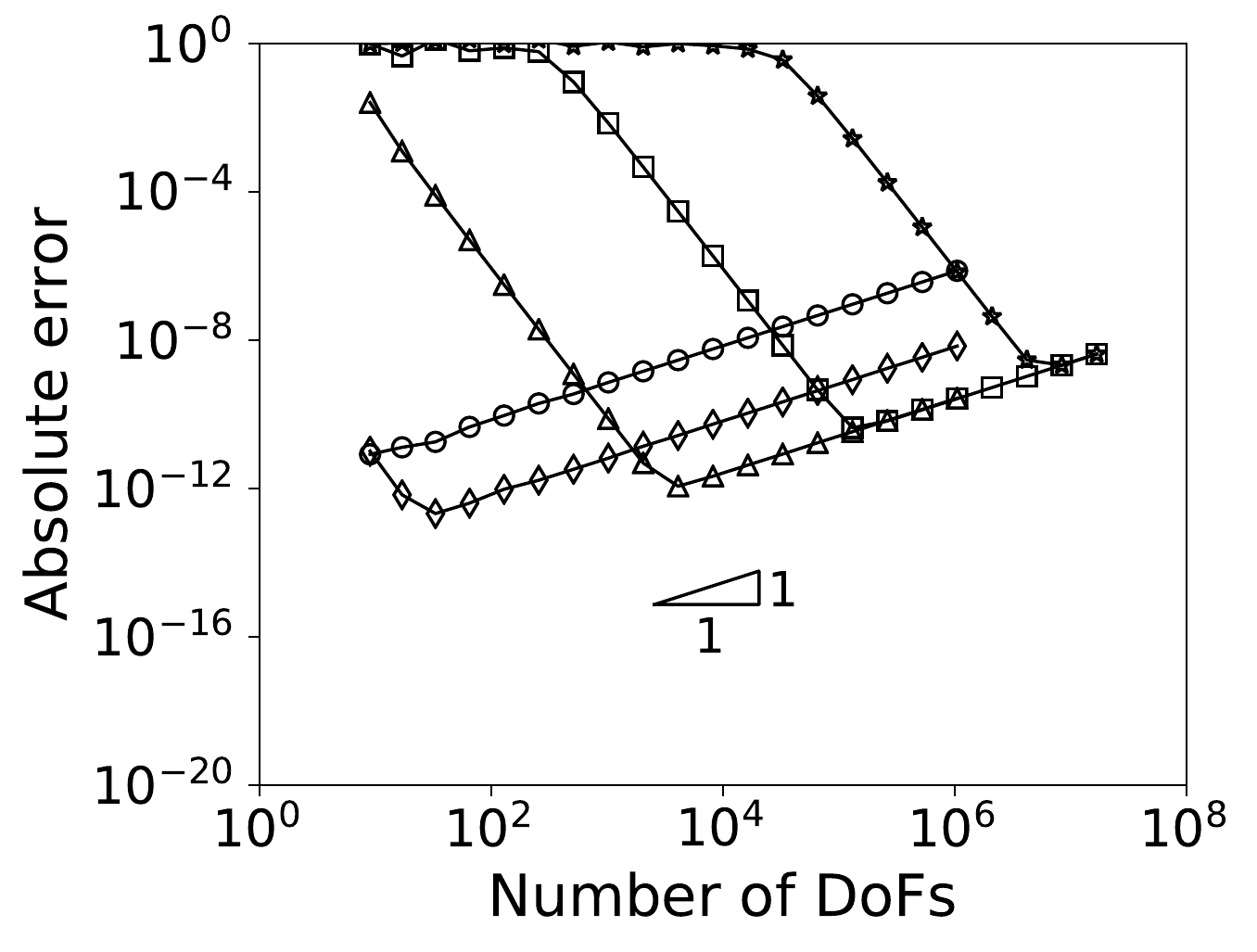}
        \caption{Second derivative}
        \label{py_L2_Pois3_MM_scaling_no_2ndd}
    \end{subfigure}
\caption{Absolute errors for Case 3 in Table \ref{scaling_cases_Poisson} using the mixed FEM without scaling the right-hand side.}
\label{py_L2_Pois3_MM_scaling_no}
\end{figure}

\begin{figure}[!ht]
    \begin{subfigure}{5.5cm}
        \includegraphics[width=1.0\linewidth]{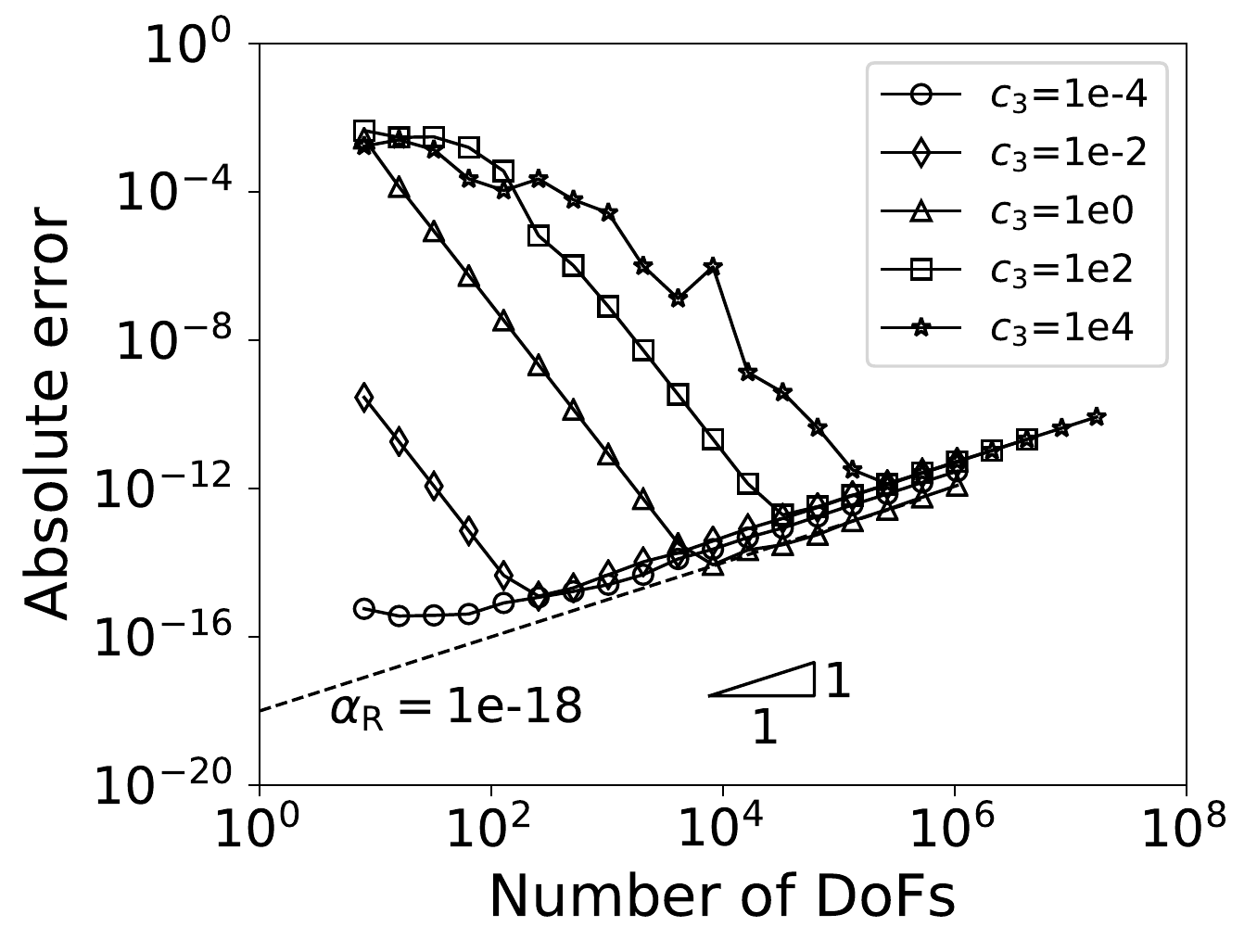}
        \caption{Solution}
        \label{py_L2_Pois3_MM_scaling_M1_solu}
    \end{subfigure}
    \hspace{-0.2cm}
    \begin{subfigure}{5.5cm}
        \includegraphics[width=1.0\linewidth]{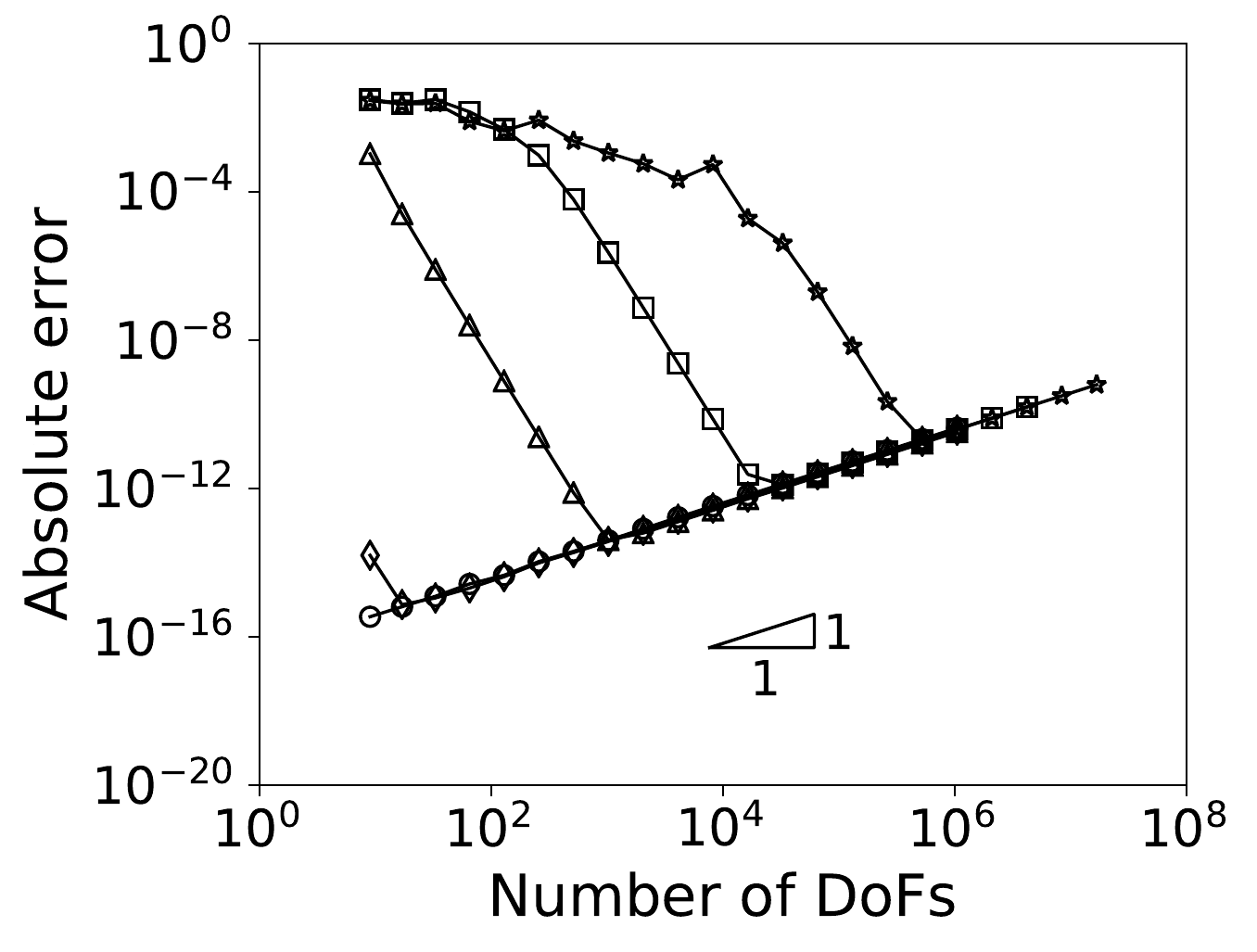}
        \caption{First derivative}
        \label{py_L2_Pois3_MM_scaling_M1_grad}
    \end{subfigure}
    \hspace{-0.2cm}
    \begin{subfigure}{5.5cm}
        \includegraphics[width=1.0\linewidth]{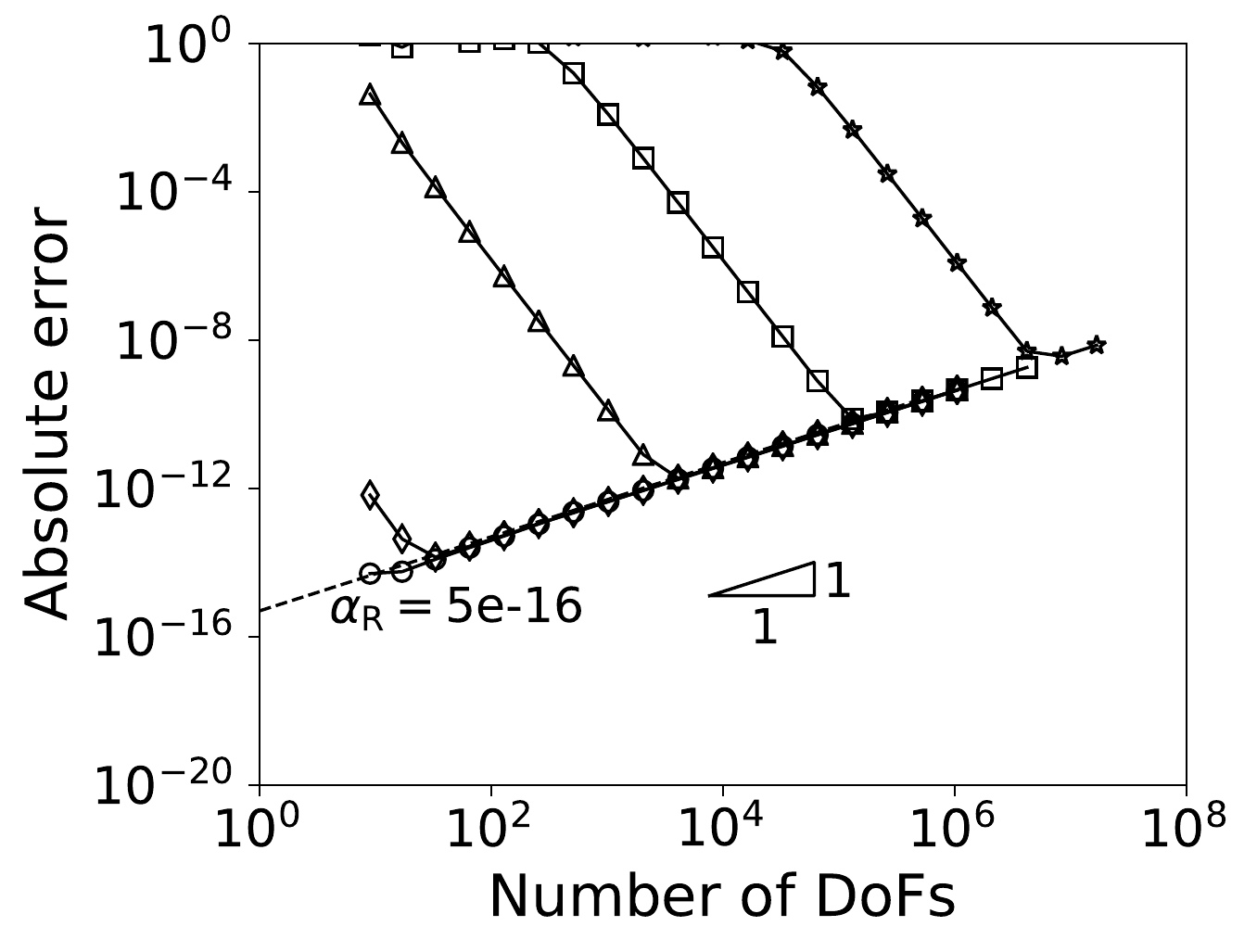}
        \caption{Second derivative}
        \label{py_L2_Pois3_MM_scaling_M1_2ndd}
    \end{subfigure}
\caption{Absolute errors for Case 3 in Table \ref{scaling_cases_Poisson} using scheme $M_1$.}
\label{py_L2_Pois3_MM_scaling_M1}
\end{figure}

\begin{figure}[!ht]
    \begin{subfigure}{5.5cm}
        \includegraphics[width=1.0\linewidth]{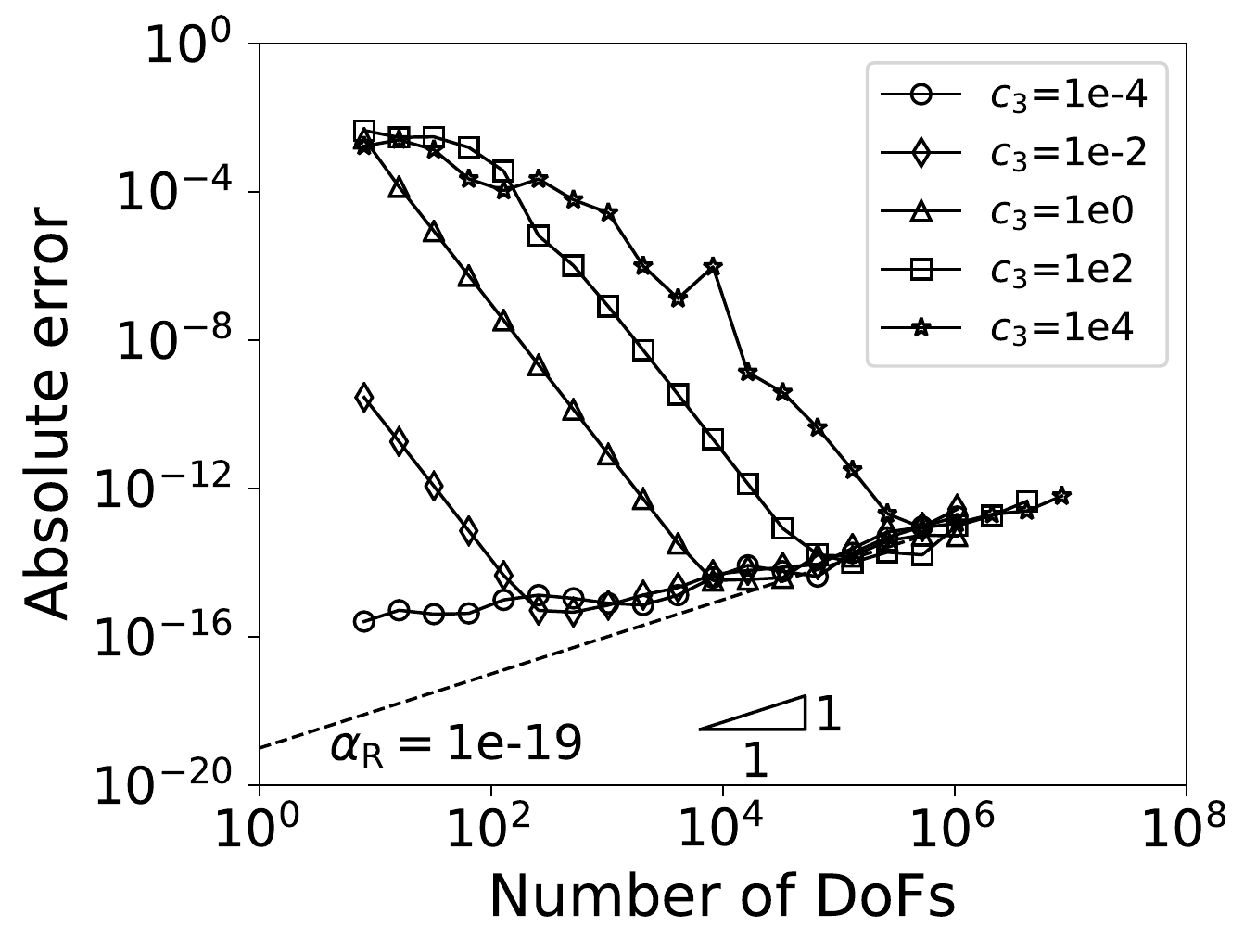}
        \caption{Solution}
        \label{py_L2_Pois3_MM_scaling_M2_solu}
    \end{subfigure}
    \hspace{-0.2cm}
    \begin{subfigure}{5.5cm}
        \includegraphics[width=1.0\linewidth]{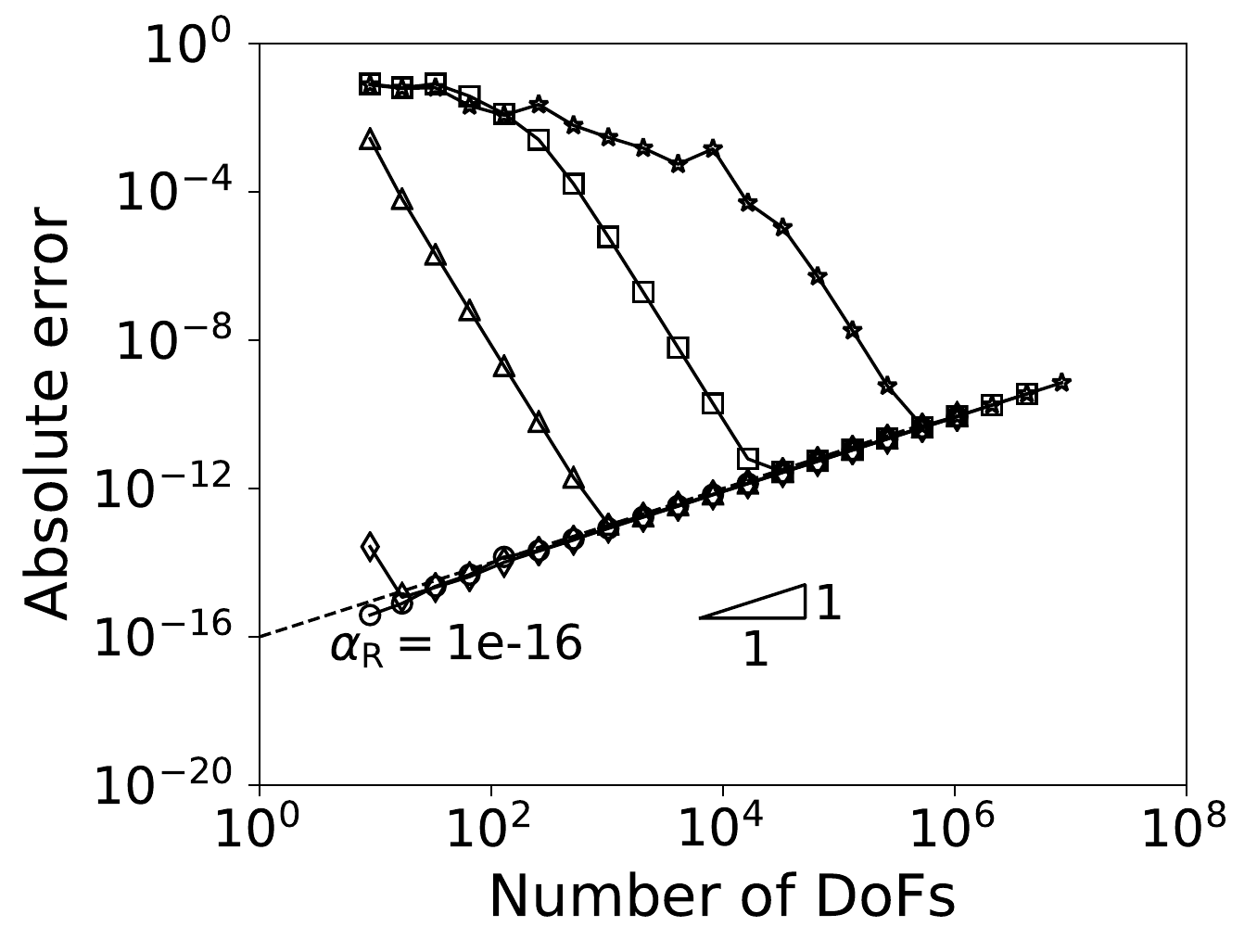}
        \caption{First derivative}
        \label{py_L2_Pois3_MM_scaling_M2_grad}
    \end{subfigure}
    \hspace{-0.2cm}
    \begin{subfigure}{5.5cm}
        \includegraphics[width=1.0\linewidth]{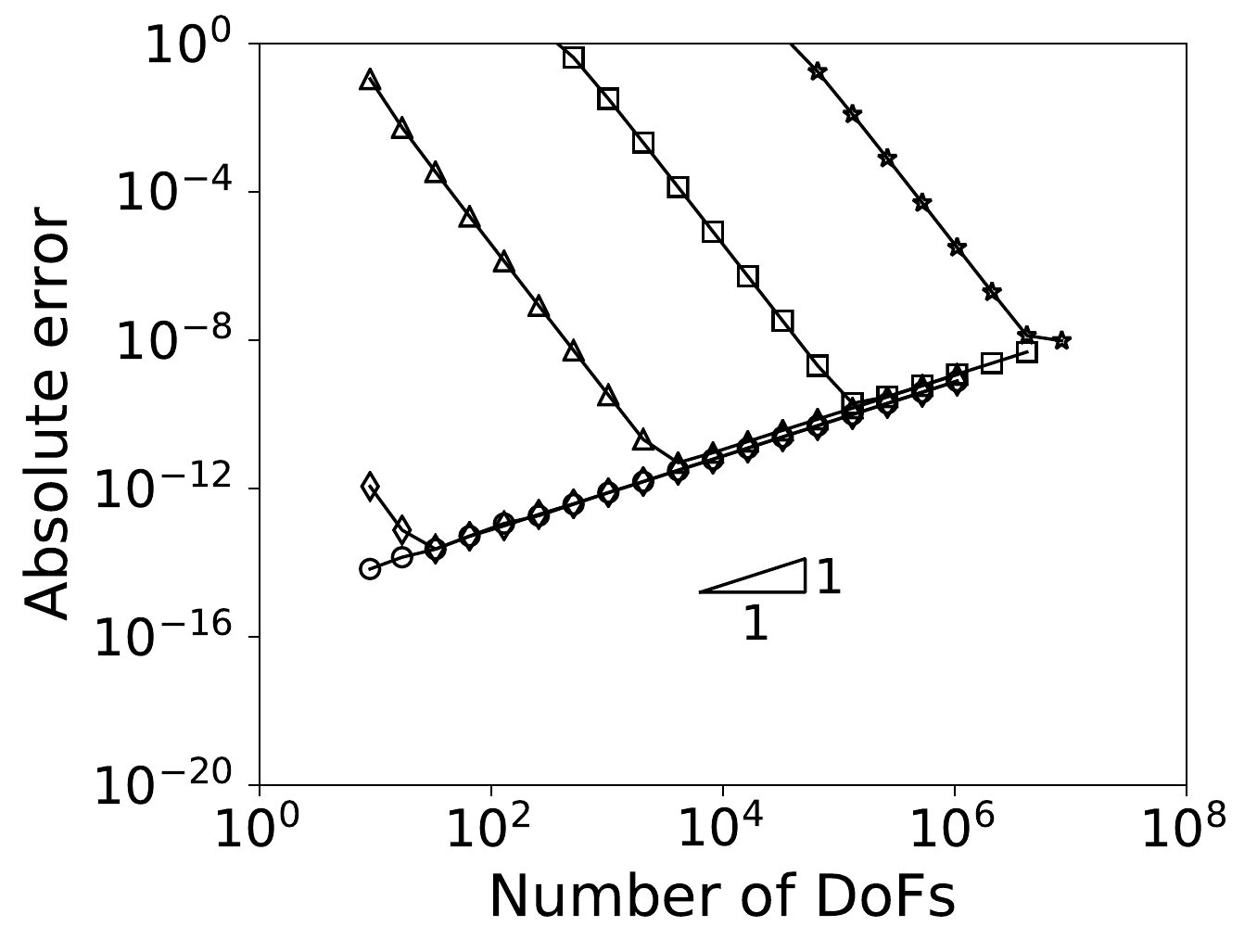}
        \caption{Second derivative}
        \label{py_L2_Pois3_MM_scaling_M2_2ndd}
    \end{subfigure}
\caption{Absolute errors for Case 3 in Table \ref{scaling_cases_Poisson} using scheme $M_2$.}
\label{py_L2_Pois3_MM_scaling_M2}
\end{figure}

\newpage
\paragraph{Case 3}
For Case 3, using the mixed FEM without scaling the right-hand side and schemes ${\rm M}_1$ and ${\rm M}_2$, the absolute errors are shown in Figs. \ref{py_L2_Pois3_MM_scaling_no}--\ref{py_L2_Pois3_MM_scaling_M2}.

\begin{figure}[!ht]
    \begin{subfigure}{5.5cm}
        \includegraphics[width=1.0\linewidth]{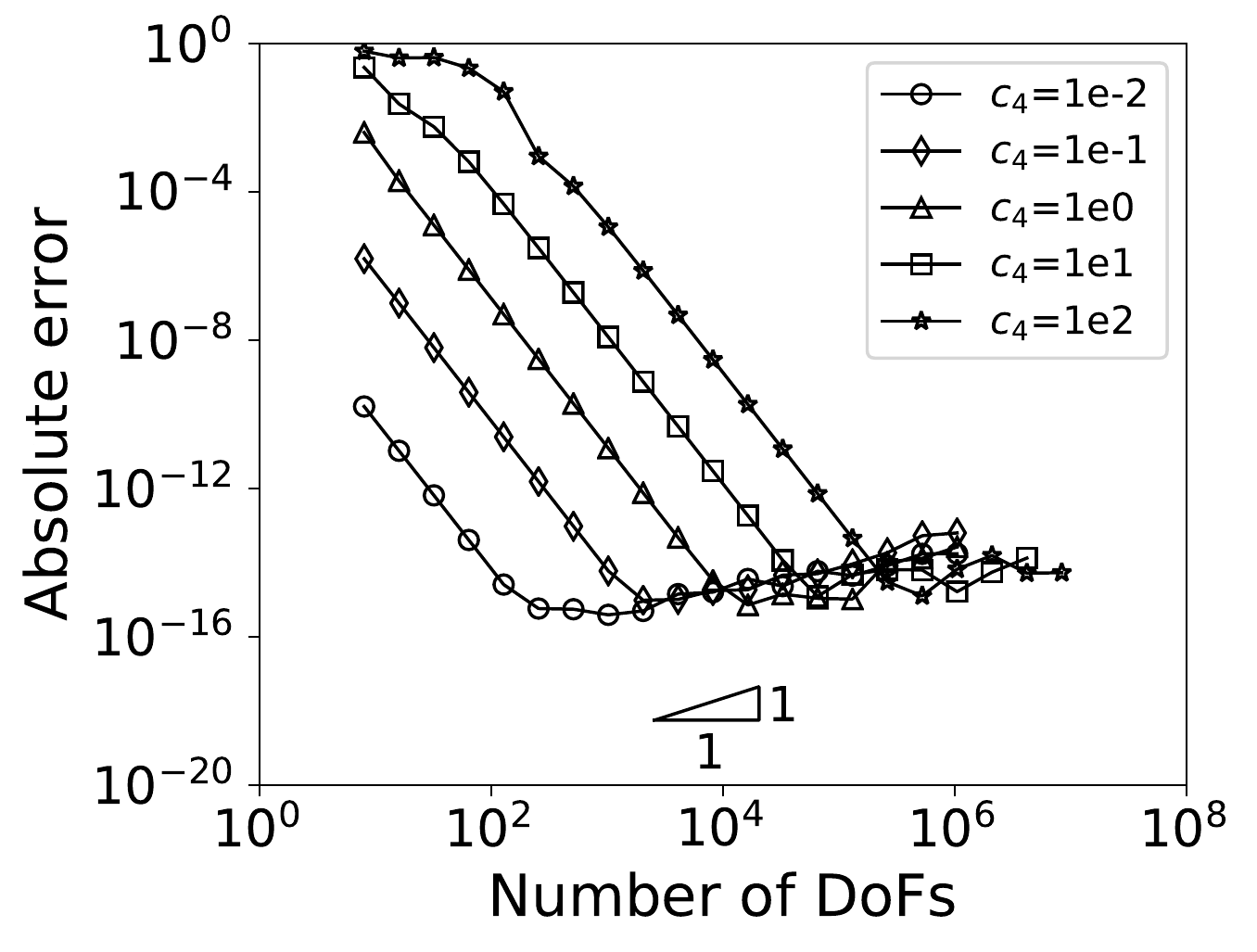}
        \caption{Solution}
        \label{py_L2_Pois4_MM_scaling_no_solu}
    \end{subfigure}
    \hspace{-0.2cm}
    \begin{subfigure}{5.5cm}
        \includegraphics[width=1.0\linewidth]{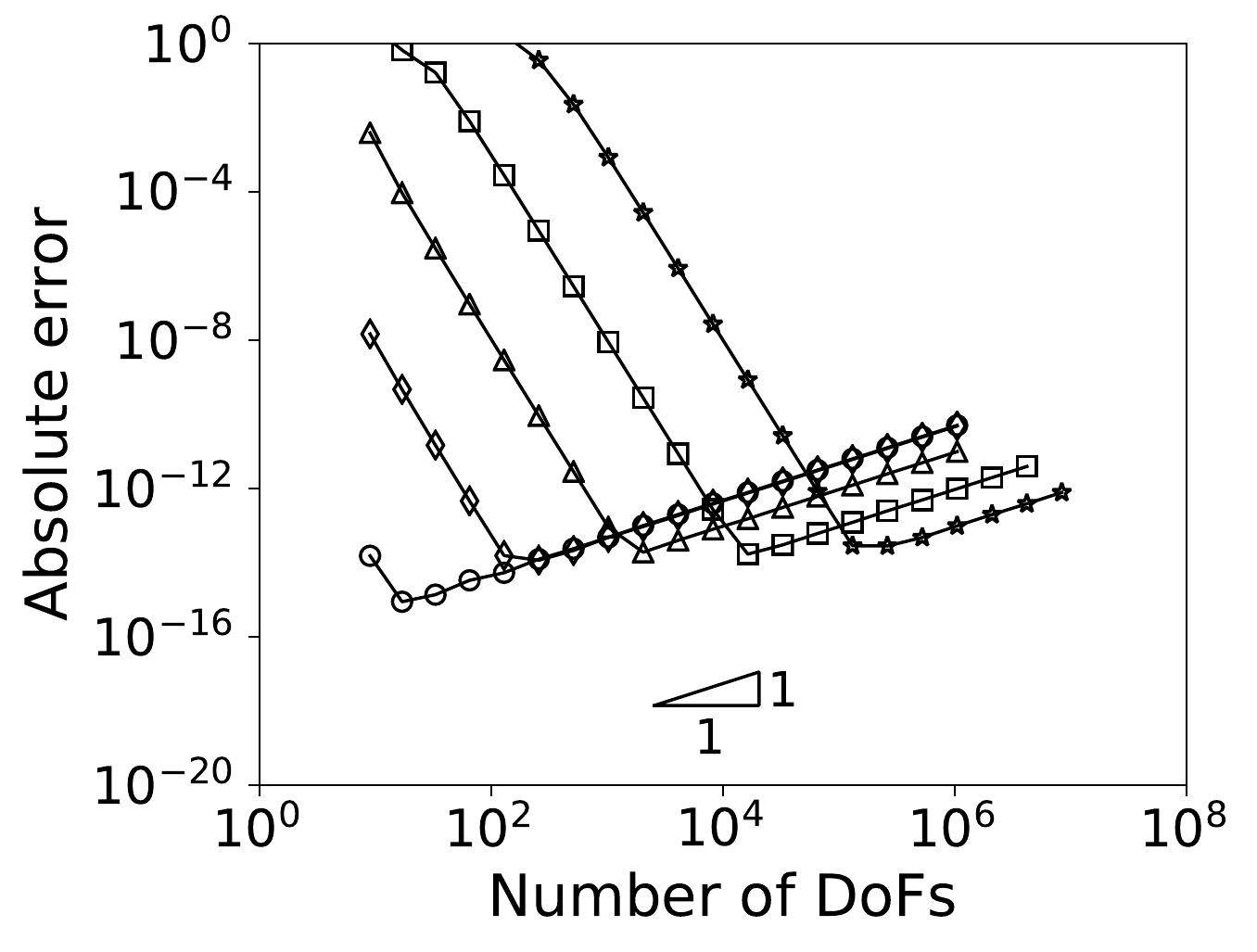}
        \caption{First derivative}
        \label{py_L2_Pois4_MM_scaling_no_grad}
    \end{subfigure}
    \hspace{-0.2cm}
    \begin{subfigure}{5.5cm}
        \includegraphics[width=1.0\linewidth]{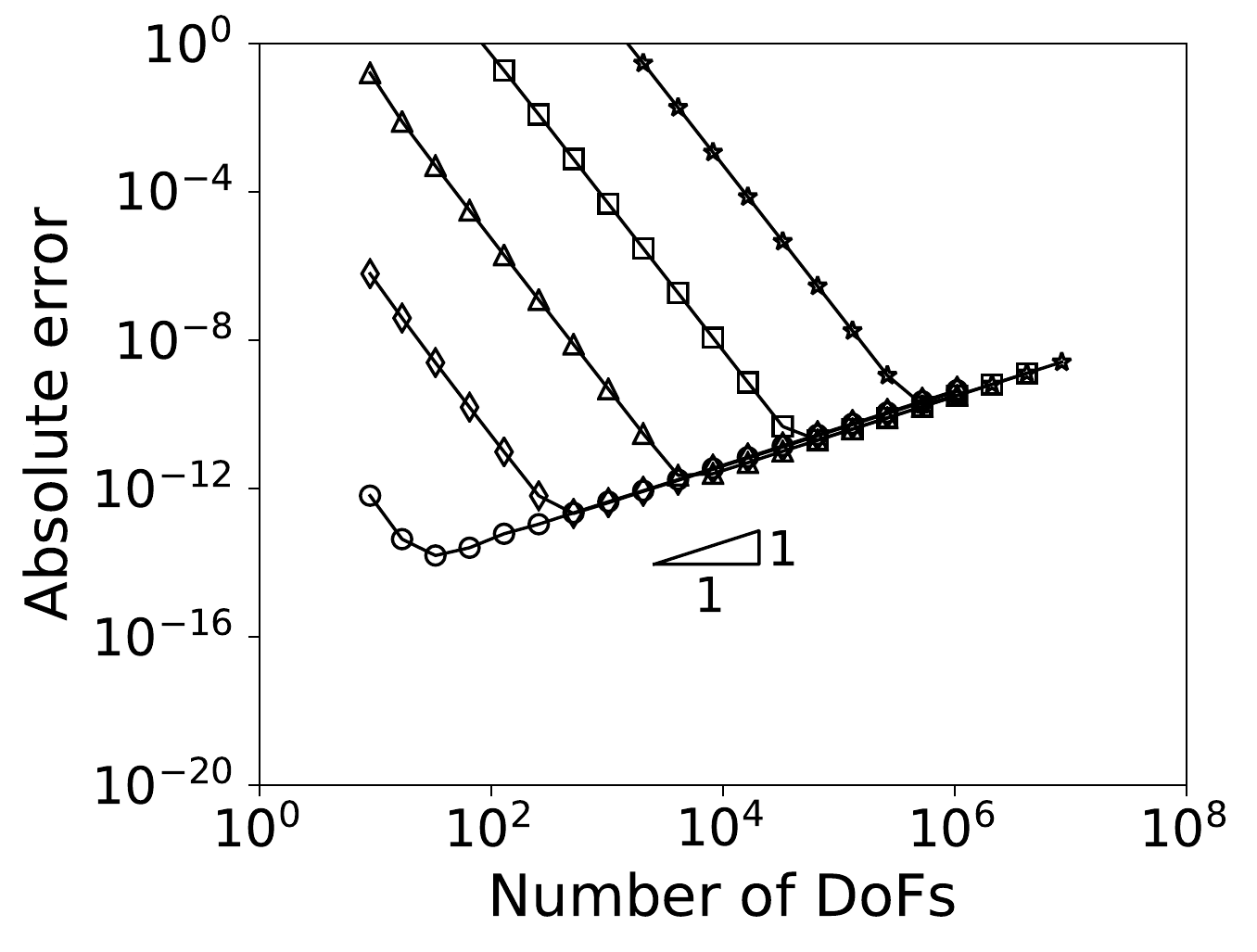}
        \caption{Second derivative}
        \label{py_L2_Pois4_MM_scaling_no_2ndd}
    \end{subfigure}
\caption{Absolute errors for Case 4 in Table \ref{scaling_cases_Poisson} using the mixed FEM without scaling the right-hand side.}
\label{py_L2_Pois4_MM_scaling_no}
\end{figure}

\begin{figure}[!ht]
    \begin{subfigure}{5.5cm}
        \includegraphics[width=1.0\linewidth]{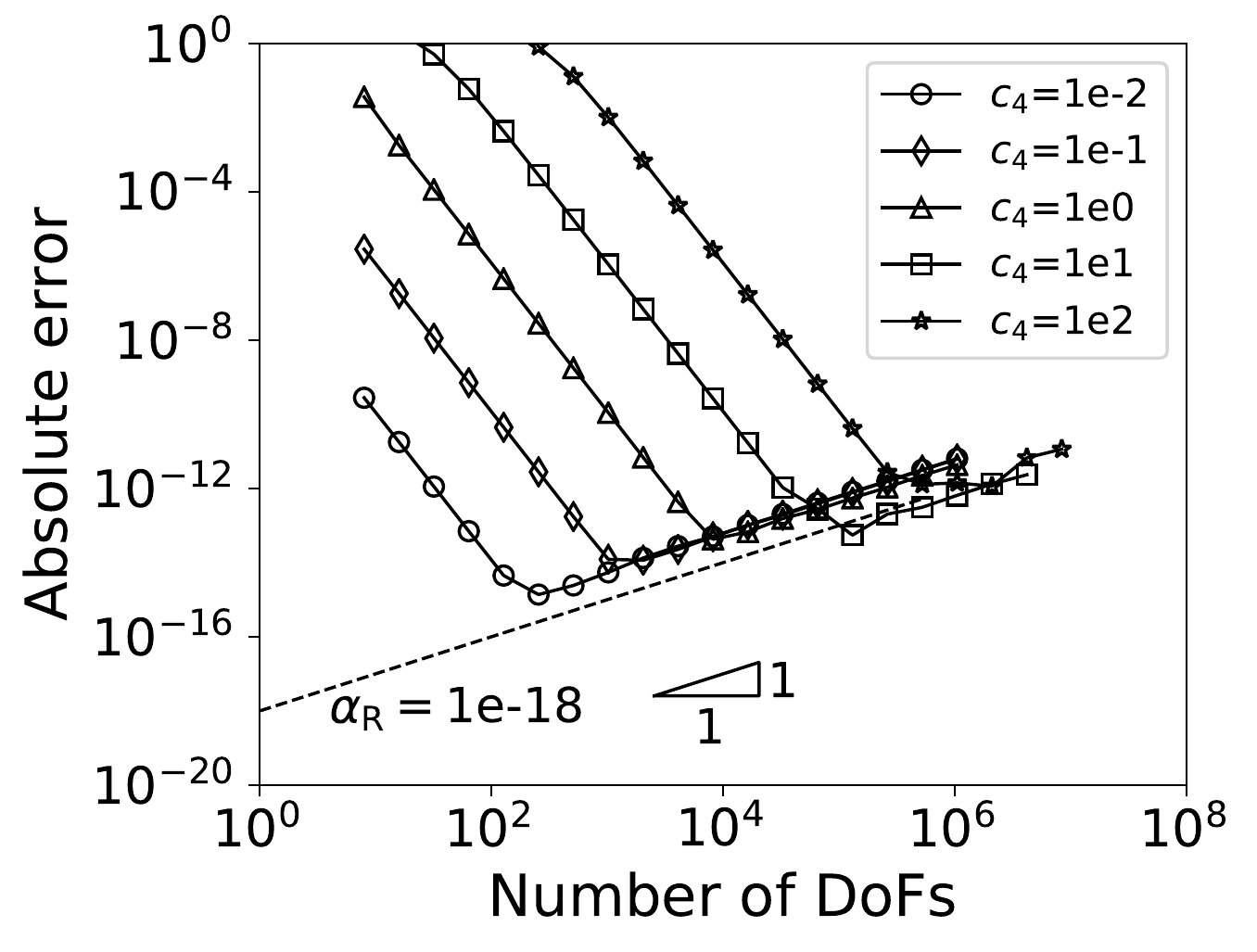}
        \caption{Solution}
        \label{py_L2_Pois4_MM_scaling_M1_solu}
    \end{subfigure}
    \hspace{-0.2cm}
    \begin{subfigure}{5.5cm}
        \includegraphics[width=1.0\linewidth]{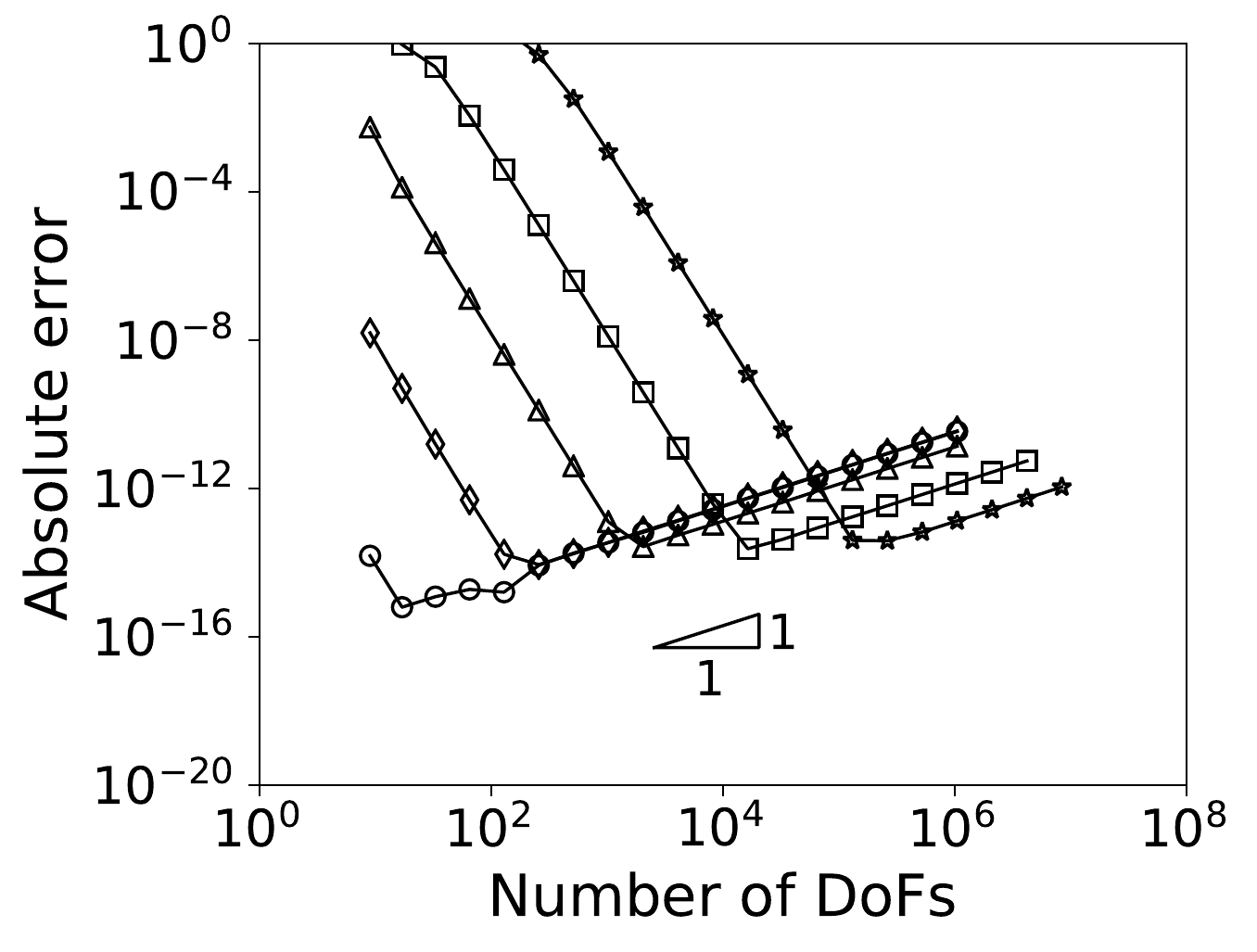}
        \caption{First derivative}
        \label{py_L2_Pois4_MM_scaling_M1_grad}
    \end{subfigure}
    \hspace{-0.2cm}
    \begin{subfigure}{5.5cm}
        \includegraphics[width=1.0\linewidth]{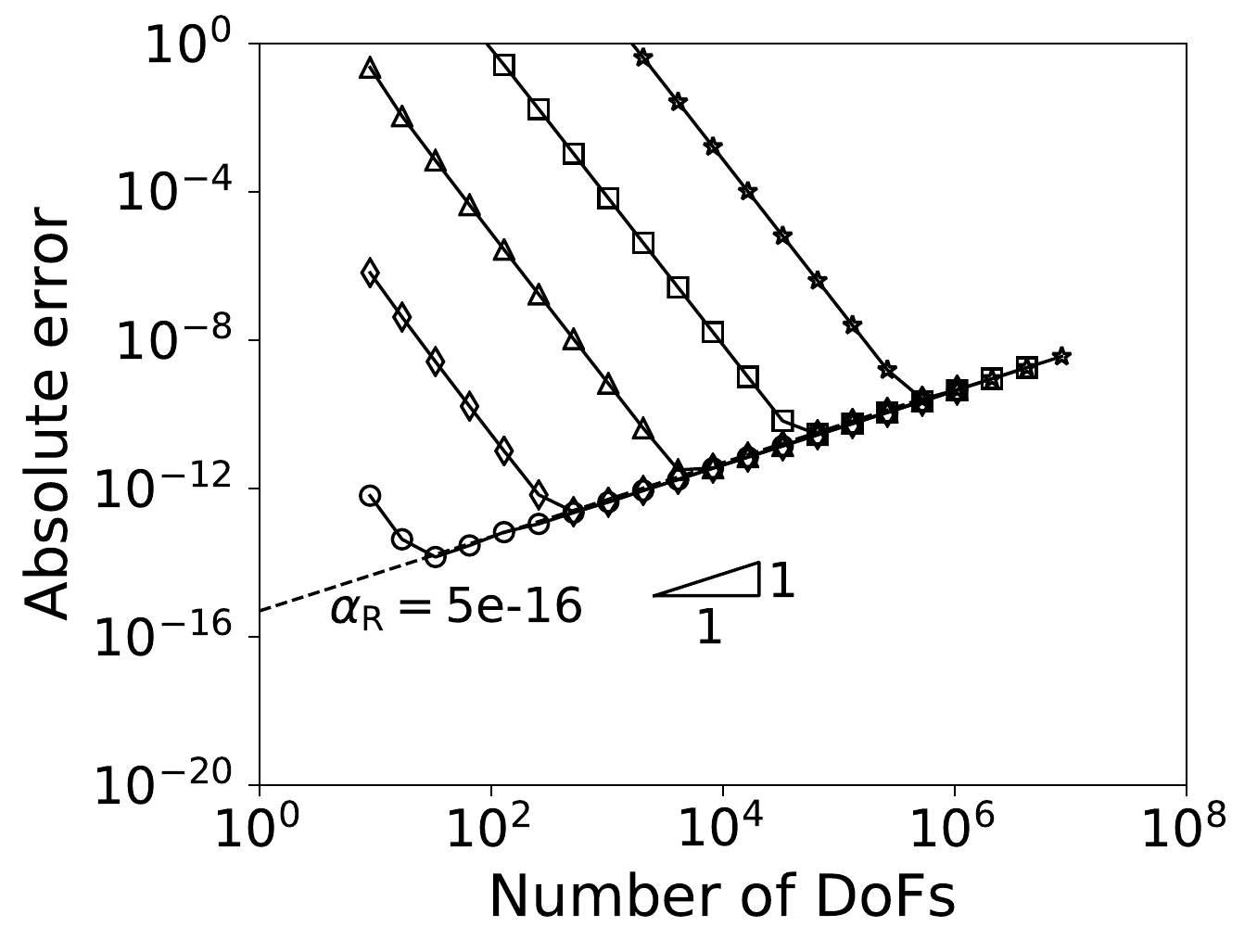}
        \caption{Second derivative}
        \label{py_L2_Pois4_MM_scaling_M1_2ndd}
    \end{subfigure}
\caption{Absolute errors for Case 4 in Table \ref{scaling_cases_Poisson} using scheme $M_1$.}
\label{py_L2_Pois4_MM_scaling_M1}
\end{figure}

\begin{figure}[!ht]
    \begin{subfigure}{5.5cm}
        \includegraphics[width=1.0\linewidth]{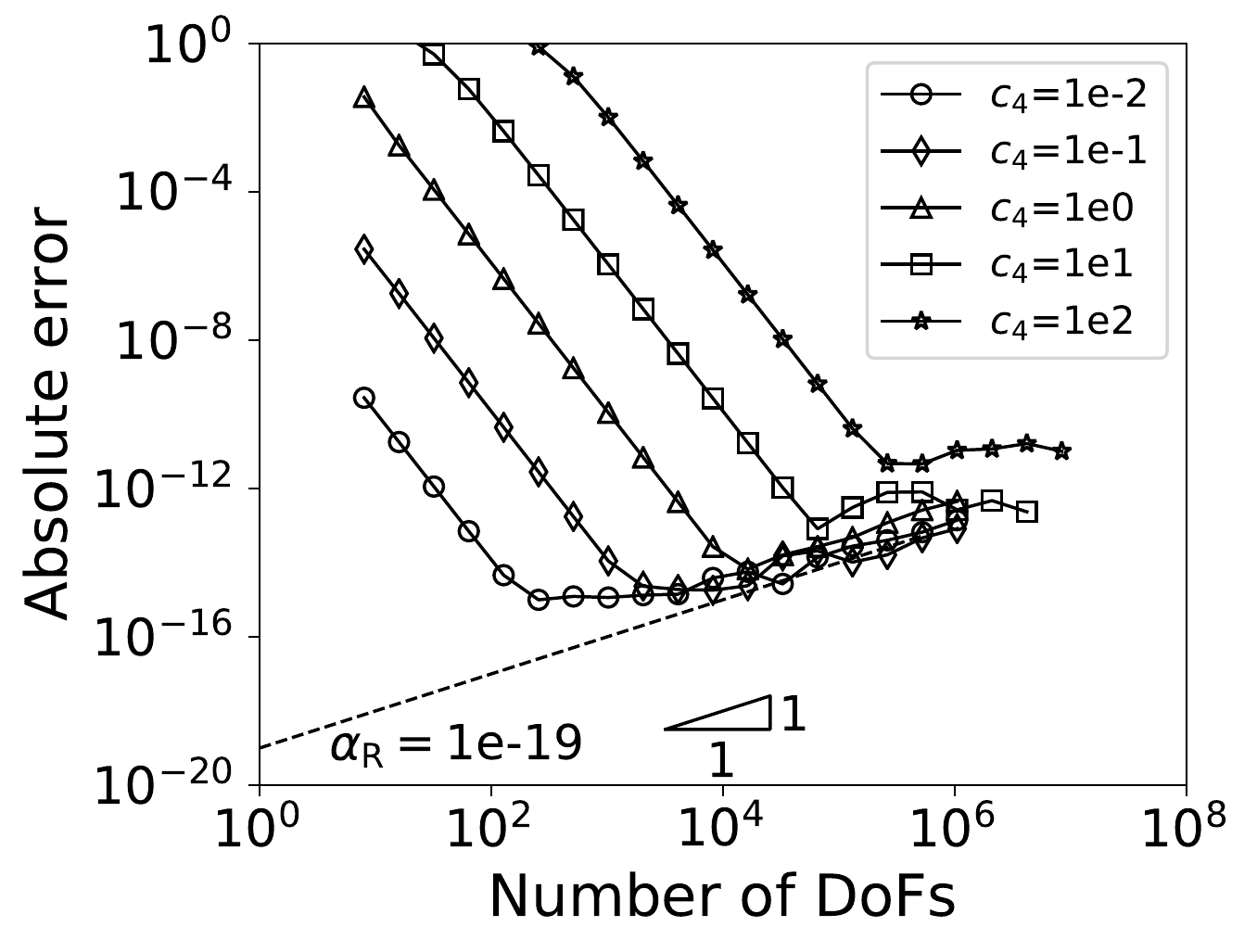}
        \caption{Solution}
        \label{py_L2_Pois4_MM_scaling_M2_solu}
    \end{subfigure}
    \hspace{-0.2cm}
    \begin{subfigure}{5.5cm}
        \includegraphics[width=1.0\linewidth]{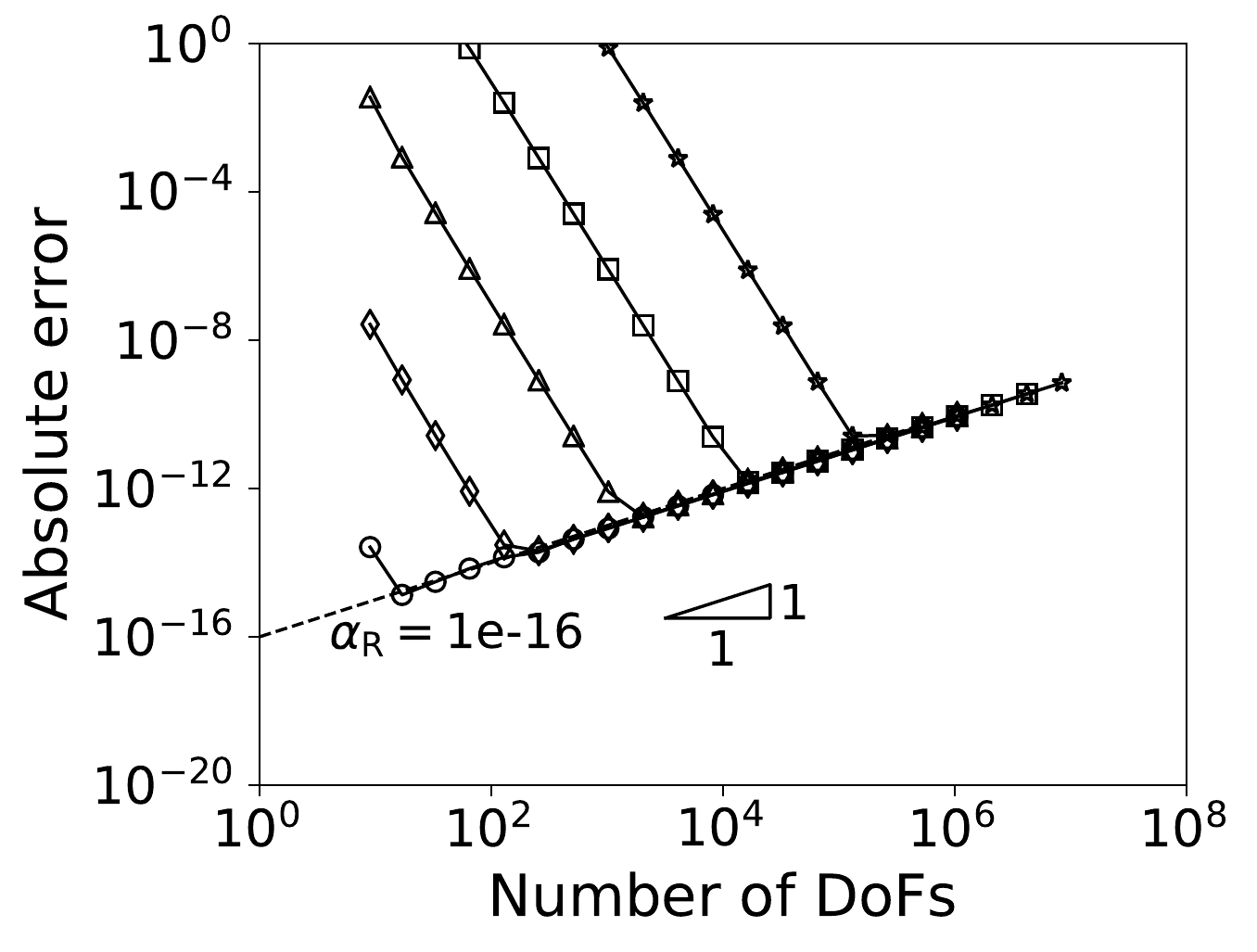}
        \caption{First derivative}
        \label{py_L2_Pois4_MM_scaling_M2_grad}
    \end{subfigure}
    \hspace{-0.2cm}
    \begin{subfigure}{5.5cm}
        \includegraphics[width=1.0\linewidth]{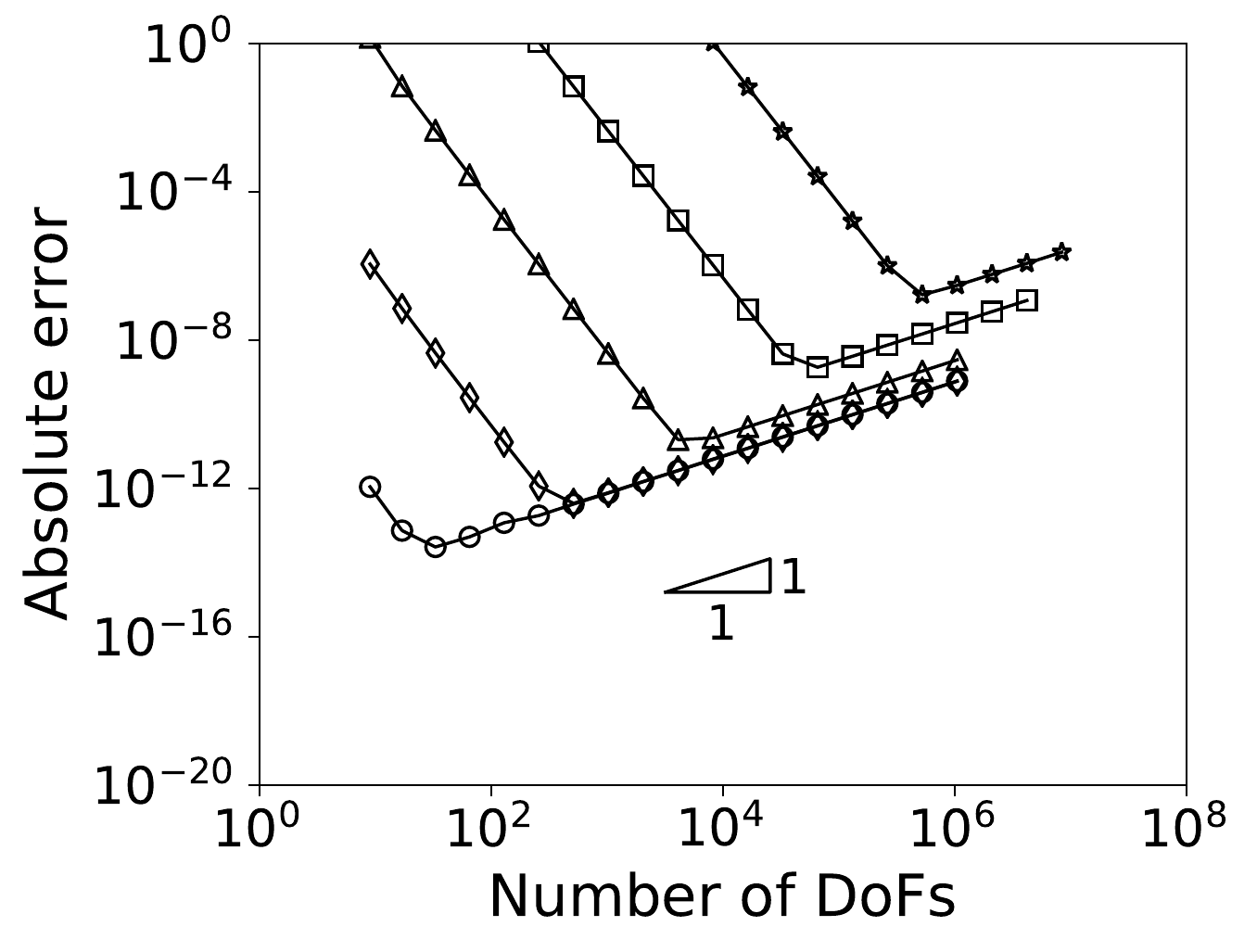}
        \caption{Second derivative}
        \label{py_L2_Pois4_MM_scaling_M2_2ndd}
    \end{subfigure}
\caption{Absolute errors for Case 4 in Table \ref{scaling_cases_Poisson} using scheme $M_2$.}
\label{py_L2_Pois4_MM_scaling_M2}
\end{figure}

\newpage
\paragraph{Case 4}
For Case 4, using the mixed FEM without scaling the right-hand side and schemes ${\rm M}_1$ and ${\rm M}_2$, the absolute errors are shown in Figs. \ref{py_L2_Pois4_MM_scaling_no}--\ref{py_L2_Pois4_MM_scaling_M2}.

\begin{figure}[!ht]
    \begin{subfigure}{5.5cm}
        \includegraphics[width=1.0\linewidth]{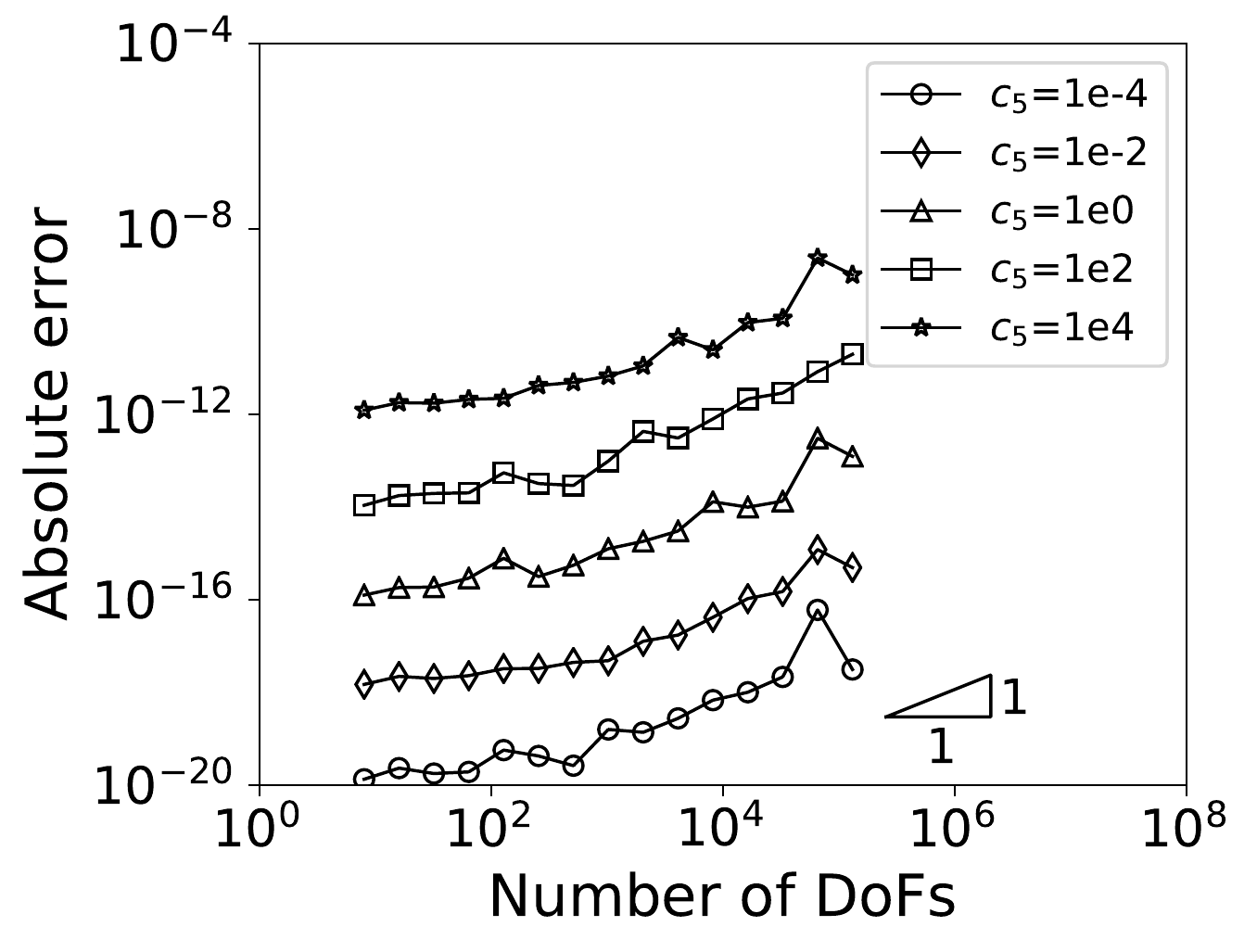}
        \caption{Solution}
        \label{py_L2_Pois5_MM_scaling_no_solu}
    \end{subfigure}
    \hspace{-0.2cm}
    \begin{subfigure}{5.5cm}
        \includegraphics[width=1.0\linewidth]{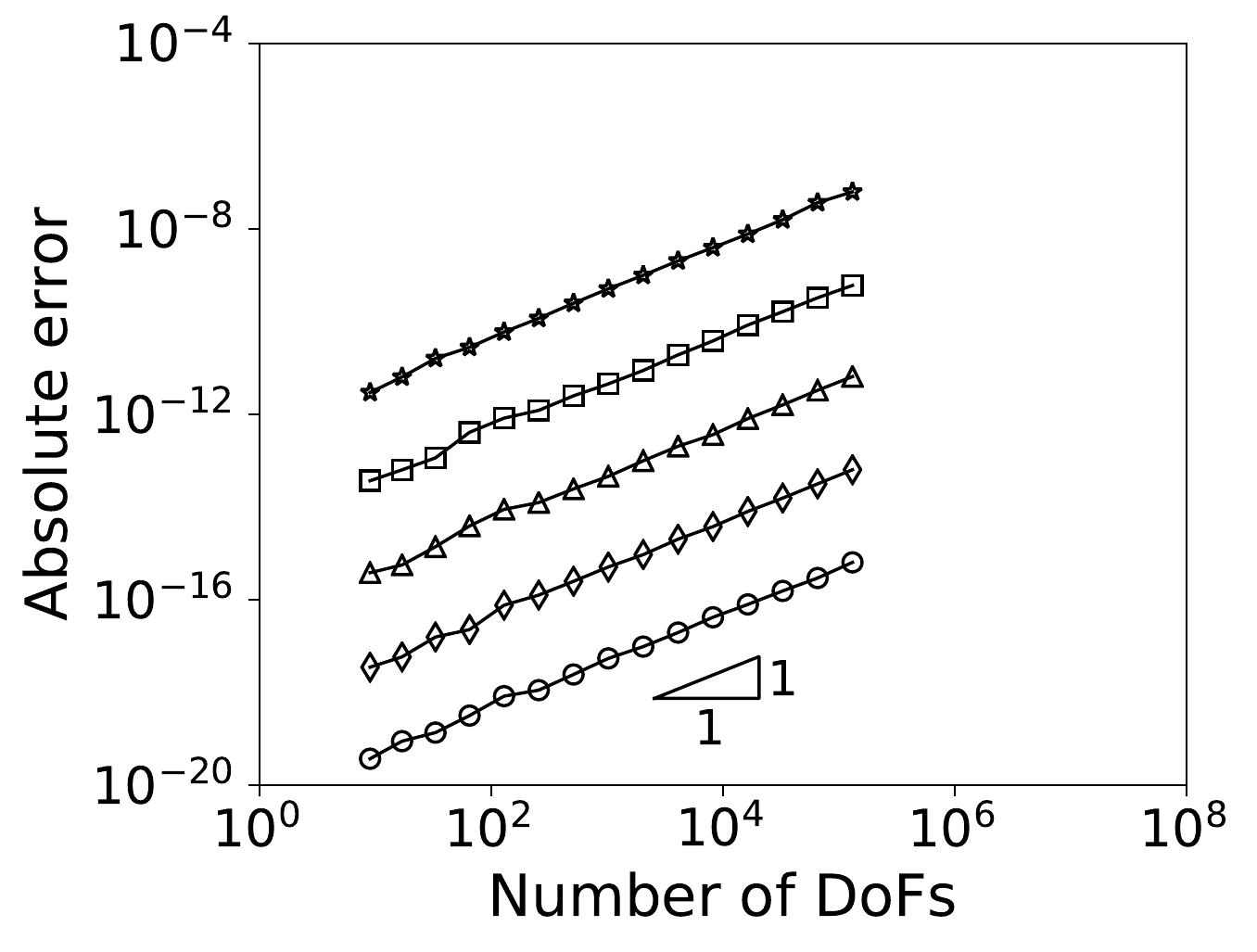}
        \caption{First derivative}
        \label{py_L2_Pois5_MM_scaling_no_grad}
    \end{subfigure}
    \hspace{-0.2cm}
    \begin{subfigure}{5.5cm}
        \includegraphics[width=1.0\linewidth]{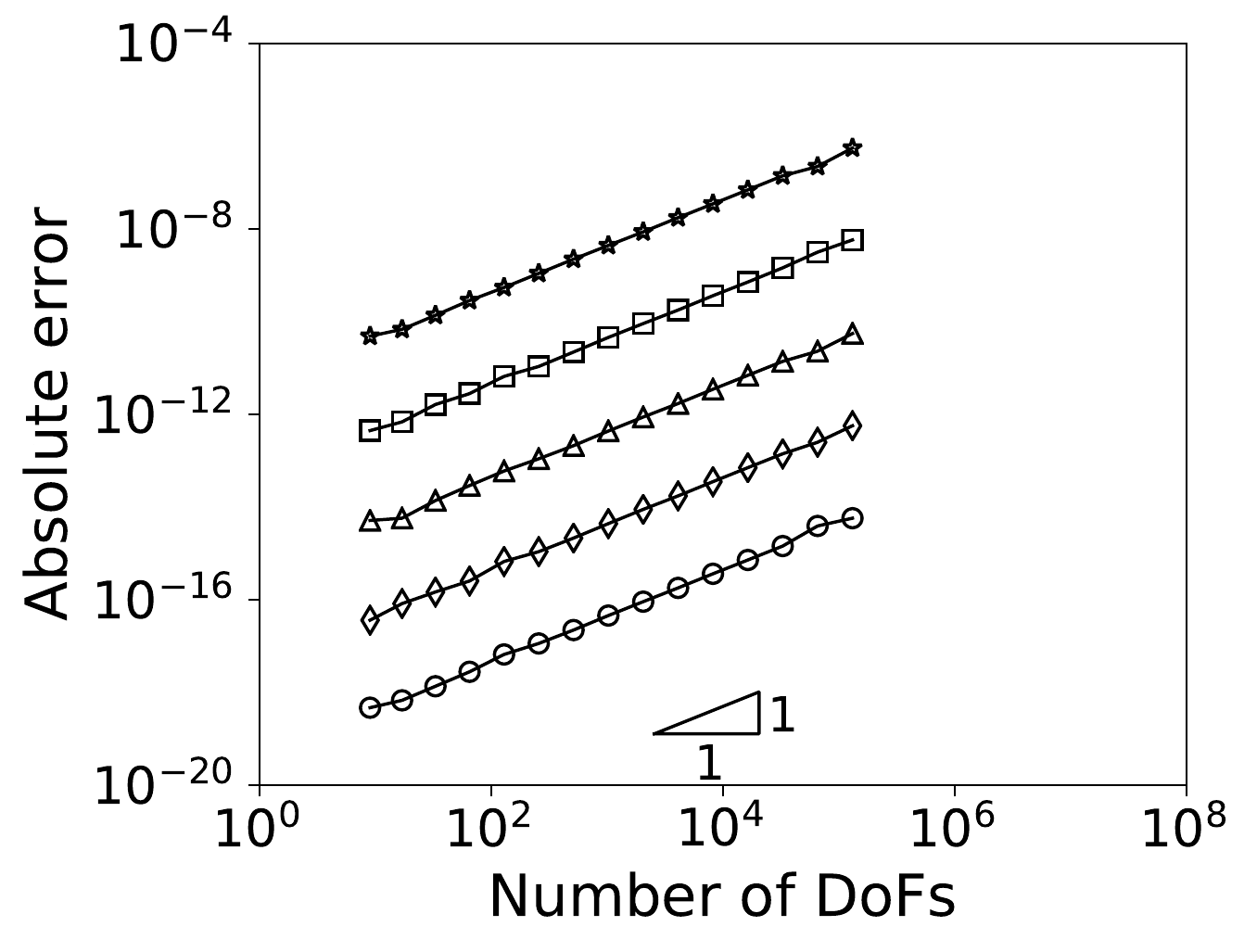}
        \caption{Second derivative}
        \label{py_L2_Pois5_MM_scaling_no_2ndd}
    \end{subfigure}
\caption{Absolute errors for Case 5 in Table \ref{scaling_cases_Poisson} using the mixed FEM without scaling the right-hand side.}
\label{py_L2_Pois5_MM_scaling_no}
\end{figure}

\begin{figure}[!ht]
    \begin{subfigure}{5.5cm}
        \includegraphics[width=1.0\linewidth]{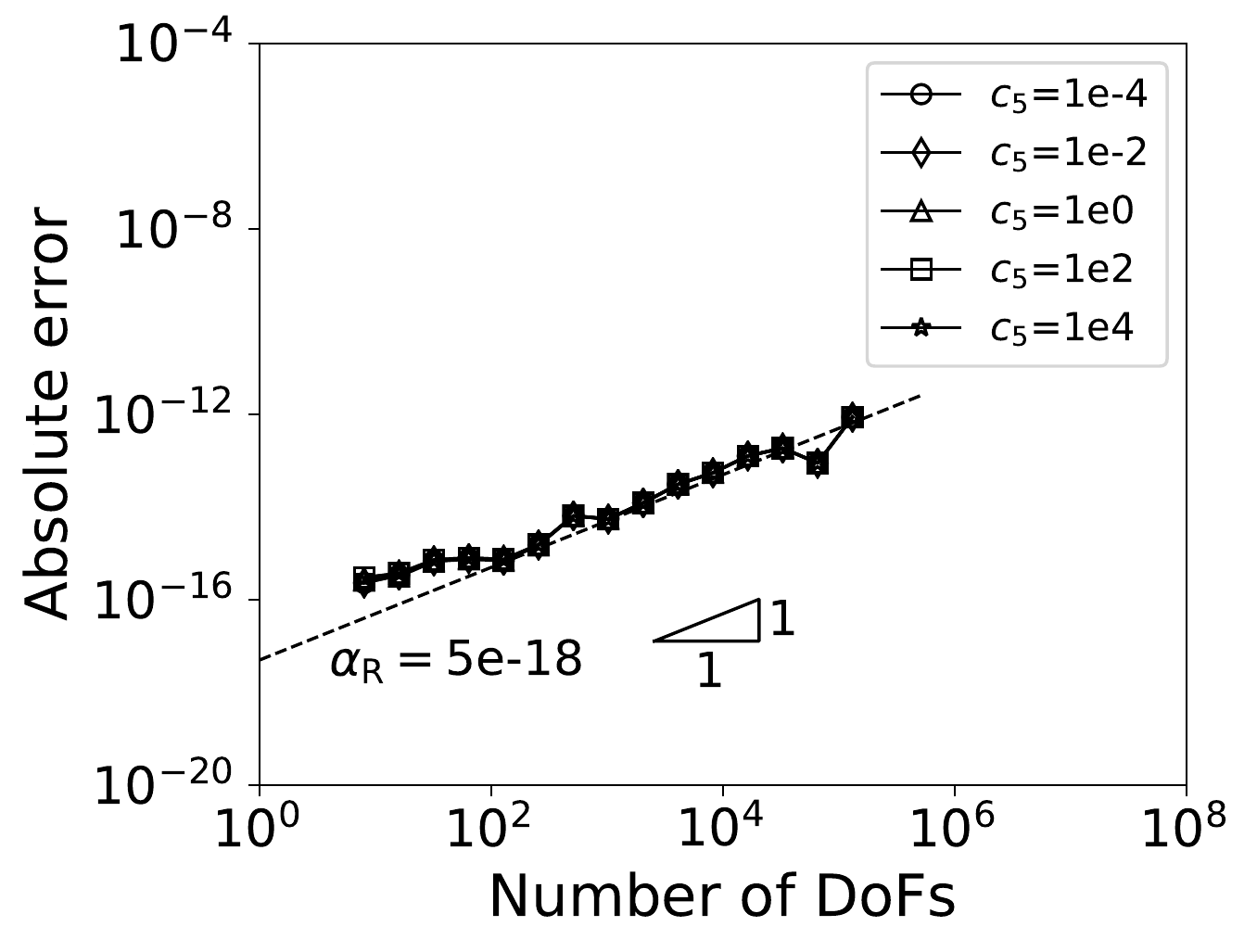}
        \caption{Solution}
        \label{py_L2_Pois5_MM_scaling_M1_solu}
    \end{subfigure}
    \hspace{-0.2cm}
    \begin{subfigure}{5.5cm}
        \includegraphics[width=1.0\linewidth]{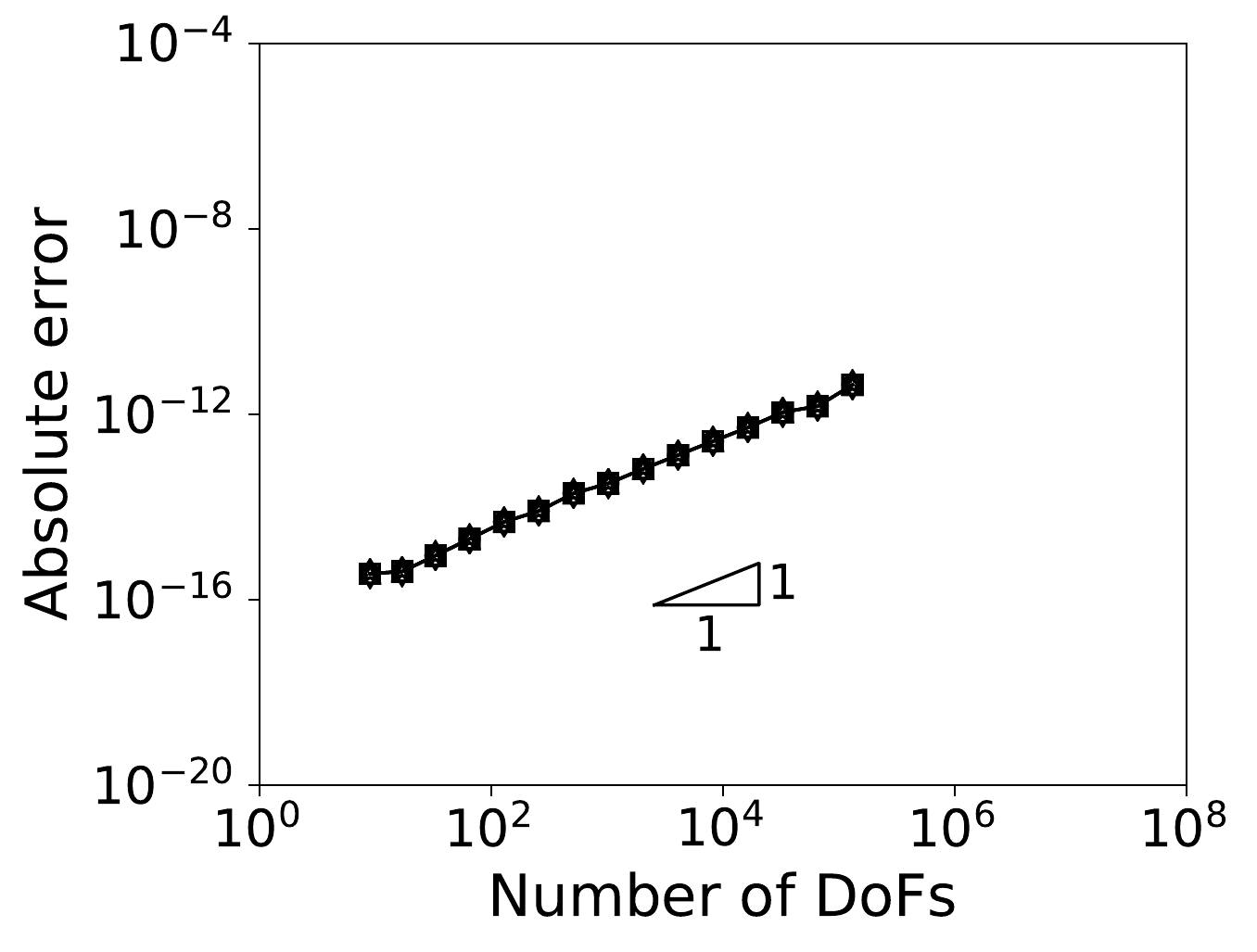}
        \caption{First derivative}
        \label{py_L2_Pois5_MM_scaling_M1_grad}
    \end{subfigure}
    \hspace{-0.2cm}
    \begin{subfigure}{5.5cm}
        \includegraphics[width=1.0\linewidth]{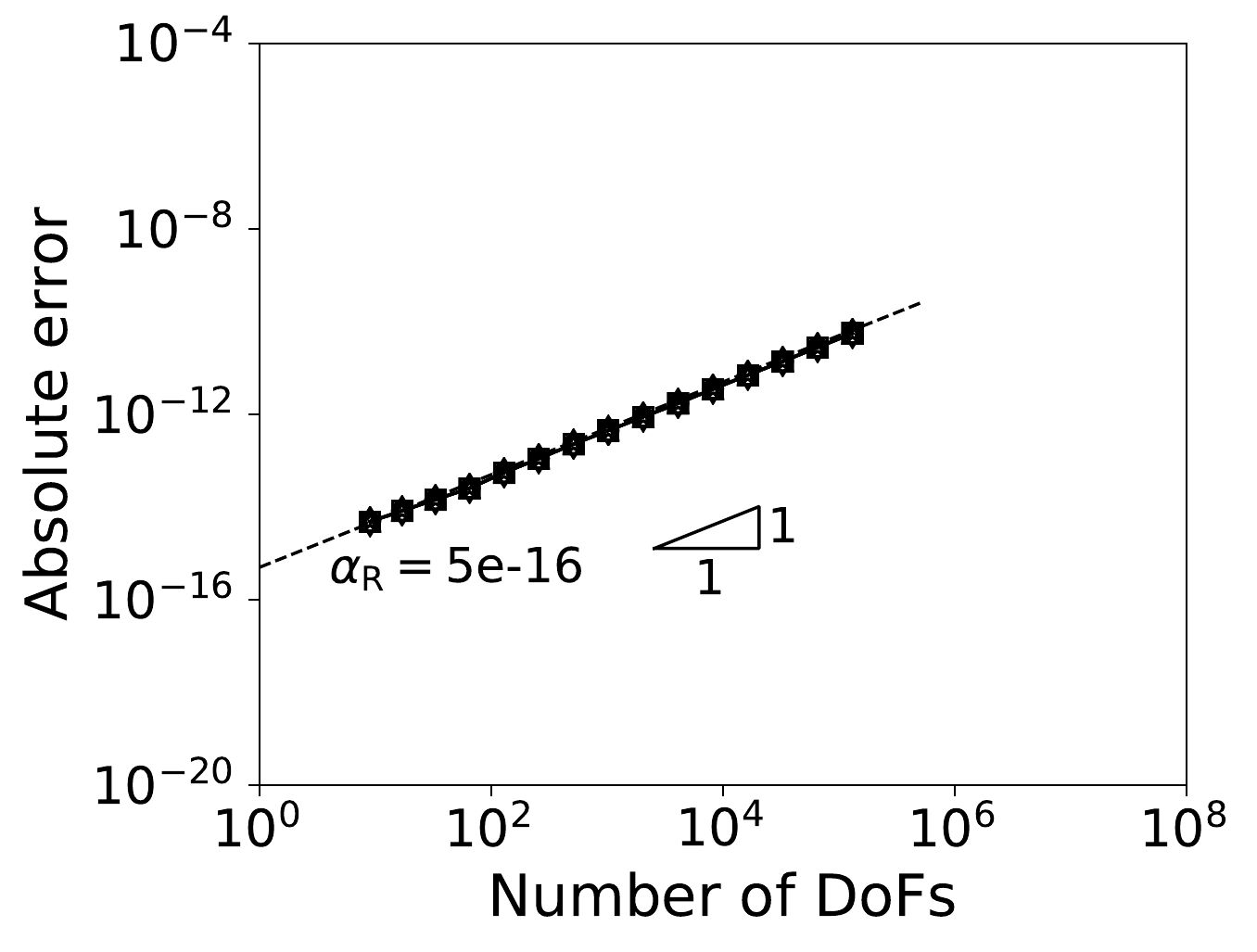}
        \caption{Second derivative}
        \label{py_L2_Pois5_MM_scaling_M1_2ndd}
    \end{subfigure}
\caption{Absolute errors for Case 5 in Table \ref{scaling_cases_Poisson} using scheme $M_1$.}
\label{py_L2_Pois5_MM_scaling_M1}
\end{figure}

\newpage
\paragraph{Case 5}
For Case 5, using the mixed FEM without scaling the right-hand side and schemes ${\rm M}_1$, the absolute errors are shown in Figs. \ref{py_L2_Pois5_MM_scaling_no}--\ref{py_L2_Pois5_MM_scaling_M1}.

\begin{figure}[!ht]
    \begin{subfigure}{5.5cm}
        \includegraphics[width=1.0\linewidth]{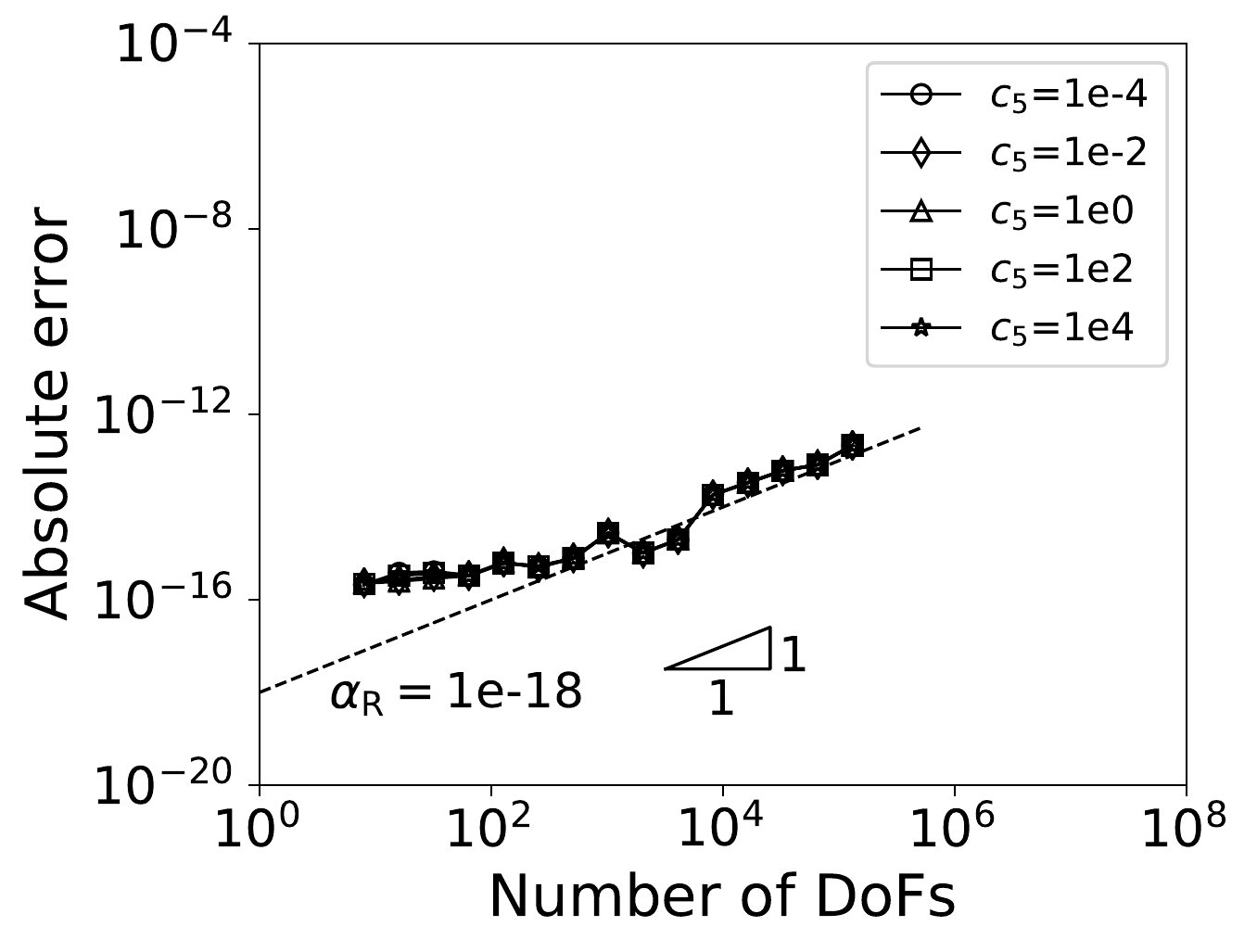}
        \caption{Solution}
        \label{py_L2_Pois5_MM_scaling_M2_solu}
    \end{subfigure}
    \hspace{-0.2cm}
    \begin{subfigure}{5.5cm}
        \includegraphics[width=1.0\linewidth]{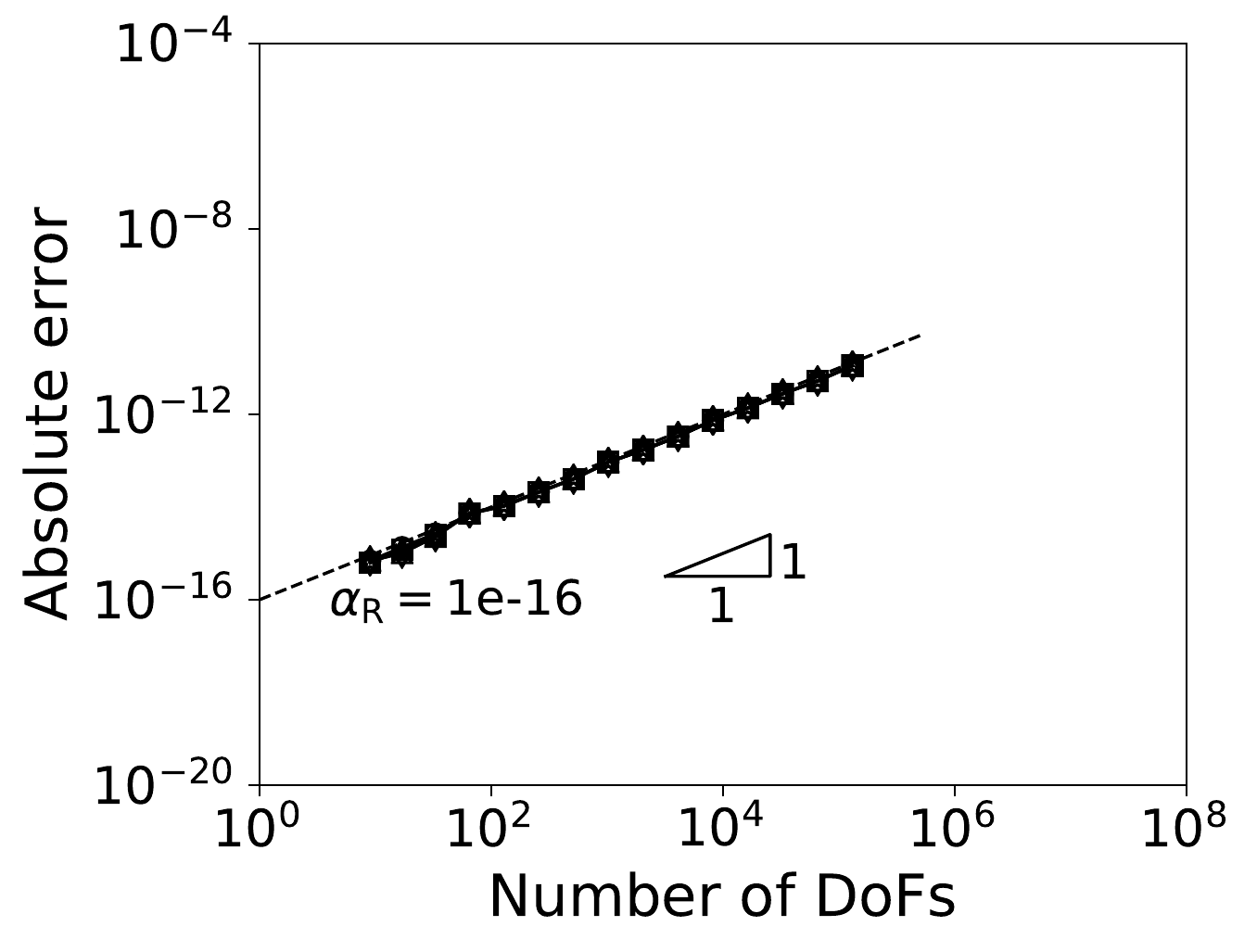}
        \caption{First derivative}
        \label{py_L2_Pois5_MM_scaling_M2_grad}
    \end{subfigure}
    \hspace{-0.2cm}
    \begin{subfigure}{5.5cm}
        \includegraphics[width=1.0\linewidth]{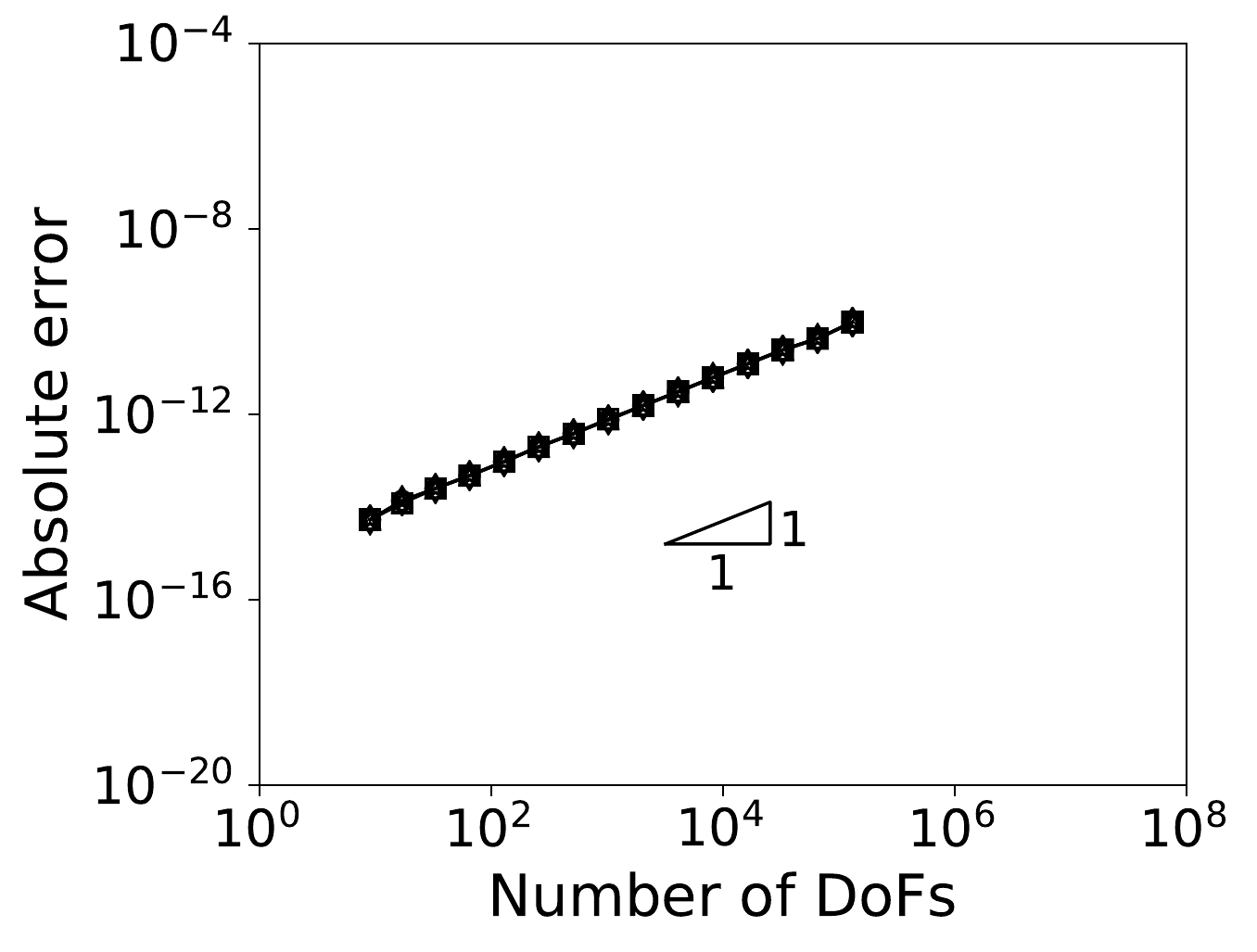}
        \caption{Second derivative}
        \label{py_L2_Pois5_MM_scaling_M2_2ndd}
    \end{subfigure}
\caption{Absolute errors for Case 5 in Table \ref{scaling_cases_Poisson} using scheme $M_2$.}
\label{py_L2_Pois5_MM_scaling_M2}
\end{figure}

\newpage

\bibliographystyle{unsrt}  
\bibliography{mybibfile}  

\begin{thebibliography}{10}

\bibitem{Kumar2016}
Mohit Kumar, Henk~M. Schuttelaars, Pieter~C. Roos, and Matthias M{\"o}ller.
\newblock Three-dimensional semi-idealized model for tidal motion in tidal
  estuaries.
\newblock {\em Ocean Dynamics}, 66(1):99--118, 2016.

\bibitem{carey1982derivative}
GF~Carey.
\newblock Derivative calculation from finite element solutions.
\newblock {\em Computer Methods in Applied Mechanics and Engineering},
  35(1):1--14, 1982.

\bibitem{ferziger2012computational}
Joel~H Ferziger and Milovan Peric.
\newblock {\em Computational methods for fluid dynamics}.
\newblock Springer Science \& Business Media, 2012.

\bibitem{zuras2008ieee}
Dan Zuras, Mike Cowlishaw, Alex Aiken, Matthew Applegate, David Bailey, Steve
  Bass, Dileep Bhandarkar, Mahesh Bhat, David Bindel, Sylvie Boldo, et~al.
\newblock {IEEE} standard for floating-point arithmetic.
\newblock {\em IEEE Std 754-2008}, pages 1--70, 2008.

\bibitem{gockenbach2006understanding}
Mark~S Gockenbach.
\newblock {\em Understanding and implementing the finite element method},
  volume~97.
\newblock Siam, 2006.

\bibitem{guo1986hp}
B~Guo and I~Babu{\v{s}}ka.
\newblock The hp version of the finite element method.
\newblock {\em Computational Mechanics}, 1(1):21--41, 1986.

\bibitem{Babuska2018Roundoff}
Ivo Babuska and Gustaf S\"oderlind.
\newblock On roundoff error growth in elliptic problems.
\newblock {\em ACM Transactions on Mathematical Software}, 44(3):1--22, 2018.

\bibitem{ling1984numerical}
Fuyun Ling and J~Proakis.
\newblock Numerical accuracy and stability: Two problems of adaptive estimation
  algorithms caused by round-off error.
\newblock In {\em Acoustics, Speech, and Signal Processing, IEEE International
  Conference on ICASSP'84.}, volume~9, pages 571--574. IEEE, 1984.

\bibitem{mou2017example}
Shan-Cong Mou, Yu-Xuan Luan, Wen-Tao Ji, Jian-Fei Zhang, and Wen-Quan Tao.
\newblock An example for the effect of round-off errors on numerical heat
  transfer.
\newblock {\em Numerical Heat Transfer, Part B: Fundamentals}, 72(1):21--32,
  2017.

\bibitem{ainsworth1992procedure}
Mark Ainsworth and J~Tinsley Oden.
\newblock A procedure for a posteriori error estimation for hp finite element
  methods.
\newblock {\em Computer Methods in Applied Mechanics and Engineering},
  101(1-3):73--96, 1992.

\bibitem{kelly1983posteriori}
DW~Kelly, De~SR Gago, OC~Zienkiewicz, I~Babuska, et~al.
\newblock A posteriori error analysis and adaptive processes in the finite
  element method: Part {I} -- error analysis.
\newblock {\em International Journal for Numerical Methods in Engineering},
  19(11):1593--1619, 1983.

\bibitem{boffi2013mixed}
Daniele Boffi, Franco Brezzi, Michel Fortin, et~al.
\newblock {\em Mixed finite element methods and applications}, volume~44.
\newblock Springer, 2013.

\bibitem{alzetta2018deal}
Giovanni Alzetta, Daniel Arndt, Wolfgang Bangerth, Vishal Boddu, Benjamin
  Brands, Denis Davydov, Rene Gassm{\"o}ller, Timo Heister, Luca Heltai,
  Katharina Kormann, et~al.
\newblock The deal.\rom{2} library, version 9.0.
\newblock {\em Journal of Numerical Mathematics}, 26(4):173--183, 2018.

\bibitem{davis2004algorithm}
Timothy~A Davis.
\newblock Algorithm 832: {UMFPACK} {V}4.3 -- an unsymmetric-pattern
  multifrontal method.
\newblock {\em ACM Transactions on Mathematical Software (TOMS)},
  30(2):196--199, 2004.

\bibitem{Runborg2012VerifyingNC}
Olof Runborg.
\newblock Lecture notes in numerical solutions of differential equations
  (dn2255): Verifying numerical convergence rates, 2012.

\bibitem{WalterFrei}
Meshing considerations for linear static problems.
\newblock
  \url{https://www.comsol.com/blogs/meshing-considerations-linear-static-problems/}.
\newblock Accessed: 2019-12-9.

\bibitem{kahan2013floating}
W~Kahan.
\newblock Floating-point tricks to solve boundary-value problems faster.
\newblock {\em University of California@ Berkeley}, 2013.

\bibitem{ginsburg1963cg}
Theo Ginsburg.
\newblock The conjugate gradient method.
\newblock {\em Numer. Math.}, 5(1):191--200, December 1963.

\bibitem{bazilevs2007weak}
Yuri Bazilevs and Thomas~JR Hughes.
\newblock Weak imposition of {D}irichlet boundary conditions in fluid
  mechanics.
\newblock {\em Computers \& Fluids}, 36(1):12--26, 2007.

\end{thebibliography}

\end{document}